\documentclass[11pt,letterpaper]{article}
\usepackage[margin=0.91in]{geometry}
\usepackage[T1]{fontenc}
\usepackage{lmodern,microtype}
\usepackage{amsmath,amssymb,amsthm,mathtools,mathrsfs}
\usepackage{enumitem,booktabs,array,longtable,tabularx}
\usepackage[dvipsnames]{xcolor}
\usepackage[colorlinks=true,linkcolor=NavyBlue,urlcolor=NavyBlue,citecolor=NavyBlue]{hyperref}
\usepackage{xurl}

\hypersetup{
	colorlinks=true,
	linkcolor=blue,
	citecolor=blue,
	urlcolor=blue,
}

\newtheorem{theorem}{Theorem}[section]
\newtheorem{proposition}[theorem]{Proposition}
\newtheorem{lemma}[theorem]{Lemma}
\newtheorem{corollary}[theorem]{Corollary}
\theoremstyle{definition}\newtheorem{definition}[theorem]{Definition}
\newtheorem{example}[theorem]{Example}
\theoremstyle{remark}\newtheorem{remark}[theorem]{Remark}
\newcommand{\Pkg}{\mathfrak P}
\newcommand{\HH}{\mathcal H}
\newcommand{\Obs}{\mathcal O}
\newcommand{\Tr}{\mathsf T}

\newcommand{\tw}{\omega}
\newcommand{\ipr}{\operatorname{ipr}}

\newcommand{\supp}{\operatorname{supp}}

\newcommand{\genlen}{\ell_{\mathrm{gen}}}

\theoremstyle{definition}\newtheorem{ddefinition}[theorem]{Definition}
\theoremstyle{plain}\newtheorem{dtheorem}[theorem]{Theorem}
\newtheorem{dlemma}[theorem]{Lemma}
\newtheorem{dproposition}[theorem]{Proposition}
\newtheorem{dcorollary}[theorem]{Corollary}
\theoremstyle{definition}\newtheorem{dexample}[theorem]{Example}
\theoremstyle{remark}\newtheorem{dremark}[theorem]{Remark}

\title{
	Dynamical Packaging Systems: Evolution, Invariants, and Local-to-Global Realization
}

\author{
	Rongchao Ma \\
	\textit{Department of Physics, University of Alberta} \\
	\textit{Edmonton, Alberta, Canada} \\
}

\date{}

\begin{document}
	\maketitle

	\begin{abstract}
		We couple deterministic, relational, and stochastic evolution to realized packaging systems through observations, interactions and admission.
		We characterize autonomous observation, exact formation descent, and complete temporal records for finite-group subgroup-matching interconnections.
		Under linear synchronization on a $D$-invariant subspace $W\le A$, conserved linear readouts modulo inherited ones form the dual of $(W\cap DA)/DW$, with $D$ the difference operator or generator.
		For a finite cover of obstruction rank $\nu$, sharp universal extension thresholds are $2k+1<\nu$ for stationary persistent empty-preserving causal closures of uniform prefix Carath\'eodory bound $k\ge1$; $1+q(r-1)<\nu$ for first-order relational atlases with projectively $r$-decomposable path relations at every horizon ($q,r\ge2$, alphabet size at most $q$); and $\max\{2,2k-1\}<\nu$ for independent coordinate-history blocks of size at most $k\ge1$.
		We characterize block-bound preservation and decide finiteness and compute finite temporal witness width for finite-state systems with finitely many observation coordinates and finite alphabets.
		Compatible initialized finite hidden Markov laws on join trees admit common finite realizations, with constructive bounds and rational solutions for rational input.
		Every non-join-tree cover admits probability-extendible finite hidden Markov local laws forcing infinite Hankel rank in every extension.
		Binary stationary ergodic obstructions persist, and structural admission can obstruct finite realization even on join trees.
	\end{abstract}

	\tableofcontents

	\section{Introduction}
	\label{dyn45:sec:intro}

	\subsection{Evolution on assembled systems}

	How does the structure of an assembled system govern its evolution, and how much of that evolution can be recovered from its parts?
	Static packaging~\cite{MaPackaging2026} describes native realization, while structural packaging~\cite{MaStructural2026} distinguishes interface matching from operational admission.
	Adding a temporal law raises three distinct problems.
	First, a global evolution may fail to induce autonomous laws on observed parts.
	The stochastic version of this problem is related to Markov lumpability~\cite{KemenySnell}.
	Second, compatible local laws may fail to extend globally, as in classical marginal extension and relational join theory~\cite{Vorobev1962,wn09:BeeriFaginMaierYannakakis1983}.
	Third, classical stochastic realization theory distinguishes the existence of a probability law from its finite hidden Markov realizability~\cite{dyn49:Heller,Anderson1999,Vidyasagar2011}.
	In the local-to-global setting, we ask whether prescribed local laws admit a common finite realization.
	These problems concern different aspects of the same system:
	its observations, structural admission, and temporal law.

	The distinction is visible in a two-coordinate system. For the binary swap $F(a,b)=(b,a)$ and the observation $r(a,b)=a$, there is no autonomous first-coordinate update:
	$r(Fx)$ is not constant on the fibers of $r$.
	Yet two successive first-coordinate observations determine the initial state (Theorem~\ref{thm:exact-deterministic-chart-factor-criterion} and Example~\ref{wn03:ex:swap}).
	For an assembly of two copies of a set, diagonal admission is preserved by the same permutation acting on both copies, but may fail under different permutations. Autonomy, recovery and preservation of admission are therefore distinct properties (Propositions~\ref{dyn46:prop-admission} and~\ref{dyn46:prop-stochastic-admission}).

	Dynamical packaging provides a setting in which these properties can be formulated together.
	A realized static system records a native context $\Xi$, a formation $P$, an output $\Pkg$, and the realization judgment $\Xi\vdash P\rightsquigarrow\Pkg$.
	We extend these data by the coupled object
	\[
	\mathbf D=(\mathbb P,X,\mathcal T,\mathcal C),\qquad
	\mathbb P=(\Xi,P,\Pkg),
	\]
	where $X$ is a state sort, $\mathcal T$ a typed temporal law, and $\mathcal C$ a structural--temporal coupling presentation (Definition~\ref{dyn46:def-coupling}).
	The coupling assigns observations and interaction supports to the realized base and specifies the declared constraints on temporal composition, admission and reconfiguration.
	State interpretation, native formation and operational admission thus retain their own mathematical roles.
	The raw temporal extension is conservative over the static axioms (Proposition~\ref{dyn45:prop-conservative-axioms}).
	Coupling also refines the comparison of systems:
	the same static base and conjugate unobserved time laws can give inequivalent observed dynamics (Example~\ref{dyn46:ex-coupled-noncollapse}).

	The theory proceeds from the comparison of dynamical systems and the recovery of their histories to exact extension and finite autonomous realization, and then to the invariants of structural interconnections.
	The atlas results use the declared observation data.
	The formation, admission and compositional results retain the additional native structure. In particular, exact extension of local laws, finite realization of those laws, and preservation of structural constraints are separate mathematical problems.

	\subsection{Principal results}

	Let $\mathcal H$ be a finite indexed observation cover of a finite coordinate set $V$.
	The geometric obstruction rank $\nu(\mathcal H)$, inherited from~\cite{MaPackaging2026} (Definition~\ref{def:nu}), equals $3$ on a nonchordal primal graph and otherwise equals the smallest unguarded clique size, or $\infty$ if every clique is guarded.
	A clique is guarded if it lies in an original chart.
	Thus for $m\ge3$, $m<\nu(\mathcal H)$ means chordality together with guarding of every clique of size at most $m$.
	The condition $\nu=\infty$ is exactly the classical join-tree condition~\cite{wn09:BeeriFaginMaierYannakakis1983}. Universal extension means extension of every compatible atlas in the specified class on the fixed cover.

	\paragraph{1. Temporal laws, observations and structural evolution.}
	Autonomous factors characterize the observations that support laws of their own.
	The greatest relation respecting the local rules on a prescribed admissible state space determines exact realization of all finite paths (Theorem~\ref{thm:greatest-relation}).
	For first-order chart systems on a join tree, agreement of complete separator path languages gives a single exact global relation (Theorem~\ref{thm:acyclic-temporal-exactness-criterion}); finite-state tests detect disagreement.
	For a possibly partial and multivalued native realization relation $J\subseteq S\times Y$ and source law $E$, Theorem~\ref{dyn47:thm:rel-descent} characterizes every output law $D$ satisfying
	\[
	J\circ E=D\circ J.
	\]
	It identifies existence, the greatest solution, forced transitions and finite minimum costs.
	Theorem~\ref{dyn48:thm:reconfiguration} extends the result to changing state types and masks.

	\paragraph{2. Intrinsic temporal recovery and sharp extension.}
	Observation depth is the least horizon at which the observed history determines the initial state, or infinity if none exists. Residual languages describe future continuations, and temporal witness width measures the number of specified coordinate--time cells needed to refute a partial history.
	The first two notions draw on classical sequential-machine and residual-language methods~\cite{Moore1956,Nerode1958}. Theorem~\ref{wn15:thm:effective-exact-rank} recognizes minimal impossible partial histories by a finite automaton and decides whether witness width is finite, computing the exact finite value.
	Uniform time expansion preserves $\nu(\mathcal H)$ (Theorem~\ref{thm:obstruction-rank-invariance-expansion}).
	The following sharp thresholds govern three recovery mechanisms:
	\begin{enumerate}[label=(\roman*),leftmargin=*,nosep]
		\item \emph{Causal closures.}
		Extending the static classification of~\cite{MaPackaging2026}, we show that for compatible stationary, persistent, empty-preserving causal closures with uniform finite-horizon Carath\'eodory bound $k\ge1$, universal exact extension holds if and only if $2k+1<\nu(\mathcal H)$.
		The least extension is coherent under truncation (Theorem~\ref{thm:all-horizon-causal-Horn-globalization}).
		
		\item \emph{Relational histories.}
		For $q,r\ge2$, compatible nonempty relations over alphabets of size at most $q$, every projection of which is determined by its at-most-$r$ coordinate projections, have universal exact extension if and only if $1+q(r-1)<\nu(\mathcal H)$ (Theorem~\ref{thm:pbr}).
		The numerical bound is the classical strong-consistency bound~\cite{Dechter1992}.
		The same cover criterion holds for first-order chart dynamics whose complete path relations satisfy these hypotheses at every horizon and agree on overlaps, on the compatible-state envelope (Theorem~\ref{wn13:thm:temporal-pbr}).
		
		\item \emph{Weighted histories.}
		For compatible initialized processes with independent coordinate-history blocks of size at most $k\ge1$, allowing binary alphabets, the universal probability extension criterion is $\max\{2,2k-1\}<\nu(\mathcal H)$ (Theorem~\ref{dyn48:thm:sharp-wug}).
		This is also the criterion for common finite hidden Markov realizability within the finite-input class (Corollary~\ref{dyn53:cor:finite-wug}).
		To preserve the block bound in the extension, for $k\ge2$ one additionally requires every primal connected component of size greater than $k$ to be complete; $k=1$ requires no restriction (Theorem~\ref{dyn53:thm:block-preserving}).
	\end{enumerate}
	The three criteria concern respectively implication determination, recovery by witnesses, and independent probability blocks.
	In each case counterexamples on every failing cover establish the converse.

	\paragraph{3. Finite autonomous realization and its geometric boundary.}
	On a join tree, agreement of separator word laws is necessary and sufficient for a common finite realization of prescribed initialized finite hidden Markov laws (Theorem~\ref{dyn53:thm:universal}).
	The constructive proof combines positive separator refinement, reversal of stochastic links, a shared driver, and reset augmentation with harmonic extraction.
	It supplies state bounds and preserves rationality.
	The converse is a sharp classification: allowing the finite alphabets to vary, every non-join-tree cover admits probability-extendible finite rational local laws whose every global extension has infinite Hankel rank (Theorem~\ref{dyn54:thm:conditional-spectrum}).
	The obstruction to finite hidden Markov realization persists for binary stationary ergodic inputs (Corollary~\ref{dyn54:cor:binary-stationary}). The distinction is between extending probability laws and realizing them by one finite mechanism.
	The realization here preserves the prescribed word laws, and may introduce hidden states.
	Requiring additional operational admission can obstruct finite realization even on a join-tree observation cover (Corollary~\ref{dyn54:cor:regular-mask}).

	\paragraph{4. Compositional invariants, conservation and resilience.}
	For finite-group dynamics, temporal records compose exactly and classify the declared subgroup-matching continuations (Theorem~\ref{dinv:thm:complete}).
	Theorem~\ref{dinv:thm:quotient} gives a canonical observed quotient and a universal conserved homomorphic quotient.
	In linear synchronization, let $W\le A$ be invariant under the ambient difference operator or generator $D$.
	The space of conserved linear readouts on $W$ modulo those inherited from $A$ is
	\[
	\mathcal E(W)\cong ((W\cap DA)/DW)^*.
	\]
	Theorem~\ref{dinv:thm:emergence} identifies its dimension, observable and invisible components, and the condition for its universal vanishing. Synchronization can create arbitrarily many new conserved readouts even for invertible discrete evolution.
	Finally, a monotone failure model illustrates the fixed-point principles through an exact greatest-shield characterization and minimum-cost intervention formula (Theorems~\ref{thm:shield} and~\ref{thm:rescue}).

	\subsection{Relation to earlier work}

	\paragraph{Static and relational extension.}
	The geometry of covers comes from Vorob'ev's marginal extension theorem and the acyclic join theory of Beeri, Fagin, Maier and Yannakakis~\cite{Vorobev1962,wn09:BeeriFaginMaierYannakakis1983,Yannakakis1981}.
	Recent annotated-relation consistency results and closure-atlas results are developed in~\cite{AK2025,AKWitness2025,AK2026,Kim2026}.
	Our causal closure theorem uses the sharp static bounded-Carath\'eodory theorem of~\cite{MaPackaging2026} to obtain coherent finite-prefix extensions and preserve the temporal axioms.
	Dechter's strong-consistency theorem already gives the bound $1+q(r-1)$, with a binary lower bound in Theorem~3.4 of~\cite{Dechter1992}.
	The related mechanisms of elimination, tightness, set intersection and algebraic decomposability are studied in~\cite{DechterVanBeek1997,Freuder1982,VanBeekDechter1997,ZhangYap2006,JeavonsCohenCooper1998}; Gyssens' degree of cyclicity gives another hierarchy for join dependencies~\cite{Gyssens1986}.
	The present relational result is an exact fixed-cover classification for hereditary projective decomposability, with original-chart recovery and a single law for all finite horizons.
	The failing-cover reduction follows the classical chordless-cycle and unguarded-clique strategy~\cite{wn09:BeeriFaginMaierYannakakis1983,wn11:AtseriasKolaitis2020}.
	The weighted threshold instead uses independent history blocks.
	The witness-width algorithm uses regular and partial-word language methods~\cite{Nerode1958,wn09:MyhillNerode,wn15:Dassow2014Partial,wn15:BlanchetSadri2017Minimal} to count specified observable cells.

	\paragraph{Stochastic observation and realization.}
	Lumpability characterizes exact autonomous Markov observation~\cite{KemenySnell}.
	Classical positive and linear realization theory treats individual word laws~\cite{dyn49:Heller,CarlylePaz1971,Anderson1999,Vidyasagar2011}; related work treats sofic measures and exact hidden Markov reduction~\cite{dyn52:BoylePetersen,BealJugeMairessePerrin2026,GT2023}.
	The common positive refinement is related to span-shaped simulations of weighted automata~\cite{dyn53:SokolovaWoracek}.
	Crane, Holroyd and Russell~\cite{CraneHolroydRussell2026} construct Markov couplings of given chains with a common Markov image; their Example~1 exhibits failure when the image is not Markov.
	Our unmarked realization contract allows additional hidden states and arbitrary finite hidden Markov separator laws.
	Their construction also uses harmonic normalization; combined with a finite common positive refinement, their Markov-image theorem implies qualitative finite joinability by exposing common driver edges (comparison following Theorem~\ref{dyn53:thm:universal}).
	The proof here is constructive with explicit state bounds and rational realizability.
	The sharp geometric converse uses classical coupling bounds~\cite{Griffeath1975,Goldstein1979}, forced to saturation by a third marginal, to obtain infinite Hankel rank in every probability extension.

	\paragraph{Structural evolution and invariance.}
	Relational descent is related to weakest prespecification~\cite{HoareHe1987}; manifest and latent behavior under interconnection is central to behavioral systems theory~\cite{Willems2000}.
	The observation quotient, temporal group records and conservation defect use classical sequential-machine, group-code and linear-observability methods~\cite{Moore1956,dinv:ForneyTrott2004,dinv:ForneyGL2012,dinv:Sontag1998}.
	The structural results characterize descent, exact temporal composition and conservation relative to declared admission and comparison data.
	They apply established mathematical tools to the coupled objects, while the recovery and realization results supply quantitative classifications and obstruction constructions.

	\paragraph{Static foundations and collective extension.}
	We refer to~\cite{MaPackaging2026} and~\cite{MaStructural2026} as Papers~I and~II.
	Paper~I~\cite{MaPackaging2026} supplies the static realization axioms and the static recovery theorem recalled in Section~\ref{sec:static}. \phantomsection\label{comp:structural}Paper~II~\cite{MaStructural2026} supplies the structural admission and continuation language; the required definitions and dynamical arguments are stated here.
	\phantomsection\label{comp:collective}Collective dynamics concerns size-indexed interactions, macroscopic limits and phase behavior.
	Its finite-system interface is described in Section~\ref{dyn46:subsec-collective-interface}.

	\subsection{Organization}

	Sections~\ref{sec:fixed}--\ref{sec:morphisms-invariants} establish temporal foundations, comparison and invariance. Sections~\ref{sec:local-dynamical-atlases}--\ref{sec:global-theorems} develop local transition laws, their recovery and exact extension. Section~\ref{sec:temporal-formation} treats native formation descent and reconfiguration. Sections~\ref{sec:gluing-probability-laws}--\ref{sec:structural-admission-obstructions} develop probability gluing, finite autonomous realization and the obstructions imposed by cover geometry and admission.
	Section~\ref{sec:compositional-invariants-conservation} establishes compositional records, conservation and structural resilience. 
	Section~\ref{sec:discussion} identifies further mathematical directions.
	For a direct route to the finite hidden Markov results, read Definition~\ref{def:nu} and Sections~\ref{sec:gluing-probability-laws}--\ref{sec:structural-admission-obstructions}.
	For relational recovery, start with Definition~\ref{def:nu} and Sections~\ref{sec:local-dynamical-atlases}--\ref{sec:global-theorems}.

	\section{Foundations and non-collapse of dynamical packaging}
	\label{sec:fixed}

	A time law on a fixed package may be specified independently of incidence.
	Coupling enters when a local observation is to evolve autonomously or a joint update is to respect an admitted assembly.
	The factor criterion (Theorem~\ref{thm:exact-deterministic-chart-factor-criterion}) and the admission test (Proposition~\ref{dyn46:prop-admission}) give these requirements explicitly.
	Channel support identifies the coordinates on which they must be checked.

	\subsection{Static packaging data and specialization principles}
	\label{sec:static}

	First, we recall the static data needed below.
	The static axioms and the bounded-Carath\'eodory theorem are imported from Paper~I~\cite{MaPackaging2026}.
	The distinction between interface matching and operational admission is imported from Paper~II~\cite{MaStructural2026}.
	The dynamical definitions that follow add state sorts, observation laws and temporal comparisons to these data.
	A packaged output is
	\begin{equation}\label{eq:five-tuple}
		\Pkg=(\operatorname{Pkg}_{\Pkg},E_{\Pkg},\triangleleft_{\Pkg},\mathsf L_{\Pkg},\mathsf R_{\Pkg}).
	\end{equation}
	Its entries are the carriers, internal units, incidence, candidate sets and selected families, respectively.
	A native context $\Xi=(\mathscr C_\Xi,\operatorname{Car}_\Xi,\operatorname{Form}_\Xi,\mathsf J_\Xi)$ specifies the structural domain, admissible carriers, formations and realization relation.
	A realized system $\mathbb P=(\Xi,P,\Pkg)$ retains $(P,\Pkg)\in\mathsf J_\Xi$, written $\Xi\vdash P\rightsquigarrow\Pkg$.
	Realization may be many-to-many, so identifying formation states with endpoints requires a declared map or relation.
	In a structural presentation, $\operatorname{Match}_\alpha$ tests interface equality and $\mathsf M_\alpha$ selects operationally admitted tuples; interface soundness (H2) is $\mathsf M_\alpha\subseteq\operatorname{Match}_\alpha$.
	A native realization judgment and an allowed temporal step must be checked separately.

	For reference, the inherited conditions have the following content.
	Fix the small-set and ambient-set conventions of Paper~I.
	Conditions (S1)--(S5) require, respectively, carriers in $\mathscr C_\Xi$, a small internal-unit set, incidence in $E_{\Pkg}\times\operatorname{Pkg}_{\Pkg}$, a small candidate set for each carrier, and $\mathsf R_{\Pkg}(K)\subseteq\mathsf L_{\Pkg}(K)$.
	The context clauses (P1)--(P3) of Paper~I require $\operatorname{Car}_\Xi\subseteq\mathscr C_\Xi$, a set $\operatorname{Form}_\Xi$, and a relation $\mathsf J_\Xi\subseteq\operatorname{Form}_\Xi\times\mathbf{PkgDat}_\Xi$, where $\mathbf{PkgDat}_\Xi$ consists of packaged outputs over $\mathscr C_\Xi$ whose carriers are admitted.
	Conditions (Y1)--(Y4) require such a context, a formation in its formation-data space, an admissible packaged output, and the actual realization pair in $\mathsf J_\Xi$.
	No deterministic or total realization relation is assumed.

	In the structural enrichment of Paper~II,
	(A0) supplies a nonempty indexed family of realized bases,
	(A1) an observation-site poset with candidate sets, restriction maps and type-indexed selected families,
	(A2) typed formation signatures and their specified interface comparisons, and
	(A3) admission relations in the products of the selected families of the inputs and output.
	Original types retain their base candidates and selections; prospective output types have declared candidates.
	Condition
	(H1) requires identity and composition for restriction maps;
	(H2) requires interface soundness; the optional
	(H3) assigns natural-number ranks to types strictly increasing from each input to the output of a primitive formation.
	The rank concerns formation rather than time.
	These are the inherited typing requirements.
	The finite observable specialization and the incidence construction below link these data to the observation cover, while the imported static theorem supplies the geometric recovery criterion used in Section~\ref{sec:global-theorems}.

	For atlas specializations, let $\HH=(I_j)$ be a finite indexed cover of a descriptor-label set $V$.
	This cover need not be the family of primitive package interiors.
	Its geometric obstruction rank is recalled below.

	\begin{definition}[Typed finite observable specialization]
		\label{def:typed-specialization-v05}
		A \textbf{finite observable specialization} of a realized packaged system	$\mathbb P=(\Xi,P,\Pkg)$ consists of a finite descriptor set $V$, nonempty finite alphabets $(A_v)_{v\in V}$, a nonempty admissible state set $X\subseteq\prod_{v\in V} A_v$, and a specified finite cover $\HH$ of $V$.
		For $U\subseteq V$, write $X_U=\pi_U X$ and $r_U=\pi_U|_X:X\to X_U$.
		For $S\subseteq U$, the restriction map $\rho^U_S:X_U\to X_S$ is the coordinate projection.
		The maps obey $\rho^U_U=\mathrm{id}$ and $\rho^S_R\rho^U_S=\rho^U_R$ whenever $R\subseteq S\subseteq U$.
		The cover, alphabets, states, and coordinate identifications are additional data, not consequences of axioms (S1)--(S5).
	\end{definition}

	\begin{remark}[Two levels of restriction]\label{rem:restriction-level-v05}
		State restriction $r_U$ is always defined for a typed coordinate presentation.
		A transition relation need not admit an autonomous local transition system:
		projection of one-step transitions may create spurious longer trajectories.
		Accordingly, local trajectory semantics must be defined by projection of complete trajectories, unless a separate factorization or exact-lifting hypothesis is proved.
		The terms ``chart overlap'' and ``primitive package interior'' are not identified without an explicit incidence-resolved specialization.
	\end{remark}

	\begin{definition}[Geometric obstruction rank of a cover]
		\label{def:nu}
		Let $\HH$ be a finite cover of a finite set $V$.
		\begin{itemize}
			\item The \textbf{primal graph} of $\HH$ is the graph
			$G_{\HH}=(V,E_{\HH})$, where
			\[
			E_{\HH}
			=\bigl\{\{u,v\}:u\ne v,\ 
			\{u,v\}\subseteq I\text{ for some }I\in\HH\bigr\}.
			\]
			A nonempty clique $C$ of $G_{\HH}$ is \textbf{guarded} if $C\subseteq I$ for some $I\in\HH$, and \textbf{unguarded} otherwise.
			A graph is \textbf{chordal} if it has no induced cycle of length at least four.
			
			\item The \textbf{geometric obstruction rank} of $\HH$ is
			\[
			\nu(\HH):=
			\begin{cases}
				3, & G_{\HH}\text{ is nonchordal},\\
				\min\{|C|:C\text{ is an unguarded clique}\},
				& G_{\HH}\text{ is chordal},
			\end{cases}
			\]
			with the convention $\min\varnothing=\infty$.
		\end{itemize}
	\end{definition}

	The static rank $\ipr(\mathscr C)$ of Paper~I is a classwise rooted-circuit original-chart recovery rank for finite, empty-preserving closure systems.
	We retain the notation to distinguish this imported invariant from the temporal and weighted recovery parameters defined below.

	\begin{theorem}[Static globalization for bounded Carath\'eodory number, Ma~\cite{MaPackaging2026}]
		\label{thm:source}
		For finite compatible closure atlases with empty-preserving local Carath\'eodory number at most $k\ge1$, universal exact globalization holds exactly when
		\[
		2k+1<\nu(\HH).
		\]
		Under this condition, the closure generated by the union of the local implicational bases is the least exact extension.
		The extension is taken in the ambient class of finite empty-preserving closures, not necessarily in the class with Carath\'eodory number at most $k$.
	\end{theorem}

	This is Theorem~4.22 and the closure case of Theorem~8.52(i) of \cite{MaPackaging2026}.
	The temporal Horn correspondence uses this static theorem at each finite horizon and proves coherence of the resulting closures under time restriction.

	\subsection{Incidence-generated observations and the two modes of atlas data}

	The incidence relation determines which units a carrier can observe.
	Selecting the carriers gives the chart cover.
	One can then either prescribe a system and project it or prescribe local data and ask for a realization.

	\begin{definition}[Incidence supports and incidence-generated chart families]
		\label{def:incidence-charts-v06}
		Fix the packaged output $\Pkg$ of a realized system 
		as in~\eqref{eq:five-tuple}.
		\begin{enumerate}[label=(\roman*),leftmargin=*,nosep]
			\item \textbf{Incidence support.}
			The incidence support of a carrier 
			$K\in\operatorname{Pkg}_{\Pkg}$ is
			\[
			I_K:=\{e\in E_{\Pkg}:e\triangleleft_{\Pkg}K\}.
			\]
			
			\item \textbf{Incidence-generated chart family.}
			For a chosen nonempty family of observed carriers
			$\mathcal K\subseteq\operatorname{Pkg}_{\Pkg}$, 
			the incidence-generated chart family is
			\[
			\mathcal H_{\mathcal K}:=(I_K)_{K\in\mathcal K},
			\qquad
			V_{\mathcal K}:=\bigcup_{K\in\mathcal K}I_K.
			\]
			Carrier indices are retained even when supports coincide.
			Empty supports may be discarded only if they carry
			no additional chart data.
			
			\item \textbf{Full and partial incidence observation.}
			The incidence specialization is fully observed if
			$V_{\mathcal K}=E_{\Pkg}$, and partially observed
			otherwise, with unobserved units
			$E_{\Pkg}\setminus V_{\mathcal K}$ recorded explicitly.
		\end{enumerate}
	\end{definition}

	\begin{remark}[Intrinsic versus auxiliary charts]\label{rem:charts-source-v06}
		An incidence-generated cover is canonical \emph{relative to the selected observed carriers $\mathcal K$}, not canonical from the five-tuple
		without that selection.
		A descriptor cover on unrelated labels $V$ remains an admissible \emph{auxiliary atlas specialization}, but it is not called intrinsic.
		The candidate spaces $\mathsf L_{\Pkg}$ and	selected families $\mathsf R_{\Pkg}$ do not automatically determine state alphabets or chart state relations.
		If their contents are used to select observed carriers or constrain states, that use is an additional specified construction and must be verified for the intended example.
	\end{remark}

	\begin{definition}[Atlas-first and realized presentations]
		\label{def:two-modes-v06}
		Fix a cover $\mathcal H$ of $V$, coordinate alphabets $A_v$, and coherent restriction maps.
		Let $\pi_I$ denote the coordinate projection onto $I$.
		\begin{enumerate}[label=(\roman*),leftmargin=*,nosep]
			\item \textbf{Atlas-first presentation.}
			An atlas-first presentation is a family $(Z_I,R_I)_{I\in\mathcal H}$, where
			\[
			\varnothing\ne Z_I\subseteq\prod_{v\in I}A_v,
			\qquad R_I\subseteq Z_I^2,
			\]
			subject to the stated overlap-compatibility conditions.
			No global realization is assumed.
			
			\item \textbf{Realized presentation.}
			A realized presentation consists of a global admissible
			state set $X\subseteq\prod_{v\in V}A_v$ and a transition
			relation $E\subseteq X^2$.
			Its chart state sets, one-step shadows, and trajectory
			shadows are, respectively,
			\[
			X_I:=\pi_I X,\qquad
			(\pi_I\times\pi_I)(E),\qquad
			\pi_{I\times[0,m]}\mathsf T_E(m).
			\]
			For $m>1$, the trajectory shadow may be strictly smaller
			than the trajectory language of the one-step shadow.
		\end{enumerate}
		A realized presentation $(X,E)$ is an exact globalization of an atlas-first presentation $(Z_I,R_I)_{I\in\mathcal H}$ if
		\[
		\pi_{I\times[0,m]}\mathsf T_E(m)
		=\mathsf T_{R_I}(m)
		\quad\text{for all }I\in\mathcal H,\ m\geq0,
		\]
		where $\mathsf T_{R_I}(0)=Z_I$.
	\end{definition}

	\begin{proposition}[Faithfulness of incidence observations]\label{prop:incidence-faithfulness-v06}
		Put $U=\bigcup_{K\in\mathcal K}I_K$.
		For any set of states $X\subseteq\prod_{e\in E_{\Pkg}} A_e$, the joint observation $x\mapsto(\pi_{I_K}x)_{K\in\mathcal K}$ is injective exactly when $\pi_U$ is injective on $X$.
		In particular coverage $U=E_{\Pkg}$ suffices.
	\end{proposition}

	\begin{proof}
		Two observations agree precisely when the states agree on every coordinate of each $I_K$, equivalently on $U$.
		For the full product, an uncovered coordinate with two possible values therefore prevents injectivity.
		For a constrained $X$, the values on $U$ can determine all remaining coordinates.
		Even with coverage, compatible local tuples can lie outside the image.
		For example, take two singleton charts and restrict the joint binary states to the diagonal.
	\end{proof}

	\subsection{Primitive dynamical axioms and their models}
	\label{dyn45:subsec-axioms}

	The raw dynamical expansion adds a state sort and an evolution law to the realized packaging theory of Section~\ref{sec:static}.
	The coupled theory then specifies the additional diagrams and law data in Definition~\ref{dyn46:def-coupling}.
	Its primitive relational core is deliberately weaker than the finite-state, nonblocking, autonomous-factor, or exact-globalization specializations
	used later.
	The labels below specify which conditions belong to the definition of a model and which are additional hypotheses.

	\begin{definition}[Primitive fixed-support dynamical axioms]
		\label{dyn45:def-core-axioms}
		A \textbf{fixed-support relational dynamical packaged object} is a triple
		\[
		\mathbf D=(\mathbb P,X,E)
		\]
		satisfying the following axioms:
		\begin{description}[style=nextline,leftmargin=3em]
			\item[(D0) Realized base.]
			$\mathbb P=(\Xi,P,\mathfrak P)$ is a realized packaging system satisfying (S1)--(S5) and (Y1)--(Y4), with
			$\Xi\vdash P\rightsquigarrow\mathfrak P$.
			Its structural domain and declared native markings are retained.
			
			\item[(D1) State space.]
			$X$ is a specified nonempty small set associated with $\mathbb P$.
			Its elements are called states.
			No identification of $X$ with a native candidate or selected family is required.
			Any such identification is specified separately.
			
			\item[(D2) Temporal law.]
			$E\subseteq X\times X$ is the relation of permitted one-step transitions: $(x,y)\in E$ means that state $x$ may evolve to state $y$ in one step.
			No totality, determinism, or reversibility is required, and $E$ may be empty.
		\end{description}
		The base $\mathbb P$, state space $X$, and temporal law $E$ remain fixed throughout the evolution.
		The primitive dynamical axiom system is therefore
		\[
		\mathsf{Ax}_{\mathrm{Dyn}}^{\mathrm{raw}}
		=\mathsf{Ax}_{\mathrm{Pkg}}+(\mathrm{D1})+(\mathrm{D2}),
		\]
		where (D0) specifies the inherited static model.
	\end{definition}

	\begin{definition}[Trajectories and initialization]
		\label{dyn45:def-initialized}
		Let $X$ be a state space and $E\subseteq X^2$ a transition relation.
		For $m\in\mathbb N_0$:
		\begin{enumerate}[label=(\roman*),leftmargin=*,nosep]
			\item \textbf{Trajectories.}
			The set of trajectories of length $m$ is
			\[
			\mathsf T_E(m):=
			\{(x_0,\ldots,x_m)\in X^{m+1}:
			(x_t,x_{t+1})\in E\text{ for }0\le t<m\}.
			\]
			In particular, $\mathsf T_E(0)=X$.
			
			\item \textbf{Initialization.}
			An initialization is a specified initial state set
			$X_{\mathrm{init}}\subseteq X$.
			The corresponding initialized trajectory set is
			\[
			\mathsf T_{E,X_{\mathrm{init}}}(m):=
			\{(x_0,\ldots,x_m)\in\mathsf T_E(m):
			x_0\in X_{\mathrm{init}}\}.
			\]
			In particular,
			$\mathsf T_{E,X_{\mathrm{init}}}(0)=X_{\mathrm{init}}$.
		\end{enumerate}
		Unless otherwise specified, $X_{\mathrm{init}}=X$.
		Empty initial sets are allowed unless nonemptiness
		is required.
		We also write $X_0$ for $X_{\mathrm{init}}$ in initialized specializations.
	\end{definition}

	\begin{proposition}[Temporal consistency of the relational core]
		\label{dyn45:prop-path-semantics}
		Every contiguous subpath of a member of $\mathsf T_E(m)$ is a path.
		Two paths can be concatenated when the final state of the first is the initial state of the second.
		Consequently, if $E^{[m]}$ denotes the relation of endpoints of length-$m$ paths, then $E^{[0]}=\Delta_X$ and
		\[
		E^{[m+n]}=E^{[m]}\,;\,E^{[n]},\qquad
		R\,;\,S:=\{(x,z):\exists y\ ((x,y)\in R,(y,z)\in S)\}.
		\]
		For initialized paths, prefixes retain their initialization; a suffix need not start in $X_{\mathrm{init}}$.
	\end{proposition}

	\begin{proof}
		Every retained consecutive pair still belongs to $E$. Concatenating
		at the same state preserves this condition at the junction. Cutting
		a path after its $m$th edge gives one inclusion in the displayed
		identity, and concatenating witnesses gives the other. Initialization
		is a condition on the first state, which explains its prefix and
		suffix behavior.
	\end{proof}

	\begin{definition}[Finite-state, nonblocking, deterministic, and reversible dynamics]
		\label{dyn45:def-specializations}
		Let $\mathbf D=(\mathbb P,X,E)$ be a fixed-support
		relational dynamical packaged object. For each $x\in X$,
		define its successor set and the domain of $E$ by
		\[
		E(x):=\{y\in X:(x,y)\in E\},
		\qquad
		\operatorname{dom}(E):=\{x\in X:E(x)\ne\varnothing\}.
		\]
		The object $\mathbf D$ is called:
		\begin{enumerate}[label=(\roman*),leftmargin=*,nosep]
			\item \textbf{Finite-state} if $|X|<\infty$.
			
			\item \textbf{Nonblocking} if $E(x)\ne\varnothing$
			for every $x\in X$.
			
			\item \textbf{Partial deterministic} if $|E(x)|\leq 1$
			for every $x\in X$.
			Then $E=\operatorname{Graph}(F)$ for a unique map
			$F:\operatorname{dom}(E)\to X$.
			
			\item \textbf{Total deterministic} if $|E(x)|=1$
			for every $x\in X$.
			Then $E=\operatorname{Graph}(F)$ for a unique map
			$F:X\to X$, and we write $\mathbf D=(\mathbb P,X,F)$.
			
			\item \textbf{Reversible deterministic} if it is
			total deterministic and $F:X\to X$ is bijective.
		\end{enumerate}
		These conditions are additional to the primitive axioms.
		The primitive axioms also impose no requirement of
		trajectory existence at every horizon, joint observational
		faithfulness, conservation, fixed-point existence,
		or exact globalization.
	\end{definition}

	In the total deterministic case the equations $F^0=\mathrm{id}_X$ and $F^{m+n}=F^mF^n$ follow by iteration;
	they are not additional independent axioms.
	A specified time monoid can instead act by maps $\Phi_t$ satisfying $\Phi_0=\mathrm{id}_X$ and $\Phi_{s+t}=\Phi_s\Phi_t$.
	For a continuous-time specialization, topology, continuity or smoothness, and the time domain are part of that additional data.
	Commuting with one sampled time map is weaker than commuting with the whole action.
	Nonautonomous and reconfigurable laws are separately typed in Section~\ref{sec:morphisms-invariants}.
	They are not required to be iterates of one autonomous update.

	\begin{proposition}[Conservative extension and proper dynamical subclasses]
		\label{dyn45:prop-conservative-axioms}
		Every realized packaging system admits a model of $\mathsf{Ax}_{\mathrm{Dyn}}^{\mathrm{raw}}$ over that same base, even with finite, total deterministic and reversible dynamics.
		Forgetting the new data recovers the original realized system exactly.
		Thus this bare expansion is conservative for sentences in the original static signature.
		In contrast, imposing a condition on already specified dynamical data can select a proper subclass of dynamical models.
	\end{proposition}

	\begin{proof}
		Use $X=\{*\}$ and $E=\{(*,*)\}$ while keeping $\Xi,P,\mathfrak P$
		and its realization judgment unchanged. Every static model therefore
		has such an expansion. A static countermodel remains a countermodel
		after expansion, establishing conservativity. Properness at the
		dynamical level is different: the same base with a singleton state
		set and $E=\varnothing$ is a raw model but is not nonblocking.
		Neither construction changes the static base.
	\end{proof}

	\begin{remark}[Empty native data and three different relations]
		\label{dyn45:rem-native-typing}
		The static axioms permit empty carrier, candidate and selected
		families. The nonempty dynamic state sort in (D1) does not exclude
		them because it is not identified with those families. A model that
		requires a state in each selected family is a further specialization
		whose conservativity must be checked separately. Likewise,
		$\mathsf J_\Xi$ relates formation data to packaged outputs;
		$\mathsf M_\alpha$ admits typed operational assembly tuples in
		Structural Packaging; and $E$ relates temporal states. These relations
		have different domains and meanings.
		No equality among them follows from (D0)--(D2).
	\end{remark}

	\subsection{Observation axioms and structural enrichments}
	\label{dyn45:subsec-observation-axioms}

	An observation assigns compatible readings to one or more sites.
	Projecting an edge is immediate.
	Composing two projected edges requires an intermediate state that realizes both.

	\begin{definition}[Coherent observed dynamical axioms]
		\label{dyn45:def-observed-axioms}
		Let $\mathbf D=(\mathbb P,X,E)$ be a fixed-support
		relational dynamical packaged object.
		A \textbf{coherent observation structure} on $\mathbf D$
		satisfies the following axioms:
		\begin{description}[style=nextline,leftmargin=3em]
			\item[(O1) Observation sites.]
			Let $D$ be a small partially ordered set of sites.
			Each site $U\in D$ has a small observation set $X_U$.
			For every $U\leq V$, there is a restriction map
			$\rho_{V,U}:X_V\to X_U$.
			When original charts are distinguished, their indexed
			family $\mathcal H$ and site assignment
			$j:\mathcal H\to D$ are specified.
			Distinct chart occurrences retain their indices
			even when assigned to the same site.
			
			\item[(O2) Coherent restrictions.]
			The restriction maps satisfy
			\[
			\rho_{U,U}=\mathrm{id}_{X_U},
			\qquad
			\rho_{W,U}=\rho_{V,U}\circ\rho_{W,V}
			\]
			for every $U\in D$ and $U\leq V\leq W$,
			respectively.
			
			\item[(O3) Compatible observations.]
			Each site $U\in D$ has an observation map
			$r_U:X\to X_U$ satisfying
			\[
			r_U=\rho_{V,U}\circ r_V
			\qquad\text{whenever }U\leq V.
			\]
		\end{description}
		A dynamical packaged object equipped with such a structure is called a \textbf{coherent observed dynamical packaged object}.
		Its axiom system is
		\[
		\mathsf{Ax}_{\mathrm{Dyn}}^{\mathrm{obs}}
		=\mathsf{Ax}_{\mathrm{Dyn}}^{\mathrm{raw}}
		+(\mathrm{O1})+(\mathrm{O2})+(\mathrm{O3}).
		\]
		Surjectivity of the observation maps, injectivity of the
		joint observation map
		$r=(r_U)_{U\in D}:X\to\prod_{U\in D}X_U$,
		and autonomous local evolution are additional properties.
		By (O3), the restriction maps restrict to the attained
		observation sets $r_U(X)$.
		Observations at two sites can be compared on any specified common lower site; binary meets in $D$ are not required.
	\end{definition}

	\begin{proposition}[Observation shadows and conservative enrichment]
		\label{dyn45:prop-observation-shadows}
		For an observed relational model put $E_U=(r_U\times r_U)(E)$.
		If $U\le V$, then
		$
		(\rho_{V,U}\times\rho_{V,U})(E_V)=E_U.
		$
		At every horizon the analogous identity holds for the \emph{projected global path relations}.
		It need not hold after replacing those path relations by all paths of their one-step shadows.
		Every raw model admits an observed enrichment satisfying (O1)--(O3), so the unrestricted observation enrichment is conservative over the raw dynamical theory.
	\end{proposition}

	\begin{proof}
		Apply $r_U=\rho_{V,U}r_V$ to the endpoints of each edge, or to every coordinate of a global path.
		These image identities do not assert that independently liftable local edges share intermediate state witnesses.
		For the swap in Example~\ref{wn03:ex:swap}, take the
		full state site and its first-coordinate lower site.
		The projected length-two paths from the full site satisfy $a_2=a_0$, whereas the one-step shadow on the lower site allows all binary triples.
		For conservativity choose one site with $X_U=\{*\}$, its identity restriction, and the unique observation from $X$.
		Forgetting this	enrichment leaves every raw model unchanged.
	\end{proof}

	\paragraph{Coordinate and incidence specializations.}
	Definitions~\ref{def:typed-specialization-v05} and \ref{def:incidence-charts-v06} specialize the site language to coordinate supports.
	In this case common lower sites may be actual intersections and restrictions are coordinate projections.
	The	chosen observable atoms are retained structure; a general bijection of instantaneous state sets does not preserve their erasure tests.
	The terminal empty-support observation $X_\varnothing=\{*\}$ is a coordinate convention, not a requirement on arbitrary native sites.
	These observation sets are not automatically the static candidate sets $\mathsf L_{\mathfrak P}(K)$.

	\paragraph{Structural dynamical models.}
	One may retain a raw, coherent, or hierarchical Level-II structural presentation of Structural Packaging together with the dynamical data above.
	The retained structural axioms are declared separately:
	(A0)--(A3) specify assembly data, (H1)--(H2) require coherent restrictions and interface soundness, and (H3) adds a finite-valued formation rank.
	To identify its operational states or observation sites with dynamical ones requires a specified interpretation, not merely coincident notation. For example, when deterministic updates $F_K$ are declared on the operational states of a signature $\alpha:(K_1,\ldots,K_m)\to K$, temporal preservation of admission is the additional condition
	\[
	(F_{K_1}\times\cdots\times F_{K_m}\times F_K)
	(\mathsf M_\alpha)\subseteq\mathsf M_\alpha.
	\]
	Iteration then preserves admitted tuples by induction.
	This condition is not implied by interface soundness.
	H3 concerns the	shape of formation expressions and does not prohibit temporal cycles.
	Structural and dynamical data can be combined, but neither extension forces all of the additional data of the other.

	More explicitly, a coherent structural--dynamical specialization is specified by the schema
	\[
	\mathsf{Ax}_{\mathrm{StrDyn}}^{\mathrm{coh}}
	=\mathsf{Ax}_{\mathrm{II}}^{\mathrm{coh}}
	+\mathsf{Ax}_{\mathrm{Dyn}}^{\mathrm{raw}}
	+\Sigma_{\mathrm{interp}}.
	\]
	Here $\mathsf{Ax}_{\mathrm{II}}^{\mathrm{coh}}$ denotes the retained
	Level-II coherent assembly axioms, and $\Sigma_{\mathrm{interp}}$
	declares the sorts, maps, and equations connecting the two
	presentations, identifying their common inherited base symbols.
	It does not automatically identify the dynamical state sort with
	the native realization or operational-state sorts. To require
	evolution to preserve assembly, one further adds an explicitly
	chosen $\Sigma_{\mathrm{evol}}$, such as the admission-preservation
	condition above. Adding (H3) yields a hierarchical structural
	specialization independently of those temporal conditions. This
	is a schema for declared enrichments; no undeclared interpretation
	or universal assembly category is asserted.

	\paragraph{Atlas-first data and derived existence questions.}
	Independently specified local states and transitions form a raw	atlas-first presentation as in Definition~\ref{def:two-modes-v06}.
	No global state sort is postulated in that presentation.
	Its compatible-state envelope may be empty; such emptiness records a failed realization problem and is not a nonempty-state model of (D1).
	State compatibility, one-step compatibility, all-horizon compatibility, and exactness remain the separate conditions of Definition~\ref{def:v07-exact}.
	Exactness, nonblocking and fixed-point existence are conclusions to test or hypotheses of named subclasses, not axioms imposed on all raw atlases.

	\subsection{Conservativity and non-collapse under forgetting coupling}
	\label{dyn46:subsec-noncollapse}

	The axioms (D0)--(D2) describe the raw relational expansion. Forgetful comparisons distinguish the coupled theory from the independent packaging and temporal components.

	\begin{proposition}[Product decomposition of the raw theory]
		\label{dyn46:prop-product}
		Fix a category $\mathbf B$ of realized bases and let
		$\mathbf{RelDyn}$ have nonempty small state sets with one-step relations as objects and edge-preserving functions as arrows.
		Equip the raw expansion with componentwise arrows $(\alpha,f)$, with no extra preservation condition coupling the components.
		Then
		\[
		\mathbf{DynPkg}_{\rm raw}\cong\mathbf B\times\mathbf{RelDyn}.
		\]
		The analogous assertion holds for any fixed temporal-law category with independent arrows.
		In particular a property of the temporal component alone need not be a packaging-dependent property.
	\end{proposition}

	\begin{proof}
		The object correspondence sends $(\mathbb P,X,E)$ to
		$(\mathbb P,(X,E))$.
		The arrow sets, identities and compositions	are exactly the componentwise ones.
		This gives inverse functors.
		The hypothesis on arrows is essential; the statement is not a product decomposition of the observed or interaction-marked category.
	\end{proof}

	\begin{proposition}[Temporal information is not determined by a static base]
		\label{dyn46:prop-time-independence}
		For every realized base there are total deterministic reversible expansions with the same two-element state set that are not conjugate, even after forgetting the packaging.
		Consequently the static base alone does not determine the dynamical isomorphism class in the raw theory.
	\end{proposition}

	\begin{proof}
		On $\{0,1\}$ take the identity and the transposition. The former has two fixed states and the latter none. Conjugacy preserves fixed states.
		Both expansions satisfy (D0)--(D2) over the unchanged base.
	\end{proof}

	\begin{example}[Coupling is not determined by the two separate isomorphism classes]
		\label{dyn46:ex-coupled-noncollapse}
		Take a realized base with three internal units $a,b,c$ and three observed singleton carriers (with the indicated incidence and the same native data in both models).
		Such an output satisfies (S1)--(S5), for example using singleton candidate and selected sets and a domain
		containing the three carriers.
		Choose a context admitting those carriers, a singleton formation-data set and its realization pair with this output;
		(P1)--(P3) and (Y1)--(Y4) allow that context.
		On the full binary configuration space $X=\mathbb F_2^3$, consider
		\[
		F(a,b,c)=(a+1,b,c),\qquad G(a,b,c)=(a+1,b+c,c).
		\]
		They are conjugate after forgetting the observations: the involution $h(a,b,c)=(a,b+ac,c)$ satisfies $hF=Gh$.
		But the $b$-observation is an autonomous factor of $F$ and not of $G$.
		In fact every singleton observation of $F$ is autonomous, whereas only the $a$- and $c$-observations of $G$ are autonomous.
		For $G$ the states $(0,0,0)$ and $(0,0,1)$ have the same $b$ value but different next $b$ values.
		Any observed isomorphism with inverse transports the
		existence of a local autonomous factor, by conjugating that factor with the corresponding observation bijection.
		Thus no isomorphism	preserving the singleton chart family, even allowing chart permutations, exists between these two models.
		The bare base and unobserved time-law isomorphism classes coincide but the coupled classes differ.
	\end{example}

	These examples prove non-recoverability relative to the indicated forgetful comparisons.
	The coupled theory studies the interaction of native types, observations, interfaces and temporal operators.

	\subsection{Typed temporal laws and their reducts}
	\label{dyn46:subsec-time-types}

	A path relation records which steps are possible.
	A stochastic update also records their weights, and an evolution with an input records which steps an input permits.
	Each comparison below specifies which of these data it retains.

	\begin{definition}[Types of temporal laws]
		\label{dyn46:def-time-types}
		A \textbf{typed temporal law} specifies a state space,
		an evolution law, and the auxiliary structures required
		by that law.
		The principal types are:
		\begin{enumerate}[label=(\roman*),leftmargin=*,nosep]
			\item \textbf{Discrete relational law.}
			A relation $E\subseteq X^2$ specifies the permitted
			one-step transitions. Partial and total deterministic
			laws are the specializations in
			Definition~\ref{dyn45:def-specializations}.
			
			\item \textbf{Deterministic time action.}
			A time monoid $T$, written additively with identity $0$,
			acts on $X$ through maps $\Phi_t:X\to X$ satisfying
			\[
			\Phi_0=\mathrm{id}_X,\qquad
			\Phi_{s+t}=\Phi_s\circ\Phi_t
			\quad(s,t\in T).
			\]
			Continuity, differentiability, and completeness
			are imposed when required by the chosen class.
			
			\item \textbf{Discrete Markov law.}
			On a measurable state space $(X,\mathcal F)$,
			a Markov kernel $P(x,dy)$ specifies the conditional
			distribution of the next state given the current state.
			An initial probability measure $\mu_0$ determines
			the corresponding finite-history laws.
			
			\item \textbf{Markov semigroup.}
			On a measurable state space $(X,\mathcal F)$,
			a family of Markov kernels $(P_t)_{t\in T}$
			satisfies
			\[
			P_0(x,\cdot)=\delta_x,\qquad
			P_{s+t}(x,B)=\int_X P_t(y,B)\,P_s(x,dy)
			\]
			for $s,t\in T$ and $B\in\mathcal F$.
			For finite-state continuous-time dynamics, a
			conservative rate matrix $Q$ satisfies
			\[
			Q(x,y)\geq0\quad(x\ne y),\qquad
			Q(x,x)=-\sum_{y\ne x}Q(x,y),
			\]
			and determines the semigroup $P_t=\exp(tQ)$.
			
			\item \textbf{Controlled or scheduled law.}
			A family of evolution laws is indexed by controls
			or schedules, together with conditions specifying
			the admissible choices.
			Nonautonomous and reconfigurable laws additionally
			specify their time dependence, changing types
			when applicable, and the corresponding transports.
		\end{enumerate}
		These types specify alternative dynamical structures.
		An energy function, potential, or stationary measure alone does not specify a temporal law.
		An equilibrium model without an evolution law remains a static enrichment.
	\end{definition}

	\noindent For each type, the underlying static base is unchanged.
	The singleton expansion still demonstrates conservativity for the unrestricted types allowing a trivial evolution; imposing nontrivial rates, faithful observations, or prescribed state interpretations requires a separate expansion argument.

	\begin{example}[A law's support does not retain its weights]
		\label{dyn46:ex-stochastic-reduct}
		On $\{0,1\}$ let the rates $0\to1$ and $1\to0$ be $\alpha,\beta>0$.
		The off-diagonal support is the same for all such rates, whereas the stationary law is
		\[
		\pi=\left(\frac{\beta}{\alpha+\beta},
		\frac{\alpha}{\alpha+\beta}\right).
		\]
		This follows by solving $\pi Q=0$ with $\pi\mathbf 1=1$. For instance $(\alpha,\beta)=(1,1)$ and $(3,1)$ have different stationary-law weight multisets.
		The off-diagonal support reduct cannot determine
		this invariant.
		It is not the time-$t$ transition support of the
		semigroup, which for positive $t$ also contains both diagonal pairs.
		A relational path theorem for either reduct says nothing by itself about weighted path laws, stochastic stationarity or scaling limits.
	\end{example}

	\subsection{Structural--temporal coupling data and axiomatic enrichments}
	\label{dyn46:subsec-coupling}

	We now specify how the temporal law acts on the package.
	A read set identifies the inputs of an update, a write set identifies the coordinates it can change, and admission determines which assemblies survive the update.

	\begin{definition}[Coupled dynamical packaged object and coupling axioms]
		\label{dyn46:def-coupling}
		A \textbf{coupled dynamical packaged object} with fixed support	is a quadruple
		$
		\mathbf D=(\mathbb P,X,\mathcal T,\mathcal C),
		$
		where $\mathbb P$ is a realized packaging system satisfying (D0), $X$ is a state space satisfying (D1),
		$\mathcal T$ is a temporal law of a type specified in
		Definition~\ref{dyn46:def-time-types}, and $\mathcal C$
		is a coupling presentation specifying which of the following clauses are imposed, together with their data and conditions.		
		\begin{description}[style=nextline,leftmargin=3em]
			\item[(C1) State interpretation.]
			Specify the native sorts or coordinates represented by the states in $X$ and the maps or relations connecting them.
			A coordinate realization consists of a descriptor set $V$, alphabets $A_v$, and a state space
			\[
			X\subseteq\prod_{v\in V}A_v.
			\]
			Any identification of coordinates or states with native units or operational candidates must be specified explicitly.
			
			\item[(C2) Observation provenance.]
			Specify a coherent observation structure satisfying
			(O1)--(O3), together with the association of its chart occurrences to carriers, interfaces, or auxiliary sites.
			Incidence-generated observations are given by
			Definition~\ref{def:incidence-charts-v06}.
			Observation autonomy is a separate property.
			
			\item[(C3) Interaction provenance.]
			Specify an indexed interaction family and its supports as in Definition~\ref{dyn46:def-interaction}.
			Distinguish carrier-supported interactions from auxiliary interactions.
			Interaction supports need not be observation
			charts.
			
			\item[(C4) Temporal interaction.]
			Specify the channel laws, their read and write supports when applicable, and the rule composing them into $\mathcal T$.
			The resulting temporal law must satisfy the declared
			composition and compatibility conditions.
			
			\item[(C5) Structural compatibility.]
			For a retained Level-II presentation, specify the
			interpretation of its realized bases, operational types, signatures, admission relations, and interface observations.
			Declare the required compatibility and preservation conditions.
			Admission preservation, autonomous evolution, and exact lifting are tested separately.
			
			\item[(C6) Reconfiguration.]
			For changing base types, replace $(\mathbb P,X)$ by
			an indexed family $(\mathbb P_i,X_i)$ of realized bases and state spaces.
			Specify the allowed temporal laws between state spaces, persistent-unit identifications, and structural transports.
		\end{description}		
		Each declared clause specifies its domains, markings, and
		required equations.
		The clauses are optional enrichments, not universal existence axioms.
		A model without coupling enrichments is uncoupled.
		The observed relational specialization consists of
		(D0)--(D2), (O1)--(O3), and its declared provenance.		
		The shared dynamical axiom schema is
		\[
		\boxed{\mathsf{Ax}_{\rm Dyn}
		=\mathsf{Ax}_{\rm Pkg}
		+\Sigma_{\rm Time}
		+\Sigma_{\rm Coupling}},
		\]
		where $\Sigma_{\rm Time}$ specifies the state space and
		temporal-law type, and $\Sigma_{\rm Coupling}$ specifies
		the declared coupling clauses.
		In structural specializations, $\mathsf{Ax}_{\rm II}$
		is inherited on the retained realized bases.
	\end{definition}

	\begin{definition}[Interaction assignments and their supports]
		\label{dyn46:def-interaction}
		Fix a coordinate realization on a descriptor set $V$.
		An \textbf{interaction assignment} consists of a small
		index set $\mathcal J$, a support map
		$
		S:\mathcal J\longrightarrow 2^V,
		$
		and its declared association with native package incidence.
		For each interaction occurrence $J\in\mathcal J$,
		$S(J)$ is its \textbf{interaction support}.
		\begin{enumerate}[label=(\roman*),leftmargin=*,nosep]
			\item \textbf{Carrier-supported assignment.}
			An assignment is carrier-supported if it specifies
			a carrier map
			$\iota:\mathcal J\to\operatorname{Pkg}_{\mathfrak P}$
			satisfying
			\[
			S(J)=I_{\iota(J)}^V
			\qquad(J\in\mathcal J),
			\]
			where $I_K^V\subseteq V$ is the incidence support of
			carrier $K$ under the declared native-unit interpretation.
			In the native-unit specialization $V=E_{\mathfrak P}$, $I_K^V=I_K$.
			The carrier map may be restricted to selected observed carriers.
			
			\item \textbf{Auxiliary assignment.}
			An assignment is auxiliary if its supports are
			specified without requiring a carrier map.
		\end{enumerate}
		The support map $S$ need not be injective, and its supports need not coincide with observation charts.
		Interaction indices are retained even when supports
		are empty or identical.
		This assignment provides the microscopic interaction data
		for the collective families in
		Definition~\ref{dyn46:def-collective-schema}.
	\end{definition}

	\begin{remark}[Four different structures]
		\label{dyn46:rem-four-structures}
		Incidence specifies membership, observations specify accessible
		information, interaction assignments specify participating
		coordinates, and temporal laws specify the updates or their rates.
		For an interaction $J$ write $R_J,W_J\subseteq S(J)$ for declared
		read and write sets when a law has these roles. They are not inferred
		from an observation chart. The same carrier can support several
		channels, and the same aggregate law can have several decompositions.
		The interaction assignment is interpreted together with its declared law and composition.
	\end{remark}

	\subsection{Native interaction laws: support and compatibility tests}
	\label{dyn46:subsec-native-interaction}

	The following tests concern packaging-dependent
	coupling.
	They are finite-coordinate or relational specializations of
	(C1)--(C4), not requirements on all observed models.

	\begin{definition}[Local relational channels, lifts, and composition]
		\label{dyn46:def-channel-lift}
		Fix a state space $X\subseteq\prod_{v\in V}A_v$
		and an indexed interaction family $\mathcal J$.
		For each $J\in\mathcal J$, let $S_J=S(J)$ be its
		interaction support and $W_J\subseteq S_J$ its
		write set, and put $X_{S_J}:=\pi_{S_J}X$.
		\begin{enumerate}[label=(\roman*),leftmargin=*,nosep]
			\item \textbf{Local relational channel.}
			A local relational channel on $S_J$ is a relation
			\[
			L_J\subseteq X_{S_J}\times X_{S_J},
			\]
			specifying the permitted local transitions.
			
			\item \textbf{Global lift.}
			The global lift of $L_J$ to $X$ is the relation
			\[
			\operatorname{Lift}_{X,J}(L_J):=
			\left\{(x,y)\in X^2:
			\begin{array}{l}
				(\pi_{S_J}x,\pi_{S_J}y)\in L_J,\\
				\pi_{V\setminus W_J}x=\pi_{V\setminus W_J}y
			\end{array}
			\right\}.
			\]
			Thus both endpoints are admissible, and only
			coordinates in $W_J$ may change.
			
			\item \textbf{Channel composition.}
			The asynchronous single-channel law is
			\[
			E=\bigcup_{J\in\mathcal J}
			\operatorname{Lift}_{X,J}(L_J).
			\]
			For an ordered schedule $(J_1,\ldots,J_n)$,
			the corresponding transition relation is
			\[
			E_{J_1,\ldots,J_n}
			=
			\operatorname{Lift}_{X,J_n}(L_{J_n})
			\circ\cdots\circ
			\operatorname{Lift}_{X,J_1}(L_{J_1}).
			\]
			Synchronous admission uses a separate construction.
		\end{enumerate}
		The composition rule is part of the declared temporal law.
	\end{definition}

	\begin{proposition}[Exact channel-locality test and canonical support rule]
		\label{dyn46:prop-channel-saturation}
		Let $E_J\subseteq X^2$ be a marked channel. Define
		$B_J=\{(x,y)\in X^2:x|_{V\setminus W_J}=y|_{V\setminus W_J}\}$
		and $p_J:B_J\to X_{S_J}^2$ by the two support projections.
		\begin{enumerate}[label=(\roman*),leftmargin=*,nosep]
			\item \emph{Existence.} There exists a local relation whose lift equals $E_J$ if and only if
			$E_J\subseteq B_J$ and
			\[
			u,v\in B_J,\quad p_J(u)=p_J(v),\quad u\in E_J
			\ \Longrightarrow\ v\in E_J.
			\]
			
			\item \emph{Uniqueness on attained pairs.} When this holds, $p_J(E_J)$ is the unique rule on attained pairs
			$p_J(B_J)$, and its lift is $E_J$. On unattained pairs a rule is
			not identifiable.
			
			\item \emph{Least local enlargement.} For arbitrary $E_J\subseteq B_J$ the lift of
			$p_J(E_J)$ is the least support-local superset of $E_J$ in $B_J$.
		\end{enumerate}
	\end{proposition}

	\begin{proof}
		\leavevmode
		\begin{enumerate}[label={(iii)},leftmargin=*,nosep]
			\item[(i)] Every lift is a preimage under $p_J$ and hence is saturated in
			$B_J$. Conversely saturation gives $E_J=p_J^{-1}p_J(E_J)$.
			
			\item[(ii)] If a rule has that inverse image, membership on attained pairs
			is forced.
			
			\item[(iii)] The preimage of the image is contained in every saturated
			superset, proving minimality. \qedhere
		\end{enumerate}
	\end{proof}

	\begin{proposition}[Exact native support test for a deterministic coordinate]
		\label{dyn46:prop-native-read}
		Let $V$ be finite, $X=\prod_{v\in V}A_v$ with all $A_v$ nonempty,
		and $F:X\to X$. For $w\in V$ define the essential read set
		\[
		D_w(F)=\{v:\exists x,x'\in X,\ x|_{V\setminus\{v\}}=x'|_{V\setminus\{v\}},
		\ F(x)_w\ne F(x')_w\}.
		\]
		\begin{enumerate}[label=(\roman*),leftmargin=*,nosep]
			\item \emph{Read set.} Then $F(x)_w$ factors through $\pi_S$ if and only if $D_w(F)\subseteq S$.
			
			\item \emph{Single-carrier update.} In particular, if a local mechanism updating $w$ is required to use
			one carrier containing $w$, an admissible carrier exists exactly when
			$D_w(F)\cup\{w\}\subseteq I_K$ for some allowed carrier $K$.
			
			\item \emph{Synchronous realization.} If
			such a condition is imposed independently for every $w$ and simultaneous
			coordinate updates are permitted, these tests are necessary and sufficient
			for a carrier-supported synchronous realization of $F$.
		\end{enumerate}
	\end{proposition}

	\begin{proof}
		\leavevmode
		\begin{enumerate}[label={(iii)},leftmargin=*,nosep]
			\item[(i)] If $F_w$ factors through $\pi_S$, no change outside $S$ can affect
			it. Conversely, for two states agreeing on $S$, change the finitely
			many outside coordinates one at a time. Every intermediate state
			belongs to the full product, and no changed coordinate belongs to
			$D_w(F)$, so the output coordinate remains unchanged. This proves
			the fiber condition and hence the factorization.
			
			\item[(ii)] By (i), a carrier $K$ containing $w$ supports the required update exactly when $D_w(F)\subseteq I_K$.
			
			\item[(iii)] Choosing one allowed carrier for each output gives local rules; taking all outputs simultaneously recovers $F$. Conversely each rule must satisfy the carrier condition in (ii). \qedhere
		\end{enumerate}
	\end{proof}

	\begin{remark}[Constrained state spaces and a relative coupling invariant]
		\label{dyn46:rem-read-invariant}
		For general $X\subsetneq\prod_v A_v$, use the fiber condition
		$x|_S=x'|_S\Rightarrow F(x)_w=F(x')_w$ on $X$ instead; the
		single-coordinate proof can fail because intermediate configurations
		need not be admissible. In the full-product model the family
		$(D_w(F))_w$ transports under native unit permutations and coordinatewise
		alphabet bijections. In particular, when the allowed carrier library is transported, the
		minimum allowed carrier size containing $D_w(F)\cup\{w\}$ (infinity
		when none exists) is a comparison-relative coupling invariant, not
		an invariant of $F$ up to arbitrary conjugacy. In Example~\ref{dyn46:ex-coupled-noncollapse},
		all outputs of $F$ pass the singleton-carrier test, but the $b$ output
		of $G$ requires both $b$ and $c$. Enlarging the static carrier library
		changes the test without changing $G$.
	\end{remark}
	The minimum native carrier size is a separate realization problem for the declared support library.

	\begin{proposition}[Local rate representation and aggregation of finite jump generators]
		\label{dyn46:prop-jump-support}
		Let $X$ and $\mathcal J$ be finite.
		For each $J\in\mathcal J$, prescribe a nonnegative rate function $q_J:p_J(B_J)\to[0,\infty)$ and define, for $x\ne y$,
		\[
		Q_J(x,y):=
		\begin{cases}
			q_J(p_J(x,y)), &(x,y)\in B_J,\\
			0, &\text{otherwise}.
		\end{cases}
		\]
		Set $Q_J(x,x):=-\sum_{y\ne x}Q_J(x,y)$.
		\begin{enumerate}[label=(\roman*),leftmargin=*,nosep]
			\item \textbf{Generator aggregation.}
			Each $Q_J$ is a conservative generator, and so is
			\[
			Q:=\sum_{J\in\mathcal J}Q_J.
			\]
			The off-diagonal support of $Q$ is the union of
			the off-diagonal supports of the $Q_J$.
			
			\item \textbf{Local rate representation.}
			A conservative generator on $X$ admits the above
			representation for channel $J$ if and only if
			its off-diagonal rates vanish outside $B_J$ and
			are constant on the fibers of $p_J$ restricted
			to off-diagonal pairs in $B_J$.
		\end{enumerate}
	\end{proposition}

	\begin{proof}
		\leavevmode
		\begin{enumerate}[label=(\roman*),leftmargin=*,nosep]
			\item Each $Q_J$ has nonnegative off-diagonal entries
			and zero row sums, and hence so does $Q$.
			Since the rates are nonnegative, $Q(x,y)>0$ for
			$x\ne y$ if and only if $Q_J(x,y)>0$ for some $J$.
			
			\item Necessity follows directly from the defining
			formula.
			Conversely, constancy on the specified fibers defines $q_J$ on support pairs attained
			off-diagonally.
			Set $q_J=0$ on the remaining support pairs.
			The diagonal entries are uniquely determined by the zero row-sum condition.
		\end{enumerate}
	\end{proof}

	\noindent This is a support-local specialization allowing reads on $S_J$ and writes on $W_J$.
	More general declared read/write laws may use separate rate arguments and require their own fiber test.
	The support condition does not imply autonomous observation; admissibility of a global target may itself depend on unobserved coordinates.
	The	aggregate generator does not remember how rates are distributed among repeated occurrences:
	two channels with the same write and support sets can split a positive rate in different proportions while leaving
	$Q$ unchanged.
	A comparison retaining marked channel rates is thus strictly finer in general than aggregate-generator comparison.
	The aggregate-law realization problem under a prescribed carrier library is a separate comparison contract.

	\subsection{Temporal realization and structural admission}
	\label{dyn46:subsec-structural-coupling}

	The following specializations use the native distinction between the realization relation $\mathsf J_\Xi$, the interface matching relation, and operational admission $\mathsf M_\alpha$.
	The latter is not a transition relation.
	A Level-II type may refer to an occurrence in an indexed family of realized bases or to a prospective operational
	type.
	Empty admission relations are permitted by (A3).

	\paragraph{Functional realization descent.}
	For a declared surjection $q:S\twoheadrightarrow Y$ from formation states to realized output states and a total formation update $F_S$, Theorem~\ref{thm:exact-deterministic-chart-factor-criterion} gives the unique output law $F_Y$ with $qF_S=F_Yq$ precisely when $q(s)=q(s')$ implies $q(F_Ss)=q(F_Ss')$.
	This specialization does not replace a general relational realization judgment by a function.

	\noindent For a relational realization $J\subseteq S\times Y$
	and declared time relations on $S,Y$, one can instead test the condition that every compatible pair and every permitted source step admit a matching output step to another compatible pair.
	This is an existence-of-lift condition on $J$, not the functional uniqueness statement above.
	Both directions, totality and exact trajectory matching are independent further tests.
	The exact relational square, its greatest solution, composition, and rigidity are treated in
	Section~\ref{dyn47:subsec:relational-descent}.
	This is a distinct comparison from representative-wise forward lifting.

	\begin{definition}[Temporal assembly and synchronized laws]
		\label{dyn46:def-admission-dynamics}
		Fix a structural signature $\alpha:(K_1,\ldots,K_m)\to K$.
		Put $Y_i=\mathsf R(K_i)$ for $1\leq i\leq m$ and
		$Y_{m+1}=\mathsf R(K)$, and let
		\[
		A:=\mathsf M_\alpha\subseteq\prod_{i=1}^{m+1}Y_i
		\]
		be the native admission set.
		For local transition relations $E_i\subseteq Y_i^2$,
		define:
		\begin{enumerate}[label=(\roman*),leftmargin=*,nosep]
			\item \textbf{Maximal synchronized law.}
			The maximal synchronized law on $A$ is
			\[
			E_\alpha^{\max}:=
			\{(u,v)\in A^2:(u_i,v_i)\in E_i
			\text{ for all }i\}.
			\]
			
			\item \textbf{Admissible synchronized law.}
			A relation $E_\alpha\subseteq A^2$ is an admissible
			synchronized law if
			$E_\alpha\subseteq E_\alpha^{\max}$.
			
			\item \textbf{Temporal assembly.}
			A temporal assembly is the pair $(A,E_\alpha)$,
			where $E_\alpha$ is an admissible synchronized law.
		\end{enumerate}
		This specialization uses operational states as temporal
		states.
		Other state spaces require an explicit interpretation.
		The construction is defined even when $A=\varnothing$,
		but yields a dynamical model satisfying (D1) only when
		$A\ne\varnothing$.
	\end{definition}

	\begin{proposition}[Maximal synchronized evolution, fixed states, and viable domains]
		\label{dyn46:prop-admission}
		\leavevmode
		\begin{enumerate}[label=(\roman*),leftmargin=*,nosep]
			\item \emph{Maximality and paths.}
			The relation $E_\alpha^{\max}$ contains every temporal relation on $A$ preserving all the coordinate transitions.
			Its length-$n$ paths are exactly the synchronized tuples of length-$n$ local paths whose assemblies belong to $A$ at every time.
			
			\item \emph{Fixed states.}
			Moreover,
			\[
			\operatorname{Fix}(E_\alpha^{\max})
			=A\cap\prod_i\operatorname{Fix}(E_i).
			\]
			
			\item \emph{Deterministic assembly.}
			For total local maps $F_i$ write $G=\prod_iF_i$.
			The synchronized law is a total deterministic law on $A\ne\varnothing$ if and only if $G(A)\subseteq A$.
			In the absence of this condition its largest
			forward-viable domain is
			$A_\infty=\bigcap_{n\ge0}G^{-n}(A)$.
			If the ambient product is finite the decreasing iteration $A_0=A$, $A_{n+1}=A\cap G^{-1}(A_n)$ stabilizes after at most $|A|$ strict decreases.
		\end{enumerate}
	\end{proposition}

	\begin{proof}
		\leavevmode
		\begin{enumerate}[label={(iii)},leftmargin=*,nosep]
			\item[(i)] Maximality is the defining conjunction. Applying it at each step proves the path identity.
			
			\item[(ii)] Applying it to $(u,u)$ proves the stationary identity.
			
			\item[(iii)] In the deterministic case the only possible successor is $G(u)$, which exists in the admitted set exactly when $G(u)\in A$.
			A state belongs to $A_\infty$ precisely when its entire forward orbit lies in $A$.
			Thus $A_\infty$ is invariant and contains every invariant subset of $A$.
			Induction identifies $A_n$ with the finite intersection through time $n$.
			A strictly decreasing finite chain has at most
			$|A|$ deletions before stabilization. \qedhere
		\end{enumerate}
	\end{proof}

	\begin{example}[Admission preservation is a coupled property]
		\label{dyn46:ex-admission-noncollapse}
		Take two occurrence-labeled operational sorts equal to $\{0,1,2\}$	with admission $A=\{(i,i):i=0,1,2\}$. Choose two carrier occurrences with candidate and selected sets $\{0,1,2\}$ and native singleton
		formation witnesses as in Example~\ref{dyn46:ex-coupled-noncollapse}.
		Use a unary signature with a common lower observation site of that same set and identity restriction maps.
		Equality is the interface matching relation and $A$ is contained in it.
		Thus (A0)--(A3), (H1)--(H2) allow this choice.
		Compare local permutations $((01),(01))$ and $((01),(12))$.
		Their separate local conjugacy classes coincide, and the admission relation is identical.
		In the first case the product update preserves $A$; in the second case $(0,0)\mapsto(1,0)$ leaves $A$.
		Thus separate local conjugacy classes and static admission data do not determine total synchronized
		evolution on the admitted set.
		Any comparison transporting both admission and the two local laws preserves this property.
	\end{example}

	\begin{proposition}[Admission preservation for finite jump laws]
		\label{dyn46:prop-stochastic-admission}
		Let $Z$ be a finite nonempty ambient operational state set, $\varnothing\ne A\subseteq Z$ its admitted states, and $Q=\sum_JQ_J$ a finite sum of conservative generators.
		The original Markov evolution started at any state of $A$ remains in $A$ almost surely for all times
		if and only if
		\[
		Q(x,Z\setminus A):=\sum_{y\notin A}Q(x,y)=0\qquad(x\in A).
		\]
		This holds exactly when every marked summand has zero such exit rate from every $x\in A$.
		In that case $Q|_{A\times A}$ is a conservative
		generator whose semigroup is the restriction of the original one.
	\end{proposition}

	\begin{proof}
		If a rate from $x\in A$ to $y\notin A$ is positive, then
		$(e^{tQ})(x,y)=tQ(x,y)+O(t^2)>0$ for small positive $t$.
		Conversely if every such rate vanishes, ordering states by $A$ and its complement gives a block triangular matrix with zero upper-right block, a property preserved by powers and the exponential.
		Hence all one-time probabilities of exit vanish; finite-state jump paths are right-continuous with finitely many jumps on bounded intervals, so checking rational times excludes any exit almost surely. Nonnegative summand rates prove the channel criterion, and zero exit gives the claimed restricted diagonal row sums and semigroup.
	\end{proof}

	\begin{example}[Synchronization is not independent clock evolution]
		\label{dyn46:ex-clock-noncollapse}
		On $\{0,1\}^2$ with diagonal admission, the simultaneous flip $(a,b)\mapsto(1-a,1-b)$ preserves admission.
		A continuous-time channel flipping both bits at rate one also preserves it.
		The sum of two rate-one generators flipping the bits independently does not:
		from $(0,0)$ there are positive rates to $(1,0)$ and $(0,1)$.
		Each coordinate separately has the same rate-one flip law in the simultaneous and independent continuous-time models.
		Hence even identical marginal dynamics do not determine admission preservation.
		The joint interaction rule is essential. Suppressing exiting rates and recomputing diagonals would define a different kinetic law.
	\end{example}

	\subsection{Interface to collective dynamics}
	\label{dyn46:subsec-collective-interface}

	At a fixed size the collective construction supplies a coupled system of the type specified above.
	Passage to a macroscopic observation requires a family
	of such systems with the following additional data.

	\begin{definition}[Collective dynamical families]
		\label{dyn46:def-collective-schema}
		Fix a parameter space $\Theta$ and specified temporal-law
		and coupling types.
		A \textbf{collective dynamical family} consists of:
		\begin{enumerate}[label=(\roman*),leftmargin=*,nosep]
			\item \textbf{Microscopic dynamics.}
			A size- and parameter-indexed family
			\[
			(\mathbf D_{N,\lambda})_{N\geq1,\;\lambda\in\Theta},
			\]
			with specified interaction assignments, normalized
			kinetic parameters, and initial states or laws.
			
			\item \textbf{Macroscopic observations.}
			For each size $N$, a macroscopic observation map
			\[
			M_N:X_N\longrightarrow Y_N,
			\]
			where $X_N$ is the microscopic state space and $Y_N$
			is the macroscopic observation space.
			Parameter dependence is suppressed in this notation.
			
			\item \textbf{Collective comparison.}
			Specified markings, parameter comparisons, time
			conventions, and stationary or criticality tests.
			When a limit is studied, a common comparison space
			and a mode of convergence are also specified.
		\end{enumerate}
		Limiting equations, autonomous macroscopic laws, equilibrium multiplicity, large-deviation principles, and phase transitions are properties to be established, not defining axioms.
		A family specified only by packaged outputs is an
		output-level collective specialization.
		A realized-base family additionally requires native contexts, formations, and realization witnesses.
	\end{definition}

	The interaction family may be extensive, local or all-to-all.
	The normalization is not inferred from incidence alone.
	Exact factors of one temporal system, approximate reductions with quantitative error, and asymptotic autonomy of a family are distinct.
	The finite-generator comparison in Section~\ref{dyn48:subsec:continuous} concerns two declared finite systems.
	Collective reduction, stationary transfer and model-specific phase selection belong to the collective extension.

	\subsection{Deterministic evolution, autonomous factors, observations, and changes of type}
	\label{subsec:deterministic-evolution}

	\paragraph{Autonomous factors.}	
	An observation has its own deterministic evolution precisely when all states with the same present reading give the same next reading.
	Surjectivity onto the attained observations makes that evolution unique.

	\begin{dtheorem}[Exact deterministic chart-factor criterion]
		\label{thm:exact-deterministic-chart-factor-criterion}
		Let \(r:X\twoheadrightarrow Y\) be surjective and \(F:X\to X\) deterministic.
		The following are equivalent:
		\begin{enumerate}[label=(\alph*),leftmargin=2.2em]
			\item There exists \(f:Y\to Y\) with \(rF=fr\).
			
			\item The kernel equivalence of \(r\) is forward invariant:
			\[
			\boxed{r(x)=r(y)\Longrightarrow r(Fx)=r(Fy).}
			\]
			
			\item The image of the update graph, or transition shadow, \(R_r(F)=\{(r(x),r(Fx)):x\in X\}\) is the graph of a total single-valued map on \(Y\).
		\end{enumerate}
		If these hold, \(f\) is unique, and \(rF^t=f^tr\) for all \(t\in\mathbb N_0\).
	\end{dtheorem}

	\begin{proof}
		\leavevmode
		\begin{enumerate}[label={},leftmargin=0pt,nosep]
			\item \textnormal{(a)$\Rightarrow$(b).}\enspace If \(rF=fr\), two states in one \(r\)-fiber have the same next observation.
			
			\item \textnormal{(b)$\Rightarrow$(a).}\enspace Conversely, if \(r\)-fibers have uniquely determined next observations, define \(f(r(x)):=r(Fx)\).
			The invariance condition makes this independent of the chosen preimage, and surjectivity defines it on every \(y\in Y\).
			
			\item \textnormal{(a)$\Leftrightarrow$(c).}\enspace If (a) holds, the shadow is the graph of $f$ by surjectivity.
			Conversely, if the shadow is a total single-valued graph, its map satisfies $rF=fr$.
		\end{enumerate}
		Surjectivity also gives uniqueness. Iteration proves the last identity.
	\end{proof}

	\begin{remark}
		Autonomy makes the entire observed history a function of the initial	reading.
		The precise consequence for temporal observation depth is
		Theorem~\ref{wn03:thm:no-go};
		Example~\ref{wn03:ex:swap} shows why a nonautonomous observation can recover additional information.
	\end{remark}

	For deterministic dynamics we can follow a single initial state through time.
	This makes it possible to distinguish an autonomous observation from an observation history that eventually identifies the initial state.

	\subsubsection{Static provenance and the conservative dynamical layer}
	\label{subsec:provenance-conservative-dynamical-layer}

	Retain the realized base of (D0) in Definition~\ref{dyn45:def-core-axioms}.
	We now specify the deterministic specialization of the temporal-law schema (Definition~\ref{dyn46:def-time-types}), without adding a new static axiom.
	Descriptor charts need not be primitive package interiors.

	\begin{ddefinition}[Autonomous and nonautonomous deterministic fixed-support systems]
		\label{wn03:def:dynamic-fixed}
		Fix a realized packaging system $\mathbb P$ and a
		nonempty small state space $X$.
		\begin{enumerate}[label=(\roman*),leftmargin=*,nosep]
			\item \textbf{Autonomous deterministic system.}
			An autonomous deterministic fixed-support dynamical
			packaged object is a triple
			\[
			\mathbf D=(\mathbb P,X,F),
			\qquad F:X\to X.
			\]
			The trajectory starting at $x\in X$ is
			\[
			x_t=F^t(x),\qquad t\in\mathbb N_0,
			\]
			where $F^0=\mathrm{id}_X$.
			
			\item \textbf{Nonautonomous deterministic system.}
			A nonautonomous deterministic fixed-support dynamical
			packaged object is a triple
			\[
			\mathbf D=(\mathbb P,X,(F_t)_{t\in\mathbb N_0}),
			\qquad F_t:X\to X.
			\]
			Its evolution from time $s$ to time $t>s$ is
			\[
			\Phi(t,s):=F_{t-1}\circ\cdots\circ F_s,
			\qquad \Phi(s,s):=\mathrm{id}_X.
			\]
			Thus a trajectory satisfies $x_t=\Phi(t,0)(x_0)$.
		\end{enumerate}
		In both cases, $\mathbb P$ and $X$ remain fixed
		throughout the evolution.
	\end{ddefinition}

	The temporal update is distinct from the native realization
	relation $\mathsf J_\Xi$.
	States may represent labelings, configurations, formations with memory, graph states, or algebraic objects, according to the declared interpretation.

	\begin{ddefinition}[Incidence-resolved state presentation]
		\label{wn03:def:incidence-state}
		A fixed-support dynamical packaged object is
		\textbf{incidence-resolved} if it specifies:
		\begin{enumerate}[label=(\roman*),leftmargin=*,nosep]
			\item \textbf{Coordinate state space.}
			A nonempty alphabet $A_e$ for each internal unit
			$e\in E_{\mathfrak P}$, with
			\[
			X\subseteq\prod_{e\in E_{\mathfrak P}}A_e.
			\]
			
			\item \textbf{Carrier observation.}
			For each carrier $K\in\operatorname{Pkg}_{\mathfrak P}$,
			let
			\[
			I_K:=\{e\in E_{\mathfrak P}:
			e\triangleleft_{\mathfrak P}K\}
			\]
			be its incidence support.
			Its local state space and observation map are
			\[
			X_K:=\{x|_{I_K}:x\in X\},
			\qquad
			r_K:X\to X_K,\quad r_K(x):=x|_{I_K}.
			\]
		\end{enumerate}
		Incidence determines the observed coordinates, but not
		their autonomous evolution.
		Any identification with	native candidate spaces $\mathsf L_{\mathfrak P}(K)$	or realized spaces $\mathsf R_{\mathfrak P}(K)$	requires a separate interpretation.
	\end{ddefinition}

	Proposition~\ref{prop:incidence-faithfulness-v06}
	characterizes when the carrier observations jointly
	distinguish states.	
	By Proposition~\ref{dyn45:prop-conservative-axioms},
	the deterministic extension is conservative over the
	primitive static theory.
	Additional restrictive coupling	conditions require separate justification.

	\begin{ddefinition}[Partial deterministic, relational, and controlled laws]
		\label{wn03:def:transitions}
		Let $X$ be a nonempty state space.
		The following discrete temporal laws specialize
		Definition~\ref{dyn46:def-time-types}:
		\begin{enumerate}[label=(\roman*),leftmargin=*,nosep]
			\item \textbf{Partial deterministic law.}
			A map $F:D\to X$, where $D\subseteq X$, specifies
			the unique successor $F(x)$ of each $x\in D$.
			States outside $D$ have no successor.
			
			\item \textbf{Relational law.}
			A relation $R\subseteq X^2$ specifies the
			permitted transitions, with successor set
			\[
			R(x):=\{y\in X:(x,y)\in R\}.
			\]
			
			\item \textbf{Deterministic controlled law.}
			A triple $(U,A,F)$ consists of a control set $U$,
			an admissibility relation $A\subseteq X\times U$,
			and an update map $F:A\to X$.
			A control $u$ is admissible at $x$ if $(x,u)\in A$,
			in which case the next state is $F(x,u)$.
			Controls are unrestricted if $A=X\times U$.
		\end{enumerate}
	\end{ddefinition}

	\subsubsection{Structural evolution and changes of type}
	\label{wn03:subsec:reconfiguration}

	A fixed-support evolution changes states without changing the base.
	The following deterministic specialization of the relational reconfiguration in Definition~\ref{def:reconfig-v05} allows the realized base itself to change.

	\begin{ddefinition}[Bare and equipped reconfigurable evolutions]
		\label{wn03:def:reconfigurable}
		\leavevmode
		\begin{enumerate}[label=(\roman*),leftmargin=*,nosep]
			\item \textbf{Bare reconfigurable evolution.}
			A bare reconfigurable evolution is a pair
			\[
			\mathbf D^{\mathrm{rec}}
			:=\bigl((\mathbb P_t,X_t)_{t\geq0},
			(F_t)_{t\geq0}\bigr),
			\]
			where, for each $t\geq0$,
			$\mathbb P_t=(\Xi_t,P_t,\mathfrak P_t)$ is a
			realized packaging system, $X_t$ is a nonempty
			state space associated with $\mathbb P_t$, and
			$F_t:X_t\to X_{t+1}$ is the temporal update.
			A trajectory $(x_t)_{t\geq0}$ satisfies
			\[
			x_t\in X_t,\qquad x_{t+1}=F_t(x_t).
			\]
			
			\item \textbf{Equipped reconfigurable evolution.}
			An equipped reconfigurable evolution is a pair
			\[
			\bigl(\mathbf D^{\mathrm{rec}},(\gamma_t)_{t\geq0}\bigr),
			\]
			where $\mathbf D^{\mathrm{rec}}$ is a bare
			reconfigurable evolution and
			\[
			\gamma_t:\mathfrak P_t\dashrightarrow\mathfrak P_{t+1}
			\]
			is a structural transport of a declared type, such as
			spans of internal-unit and carrier sets with
			incidence-respecting maps.
		\end{enumerate}
		The updates $F_t$ determine state evolution, while the transports $\gamma_t$ specify structural correspondences.
		Without structural transports, no canonical identification of units or carriers across time is assumed.
	\end{ddefinition}

	A choice of the next type can itself depend on the current state.
	Write $X=\coprod_i X_i$ for the disjoint union of state fibers over the possible bases and $p:X\to I$ for the type map.
	A joint step from $X_i$ to $X_j$ is a relation $E_{ij}\subseteq X_i\times X_j$.
	If a total deterministic $F:X\to X$	implements these steps, the type sequence has an autonomous update exactly when $p(Fx)=p(Fx')$ for all $x,x'\in X_i$, by Theorem~\ref{thm:exact-deterministic-chart-factor-criterion}
	on the attained types.
	Even when this condition fails, the joint law specifies
	which changes of type can occur.
	Comparing a named carrier before and after a change additionally uses the declared transport.
	This gives a precise distinction between a changing configuration on one base and state-dependent change of the base itself.

	\begin{dremark}[Transport between base types]
		The primitive packaging axioms give no mandatory maps of carriers, internal units, or candidates.
		Consequently, reconfiguration transport requires declared typing and preservation rules. Section~\ref{sec:morphisms-invariants} now defines categories of typed observed morphisms and time-preserving comparisons, while more general spans and changing-incidence assembly remain separately specified enrichments.
		Reconfigurable evolution extends nonautonomous fixed support.
		It is not required by the elementary globalization and recovery theorems below.
	\end{dremark}

	\subsubsection{Chart observations and the meaning of restriction}
	\label{wn03:subsec:observations}

	Consider a nonempty finite descriptor-label cover \(\mathcal H\) of a finite set \(V\).
	It is extra atlas data, rather than an axiom imposed on all packaged structures.
	Suppose a nonempty global state space \(X\) is given, together with local state spaces \(X_I\) and observations
	\[
	r_I:X\longrightarrow X_I\qquad (I\in\mathcal H).
	\]
	To compare distinct charts, take the support poset of subsets $S$ contained in some chart (with a separately assigned occurrence index when supports repeat).
	Specify $X_S$ and coherent restrictions $\rho^S_T:X_S\to X_T$ for $T\subseteq S$.
	In this specialization charts with the same support use the same support observation, and their occurrence markings are retained separately.
	Choose $X_\varnothing=\{*\}$, a convention not implied by coherence.
	The chart observations satisfy
	\(\rho^I_{I\cap J}r_I=\rho^J_{I\cap J}r_J\).
	For every $S$ in this poset define $r_S=\rho^I_Sr_I$ using any chart $I\supseteq S$; overlap agreement and coherence make this independent of $I$ and give (O3).
	Observations on other sites, or distinct readings on occurrences with the same support, require their own declared site data.
	General native observation sites instead use specified common lower sites, without assuming intersections exist.

	\begin{ddefinition}[Local observation, autonomous restriction, and transition shadow]
		\label{wn03:def:chart-shadow}
		Let $F:X\to X$ be a deterministic temporal law and
		$X_I$ the observation space of a chart $I$.
		\begin{enumerate}[label=(\roman*),leftmargin=*,nosep]
			\item \textbf{Local observation.}
			A local observation at $I$ is a map
			$r_I:X\to X_I$.
			
			\item \textbf{Autonomous restriction.}
			The observation $r_I$ admits an autonomous deterministic restriction if there exists a map $f_I:X_I\to X_I$ satisfying
			\[
			r_I\circ F=f_I\circ r_I.
			\]
			Such a map $f_I$ is called an autonomous local update.
			Its restriction to $r_I(X)$ is uniquely determined
			and defines an autonomous factor on the attained image.
			
			\item \textbf{Transition shadow.}
			The transition shadow of $F$ under $r_I$ is the relation
			\[
			R_I(F):=
			\{(r_I(x),r_I(F(x))):x\in X\}
			\subseteq X_I^2.
			\]
			It is always defined, has both endpoints in $r_I(X)$,
			and may be multivalued.
		\end{enumerate}
	\end{ddefinition}

	The existence of an autonomous local map is an additional dynamical property, not a consequence of the primitive packaging data or the ability to observe a component.
	An unobserved input can change the next chart value on a fixed observation fiber.

	The factor criterion of Theorem~\ref{thm:exact-deterministic-chart-factor-criterion} applies on the attained observation set $r_I(X)$.
	The shadow is the graph of its autonomous map exactly when every observed state determines one successor.
	Such a map extends to the nominal $X_I$ by arbitrary choices of successors on the unattained values.
	The data force no such choices.

	\begin{dexample}[Coupled coordinates have no autonomous one-coordinate factor]\label{wn03:ex:swap}
		Take \(X=\{0,1\}^2\), \(F(a,b)=(b,a)\), and observe \(r(a,b)=a\).
		The states \((0,0)\) and \((0,1)\) have the same observed value but different observed successors, so no \(f:\{0,1\}\to\{0,1\}\) satisfies \(rF=fr\).
		Its transition shadow is the full relation $\{0,1\}^2$. At horizon two the shadow allows the word $(0,0,1)$, whereas every observed swap trajectory satisfies $a_{t+2}=a_t$.
		Nevertheless, successive observations of the first coordinate at times zero and one recover the entire initial state \((a,b)\).
		This example will separate autonomous restriction from temporal observation depth.
	\end{dexample}

	\subsubsection{Compatible-state envelopes and exact globalization}
	\label{wn03:subsec:envelopes}

	The finite cover admits a canonical set of compatible chart states,
	\begin{equation}\label{wn03:eq:compatible-envelope}
		M_{\mathcal H}:=\left\{(x_I)_{I\in\mathcal H}\in\prod_{I\in\mathcal H}X_I:
		\rho^I_{I\cap J}(x_I)=\rho^J_{I\cap J}(x_J)\ \text{for all }I,J\right\}.
	\end{equation}
	Assume the observation tuple \(\iota=(r_I)_{I\in\mathcal H}:X\hookrightarrow M_{\mathcal H}\) is injective (\emph{joint faithfulness}). The image \(\iota(X)\) represents the \emph{globally admissible} compatible families. We do not require every element of \(M_{\mathcal H}\) to be globally admissible. A chart is \emph{state-realized} if \(r_I(X)=X_I\); this extra condition matters when a local update is to be identified on its entire nominal chart state space.

	Let \(f_I:X_I\to X_I\) be prescribed chart updates and let \(G=\prod_{I\in\mathcal H}f_I\) be their coordinatewise product on \(\prod_I X_I\). The local updates are \emph{gluing-compatible} if \(G(M_{\mathcal H})\subseteq M_{\mathcal H}\). This is stronger than merely agreeing along one realized trajectory. One sufficient local check is that, for each pair \(I,J\) and any pair \(u\in X_I,v\in X_J\) agreeing on the overlap, the updated values also agree on that overlap.

	\begin{dtheorem}[Exact deterministic globalization on a prescribed state space]
		\label{wn03:thm:global}
		Assume \(\iota:X\hookrightarrow M_{\mathcal H}\) is injective.
		Then a global deterministic law \(F:X\to X\) satisfying
		\begin{equation}\label{wn03:eq:intertwine}
			r_I\circ F=f_I\circ r_I\qquad(I\in\mathcal H)
		\end{equation}
		exists if and only if
		\begin{equation}\label{wn03:eq:invariant-image}
			\quad G(\iota(X))\subseteq\iota(X).
		\end{equation}
		When it exists it is unique and is given by \(F=\iota^{-1}\circ G\circ\iota\).
		No condition on $G$ outside $\iota(X)$ is needed for this equivalence.
		In particular, if \(\iota:X\overset{\sim}{\to} M_{\mathcal H}\), existence is equivalent to gluing compatibility on the entire envelope.
	\end{dtheorem}

	\begin{proof}
		If \(F\) satisfies \eqref{wn03:eq:intertwine}, then \(\iota F=G\iota\) by equality of every component, so the image is \(G\)-invariant.
		Conversely, if \eqref{wn03:eq:invariant-image} holds, injectivity of \(\iota\) allows the definition \(F(x)=\iota^{-1}(G\iota(x))\).
		This is a total map to \(X\) and satisfies every required square.
		If two global maps satisfy all squares, their observations agree on every chart and joint faithfulness forces equality.
		When the image is the entire envelope, \eqref{wn03:eq:invariant-image} is exactly gluing compatibility.
	\end{proof}

	\begin{dremark}[Exactness versus unreachable chart states]
		The equation \(r_IF=f_Ir_I\) realizes each prescribed update on chart states that occur as observations of global states.
		If some \(r_I(X)\subsetneq X_I\), a global map cannot, by itself, test or recover the prescribed behavior of \(f_I\) outside \(r_I(X)\).
		For the stronger interpretation of exact extension on every nominal local state, require state-realized charts \(r_I(X)=X_I\), or specify how unrealized chart states are treated.
		This is distinct from joint faithfulness, which ensures uniqueness of the global update.
	\end{dremark}

	\begin{dcorollary}[Time-dependent globalization]
		\label{wn03:cor:time-global}
		Suppose \(\iota_t:X_t\hookrightarrow M_t\) is an injective state embedding at each time and a local time-step atlas induces \(G_t:M_t\to M_{t+1}\).
		A unique global time-step \(F_t:X_t\to X_{t+1}\) with \(\iota_{t+1}F_t=G_t\iota_t\) exists if and only if
		\(G_t(\iota_t(X_t))\subseteq\iota_{t+1}(X_{t+1})\).
		If this holds at every time, the assembled global propagators obey \(\Phi(t,s)=F_{t-1}\cdots F_s\).
	\end{dcorollary}

	\begin{proof}
		Apply the injective-image argument from Theorem~\ref{wn03:thm:global} at each time and compose the resulting maps.
	\end{proof}

	\begin{dexample}[Compatible maps can violate global state constraints]\label{wn03:ex:diagonal}
		Let \(V=\{a,b\}\), \(X_a=X_b=\{0,1\}\), and \(X=\{(0,0),(1,1)\}\subset\{0,1\}^2\). Use the two singleton charts, so their overlap restrictions are trivial and every pair of local maps is gluing-compatible. Put \(f_a(a)=1-a\) and \(f_b(b)=b\). Then \(G(0,0)=(1,0)\notin X\); no global map on \(X\) can realize both local laws. Each chart is state-realized, and the cover is acyclic. Thus static cover geometry and chartwise determinism do not by themselves ensure constrained-state dynamical globalization.
	\end{dexample}

	\paragraph{Relational atlases and a non-globalization obstruction.}

	For transition relations, projection to overlaps is natural but does not commute with joins:
	gluing an independently permitted transition on each chart can create spurious or impossible global transitions.

	\begin{dexample}[Pairwise-compatible transition relations with no global transition]
		\label{ex:pairwise-compatible-transition-relations}
		Take \(V=\{a,b,c\}\), binary coordinates, and the three charts \(\{a,b\},\{b,c\},\{a,c\}\).
		A chart relation relates a current local state \(u\) to a proposed next local state \(v\).
		Let
		\[
		R_{ab}=\{(u,v):v_a=v_b\},\qquad
		R_{bc}=\{(u,v):v_b=v_c\},\qquad
		R_{ac}=\{(u,v):v_a\ne v_c\},
		\]
		with all current-state coordinates unrestricted. On each singleton overlap, every current/next bit pair occurs, so all three relations have identical restrictions on overlaps. However, no global transition can satisfy \(v_a=v_b=v_c\) and \(v_a\ne v_c\).
		Thus a global relation \(R\subseteq(\{0,1\}^V)^2\) restricting exactly to the three nonempty chart relations cannot exist. Pairwise relational compatibility is insufficient even though every local relation is total and individually nonempty.
	\end{dexample}

	The example has compatible one-step chart shadows and an empty set of global witnesses. The same obstruction appears as failure of a global section in relational models \cite{wn11:AbramskyBrandenburger2011}.

	\subsubsection{Temporal observations and a fundamental no-go theorem}
	\label{wn03:subsec:temporal}

	To prepare Section~\ref{sec:width}, distinguish the intrinsic \emph{implication-propagation recovery rank} \(\operatorname{ipr}\) of Paper~I from a \emph{state-observation depth} of the present system.
	They are not defined by the same assertions and no equality is assumed.

	\begin{ddefinition}[Temporal observation depth]
		\label{wn03:def:observation-depth}
		Let $F:X\to X$ be a deterministic update and
		$r:X\to Y$ a state observation.
		For $m\in\mathbb N_0$, define the temporal observation map
		\[
		\mathcal O_m:X\to Y^{m+1},\qquad
		\mathcal O_m(x):=(r(F^t x))_{t=0}^{m}.
		\]
		The \textbf{temporal observation depth} of $(X,F)$
		with respect to $r$ is
		\[
		\operatorname{tod}(X,F;r):=
		\min\{m\in\mathbb N_0:\mathcal O_m
		\text{ is injective}\},
		\]
		with the convention $\min\varnothing=\infty$.
		Equivalently, $\operatorname{tod}(X,F;r)\leq m$
		if and only if
		\[
		\forall x\ne y\text{ in }X,\quad
		\exists\,t\in\{0,\ldots,m\}:
		r(F^t x)\ne r(F^t y).
		\]
	\end{ddefinition}

	\begin{dtheorem}[Autonomous observations cannot distinguish hidden states over time]
		\label{wn03:thm:no-go}
		If \(r:X\twoheadrightarrow Y\) admits an autonomous factor \(f:Y\to Y\) such that \(rF=fr\), then
		\[
		\operatorname{tod}(X,F;r)=0\ \text{if }r\text{ is injective, and }\infty\ \text{otherwise}.
		\]
		In particular, a nontrivial finite positive temporal observation depth is possible only for observations that do not define autonomous local deterministic restrictions.
	\end{dtheorem}

	\begin{proof}
		By induction, \(rF^t=f^tr\).
		Hence every observed trajectory is determined by \(r(x)\).
		If \(r\) is not injective, two states in the same fiber have identical observations at all times, and no \(\mathcal O_m\) is injective.
		If \(r\) is injective, the zero-horizon observation already determines the state.
	\end{proof}

	\begin{dproposition}[Sharp observation stabilization and depth bounds]\label{wn03:prop:finite-depth}
		Let $N=|X|\ge1$ and $b=|r(X)|$, and let $c$ be the number of classes of the all-future observation equivalence (the history quotient of Proposition~\ref{dyn45:eq-observation-quotient}). Then $b\le c\le N$.
		\begin{enumerate}[label=(\roman*),leftmargin=*,nosep]
			\item \emph{Stabilization.} Equality of two observation histories through time $c-b$ is equivalent to equality at all future times, in particular through time $N-b$.
			
			\item \emph{Observation depth.} The observation depth is finite exactly when $c=N$, and then is at most $N-b$.
			
			\item \emph{Sharpness.} For each $2\le b\le c\le N$ the stabilization bound $c-b$ is attained, even when $c<N$.
		\end{enumerate}
	\end{dproposition}

	\begin{proof}
		\leavevmode
		\begin{enumerate}[label={(iii)},leftmargin=*,nosep]
			\item[(i)] We use Moore's refinement of the state partition by observable experiments \cite{Moore1956}. Let $\Pi_m$ identify states with equal observations at times $0,\ldots,m$.
			A pair is in $\Pi_{m+1}$ exactly when its first observations agree and its images under $F$ are in $\Pi_m$. Equality of successive partitions therefore implies equality at all subsequent steps. The first partition has $b$ blocks and each strict refinement adds at least one. Every finite-history partition is coarser than the all-future partition, which has $c$ blocks. If a refinement is not strict, all later partitions coincide with it. Otherwise, after $c-b$ strict refinements there are $c$ blocks, so the partition must already equal the all-future partition. The bound for $\Pi_{N-b}$ follows from $c\le N$.
			
			\item[(ii)] Injectivity at any horizon forces $c=N$, while $c=N$ makes $\Pi_{c-b}$ discrete.
			
			\item[(iii)] For sharpness first use $c$ states: take $L=c-b+1$ states $x_1,\ldots,x_L$ with observation $0$, followed by a fixed state of observation $1$. Send $x_i$ to $x_{i+1}$ for $i<L$, and $x_L$ to the fixed state. Add $b-2$ fixed states, each with its own observation. The first visits to observation $1$ from the chain occur after $L,L-1,\ldots,1$ steps. At time $L-1$ all $c$ starting states are distinguished. If $L\ge2$, the first two starting states still agree through time $L-2$. When $L=1$ the initial observation is injective on these states. To reach any $N\ge c$, add $N-c$ states with the same observation and successor as $x_1$. Their complete histories coincide with that of $x_1$, creating no additional history classes and preserving the least stabilization horizon. Taking $c=N$ recovers the finite-depth sharpness statement. For $b=1$ all histories are constant, $c=1$, and the finite-depth case is $N=1$. \qedhere
		\end{enumerate}
	\end{proof}
	The observation partition compares two initial states. The witness widths in Section~\ref{sec:width} instead count the cells needed to refute a partial history.

	\subsubsection{Finite-support preservation certificates}
	\label{wn03:subsec:finite-support}

	The image criterion of Theorem~\ref{wn03:thm:global} characterizes invariance of the constrained state set. The following special case makes that criterion checkable from supported local rules.

	\begin{dtheorem}[Finite-window admission-preservation test]\label{wn03:thm:finite-window}
		Let \(V\) be finite, \(A\ne\varnothing\) a finite alphabet, and let global admissibility be defined by constraints \(C_\alpha\subseteq A^{S_\alpha}\) for supports \(S_\alpha\subseteq V\):
		\[
		X=\{x\in A^V:x|_{S_\alpha}\in C_\alpha\ \text{for all }\alpha\}\ne\varnothing.
		\]
		For each \(i\in V\), let \(D_i\subseteq V\) and \(\phi_i:A^{D_i}\to A\) define the coordinate update \((\widehat Fx)_i=\phi_i(x|_{D_i})\) on \(A^V\).
		Set \(U_\alpha=\bigcup_{i\in S_\alpha}D_i\).
		Then \(\widehat F\) restricts to a map \(X\to X\) if and only if, for every constraint index \(\alpha\) and every globally realizable pattern \(p\in\pi_{U_\alpha}(X)\),
		\begin{equation}\label{wn03:eq:finite-window}
			(\phi_i(p|_{D_i}))_{i\in S_\alpha}\in C_\alpha.
		\end{equation}
		If \(|S_\alpha|\le k\) and \(|D_i|\le d\), every required test depends on at most \(kd\) input coordinates.
		This is a support bound, not a computational-complexity bound for finding \(\pi_{U_\alpha}(X)\).
		For an effective finite test the constraint data must also be given effectively (for example as a finite list).
	\end{dtheorem}

	\begin{proof}
		For \(x\in X\), the output on \(S_\alpha\) is \((\phi_i((x|_{U_\alpha})|_{D_i}))_{i\in S_\alpha}\); it depends only on \(p=x|_{U_\alpha}\). 
		Invariance is exactly the requirement that this output belong to \(C_\alpha\) for every allowed \(x\) and \(\alpha\).
		Replacing \(x\) with all its projections \(p\in\pi_{U_\alpha}(X)\) gives the stated equivalent condition.
		Finally, \(|U_\alpha|\le\sum_{i\in S_\alpha}|D_i|\le kd\).
	\end{proof}

	With the global state set fixed, invariance can be tested on the inputs of each constraint.
	We next compare evolutions that need not have the same state set.

	\section{Morphisms, equivalence, invariants, and fixed points}
	\label{sec:morphisms-invariants}

	A comparison may preserve the evolution and still lose information about individual states, complete paths or recurrence.
	Proposition~\ref{dyn45:eq-observation-quotient} constructs the deterministic quotient that preserves all future observations.
	Proposition~\ref{dyn45:fp-kernels} identifies the largest safe sets for the chosen successor rule.
	Between these results we specify the morphisms and markings under which each construction is invariant.

	\subsection{Typed morphisms of observed dynamics}
	\label{dyn45:subsec-morphism-category}

	Fix a category $\mathbf B$ of permitted comparisons of realized packaging bases and a common initialization convention.
	Its arrows and their preservation requirements are specified by the native specialization; the identity-only category
	is always available.
	For structural bases these comparisons also state which observation sites, signatures, interface comparisons,
	operational admissions, and markings are retained.
	The static axioms alone do not select one category of all such maps.

	\begin{definition}[Observed relational morphism]
		\label{def:observed-morphism-v05}
		Let $\mathbf D$ and $\mathbf D'$ be observed relational
		dynamical packaged objects satisfying
		Definition~\ref{dyn45:def-observed-axioms}, with site
		posets $D$ and $D'$, respectively.
		An \textbf{observed relational morphism} is a tuple
		\[
		\Psi=(\alpha,f,\phi,(f_U)_{U\in D})
		:\mathbf D\longrightarrow\mathbf D',
		\]
		consisting of a permitted base arrow
		$\alpha:\mathbb P\to\mathbb P'$, a state map $f:X\to X'$,
		an order-preserving site map $\phi:D\to D'$,
		and local observation maps $f_U:X_U\to X'_{\phi(U)}$,
		satisfying the following conditions:
		\begin{enumerate}[label=(\roman*),leftmargin=*,nosep]
			\item \textbf{Observation compatibility.}
			The maps commute with observations and restrictions:
			\begin{align}
				r'_{\phi(U)}\circ f
				&=f_U\circ r_U
				&& (U\in D),\label{eq:observed-square-v05}\\
				f_U\circ\rho_{V,U}
				&=\rho'_{\phi(V),\phi(U)}\circ f_V
				&& (U\leq V).\label{dyn45:eq-restriction-naturality}
			\end{align}
			
			\item \textbf{Structural and occurrence compatibility.}
			The site map $\phi$ respects the declared chart types, and $(\alpha,\phi)$ preserves the declared associations between native types and observation sites.
			For separately indexed chart families, the morphism
			also includes a permitted occurrence map
			$\theta:\mathcal H\to\mathcal H'$ satisfying
			\[
			j'\circ\theta=\phi\circ j
			\]
			and preserving the required occurrence markings.
			The map $\theta$ may be omitted from the notation
			when determined by the other data.
			
			\item \textbf{Temporal preservation.}
			The state map preserves permitted transitions:
			\[
			(x,y)\in E
			\quad\Longrightarrow\quad
			(f(x),f(y))\in E'.
			\]
			
			\item \textbf{Initialization preservation.}
			For initialized systems,
			\[
			f(X_{\mathrm{init}})\subseteq X'_{\mathrm{init}}.
			\]
		\end{enumerate}
		Permitted type and occurrence maps contain identities
		and are closed under composition.
		For observations without a site-poset enrichment,
		use the discrete poset of chart occurrences.
		Temporal preservation is forward preservation;
		backward lifting and exactness are additional properties.
	\end{definition}

	\begin{proposition}[Identity, composition, and observation naturality]
		\label{prop:morphism-composition-v05}
		Typed observed morphisms form a category for each fixed compatible
		choice of native comparisons and observation types. Every morphism
		sends each finite source trajectory to a finite target trajectory;
		it preserves the stipulated initial condition. Forgetting states,
		evolution and observations gives a functor to $\mathbf B$.
	\end{proposition}

	\begin{proof}
		Identities satisfy all conditions. For composable morphisms
		$(\alpha,f,\phi,(f_U))$ and $(\beta,g,\psi,(g_W))$, define
		\[
		(\beta\alpha,gf,\psi\phi,(g_{\phi(U)}f_U)_{U\in D}).
		\]
		The observation squares compose because
		$r''_{\psi\phi(U)}gf=g_{\phi(U)}r'_{\phi(U)}f
		=g_{\phi(U)}f_Ur_U$. The two restriction-naturality equations
		likewise compose, and order preservation is closed under composition.
		Occurrence maps compose as functions: if $j'\theta=\phi j$ and
		$j''\theta'=\psi j'$, then $j''\theta'\theta=\psi\phi j$.
		Each source edge is preserved twice, and initial-set inclusions
		compose. Applying the edge condition at every time proves the path
		assertion. Associativity and identity laws follow componentwise from
		function composition and from $\mathbf B$. The projection onto the
		base arrow therefore preserves identities and composition.
	\end{proof}

	\subsection{Factors, lifting and exact comparison}
	\label{dyn45:subsec-factors-lifting}

	A factor sends fine states to coarse states and intertwines their updates.
	Recovering a prescribed coarse path also requires its successive steps to have compatible fine witnesses.

	\begin{proposition}[Deterministic morphisms and factors]
		\label{dyn45:prop-conjugacy}
		For total deterministic source and target laws, the transition
		condition in Definition~\ref{def:observed-morphism-v05} is equivalent to
		\[
		fF=F'f.
		\]
		It implies $fF^m=(F')^mf$ for every $m\ge0$. A surjective state map
		with this property is a \emph{deterministic factor map}; it realizes
		every target finite trajectory when all states are initially
		permitted. For initialized systems the same conclusion holds when
		$f(X_{\mathrm{init}})=X'_{\mathrm{init}}$.
	\end{proposition}

	\begin{proof}
		The unique target successor of $f(x)$ is $F'(f(x))$, so edge
		preservation is exactly the displayed equation. Induction proves its
		iterated version. Lift the initial state of any target path to a
		source initial state and iterate $F$; the iterated identity gives
		the prescribed path at every time.
	\end{proof}

	\begin{definition}[State-path exactness, chart-trace exactness, and prescribed-source lifting]
		\label{dyn45:def-path-exact}
		Let $\Psi:\mathbf D\to\mathbf D'$ be an observed
		relational morphism with state map $f$, site map $\phi$,
		and local maps $f_U$.
		For any map $g:A\to B$ and $m\in\mathbb N_0$, define
		\[
		g^{[m]}:A^{m+1}\to B^{m+1},\qquad
		g^{[m]}((a_t)_{t=0}^{m}):=(g(a_t))_{t=0}^{m}.
		\]
		\begin{enumerate}[label=(\roman*),leftmargin=*,nosep]
			\item \textbf{State-path exactness.}
			The morphism $\Psi$ is state-path exact if
			\[
			f^{[m]}(\mathsf T_E(m))=\mathsf T_{E'}(m)
			\qquad(m\geq0).
			\]
			
			\item \textbf{Chart-trace exactness.}
			The morphism $\Psi$ is chart-trace exact if
			\[
			f_U^{[m]}\bigl(r_U^{[m]}(\mathsf T_E(m))\bigr)
			=(r'_{\phi(U)})^{[m]}(\mathsf T_{E'}(m))
			\]
			for every matched chart site $U$ and every $m\geq0$. For a matched occurrence $I\in\mathcal H$, use $U=j(I)$ and the target site $j'(\theta(I))=\phi(U)$. If target trace languages are prescribed
			externally, use those specified languages on the right.
			
			\item \textbf{Prescribed-source forward lifting.}
			The morphism $\Psi$ has prescribed-source forward
			lifting if, for every $x\in X$ and $y'\in X'$,
			\[
			(f(x),y')\in E'
			\Longrightarrow
			\exists y\in X:\ (x,y)\in E,\quad f(y)=y'.
			\]
		\end{enumerate}
		Initialized state-path and chart-trace exactness are
		defined by replacing the trajectory sets on both sides
		with their initialized versions.
	\end{definition}

	\begin{proposition}[Exactness and a sufficient lifting condition]
		\label{dyn45:prop-path-lifting}
		\leavevmode
		\begin{enumerate}[label=(\roman*),leftmargin=*,nosep]
			\item \emph{Composition and projection.} State-path exact morphisms contain identities and are closed under
			composition; they are chart-trace exact on every matched chart.
			
			\item \emph{Lifting criterion.} A surjective state morphism with prescribed-source forward lifting
			is state-path exact. In the initialized case, it suffices to replace
			surjectivity by $f(X_{\mathrm{init}})=X'_{\mathrm{init}}$.
		\end{enumerate}
	\end{proposition}

	\begin{proof}
		\leavevmode
		\begin{enumerate}[label={(iii)},leftmargin=*,nosep]
			\item[(i)] Images of path sets compose, proving the first assertion. Apply an
			observation map to the path-image equality and use the observation
			square for the second.
			
			\item[(ii)] For lifting, choose a preimage of the target
			path's first state and successively lift each of its edges from the
			source state already chosen. This yields the required finite path;
			forward preservation gives the reverse inclusion. The same argument
			starts in $X_{\mathrm{init}}$ under the initialized hypothesis. \qedhere
		\end{enumerate}
	\end{proof}

	\begin{example}[One-step exactness need not lift paths]
		\label{dyn45:ex-edge-exact-not-path}
		Take $X=\{a,b,c,d\}$ with only the edges $a\to b$ and $c\to d$.
		Let $X'=\{0,1\}$ with the two-cycle $0\to1\to0$, and put
		$f(a)=f(d)=0$, $f(b)=f(c)=1$. Then $f$ is surjective and
		$(f\times f)(E)=E'$, but the target path $(0,1,0)$ has no source
		lift. This uses only the relational core; nonblocking and a
		prescribed-source lifting condition are additional restrictions.
	\end{example}

	\subsection{Strict structure transport and dynamical isomorphism}
	\label{subsec:strict-morphisms-v06}

	An isomorphism must transport every piece of data that the comparison marks.
	We start with the endpoint, add the observation maps, and then check how their incidence supports transform.

	\begin{definition}[Strict packaged-output and realized-system isomorphisms]
		\label{def:strict-static-isomorphism-v06}
		\leavevmode
		\begin{enumerate}[label=(\roman*),leftmargin=*,nosep]
			\item \textbf{Strict packaged-output isomorphism.}
			A strict packaged-output isomorphism from
			$\mathfrak P$ to $\mathfrak P'$ consists of bijections
			\[
			\begin{aligned}
				\alpha_C&:\operatorname{Pkg}_{\mathfrak P}
				\to\operatorname{Pkg}_{\mathfrak P'},\\
				\alpha_E&:E_{\mathfrak P}\to E_{\mathfrak P'},\\
				\alpha_{L,K}&:\mathsf L_{\mathfrak P}(K)
				\to\mathsf L_{\mathfrak P'}(\alpha_C(K))
			\end{aligned}
			\]
			for every carrier $K$, satisfying
			\[
			e\triangleleft_{\mathfrak P}K
			\iff
			\alpha_E(e)\triangleleft_{\mathfrak P'}\alpha_C(K)
			\]
			and
			\[
			\alpha_{L,K}(\mathsf R_{\mathfrak P}(K))
			=\mathsf R_{\mathfrak P'}(\alpha_C(K)).
			\]
			All bijections respect the declared ambient types
			and permitted interpretations.
			When candidate-space transports are prescribed,
			the maps $\alpha_{L,K}$ must coincide with them.
			They are not determined by $\alpha_E$ alone.
			Identity transports may be used when no other
			comparison is declared.
			
			\item \textbf{Realized-system isomorphism.}
			A realized-system isomorphism between
			$\mathbb P=(\Xi,P,\mathfrak P)$ and
			$\mathbb P'=(\Xi',P',\mathfrak P')$ consists of
			a strict packaged-output isomorphism together with
			a permitted invertible native transport carrying
			$\Xi$, $P$, and the distinguished realization
			to their counterparts in $\mathbb P'$.
			If the native law is retained, the induced formation
			and output transports $\alpha_F$ and $\alpha_O$
			satisfy
			\[
			\Xi\vdash Q\rightsquigarrow\mathfrak Q
			\iff
			\Xi'\vdash\alpha_F(Q)
			\rightsquigarrow\alpha_O(\mathfrak Q)
			\]
			on their declared domains.
		\end{enumerate}
		An output isomorphism alone does not determine
		a realized-system isomorphism.
	\end{definition}

	\begin{definition}[Observed dynamical isomorphism]
		\label{dyn45:def-observed-isomorphism}
		An \emph{observed dynamical isomorphism} is an observed
		relational morphism
		\[
		\Psi=(\alpha,f,\phi,(f_U)_{U\in D})
		:\mathbf D\longrightarrow\mathbf D'
		\]
		as in Definition~\ref{def:observed-morphism-v05},
		satisfying the following conditions:
		\begin{enumerate}[label=(\roman*),leftmargin=*,nosep]
			\item \textbf{Base isomorphism.}
			The base arrow $\alpha:\mathbb P\to\mathbb P'$
			is an isomorphism in $\mathbf B$.
			
			\item \textbf{State and observation isomorphisms.}
			The state map $f:X\to X'$ and all local maps
			$f_U:X_U\to X'_{\phi(U)}$ are bijections.
			The site map $\phi:D\to D'$ is an order isomorphism.
			All observation and restriction squares commute.
			
			\item \textbf{Structural and occurrence isomorphisms.}
			The native and site type transports are permitted
			in both directions.
			When chart occurrences are separately indexed,
			the occurrence map $\theta:\mathcal H\to\mathcal H'$
			is a bijection with permitted inverse.
			All retained markings, value labels, and costs are
			transported according to the comparison contract
			in both directions.
			
			\item \textbf{Temporal and initialization exactness.}
			The transition relations satisfy
			\[
			E'=(f\times f)(E).
			\]
			When initialization is retained,
			$X'_{\mathrm{init}}=f(X_{\mathrm{init}})$.
		\end{enumerate}
		Equivalently, $\Psi$ is an invertible observed relational
		morphism. For total deterministic laws, temporal
		exactness is equivalent to conjugacy
		\[
		F'=f\circ F\circ f^{-1}.
		\]
	\end{definition}

	\begin{proposition}[Characterization and invariance under isomorphism]
		\label{dyn45:prop-isomorphism-characterization}
		\leavevmode
		\begin{enumerate}[label=(\roman*),leftmargin=*,nosep]
			\item \emph{Invertibility.} The concrete conditions in Definition~\ref{dyn45:def-observed-isomorphism}
			are equivalent to categorical invertibility.
			
			\item \emph{Invariant properties.} They induce bijections
			of all corresponding finite path sets and preserve raw and coherent
			observed axioms, determinism, nonblocking, and reversibility when
			these are imposed.
			
			\item \emph{Exactness.} For coordinate atlases, if the comparison also transports the
			prescribed local state and law data by atom-preserving bijections,
			state-, one-step-, and all-horizon exactness are preserved.
		\end{enumerate}
	\end{proposition}

	\begin{proof}
		\leavevmode
		\begin{enumerate}[label={(iii)},leftmargin=*,nosep]
			\item[(i)] A categorical inverse forces componentwise inverse functions and
			inverse base, site and occurrence arrows. Transition preservation in both
			directions yields equality of the edge relations, and initialized
			inclusions in both directions yield equality of initial sets.
			Conversely, invert all the component bijections. The commuting
			squares can be rearranged into the inverse observation, restriction
			and occurrence squares; exact edge and initial-set transport supplies
			the remaining conditions.
			
			\item[(ii)] Applying $f$ or $f^{-1}$ at every time
			gives mutually inverse maps of path sets. The axiom identities and
			the listed properties transport through those bijections.
			
			\item[(iii)] For an	atom-preserving atlas comparison, restriction of a transported path
			is the transported restriction; applying this to each defining
			projection equality proves exactness invariance. \qedhere
		\end{enumerate}
	\end{proof}

	\begin{definition}[Incidence-observed dynamical isomorphism]
		\label{def:incidence-dyn-iso-v06}
		Let $\mathbf D$ and $\mathbf D'$ be incidence-resolved
		dynamical packaged objects.
		An \emph{incidence-observed dynamical isomorphism}
		consists of the following data and conditions:
		\begin{enumerate}[label=(\roman*),leftmargin=*,nosep]
			\item \textbf{Structural isomorphism.}
			A strict packaged-output isomorphism as in
			Definition~\ref{def:strict-static-isomorphism-v06},
			whose carrier bijection satisfies
			$\alpha_C(\mathcal K)=\mathcal K'$.
			
			\item \textbf{Coordinate isomorphism.}
			For each internal unit $e\in E_{\mathfrak P}$,
			a bijection
			$a_e:A_e\to A'_{\alpha_E(e)}$.
			These bijections induce a product bijection
			$\widehat f$ defined by
			\[
			\widehat f(x)_{\alpha_E(e)}:=a_e(x_e),
			\]
			with $\widehat f(X)=X'$.
			Put $f:=\widehat f|_X:X\to X'$.
			
			\item \textbf{Observation isomorphism.}
			For each $K\in\mathcal K$, the induced chart
			bijection $f_K:X_K\to X'_{\alpha_C(K)}$ is defined by
			\[
			f_K(u)_{\alpha_E(e)}:=a_e(u_e)
			\qquad(e\in I_K).
			\]
			Consequently,
			$r'_{\alpha_C(K)}\circ f=f_K\circ r_K$.
			
			\item \textbf{Temporal isomorphism.}
			The state bijection satisfies
			\[
			E'=(f\times f)(E),
			\]
			and $X'_{\mathrm{init}}=f(X_{\mathrm{init}})$
			when initialization is retained.
		\end{enumerate}
		For realized dynamical systems, the strict output
		isomorphism must additionally be accompanied by the
		native-base transport of
		Definition~\ref{def:strict-static-isomorphism-v06}.
		With the required markings and transports retained,
		this is a specialization of
		Definition~\ref{dyn45:def-observed-isomorphism}.
		Otherwise, the comparison is at the output level.
	\end{definition}

	\begin{proposition}[Functoriality of the incidence support construction]
		\label{prop:incidence-functoriality-v06}
		Under strict isomorphisms preserving the observed carrier family,
		$\alpha_E(I_K)=I'_{\alpha_C(K)}$. Hence the incidence-generated
		chart family, its occurrence labels and its coordinate-restriction
		diagram are preserved under identities and composition.
		The faithfulness assertion of
		Proposition~\ref{prop:incidence-faithfulness-v06} is invariant under
		such isomorphisms.
	\end{proposition}

	\begin{proof}
		The support equality is preservation and reflection of incidence.
		Coordinate relabeling commutes with restriction, and composition
		preserves all the displayed equalities. Joint observation maps
		are transported by bijections, preserving injectivity.
	\end{proof}

	\begin{remark}[Coordinate atomization is structural]
		\label{rem:atomization-v05}
		The incidence-observed category retains independently specifiable
		coordinate atoms. Coordinatewise bijections carry partial patterns
		and one-cell erasures bijectively, preserving their feasibility and
		minimal-nogood support sizes. Merging several observables into one
		atomic coordinate, or splitting an atom, is not such an isomorphism.
		Section~\ref{wn16:sec:composition-rank} gives a four-state example
		whose witness width changes under that recoding.
	\end{remark}

	\begin{remark}[Limits of the adopted convention]
		\label{rem:strict-scope-v06}
		An arbitrary state-graph isomorphism preserves that graph but need
		not preserve an observation decomposition, a native realization,
		or a structural assembly language. Likewise, a general forward
		morphism need not preserve numerical recovery ranks, reflect
		nonblocking, or lift trajectories. Every invariance statement
		therefore specifies its comparison class.
	\end{remark}

	\subsection{Coupled morphisms and temporal-law comparisons}
	\label{dyn46:subsec-coupled-morphisms}

	For the enriched theory a morphism must specify which coupling data
	it preserves, not just the underlying state dynamics. In particular
	an observed isomorphism of Definition~\ref{dyn45:def-observed-isomorphism}
	need not transport a marked channel decomposition. We give explicit
	categories for the fixed-support specializations used as the shared
	interface. Reconfiguration retains its independently declared
	transport category below.

	\begin{definition}[Temporal morphisms of fixed law type]
		\label{dyn46:def-typed-morphism}
		Fix source and target systems of the same temporal-law type, time domain
		and initialization convention. A \emph{temporal morphism} consists of a
		permitted base arrow and state map $f:X\to X'$, with the retained observation, restriction and occurrence data and compatibility conditions of Definition~\ref{def:observed-morphism-v05}. Its temporal condition is as follows:
		\begin{enumerate}[label=(\roman*),leftmargin=*,nosep]
		 \item \textbf{Deterministic actions.} $f\circ\Phi_t=\Phi'_t\circ f$
		 for every retained time $t$ (in particular $f\circ F=F'\circ f$
		 for discrete total maps).
		 \item \textbf{Markov kernels.} The map $f$ is measurable and
		 $f_*P(x,\cdot)=P'(f(x),\cdot)$ for every $x\in X$, where
		 $(f_*P(x,\cdot))(B)=P(x,f^{-1}(B))$ for measurable $B\subseteq X'$.
		 \item \textbf{Markov semigroups.} The kernel condition holds for every
		 time $t$, with $P,P'$ replaced by $P_t,P'_t$. For finite generators it
		 is equivalent to $Q(f^*g)=f^*(Q'g)$ for all real-valued target
		 functions $g$, where $f^*g=g\circ f$.
		\end{enumerate}
		Relational morphisms preserve edges; partial deterministic laws use this
		condition on their graphs, so defined steps map to defined steps with
		$f(F(x))=F'(f(x))$. An isomorphism also reflects definedness.
		Continuous or smooth comparisons require the declared regularity.
		Initialized Markov morphisms satisfy $f_*\mu_0=\mu'_0$, or the declared
		transport condition on permitted initial laws, closed under identities
		and composition. Controlled comparisons specify control maps and the
		laws they intertwine. An isomorphism has an inverse in the same category
		and preserves the full law, not only its support. Time rescalings and
		state-dependent orbital time changes are distinct comparison types.
	\end{definition}

	\begin{definition}[Marked coupling morphisms in a coordinate specialization]
		\label{dyn46:def-coupled-morphism}
		Fix two coordinate-resolved objects of the same fixed-support law type and the permitted native comparisons.
		A \textbf{marked coupling morphism} consists of the following data and conditions:
		\begin{enumerate}[label=(\roman*),leftmargin=*,nosep]
			\item \emph{Coordinate and native data.} A coordinate coupling comparison
			retains a native unit bijection $\sigma:V\to V'$ (with its compatible
			base transport), coordinate maps $f_v:A_v\to A'_{\sigma(v)}$ defining $f(x)_{\sigma(v)}=f_v(x_v)$ with $f(X)\subseteq X'$, and the typed observation comparison if present. Retained native interpretation maps commute, and retained
			interpretation relations and admitted tuples are carried forward
			by their typed maps.
			
			\item \emph{Interaction data.} For an interaction-marked comparison additionally
			specify an occurrence map $\psi:\mathcal J\to\mathcal J'$ with
			\[
			S'(\psi J)=\sigma(S(J)),\qquad W'_{\psi J}=\sigma(W_J)
			\]
			when write sets are retained, and the analogous equality for retained
			read sets. For carrier-supported assignments require
			$\iota'(\psi J)=\alpha_C(\iota(J))$, with $\alpha_C$ the carrier
			transport. Require preservation of each marked channel in its
			own law type: image inclusion for relations, intertwining for maps,
			and the pullback generator equation for finite jump channels.
			
			\item \emph{Aggregate law and other markings.} Require the aggregate temporal-law condition of Definition~\ref{dyn46:def-typed-morphism} independently of the channel conditions. This
			last condition matters for occurrence maps that merge or omit channels.
			Initialization, interface comparisons, clock or composition markings
			and other retained sorts have their stipulated preservation equations.
			
			\item \emph{Composition and isomorphism.} Choose the allowed typed maps to include identities and be closed
			under composition. A strict coupled isomorphism is an invertible
			arrow of this category. In particular it has a bijection of retained
			channel occurrences, carries all retained supports and native
			interpretations exactly and transports the aggregate and channel laws.
		\end{enumerate}
		More general site-based comparisons use their declared typed maps and the same transport principles.
	\end{definition}

	\begin{proposition}[Coupling categories, reduct functors, and non-collapse]
		\label{dyn46:prop-coupled-category}
		\leavevmode
		\begin{enumerate}[label=(\roman*),leftmargin=*,nosep]
			\item \emph{Categories.} With the above choices, same-type temporal morphisms and marked
			coordinate coupling morphisms each form categories.
			
			\item \emph{Reducts.} Forgetting
			channel markings, forgetting observations, and then forgetting the
			law (where the corresponding reduct is defined) give functors.
			
			\item \emph{Generator comparison.} In the finite generator case the pullback equation is equivalent to
			the all-time semigroup equation.
			
			\item \emph{Non-collapse.} Forgetting observations or forgetting marked channels can identify
			distinct isomorphism classes of objects. This is a statement about
			object classes, not about reflection of invertibility of a given arrow.
		\end{enumerate}
	\end{proposition}

	\begin{proof}
		\leavevmode
		\begin{enumerate}[label={(iii)},leftmargin=*,nosep]
			\item[(i)] The observation and restriction equations compose as in
			Proposition~\ref{prop:morphism-composition-v05}. The support and
			incidence equations compose under the unit and occurrence maps.
			Map equations compose, and forward images of relations and admission
			sets compose by inclusions. Pushforward of measures is functorial;
			for generators $Qf^*=f^*Q'$ and $Q'g^*=g^*Q''$ imply
			$Q(gf)^*=(gf)^*Q''$.
			
			\item[(ii)] Identities
			satisfy every condition, and forgetting retained data respects
			composition.
			
			\item[(iii)] Iteration of the generator equation and the
			convergent exponential series give its semigroup version;
			differentiating the latter at zero gives the converse.
			
			\item[(iv)] Example~\ref{dyn46:ex-coupled-noncollapse} proves the
			first isomorphism-class separation. For the second, take two fixed
			labeled channels on a binary coordinate, with rate-one flip generator
			$Q_0$. Compare $(Q_0,Q_0)$ with $(\tfrac12 Q_0,\tfrac32 Q_0)$.
			Both have aggregate $2Q_0$, with identical supports and observations,
			but their rate multisets differ under every state relabeling and
			even under a permutation of the channels. The marked models are not
			isomorphic. \qedhere
		\end{enumerate}
	\end{proof}

	\begin{remark}[Comparison contracts across temporal types]
		An observed relational morphism may forget law weights. It is not
		a stochastic factor unless the kernel or generator equation holds.
		A graph-of-map reduction retains a deterministic step but a sampled
		flow map does not retain the time action. These reducts can give
		weaker contextual equivalences; each test must be defined on the
		retained signature and observation semantics. There is no universal
		identification of relational, stochastic and continuous-time
		comparison categories. In particular, arbitrary controlled or
		scheduled composition is not implied by a single-step morphism.
	\end{remark}

	\subsection{Exact stochastic observation}
	\label{dyn46:subsec-stochastic-factor}

	A weighted observation evolves autonomously when the total weight sent to each observed fiber depends only on the present reading.
	This is strong lumpability for a Markov kernel and has a generator counterpart \cite{KemenySnell}.

	\begin{proposition}[Strong lumpability for finite kernels and generators]
		\label{dyn46:prop-weighted-factor}
		Let $r:X\twoheadrightarrow Y$ be a surjection of finite sets.
		\begin{enumerate}[label=(\roman*),leftmargin=*,nosep]
			\item For a discrete Markov kernel $P$, there exists a single Markov kernel on $Y$ governing the projected process for every initial law (strong, or Kemeny--Snell, lumpability) if and only if for all $r(x)=r(x')$ and $y\in Y$,
			\[
			\sum_{z:r(z)=y}P(x,z)=\sum_{z:r(z)=y}P(x',z).
			\]
			
			\item For a finite homogeneous generator $Q$ the analogous condition is
			\[
			\sum_{z:r(z)=y}Q(x,z)=\sum_{z:r(z)=y}Q(x',z).
			\]
		\end{enumerate}
		Here the sum includes the diagonal when applicable. In either case
		the common sums define the unique factor kernel or generator.
		Equivalently the space of functions constant on observation fibers
		is preserved by $P$ or $Q$. For a finite generator this is equivalent
		to exact factorization of the entire semigroup.
	\end{proposition}

	\begin{proof}
		Apply the operator to the indicators of target fibers. The displayed
		condition is precisely that the result is constant on each source
		fiber. These indicators span all target functions, so the condition
		gives the unique operator on $Y$ intertwined by $r^*$.
		\begin{enumerate}[label={(iii)},leftmargin=*,nosep]
			\item[(i)] For the kernel the entries are nonnegative and sum to one. Iteration of the intertwining identity gives the all-time factor property. Choosing point-mass initial laws and examining one step proves necessity for one common factor law. Conditional on the source past, the next observed law depends only on the current observed state; conditioning again on the observed past proves the Markov property.
			
			\item[(ii)] For the generator, off-diagonal fiber sums are nonnegative and full row sums are zero. Exponentiating the intertwined generators gives the all-time factor semigroup, with the Markov property by the same conditioning argument at each future time increment. Differentiation at zero of the projected semigroup identity for point-mass initial laws proves necessity. \qedhere
		\end{enumerate}
	\end{proof}

	The strong-lumpability criterion applies to the aggregate kernel or generator \cite{KemenySnell}. A sum of generators preserves the space of fiber-constant functions whenever each summand does. Cancellation also permits the aggregate to preserve this space. For example, on $\{0,1\}^2$

	observe the first bit and use two channels flipping it at respective
	rates $b$ and $1-b$, where $b$ is the hidden second bit. Each summand
	fails the fiber condition but their aggregate flips the observed bit
	at rate one. Exact autonomy is
	neither a consequence of incidence provenance nor a prerequisite
	for an approximate or large-size autonomous law. Error estimates,
	tightness and convergence of initial laws must be supplied by the
	particular collective theorem.

	\subsection{Fixed structures and reconfiguration}

	When the formation or its endpoint changes, a step can have different source and target state spaces.
	Its type records both spaces and the transport allowed between them.

	\begin{definition}[Bare and equipped relational reconfiguration]
		\label{def:reconfig-v05}
		\leavevmode
		\begin{enumerate}[label=(\roman*),leftmargin=*,nosep]
		 \item \textbf{Bare relational reconfiguration.}
		 A bare reconfigurable evolution consists of the pair
		 \[
		 ((\mathbb P_t,X_t)_{t\in\mathbb N_0},(E_t)_{t\in\mathbb N_0}),
		 \qquad E_t\subseteq X_t\times X_{t+1},
		 \]
		 where $\mathbb P_t$ is a realized base and $X_t\ne\varnothing$ its
		 specified state set. A path through horizon $m$ is a tuple
		 $(x_0,\ldots,x_m)$ with $x_t\in X_t$ and $(x_t,x_{t+1})\in E_t$
		 for $0\le t<m$. Retained observations and site types are indexed by time.
		 \item \textbf{Equipped relational reconfiguration.}
		 An equipped reconfiguration consists of a bare reconfiguration and
		 declared native transport or correspondence data across successive times.
		 Persistent units and carriers are identified only through those data,
		 not by equality of labels.
		\end{enumerate}
	\end{definition}

	\begin{proposition}[Time-preserving comparison of reconfigurable laws]
		\label{dyn45:prop-reconfiguration-comparison}
		For two bare evolutions, suppose each $\alpha_t:\mathbb P_t\to
		\mathbb P'_t$ is a permitted native comparison and
		$f_t:X_t\to X'_t$ satisfies
		\[
		(x,y)\in E_t\Longrightarrow
		(f_t(x),f_{t+1}(y))\in E'_t.
		\]
		With any retained observations satisfying their timewise
		commuting squares, these data preserve finite time-indexed paths
		and are closed under identity and timewise composition. If
		equipped reconfigurations are compared, compatibility with their
		declared transport or correspondence data is an additional
		requirement. In particular, if the structural transport is given
		by arrows $\gamma_t:\mathbb P_t\to\mathbb P_{t+1}$ in the chosen
		base category, require
		\[
		\alpha_{t+1}\gamma_t=\gamma'_t\alpha_t.
		\]
	\end{proposition}

	\begin{proof}
		Apply the displayed implication separately at each step. For two
		composable comparisons, apply the first implication and then the
		second; this proves preservation by the timewise composite.
		Observation squares compose as in
		Proposition~\ref{prop:morphism-composition-v05}.
	\end{proof}

	\begin{remark}[Changes of state type]\label{rem:reconfig-scope-v05}
		A time-dependent base requires maps between its successive state types. Constant bases and state sets give a fixed-support time-indexed comparison. If the two step relations are also constant, this compares two autonomous laws, but it is an ordinary autonomous morphism only when the required single comparison map intertwines those laws (in particular, when all comparison maps and retained comparison data are time-independent). For example, on $\{0,1\}$ take $F=\mathrm{id}$, $F'=s$ the transposition, and $f_t=s^t$. Then $f_{t+1}F=F'f_t$ at each step, yet $F$ and $F'$ are not conjugate. The time index is therefore part of the comparison contract.
	\end{remark}

	\subsection{Comparison contracts and contextual dynamical equivalence}
	\label{dyn45:eq-contextual-section}

	An isomorphism compares all data retained by a chosen category. A behavioral
	comparison may deliberately retain fewer data. The word \emph{equivalent}
	therefore always refers to a declared comparison contract. The contract
	specifies the native base or permitted base transports, boundary types,
	observable atoms and value labels, time parameter, initialization convention,
	permitted constructions, and terminal tests. A change of these choices is a
	change of comparison problem. In particular, conjugacy does not include a change of time scale, and a smooth-flow comparison must intertwine
	every named time map with the required regularity, not just one sampled map.
	The following definition extends the continuation language of the Structural
	Packaging companion (\hyperref[comp:structural]{Paper~II}) by admitting dynamical components
	and temporal tests.

	\begin{definition}[Dynamical comparison contracts and contextual equivalence]
		\label{dyn45:eq-contract}
		\leavevmode
		\begin{enumerate}[label=(\roman*),leftmargin=*,nosep]
			\item \textbf{Dynamical comparison contract.}
			A dynamical comparison contract specifies:
			\begin{itemize}[nosep]
				\item For each marked boundary type $\tau$, a collection
				$\mathscr D_\tau$ of dynamical components.
				
				\item A family $\mathcal C_\tau$ of permitted finite
				well-typed expressions $C[-]$ with distinguished
				input type $\tau$ and output type $\sigma(C)$,
				including the identity expression.
				Each expression specifies its retained temporal laws
				and interface conditions.
				
				\item For each output type $\sigma$, a family
				$\mathcal T_\sigma$ of terminal tests.
				
				\item For each $C\in\mathcal C_\tau$ and
				$t\in\mathcal T_{\sigma(C)}$, an outcome space
				$\Omega_{C,t}$ and a total outcome map
				\[
				o_{C,t}:\mathscr D_\tau\to\Omega_{C,t},
				\qquad
				\mathbf D\mapsto t(C[\mathbf D]),
				\]
				with specified semantics for failed and
				multivalued formations.
			\end{itemize}
			A success test is included whenever definedness
			of formation is to be preserved.
			Permitted outcome renamings are incorporated through
			their declared equivalence classes.
			
			\item \textbf{Contextual dynamical equivalence.}
			Two components $\mathbf D,\mathbf D'\in\mathscr D_\tau$
			are contextually dynamically equivalent if
			\[
			\mathbf D\equiv_{\tau,\mathcal C,\mathcal T}\mathbf D'
			\iff
			o_{C,t}(\mathbf D)=o_{C,t}(\mathbf D')
			\]
			for every $C\in\mathcal C_\tau$ and
			$t\in\mathcal T_{\sigma(C)}$.
		\end{enumerate}
		Comparisons across different boundary types require
		a specified transport of their marked contracts.
		Recovery, numerical, and conservation equivalences
		are obtained from their respective test families.
	\end{definition}

	\begin{proposition}[Equivalence, refinement and substitution]
		\label{dyn45:eq-congruence}
		\leavevmode
		\begin{enumerate}[label=(\roman*),leftmargin=*,nosep]
			\item \emph{Equivalence.} For each fixed contract the relation in Definition~\ref{dyn45:eq-contract}
			is an equivalence relation.
			
			\item \emph{Refinement.} Enlarging the family of permitted
			context--test pairs while keeping the original outcomes and their
			comparisons fixed can only refine it.
			
			\item \emph{Substitution.} Suppose the permitted language is
			closed under well-typed substitution and that evaluating substituted
			expressions agrees with their declared semantics. Then replacement by an
			equivalent component preserves every outcome of every permitted enclosing
			expression.
			
			\item \emph{Partial constructors.} For a permitted partial single-valued constructor
			$\alpha:\prod_i\mathscr D_{\tau_i}\rightharpoonup\mathscr D_\tau$, if
			success of $\alpha$ is an available test, equivalent inputs give the same
			definedness, and their outputs, when defined, are contextually equivalent.
		\end{enumerate}
	\end{proposition}

	\begin{proof}
		\leavevmode
		\begin{enumerate}[label={(iii)},leftmargin=*,nosep]
			\item[(i)] The first assertion follows from reflexivity, symmetry and transitivity
			of equality of each outcome.
			
			\item[(ii)] Additional tests impose additional
			equalities, proving refinement.
			
			\item[(iii)] If $A[-]$ contains the replaced component
			and $C[-]$ is a surrounding continuation, closure under substitution
			makes $C[A[-]]$ an allowed test context for that component. Its outcomes
			agree by the definition.
			
			\item[(iv)] For a constructor with finitely many inputs,
			replace those inputs one at a time, keeping the others fixed. The
			success test preserves definedness at each replacement. When the outputs
			are defined, any continuation and terminal test on them is the same
			composed test on the replaced input, so the outputs are equivalent. \qedhere
		\end{enumerate}
	\end{proof}

	\begin{remark}[Isomorphism and admissible tests]
		\label{dyn45:eq-isomorphism-tests}
		Isomorphism implies contextual equivalence only for a contract whose
		constructions and tests respect the specified isomorphisms. Indeed, if
		each continuation carries an input isomorphism to an output isomorphism
		and each terminal test has equal outcomes on isomorphic outputs, all
		the defining equalities follow. A test that reads an untransported
		state name, a newly introduced cost, or a changed observable atom can
		violate this hypothesis. For multivalued constructors, substitution
		invariance of the chosen outcome semantics must be verified; the
		single-valued assertion above is not a theorem that an arbitrary
		admission relation descends to a partial function on equivalence classes.
	\end{remark}

	\subsection{Trace equivalence and hidden-state presentations}
	\label{dyn45:eq-trace-section}

	Fix the same marked observable set $O$ for the systems being compared.
	It may be the product of a specified finite family of chart readouts.
	Its atomization is retained whenever a support width is measured.
	Let	$r:X\to O$ be the joint readout.
	A chosen initial subset $X_{\mathrm{init}}\subseteq X$ is additional data; when no initialization is specified, use $X_{\mathrm{init}}=X$.

	\begin{definition}[Initialized observable trace equivalence]
		\label{dyn45:eq-trace}
		Let $\mathbf D$ have state space $X$, transition relation
		$E\subseteq X^2$, initial set $X_{\mathrm{init}}\subseteq X$,
		and marked observation $r:X\to Y$.
		Define its finite and infinite observable languages by
		\[
		L_{\mathbf D}(m):=
		\{(r(x_t))_{t=0}^{m}:
		(x_t)_{t=0}^{m}\in
		\mathsf T_{E,X_{\mathrm{init}}}(m)\},
		\qquad m\in\mathbb N_0,
		\]
		and
		\[
		L_{\mathbf D}(\infty):=
		\{(r(x_t))_{t\geq0}:
		(x_t)_{t\geq0}\in X^{\mathbb N_0},\
		x_0\in X_{\mathrm{init}},\
		(x_t,x_{t+1})\in E\ \forall t\geq0\}.
		\]
		Two systems $\mathbf D$ and $\mathbf D'$ of the same
		marked output type are:
		\begin{enumerate}[label=(\roman*),leftmargin=*,nosep]
			\item \textbf{Bounded-horizon trace equivalent}
			at horizon $m$ if
			$L_{\mathbf D}(m)=L_{\mathbf D'}(m)$.
			
			\item \textbf{Finite-trace equivalent} if
			$
			L_{\mathbf D}(m)=L_{\mathbf D'}(m)
			$
			for all $m \in \mathbb N_0$.
			
			\item \textbf{Infinite-trace equivalent} if
			$L_{\mathbf D}(\infty)=L_{\mathbf D'}(\infty)$.
		\end{enumerate}
		For continuous-time systems, the trace domain and
		observation schedule must also be specified.
	\end{definition}

	\begin{proposition}[Transport of observable languages]
		\label{dyn45:eq-trace-transport}
		A forward transition-preserving map $f:X\to X'$ with
		$f(X_{\mathrm{init}})\subseteq X'_{\mathrm{init}}$ and $r'f=r$
		satisfies $L_{\mathbf D}(m)\subseteq L_{\mathbf D'}(m)$ for every $m$.
		An observed dynamical isomorphism satisfying
		$f(X_{\mathrm{init}})=X'_{\mathrm{init}}$ gives equality. Thus
		finite-trace equivalence is invariant under the stated marked
		isomorphisms. A coordinatewise relabeling of $O$ transports the languages
		by that relabeling, rather than giving literal equality before transport.
	\end{proposition}

	\begin{proof}
		Apply $f$ to each state of a source path. The initial and transition
		conditions make the image a target path, and $r'f=r$ identifies its
		observations. Applying the same argument to the inverse of an
		isomorphism gives the reverse inclusion. Relabeling changes each
		entry of the observed word by the stated bijection.
	\end{proof}

	\begin{example}[Equal traces do not identify hidden states]
		\label{dyn45:eq-hidden-counterexample}
		Over the same retained realized packaging base, take one state $a$ with
		$F(a)=a$, or two states $b,c$ with $F(b)=c$ and $F(c)=b$. Observe the
		constant value $0$ in both systems and initialize all states. Their
		observable languages consist of the single word $(0,\ldots,0)$ at every
		horizon, and their infinite observable languages also agree. Their
		state sets have different cardinalities, so there is no dynamical
		isomorphism. The first system has a fixed state and the second has none.
		Consequently the cardinality of the native fixed-point set is not an
		invariant of this trace comparison. A contextual contract permitting a
		native state-cardinality or fixed-state-existence test separates them.
		This example uses the bare observed specialization; it does not impose
		joint faithfulness of the readout.
	\end{example}

	\begin{proposition}[When trace semantics suffices for contextual comparison]
		\label{dyn45:eq-trace-sufficiency}
		Suppose each permitted continuation has its output finite-trace semantics
		(the family of all finite-horizon observable languages)
		determined by the input finite-trace semantics and the fixed other inputs,
		including its failure outcome, and every terminal outcome is a function
		of that output finite-trace semantics. Then finite-trace equivalent inputs are
		contextually equivalent for this contract. If, in addition, the identity
		continuation permits membership tests for every finite observable word,
		the two equivalences coincide.
	\end{proposition}

	\begin{proof}
		Equality of input finite-trace semantics propagates through each finite
		expression by the assumed compositional rule; the test outcomes are
		therefore equal. Conversely, a word belonging to one horizon language
		and not the other is distinguished by its membership test in the
		identity continuation.
	\end{proof}
	Neither exact state coverage nor equality of one-step shadows supplies
	the compositional hypothesis of this proposition. In particular, one-step
	projections cannot replace the full languages. The non-globalization
	examples in the atlas sections explain the need for that distinction.
	An equality of all finite languages is not used as an assertion about
	infinite trajectories unless the requisite extension or compactness
	hypothesis has also been supplied.

	\begin{proposition}[Canonical deterministic observation quotient]
		\label{dyn45:eq-observation-quotient}
		For a total deterministic law $F:X\to X$ and a fixed readout $r:X\to O$,
		define
		\[
		x\sim_r y\quad\Longleftrightarrow\quad
		r(F^t x)=r(F^t y)\quad\text{for every }t\ge0.
		\]
		\begin{enumerate}[label=(\roman*),leftmargin=*,nosep]
			\item \emph{Canonical quotient.} This is the largest forward-invariant equivalence relation contained
			in $\ker r$. On $Q_r=X/{\sim_r}$ there are well-defined maps
			$\overline F[x]=[Fx]$ and $\overline r[x]=r(x)$. Writing $q_r:X\to Q_r$
			for the quotient, it satisfies $q_rF=\overline Fq_r$ and
			$\overline r q_r=r$. It is canonically isomorphic, as an observed
			deterministic system, to the set of complete histories
			$\{(r(F^t x))_{t\ge0}:x\in X\}$ with left shift and evaluation at zero.
			
			\item \emph{Universal property.} If $s:X\twoheadrightarrow Y$ is any surjective deterministic factor
			with $sF=Gs$ and $a s=r$, there is a unique surjective map
			$h:Y\to Q_r$ satisfying
			\[
			h s=q_r,\qquad hG=\overline Fh,\qquad \overline r h=a.
			\]
			Thus $Q_r$ is the coarsest deterministic state quotient retaining the
			specified readout. The statement concerns quotients of the fixed source
			system; it does not classify its hidden native structure.
		\end{enumerate}
	\end{proposition}

	\begin{proof}
		\leavevmode
		\begin{enumerate}[label={(iii)},leftmargin=*,nosep]
			\item[(i)] Equality of all history entries is an equivalence relation and includes
			equality at time zero. Shifting histories proves forward invariance and
			well-definedness of both displayed quotient maps. If an equivalence
			$\approx$ lies in $\ker r$ and is forward invariant, induction gives
			$F^t x\approx F^t y$ whenever $x\approx y$; hence $x\sim_r y$.
			The history map has kernel exactly $\sim_r$, so its induced map from
			$Q_r$ is bijective and intertwines shift and readout.
			
			\item[(ii)] For the universal
			property, $s(x)=s(y)$ implies $s(F^t x)=G^t s(x)=G^t s(y)=s(F^t y)$.
			Applying $a$ gives $x\sim_r y$. Define $h(s(x))=[x]$. This is
			well-defined, onto, and unique because $s$ is onto. The two remaining
			identities follow by evaluating at $s(x)$, which runs through $Y$. \qedhere
		\end{enumerate}
	\end{proof}

	The realized packaging base is retained in this quotient construction;
	any additional native-law requirement on the state space is checked
	separately. The group-homomorphic version is
	Theorem~\ref{dinv:thm:quotient}. Finite observation depth, the number of
	observable history classes and the trace-residual index measure different
	objects, and this proposition asserts no equality between them.

	\subsection{System invariants, conserved observables and invariant subsets}
	\label{dyn45:eq-invariant-kinds}

	Three uses of the word \emph{invariant} are kept separate.

	\begin{definition}[System invariants, conserved observables, and invariant sets]
		\label{dyn45:eq-invariant-definition}
		Fix an equivalence relation $\sim$ on a specified class $\mathscr D$.
		\begin{enumerate}[label=(\roman*),leftmargin=*,nosep]
			\item A \textbf{system invariant} is a map $I:\mathscr D\to R$ such that
			$\mathbf D\sim\mathbf D'$ implies $I(\mathbf D)=I(\mathbf D')$.
			If the values are themselves structured objects compared by isomorphism,
			use their declared isomorphism classes as $R$, or state that target
			comparison explicitly.
			
			\item For one relational system $(X,E)$, a \textbf{conserved observable} with
			fixed value set $A$ is a map $q:X\to A$ satisfying
			$q(x)=q(y)$ for every $(x,y)\in E$.
			
			\item A \textbf{forward-invariant subset}
			is $K\subseteq X$ such that $x\in K$ and $(x,y)\in E$ imply $y\in K$.
		\end{enumerate}
		For $E=\operatorname{Graph}(F)$, conservation and forward invariance are respectively $q\circ F=q$ and $F(K)\subseteq K$.
		They make no assertion that $F(x)=x$.
		In a controlled or nondeterministic viability problem, an existential choice of a continuing transition is a separate condition from this universal forward-invariance condition.
		Initialized conservation requires $q(x_t)=q(x_0)$ on every specified initialized trajectory and is stated separately whenever used.
	\end{definition}

	\begin{proposition}[Conjugacy transport and pullback]
		\label{dyn45:eq-system-transport}
		If $f:X\to X'$ preserves forward transitions, pullback by $f$ sends target conserved observables to source conserved observables and target forward-invariant subsets to source forward-invariant subsets.
		Under an observed dynamical isomorphism these assignments are bijections, with inverses induced by $f^{-1}$.
		For observables with values in a fixed field, pullback is a homomorphism of the pointwise first-integral algebras.
		It is injective when $f$ is surjective and an isomorphism under dynamical isomorphism.
		In a $C^r$ dynamical category these bijections
		restrict to $C^r$ conserved observables when the isomorphism and its inverse have the required regularity and intertwine every time map.
	\end{proposition}

	\begin{proof}
		For $(x,y)\in E$, forward preservation gives $(f(x),f(y))\in E'$.
		Thus $q'(f(x))=q'(f(y))$ for a conserved $q'$. If $f(x)\in K'$,
		invariance of $K'$ gives $f(y)\in K'$, proving invariance of
		$f^{-1}(K')$. For an isomorphism the same statements hold for its
		inverse and the composites are identities. Composition with $f$ preserves
		pointwise sums, products and constants. If $f$ is onto, equality after
		pullback implies equality at every target state. The regularity assertion follows from composition
		of $C^r$ maps and the time-intertwining identities.
	\end{proof}

	The temporal width $\operatorname{twit}$ is a system invariant for the
	fixed observable-language and atomization comparison of
	Proposition~\ref{wn15:prop:invariance}. A first integral is a function on
	the states of one system. The isomorphism type of a specified
	first-integral algebra can in turn be a system invariant, but its
	completeness requires a theorem for a declared comparison. Likewise a
	recovery rank, a minimum source size, or a dimension can be invariant
	without classifying its objects. Coordinate hiding gives a monotonicity
	statement for temporal width; it is not an invertible change of
	observation atoms.

	\subsection{Complete records, composition and classification}
	\label{dyn45:eq-classification-section}

	A complete record determines the equivalence class of a system. To use such a record in an assembly one also needs to know how its value behaves under composition.

	\begin{definition}[Complete and compositional invariants]
		\label{dyn45:eq-complete}
		\leavevmode
		\begin{enumerate}[label=(\roman*),leftmargin=*,nosep]
			\item \textbf{Completeness.}
			A system invariant $I$ for $\sim$ is \emph{complete} if
			\[
			I(\mathbf D)=I(\mathbf D')\quad\Longleftrightarrow\quad
			\mathbf D\sim\mathbf D'.
			\]
			
			\item \textbf{Compositionality for functions.}
			For a typed formation language, a family $I_\tau$ is
			\emph{compositional} if every permitted single-valued constructor $\alpha$ has a well-defined operation $\widehat\alpha$ on its attained input records, with the same definedness, satisfying
			\[
			I_\tau(\alpha(\mathbf D_1,\ldots,\mathbf D_m))
			=\widehat\alpha(I_{\tau_1}(\mathbf D_1),\ldots,
			I_{\tau_m}(\mathbf D_m)).
			\]
			
			\item \textbf{Compositionality for relations.}
			For each relational constructor $\alpha$, let $\mathcal R_\alpha(\mathbf D_1,\ldots,\mathbf D_m)$ denote its declared set of output records, with the stated failure convention.
			Compositionality requires $\mathcal R_\alpha(\mathbf D_1,\ldots,\mathbf D_m)=\mathcal R_\alpha(\mathbf D'_1,\ldots,\mathbf D'_m)$ whenever $I_{\tau_i}(\mathbf D_i)=I_{\tau_i}(\mathbf D'_i)$ for all $i$.
			It is a separate preservation property, not a consequence of having numerical invariants.
		\end{enumerate}
	\end{definition}

	\begin{proposition}[Attained-record classification and its universal property]
		\label{dyn45:eq-classification}
		Let $\mathscr D$ be a set of models of fixed comparison type in the chosen universe, let $\sim$ be an equivalence relation, and let $I:\mathscr D\to R$
		be complete.
		\begin{enumerate}[label=(\roman*),leftmargin=*,nosep]
			\item \emph{Classification.}
			Then
			\[
			\mathscr D/{\sim}\ \longrightarrow I(\mathscr D),\qquad
			[\mathbf D]\longmapsto I(\mathbf D)
			\]
			is a bijection.
			
			\item \emph{Factorization of invariants.}
			Every other $\sim$-invariant $J$ factors uniquely as
			$J=\overline J\,I$ on the attained record set.
			
			\item \emph{Minimality of complete records.}
			Conversely, if a summary $S$ merely has the property that equal summaries imply equivalence, then $I$ factors uniquely through the attained $S$-values; this does not require $S$ itself to be invariant.
			
			\item \emph{Composition.}
			For contextual equivalence under the single-valued substitution hypotheses of
			Proposition~\ref{dyn45:eq-congruence}, a complete family of invariants inherits well-defined constructor operations on attained records.
		\end{enumerate}
	\end{proposition}

	\begin{proof}
		\leavevmode
		\begin{enumerate}[label={(iii)},leftmargin=*,nosep]
			\item[(i)] Invariance makes the displayed map well-defined; completeness makes it
			injective, and attainment makes it surjective.
			
			\item[(ii)] Define
			$\overline J(I(\mathbf D))=J(\mathbf D)$; equality of the arguments implies equivalence and hence equality of the right sides.
			All values are attained, giving uniqueness.
			
			\item[(iii)] Similarly define the map from
			$S(\mathbf D)$ to $I(\mathbf D)$ using the implication in the statement.
			
			\item[(iv)] For a constructor, choose representatives of the attained input records.
			Different representatives are equivalent.
			The substitution proposition preserves definedness and the output equivalence class, so completeness
			makes the output record independent of those choices. \qedhere
		\end{enumerate}
	\end{proof}

	The restriction to a set, or to a fixed sufficiently large universe, is part of this formal classification statement. An unrestricted collection of arbitrary packaging contexts need not have a set of equivalence classes in the smaller universe.
	The proposition also gives no algorithm for computing $I$, deciding equality of its values, or constructing all
	attained records.
	Those are additional effectiveness questions.

	\begin{definition}[Three different classification tasks]
		\label{dyn45:eq-classification-tasks}
		A classification statement declares its models, retained markings, comparison relation, and output parameter space.
		The following outputs are distinguished:
		\begin{enumerate}[label=(\roman*),leftmargin=*]
			\item \textbf{Marked-realization classification:}
			a parametrization of actual marked realizations in a fixed ambient object, where distinct embedded sources remain distinct;
			
			\item \textbf{Equivalence-class classification:}
			a parametrization of classes under a specified isomorphism, conjugacy, or contextual equivalence;
			
			\item \textbf{Resource optimization:}
			an optimization value, such as minimum hidden dimension, together	with whatever existence and equality-case descriptions are proved.
		\end{enumerate}
		The value in (iii) alone need not be a complete invariant for either (i) or (ii);
		completeness requires a separate theorem.
		A count in (i) does not give a count in (ii) without identifying equivalent marked realizations; when equivalence is given by a group action, the orbit
		count must account for stabilizers.
		The word unique is qualified as unique marked realization, unique	equivalence class, or unique up to the stated isomorphism.
	\end{definition}

	There is also a different classwise classification problem: determine which covers have a property for every compatible atlas in a specified class.
	The sharp globalization criteria of
	Theorems~\ref{thm:all-horizon-causal-Horn-globalization} and~\ref{thm:pbr} classify covers in this sense, with their stated local model and recovery hypotheses.

	\begin{proposition}[A finite classification in a fixed coordinate category]
		\label{dyn45:eq-finite-classification}
		Fix a realized packaging base, finitely many named observation atoms $V$,
		finite nonempty alphabets $(A_v)_{v\in V}$, and a fixed named chart family.
		Let $A^V=\prod_v A_v$. Consider all systems specified by
		$\varnothing\ne X\subseteq A^V$, an initial subset $X_{\mathrm{init}}\subseteq X$, and
		$E\subseteq X^2$, with coordinate observations. Fix a subgroup
		$\Gamma\le\prod_v\operatorname{Sym}(A_v)$ of allowed alphabet
		relabelings, given by its permutations or a finite generating list and
		acting on all these data. The native base and the named
		atoms are fixed, and any additional finite markings must be encoded
		and transported as well, or their stabilizer must be used for $\Gamma$.
		Then this category has an effective complete isomorphism invariant.
		
		In fact, choose enumerations of the finite ambient product and encode
		$X$, $X_{\mathrm{init}}$ and $E$ by their indicator strings. The
		lexicographically least code among their $\Gamma$-translates is complete
		for the stated atom-preserving isomorphisms. Restricting $E$ to total
		single-valued graphs gives the deterministic specialization.
	\end{proposition}

	\begin{proof}
		There are finitely many ambient states, possible pairs and elements of
		$\Gamma$, so each code and its minimum are computable by finite
		enumeration. Isomorphic systems have the same orbit of codes. If two
		minimum codes agree, a relabeling carrying the first system to that
		code followed by the inverse relabeling for the second is an allowed
		isomorphism preserving every encoded datum. Graphs of total maps are
		preserved by these bijections. This proves completeness.
	\end{proof}
	The classification retains the entire native base and the declared observation granularity.

	\subsection{Coupling invariants, marked classification, and comparison contracts}
	\label{dyn46:subsec-coupling-invariants}

	The non-collapse results identify invariants of the coupled object. For
	example, for total deterministic observed dynamics the family of
	sites whose observations are autonomous is invariant up to site
	transport. Its proof is conjugation of the factor in
	Theorem~\ref{thm:exact-deterministic-chart-factor-criterion}; Example~\ref{dyn46:ex-coupled-noncollapse}
	shows this invariant does not factor through the separate base and
	unobserved evolution classes. The essential read-set family and
	carrier-support feasibility of Proposition~\ref{dyn46:prop-native-read}
	are additional comparison-relative invariants in their coordinate
	class. Admission preservation of Proposition~\ref{dyn46:prop-admission}
	and zero exit of Proposition~\ref{dyn46:prop-stochastic-admission}
	are invariants of the retained assembly and temporal law together.
	In each case isomorphism transports the native data on which the
	criterion is evaluated. These invariants may cease to be defined
	or invariant after a reduct that discards those data.

	\begin{proposition}[Complete finite marked coupling code]
		\label{dyn46:prop-coupling-code}
		Fix a finite coordinate base and finite labeled alphabets, a finite
		site and interaction library, and an effectively specified finite group $\Gamma$ of allowed
		bijective coordinate, alphabet, site and occurrence transports
		preserving the fixed native base. For finite relational coupled
		models with these libraries, encode the admissible states,
		initial set when retained, observations, interaction assignments,
		read/write markings, local channels, aggregate relation, and any
		retained finite interpretation and assembly data in a fixed ordered
		finite record. Include the composition rule if it is retained;
		assume all variable record entries have effective canonical finite
		encodings, decidable equality, and computable transport by the
		given group action. Then the lexicographically least record
		in the $\Gamma$-orbit is a complete invariant for precisely the
		isomorphisms induced by $\Gamma$. This extends
		Proposition~\ref{dyn45:eq-finite-classification} to interaction-marked
		models. For finite rational rate data the same assertion holds
		using a canonical encoding of reduced fractions, and likewise for
		finite rational Markov kernels, with all retained law data encoded.
	\end{proposition}

	\begin{proof}
		Apply every element of $\Gamma$ to the full record and take its
		minimum under the fixed total lexicographic order. The orbit is
		finite even when rational weights vary over infinitely many possible
		values. Isomorphic objects have equal record orbits. Equality of
		minimum records supplies two allowed transports to the same record;
		composing one with the inverse of the other gives an isomorphism
		because the encoding retains exactly the specified data. This gives an effective orbit test under the stated encoding hypotheses.
	\end{proof}
	A complete observed trace record cannot in general classify marked
	interaction decompositions: the equal-aggregate example in
	Proposition~\ref{dyn46:prop-coupled-category} has the same observed
	weighted process under every common initialization. A comparison
	contract that tests channel-specific interventions may distinguish
	the models, whereas tests of the aggregate process alone cannot.
	No particular intervention family is presumed admissible merely
	from the existence of the decomposition. This distinction governs later recovery and conservation results: each result must retain the law, observations and native markings used by its hypotheses. The gluing dimensions of Section~\ref{dyn47:subsec:weighted-compose} are further invariants under the specified native transports. The finest history-dependence partition
	and weighted rank in Definition~\ref{dyn48:def:wdrk} are further
	intrinsic invariants of the declared initialized observable comparison.

	\paragraph{Comparison contracts in the specialized theories.}

	The later results instantiate distinct contracts. Their completeness
	and classification assertions are read with the following comparisons.
	\begin{center}\small
		\begin{tabularx}{\linewidth}{@{}p{0.28\linewidth}X@{}}
			\toprule
			Theory and reference & Declared comparison and scope \\
			\midrule
			Temporal witness width, Section~\ref{wn09:sec:intrinsic-temporal-recovery}
			& Observable languages with fixed named atoms; the scalar width is an
			invariant, not a complete code for all such languages. \\
			Sharp globalization, Thms.~\ref{thm:all-horizon-causal-Horn-globalization}, \ref{thm:pbr}
			& Covers for which every compatible atlas in the declared class admits
			the required exact globalization; this classifies a universal cover
			property rather than individual dynamical isomorphism types. \\
			Temporal group records, Theorem~\ref{dinv:thm:complete}
			& Initialized group-native continuations and the stated recovery or
			numerical tests; the complete records retain actual boundary histories. \\
			Observed group quotient, Theorem~\ref{dinv:thm:quotient}
			& The coarsest surjective homomorphic dynamical factor retaining the
			specified readout; hidden native sources need not be isomorphic. \\
			\bottomrule
		\end{tabularx}
	\end{center}

	\subsection{Stationary states, exact periods, and their preservation}
	\label{dyn45:fp-section}

	The fixed-support axioms permit stationary states, cycles and deadlocks; they
	require none of them individually. Unless stated otherwise,
	the state space and the update law are fixed throughout time.

	\begin{definition}[Fixed states, periods, absorption, and closed walks]
		\label{dyn45:fp-stationary}
		\leavevmode
		\begin{enumerate}[label=(\roman*),leftmargin=*,nosep]
			\item \textbf{Deterministic fixed points and periods.}
			Let $F:X\to X$ be a total deterministic law.
			A state $x\in X$ is fixed if $F(x)=x$, and has exact
			period $n\geq1$ if $n$ is the least positive integer
			such that $F^n(x)=x$.
			Define
			\[
			\begin{aligned}
				\operatorname{Fix}(F)
				&:=\{x\in X:F(x)=x\},\\
				\operatorname{Per}_{\mid n}(F)
				&:=\{x\in X:F^n(x)=x\},\\
				\operatorname{Per}_{=n}(F)
				&:=\{x\in X:F^n(x)=x,\
				F^j(x)\ne x\ (1\leq j<n)\}.
			\end{aligned}
			\]
			Thus $\operatorname{Per}_{\mid n}(F)$ consists of
			states whose exact periods divide $n$.
			
			\item \textbf{Relational fixed and absorbing states.}
			Let $E\subseteq X^2$ be a transition relation and
			$E(x):=\{y\in X:(x,y)\in E\}$.
			A state $x$ is fixed if $x\in E(x)$, absorbing if
			$E(x)=\{x\}$, and deadlocked if $E(x)=\varnothing$.
			The fixed and absorbing state sets are
			\[
			\operatorname{Fix}(E):=\{x\in X:x\in E(x)\},
			\qquad
			\operatorname{Abs}(E):=\{x\in X:E(x)=\{x\}\}.
			\]
			A fixed state admits a constant infinite trajectory,
			whereas an absorbing state admits only that trajectory.
			Deadlocked states belong to neither set.
			
			\item \textbf{Pointed closed walks.}
			For $n\geq1$, the set of pointed closed walks of
			length $n$ is
			\[
			\operatorname{Cyc}_n(E):=
			\{(x_0,\ldots,x_{n-1})\in X^n:
			(x_i,x_{i+1})\in E\ (0\leq i<n)\},
			\]
			where $x_n:=x_0$.
			Repeated states are allowed, and cyclic rotations
			are not identified. In nondeterministic systems,
			a state may admit several return lengths without
			having an assigned deterministic period.
		\end{enumerate}
	\end{definition}

	\begin{proposition}[Preservation and reflection of return structures]
		\label{dyn45:fp-morphism}
		\leavevmode
		\begin{enumerate}[label=(\roman*),leftmargin=*,nosep]
			\item \emph{Relational laws.} A transition-preserving map $f:(X,E)\to(Y,E')$ sends
			$\operatorname{Fix}(E)$ into $\operatorname{Fix}(E')$ and induces
			$f^{\times n}:\operatorname{Cyc}_n(E)\to\operatorname{Cyc}_n(E')$ for every
			$n\ge1$. If $f$ is bijective and preserves and reflects transitions, these
			maps are bijective and $f(\operatorname{Abs}(E))=\operatorname{Abs}(E')$.
			In particular these assertions hold for a dynamical isomorphism with the
			additional structural and observation data retained.
			
			\item \emph{Deterministic laws.} For deterministic laws satisfying $fF=F'f$, the image of a point of exact
			period $n$ has exact period dividing $n$. If $f$ is injective, its exact
			period is $n$. Consequently $|\operatorname{Fix}(F)|$ and
			$|\operatorname{Per}_{=n}(F)|$ are dynamical-isomorphism invariants whenever
			these cardinalities are formed in the specified universe.
		\end{enumerate}
	\end{proposition}

	\begin{proof}
		\leavevmode
		\begin{enumerate}[label={(iii)},leftmargin=*,nosep]
			\item[(i)] Apply forward preservation to each edge of a self-loop or a closed walk.
			Under an isomorphism the same argument applies to $f^{-1}$, and equality of
			successor sets after transport also preserves absorption.
			
			\item[(ii)] In the deterministic
			case induction gives $fF^j=(F')^jf$. Hence an $n$-periodic point maps to a
			point returning after $n$ steps. If its least return is $d$, division of $n$
			by $d$ shows $d\mid n$. If $f$ is injective, a return after $d$ steps implies
			$F^d x=x$, so minimality gives $n=d$. \qedhere
		\end{enumerate}
	\end{proof}

	\begin{remark}[A fixed observable need not lift to a fixed state]
		\label{dyn45:fp-factor-counterexample}
		The constant map from the two-state swap to the one-state identity in
		Example~\ref{dyn45:eq-hidden-counterexample} is a surjective deterministic
		factor, exact at every trajectory horizon. Its unique fixed state has no
		fixed preimage. Thus neither trace equivalence nor surjective
		semiconjugacy reflects existence of an underlying fixed state.
	\end{remark}

	\begin{proposition}[Finite recurrence and deterministic period decomposition]
		\label{dyn45:fp-finite}
		Every trajectory of a total deterministic law on a finite nonempty $X$
		eventually enters a periodic orbit. A nonblocking relation on such an $X$
		has a directed cycle, but need not have a self-loop. For a deterministic
		law and every $n\ge1$,
		\[
		\operatorname{Per}_{\mid n}(F)
		=\coprod_{d\mid n}\operatorname{Per}_{=d}(F).
		\]
	\end{proposition}

	\begin{proof}
		Among the first $|X|+1$ states of a deterministic trajectory two coincide;
		the intervening block then repeats. In a finite nonblocking relation one
		can choose $|X|$ successive edges from any initial state, and repetition
		gives a closed walk containing a directed cycle. The swap in
		Example~\ref{dyn45:eq-hidden-counterexample} has no self-loop. The final
		identity follows from the divisibility argument in the preceding proof.
	\end{proof}

	\subsection{Local detection and exact globalization of stationary behavior}
	\label{dyn45:fp-globalization}

	A global fixed point gives a fixed point on every autonomous chart. Recovering a global fixed point from local ones requires their states to be compatible under projection.

	\begin{proposition}[Detection by jointly faithful autonomous observations]
		\label{dyn45:fp-local}
		Let $F:X\to X$ and $f_I:X_I\to X_I$ be total laws with
		$r_IF=f_Ir_I$. Suppose the tuple $(r_I)_I$ is injective. Then, for every
		$n\ge1$,
		\[
		\operatorname{Per}_{\mid n}(F)
		=\bigcap_I r_I^{-1}\bigl(\operatorname{Per}_{\mid n}(f_I)\bigr).
		\]
		In particular $x$ is fixed precisely when all its observations are fixed.
		This tests a state already in $X$; existence of arbitrary local fixed
		states does not assert that they have a common preimage in $X$.
	\end{proposition}

	\begin{proof}
		Iteration gives $r_IF^n=f_I^nr_I$. The forward implication follows
		immediately; for the reverse implication all observations of $F^n x$ and
		$x$ agree, and joint injectivity gives equality.
	\end{proof}

	\begin{proposition}[Canonical stationary and closed-walk gluing identities]
		\label{dyn45:fp-join}
		Fix a coordinate cover $\mathcal H$ of $V$, local states
		$Z_I\subseteq\prod_{v\in I}A_v$, local relations $R_I\subseteq Z_I^2$,
		and a declared global state set $X\subseteq M:=\Join_I Z_I$. Define
		\[
		E_{\max}(X)=\{(x,y)\in X^2:
		(\pi_Ix,\pi_Iy)\in R_I\text{ for every }I\}.
		\]
		\begin{enumerate}[label=(\roman*),leftmargin=*,nosep]
			\item Then
			\begin{equation}\label{dyn45:fp-join-fixed}
				\operatorname{Fix}(E_{\max}(X))
				=X\cap\Join_I\operatorname{Fix}(R_I).
			\end{equation}
			
			\item More generally, identifying a local $n$-tuple with a pattern on
			$I\times\{0,\ldots,n-1\}$, one has
			\begin{equation}\label{dyn45:fp-join-cycles}
				\operatorname{Cyc}_n(E_{\max}(X))
				=X^n\cap\Join_I\operatorname{Cyc}_n(R_I)\qquad(n\ge1).
			\end{equation}
			
			\item For any other global relation $E\subseteq X^2$ whose one-step projections
			are contained in the $R_I$, only the corresponding inclusions are automatic.
			These identities do not themselves assert existence, exact local projection,
			or nonblocking.
		\end{enumerate}
	\end{proposition}

	\begin{proof}
		\leavevmode
		\begin{enumerate}[label={(iii)},leftmargin=*,nosep]
			\item[(i)] The pair $(x,x)$ is in the canonical relation exactly when every
			$(\pi_Ix,\pi_Ix)$ is locally allowed, proving the first identity.
			\item[(ii)] For the second, apply the same defining condition to each consecutive
			pair, including $(x_{n-1},x_0)$. Because the charts cover $V$, agreement
			of local patterns on overlaps gives exactly one global $n$-tuple, which
			must additionally lie in $X^n$.
			\item[(iii)] Any relation with the stated local
			containments is a subset of $E_{\max}(X)$, giving the inclusions. \qedhere
		\end{enumerate}
	\end{proof}

	\begin{proposition}[Acyclic stationary gluing and its precise hypothesis]
		\label{dyn45:fp-tree}
		Suppose the finite coordinate cover has a join tree and
		$X=M$. Put $S_I=\operatorname{Fix}(R_I)$. If every $S_I$ is nonempty and
		the two projections of $S_I,S_J$ to $I\cap J$ agree on each tree edge,
		then $\operatorname{Fix}(E_{\max}(M))$ is nonempty and projects onto every
		$S_I$. The same conclusion holds for pointed length-$n$ closed walks if
		the nonempty relations $\operatorname{Cyc}_n(R_I)$ have equal separator
		projections on the time-expanded tree. Under an independent constraint
		$X\subsetneq M$, the additional requirement $x_t\in X$ must still be checked.
	\end{proposition}

	\begin{proof}
		Choose any local fixed state in any chosen root chart. The equality of
		separator projections permits extension to each adjacent chart, then
		successively through the tree. The running-intersection property ensures
		that a new chart agrees with all previously chosen charts wherever their
		coordinates coincide. The assembled state belongs to $\Join_I S_I$.
		The root and its chosen fixed state were arbitrary, proving exact
		projection as well as nonemptiness. Replace $I$ by
		$I\times\{0,\ldots,n-1\}$ and $S_I$ by the local closed-walk relation
		to obtain the second assertion.
	\end{proof}

	\begin{example}[Trace-compatible local equilibria can still disagree]
		\label{dyn45:fp-tree-counterexample}
		On the join-tree cover $\{\{a,s\},\{s,b\}\}$ use binary alphabets and
		all local states. Let
		\[
		F_1(a,s)=(a+s,s),\qquad
		F_2(s,b)=(s,b+1-s),
		\]
		where arithmetic is modulo two. Both separator trace languages consist
		of all constant $s$-words at every horizon, and the canonical global law
		is $F(a,s,b)=(a+s,s,b+1-s)$. Every local trajectory extends globally.
		However the first local law has fixed states only at $s=0$, and the
		second only at $s=1$. Thus both local fixed sets are nonempty and the
		atlas is all-horizon exact, yet there is no global fixed point. The
		stationary separator condition in Proposition~\ref{dyn45:fp-tree} is
		stronger than equality of ordinary separator trace languages.
	\end{example}

	\begin{example}[All-horizon exactness does not force stationary lifting]
		\label{dyn45:fp-trace-counterexample}
		Let $X=\{0,1\}^2$, use its two singleton charts, and let each local
		relation be the full binary relation. The relation
		$E=\{(x,y)\in X^2:x\ne y\}$ is state-exact and all-horizon exact:
		any specified word in one coordinate lifts by alternating the other
		coordinate at every step. It is also nonblocking and jointly observed
		faithfully. Nevertheless $\operatorname{Fix}(E)=\varnothing$, whereas
		both local fixed sets are $\{0,1\}$ and
		$E_{\max}(X)=X^2$ has every state fixed. Therefore even all-horizon
		exactness cannot replace the canonical-relation hypothesis in the
		equalities of Proposition~\ref{dyn45:fp-join}. It only guarantees lifts
		of local paths, not stationary lifts of constant local paths.
	\end{example}

	\subsection{Invariant sets as fixed points of state-set operators}
	\label{dyn45:fp-safety}

	The order-theoretic fixed-point principle is classical~\cite{Tarski1955}.
	A state-set fixed point lives in $\mathcal P(X)$, whereas an equilibrium	lives in $X$.
	Keeping these types separate is essential in the invariant core, closure, holonomy and source-selection applications below.

	\begin{definition}[Selectable viability, universal safety, and nonblocking safety]
		\label{dyn45:fp-safe-def}
		Let $E\subseteq X^2$ be a transition relation and
		$K\subseteq X$ a specified safe set.
		For each $x\in X$, let $E(x):=\{y\in X:(x,y)\in E\}$.
		Define three monotone operators on $\mathcal P(X)$ by
		\begin{align*}
			\Psi_{\exists}(S)
			&:=K\cap\{x\in X:E(x)\cap S\ne\varnothing\},\\
			\Psi_{\forall}(S)
			&:=K\cap\{x\in X:E(x)\subseteq S\},\\
			\Psi_{\forall+}(S)
			&:=K\cap\{x\in X:\varnothing\ne E(x)\subseteq S\}.
		\end{align*}
		A subset $S\subseteq K$ is called:
		\begin{enumerate}[label=(\roman*),leftmargin=*,nosep]
			\item \textbf{Viable} (or \textbf{selectably viable})
			if every $x\in S$ has at least one successor in $S$,
			equivalently, $S\subseteq\Psi_{\exists}(S)$.
			
			\item \textbf{Universally safe} if every successor
			of every $x\in S$ belongs to $S$,
			equivalently, $S\subseteq\Psi_{\forall}(S)$.
			
			\item \textbf{Universally safe and nonblocking}
			if every $x\in S$ has at least one successor and
			all its successors belong to $S$,
			equivalently, $S\subseteq\Psi_{\forall+}(S)$.
		\end{enumerate}
		Universal safety permits deadlocks, whereas viability
		and nonblocking safety exclude them within $S$.
	\end{definition}

	\begin{proposition}[Greatest safe kernels and finite computation]
		\label{dyn45:fp-kernels}
		Each of the three operators in Definition~\ref{dyn45:fp-safe-def} has a
		greatest fixed point, equal to the union of all its postfixed subsets
		$S\subseteq\Psi(S)$. When $X$ is finite, it is obtained by the descending
		iteration $S_0=K$, $S_{j+1}=\Psi(S_j)$, after at most $|K|$ strict deletions.
		For finite $X$ the three kernels have the following interpretations:
		\begin{enumerate}[label=(\roman*),leftmargin=*]
			\item the $\exists$ kernel is the set of starts of at least one infinite
			trajectory contained in $K$, equivalently states from which a directed
			cycle in the subgraph on $K$ is reachable within $K$;
			
			\item the $\forall$ kernel is the set of states for which every finite
			trajectory stays in $K$, including trajectories ending in deadlocks;
			
			\item the $\forall+$ kernel is the set of starts for which every reachable
			state (including the start itself) lies in $K$ and has a successor.
		\end{enumerate}
		For a total deterministic law all three coincide, even without finiteness,
		with the largest forward-invariant subset of $K$,
		\begin{equation}\label{dyn45:fp-deterministic-core}
			\operatorname{Core}_F(K)=\bigcap_{t\ge0}F^{-t}(K).
		\end{equation}
		None of these descriptions requires an equilibrium in $K$.
	\end{proposition}

	\begin{proof}
		For any of the three monotone operators, let $S_*$ be the union of its
		postfixed sets. Monotonicity gives $S_*\subseteq\Psi(S_*)$. Applying
		$\Psi$ shows that $\Psi(S_*)$ is itself postfixed, so maximality gives
		$\Psi(S_*)\subseteq S_*$. Every fixed point is postfixed, proving
		greatestness. Since $\Psi(K)\subseteq K$, iteration from $K$ descends;
		on a finite set each strict step deletes at least one element. Any
		postfixed set is contained in every iterate by induction, so the
		stabilized set is $S_*$.
		\begin{enumerate}[label={(iii)},leftmargin=*,nosep]
			\item[(i)] In the existential case successors can be chosen
			within this finite set forever; conversely the states of any infinite
			safe path form a postfixed set. Repetition on that path gives a reachable
			cycle.
			\item[(ii)] In the universal case close the initial state under all finite reachable paths. Safety of this reachable set is exactly the $\forall$ postfixed condition for it; conversely every state in a $\forall$ postfixed set has all reachable states in that set.
			\item[(iii)] For the nonblocking universal case the same reachable set is postfixed for $\forall+$ exactly when it is safe and each of its states has a successor. Conversely a $\forall+$ postfixed set has both properties for all starts in it.
		\end{enumerate}
		Finally a deterministic
		total state has one successor, so the three operators are identical;
		the intersection in \eqref{dyn45:fp-deterministic-core} is invariant
		and contains every invariant subset of $K$.
	\end{proof}

	\begin{remark}[Transport of safe kernels]
		\label{dyn45:fp-safe-transport}
		If a transition isomorphism $f$ carries $K$ onto $K'$, it intertwines
		each of the three state-set operators via $S\mapsto f(S)$, and hence
		carries the corresponding greatest kernels bijectively. A forward
		morphism alone does not preserve universal safety: the target may have
		additional unsafe outgoing transitions. No metric attraction,
		convergence rate, or stochastic stability follows from these order and
		reachability statements.
	\end{remark}

	\subsection{Weighted equilibria, transport, and fixed-point types}
	\label{dyn46:subsec-weighted-equilibria}
	
	For a Markov law the stationary object is a probability distribution.
	Summing that distribution over the fibers of a factor gives the candidate stationary distribution on the observed states.

	\begin{proposition}[Transport of stationary laws and stationary point masses]
		\label{dyn46:prop-stationary-transport}
		\leavevmode
		\begin{enumerate}[label=(\roman*),leftmargin=*,nosep]
			\item \emph{Transport of stationary laws.} A measurable Markov-kernel factor $f$ as in
			Definition~\ref{dyn46:def-typed-morphism} carries every stationary
			probability measure $\mu$ to a stationary probability $f_*\mu$.
			The same holds for a semigroup-stationary law, i.e. one invariant
			under every time map on measures. Isomorphisms transport the sets
			of such stationary laws bijectively.
			
			\item \emph{Stationary point masses.} For a finite jump generator,
			a point mass $\delta_x$ is stationary if and only if
			$Q(x,y)=0$ for all $y\ne x$; for a discrete Markov kernel the
			criterion is $P(x,\cdot)=\delta_x$, equivalently $P(x,\{x\})=1$ when the singleton is measurable.
		\end{enumerate}
	\end{proposition}

	\begin{proof}
		\leavevmode
		\begin{enumerate}[label={(iii)},leftmargin=*,nosep]
			\item[(i)] For each bounded measurable target test $g$, the factor identity
			implies $\int P'g\,d(f_*\mu)=\int P(g\circ f)\,d\mu
			=\int g\circ f\,d\mu$. The last integral equals
			$\int g\,d(f_*\mu)$, proving stationarity. Apply this at every time
			for a semigroup and apply the inverse for isomorphisms.
			
			\item[(ii)] In the
			finite generator case $\delta_xQ=0$ is exactly the vanishing of
			row $x$, and by the row-sum condition this is equivalent to zero
			exit rate. Exponentiation and differentiation at zero give the
			stationarity equivalence. For a kernel the image of a point mass
			is its row probability measure, giving the stated condition. \qedhere
		\end{enumerate}
	\end{proof}

	A positive self-loop probability alone is not a stationary point
	mass. A support-graph self-loop or a relational stationary state
	must therefore not be called a stochastic equilibrium without
	specifying this distinction. For a smooth complete flow generated
	by a vector field, an equilibrium is a zero of that field, hence
	fixed at all times (by uniqueness of solutions); for a general
	semiflow use the all-time fixed-point definition. For an interacting
	size-indexed family, zeros of the effective drift, invariant laws
	of the finite process, and nonanalyticities of limiting equilibrium
	quantities are different objects. The passage from microscopic
	stationary laws to macroscopic equilibria or stationary phase
	selection requires additional limit arguments. These limit problems require hypotheses beyond the fixed-point axioms and the support reduct.

	\paragraph{Fixed-point types in the specialized theories.}

	The preceding definitions distinguish several constructions used later:
	\begin{enumerate}[label=(\roman*),leftmargin=*]
		\item a dynamical equilibrium satisfies $F(x)=x$ in the state space;
		
		\item a closure-fixed history satisfies $C(H)=H$ in a space of assertion
		sets, with $C$ an inference closure, not a temporal update;
		
		\item a holonomy-fixed label is fixed by every transport permutation,
		as in \eqref{wn11:eq:common-fix}; it need not be fixed by the temporal
		update on those labels;
		
		\item an invariant native subgroup under a reversible update is a fixed object of
		the induced action on subgroups, not a subgroup consisting
		only of fixed states;
		
		\item the invariant cores of Section~\ref{dinv:sec:records} are greatest
		fixed points of a state-set operator, and may contain nonconstant orbits.
	\end{enumerate}
	For nonautonomous laws $F_t:X\to X$ on a common state set, a constant
	trajectory must start in $\bigcap_t\operatorname{Fix}(F_t)$; a periodic
	schedule instead admits a return-map analysis with its schedule explicitly
	retained. If the state types $X_t$ change, the equality $x_{t+1}=x_t$ is
	not well typed without specified identifications or structural transport.
	Finally, a continuous-time equilibrium is fixed by every time map of its
	flow; being fixed by one time map can instead signify a periodic orbit.
	No general equilibrium-existence axiom is added in any of these cases.

	\section{Local dynamical atlases and exact trajectory realization}
	\label{sec:local-dynamical-atlases}

	A transition atlas prescribes local states and allowed steps before a global evolution is chosen. Throughout the relational atlas results, the finite indexed chart family is nonempty (individual supports may be empty), and local state sets are nonempty unless expressly stated otherwise. We first define compatible-state envelopes and complete path languages, then characterize their exact global realization. Lemma~\ref{wn05-06:lem:projection-composition} identifies a basic obstruction: projected edges may require different witnesses at their shared intermediate state. Theorem~\ref{thm:greatest-relation} gives the greatest respecting relation; Theorem~\ref{thm:acyclic-temporal-exactness-criterion} reduces exactness on a join tree to agreement of separator path languages.

	\subsection{Intrinsic atlas-first data and three exactness levels}
	\label{subsec:v07-atlas-semantics}

	Fix the observed carrier family $\mathcal K$ of Definition~\ref{def:incidence-charts-v06}; its indexed supports $I_K$ form the \emph{intrinsic} cover on $V_{\mathcal K}$.
	All constructions below also apply to an explicitly declared auxiliary descriptor cover, but are not then called incidence-generated.
	To retain chart provenance, duplicate supports remain indexed, even when their restrictions coincide.

	\begin{definition}[Atlas-first presentation and compatibility conditions]
		\label{def:v07-atlas}
		Fix an indexed cover $\mathcal H$ of $V$ and finite
		nonempty alphabets $A_v$, and write
		$A^I:=\prod_{v\in I}A_v$.
		An \textbf{atlas-first presentation} is an indexed family
		\[
		\mathcal A=(Z_I,R_I)_{I\in\mathcal H},
		\qquad
		\varnothing\ne Z_I\subseteq A^I,\quad R_I\subseteq Z_I^2,
		\]
		where $Z_I$ is the prescribed local state set and $R_I$
		is the prescribed local transition relation.
		No global realization is assumed, and the relations
		$R_I$ may be empty.		
		The atlas $\mathcal A$ is called:
		\begin{enumerate}[label=(\roman*),leftmargin=*,nosep]
			\item \textbf{State compatible} if
			\[
			\pi_S Z_I=\pi_S Z_J.
			\]
			
			\item \textbf{One-step compatible} if
			\[
			(\pi_S\times\pi_S)R_I
			=(\pi_S\times\pi_S)R_J.
			\]
			
			\item \textbf{Horizon-$m$ trace compatible} if
			\[
			\pi_{S\times[0,m]}\mathsf T_{R_I}(m)
			=\pi_{S\times[0,m]}\mathsf T_{R_J}(m).
			\]
		\end{enumerate}
		Each equality is required for every pair of chart
		occurrences $I,J\in\mathcal H$, with $S=I\cap J$.
		Compatibility on a specified edge set requires the
		corresponding equality only for those chart pairs.
	\end{definition}

	\begin{definition}[Global realization, exactness, and nonblocking]
		\label{def:v07-exact}
		Let $\mathcal A=(Z_I,R_I)_{I\in\mathcal H}$ be an
		atlas-first presentation.
		Its \textbf{canonical compatible-state space} is
		\[
		M:=\Join_{I\in\mathcal H}Z_I
		=\{x\in A^V:\pi_I x\in Z_I
		\text{ for every }I\in\mathcal H\}.
		\]
		A \textbf{global realization} of $\mathcal A$ is a pair
		$(X,E)$, where $X\subseteq M$ is an admissible state
		set and $E\subseteq X^2$ is a transition relation.
		The canonical state-space choice is $X=M$.		
		The global realization $(X,E)$ is called:
		\begin{enumerate}[label=(\roman*),leftmargin=*,nosep]
			\item \textbf{State exact} if
			\[
			\pi_I X=Z_I
			\qquad(I\in\mathcal H).
			\]
			
			\item \textbf{One-step exact} if
			\[
			(\pi_I\times\pi_I)E=R_I
			\qquad(I\in\mathcal H).
			\]
			
			\item \textbf{All-horizon exact} if
			\[
			\pi_{I\times[0,m]}\mathsf T_E(m)
			=\mathsf T_{R_I}(m)
			\]
			for every $I\in\mathcal H$ and $m\in\mathbb N_0$.
			
			\item \textbf{Nonblocking} if
			\[
			\forall x\in X,\quad
			\exists y\in X:\ (x,y)\in E.
			\]
		\end{enumerate}
		All-horizon exactness includes state and one-step
		exactness at horizons $0$ and $1$, respectively.
		Nonblocking is independent of the exactness requirements.
	\end{definition}

	\begin{remark}[No circular global object]\label{rem:v07-no-circularity}
		The join $M$ is computed from local data; it is not postulated as a prior physical state space. It may be empty. The pair $(X,E)$ is a candidate realization until the nonempty-state requirement (D1) is met; an empty envelope records a realization obstruction. Independent constraints $X\subsetneq M$ must be declared before the realization problem is asked, because they can remove even locally valid one-step transitions. The general exactness definitions retain such $X$, while cover-only universal globalization statements use $X=M$ unless expressly stated otherwise.
	\end{remark}

	\noindent\textbf{Presentation convention.} Unless explicitly labeled realized,
	this section uses the atlas-first presentation in
	Definition~\ref{def:two-modes-v06}; the chart cover may be the
	incidence-generated family of Definition~\ref{def:incidence-charts-v06}
	or an explicitly declared auxiliary cover.
	
	In the atlas-first notation of the next subsection, write $X_I=Z_I$ for the prescribed local state set, not an assumed projection of a global state space. For $S\subseteq I$, its attained local restriction is $\pi_S^I X_I$. Let
	\[
	\Tr_I(m)=\{(x_0,\ldots,x_m)\in X_I^{m+1}:(x_t,x_{t+1})\in R_I\text{ for }0\le t<m\}.
	\]
	For overlaps $S=I\cap J$, define the separator trace language $L_I^S(m)=\pi_{S\times[0,m]}\Tr_I(m)$. Equality of one-step separator relations does not generally imply equality of complete separator trajectories.
	We call an atlas \emph{all-horizon trace compatible} when $L_I^S(m)=L_J^S(m)$ for every pair, every overlap, and every $m\ge0$.

	\subsection{Local atlas semantics and restriction compatibility}
	\label{wn05-06:sec:local-dynamic-atlases}

	Independent charts may agree on their states and edges but disagree on their longer paths. We use separate projection tests for the three kinds of data.

	\subsubsection{Finite supported transition data}
	\label{wn05-06:subsec:transition-data}

	We retain the provenance convention of Section~\ref{sec:fixed}: a realized packaging system
	$\mathbb P=(\Xi,P,\mathfrak P)$ satisfies the primitive packaging axioms
	(S1)--(S5) and the contextual formation requirements of
	\cite{MaPackaging2026}. The present descriptor-label cover is additional atlas
	data, not an identification of chart domains with primitive package interiors.
	We initially fix a finite set $V$ of descriptor labels, a nonempty finite
	alphabet $A_v$ for every $v\in V$, and a finite cover
	$\mathcal H=(I_j)_{j\in J}$ of $V$. We write
	$A^I=\prod_{v\in I}A_v$ and use $\pi_S^I$ for the coordinate projection
	$A^I\to A^S$ whenever $S\subseteq I$.

	\begin{ddefinition}[Nonempty-relation transition atlases, compatibility, and join trees]
		\label{wn05-06:def:transition-atlas}
		Fix a finite indexed cover $\mathcal H$ of $V$ and
		finite nonempty coordinate alphabets $A_v$.
		Write $A^I:=\prod_{v\in I}A_v$ and let
		$\pi_S^I:A^I\to A^S$ denote restriction to $S\subseteq I$.
		\begin{enumerate}[label=(\roman*),leftmargin=*,nosep]
			\item \textbf{Nonempty-relation transition atlas.}
			A nonempty-relation transition atlas on $\mathcal H$
			is an indexed family
			\[
			\mathcal A=(X_I,R_I)_{I\in\mathcal H},
			\qquad
			\varnothing\ne X_I\subseteq A^I,\quad
			\varnothing\ne R_I\subseteq X_I^2.
			\]
			Here $X_I$ is the local state space and $R_I$
			is the local transition relation.
			For $S\subseteq I$, the relational restriction is
			\[
			R_I|_S:=
			\{(\pi_S^Iu,\pi_S^Iv):(u,v)\in R_I\}
			\subseteq (A^S)^2.
			\]
			
			\item \textbf{Compatibility.}
			The atlas $\mathcal A$ is state-compatible if
			$
			\pi_S^I X_I=\pi_S^J X_J,
			$
			and one-step compatible if
			$
			R_I|_S=R_J|_S,
			$
			for every pair of charts $I,J\in\mathcal H$,
			where $S=I\cap J$.
			
			\item \textbf{Join tree and separators.}
			A join tree of $\mathcal H$ is a tree $T$
			whose nodes are the original chart occurrences
			and which satisfies the running-intersection condition:
			for every $v\in V$, the nodes
			$\{I\in\mathcal H:v\in I\}$ induce a connected subtree.
			The separator of an edge $IJ$ of $T$ is
			$S_{IJ}:=I\cap J$.
		\end{enumerate}
		Nonemptiness of the local relations is specific to this
		specialization.
		No autonomy of the intersection	observations is assumed.
	\end{ddefinition}

	A transition relation may be specified implicitly by finitely many supported assertions:
	permitted or forbidden current--next patterns on finite subsets of the chart.
	Intersecting these requirements inside a chart defines its $R_I$. Existence of a global evolution is a further question: an
	atlas of individually valid local transition rules is not by definition one global dynamical law.
	The framework therefore includes constrained systems, automata, relational control laws, higher-order incidence dynamics, and discrete models with hidden inputs.

	\begin{ddefinition}[State envelope, natural-join law, and exactness]
		\label{wn05-06:def:join-transition}
		For a transition atlas $\mathcal A=(X_I,R_I)_{I\in\mathcal H}$:\begin{enumerate}[label=(\roman*),leftmargin=*,nosep]
			\item The \textbf{state envelope} is
			\[
			M_{\mathcal H}=\Bigl\{x\in A^V:
			\pi_I^Vx\in X_I\ \text{for all }I\in\mathcal H\Bigr\}.
			\]
			
			\item The \textbf{natural-join global relation} is
			\begin{equation}\label{wn05-06:eq:global-join-relation}
				R_{\mathcal A}=\Bigl\{(x,y)\in M_{\mathcal H}^{2}:
				(\pi_I^Vx,\pi_I^Vy)\in R_I\ \text{for all }I\in\mathcal H\Bigr\}.
			\end{equation}
			
			\item An \textbf{exact one-step relational globalization} is any
			$R\subseteq M_{\mathcal H}^{2}$ satisfying $(\pi_I^V\times\pi_I^V)(R)=R_I$ for every original chart. State exactness on the envelope is the separate requirement $\pi_I^V M_{\mathcal H}=X_I$.
		\end{enumerate}
	\end{ddefinition}

	Here $M_{\mathcal H}$ is the envelope $M$ of Definition~\ref{def:v07-exact}, and $R_{\mathcal A}$ is its greatest respecting relation (Theorem~\ref{thm:greatest-relation}).
	Under independent global state constraints $X\subsetneq M_{\mathcal H}$, the greatest respecting relation on $X$ is $R_{\mathcal A,X}=R_{\mathcal A}\cap(X\times X)$.
	Exactness of its chart projections remains to be tested.

	\subsubsection{Temporal assertions, projection, and trace defects}
	\label{wn05-06:subsec:time-assertions}

	An observed path records every time slice.
	The intermediate state must witness both adjacent steps, which is the constraint lost when the steps are projected separately.

	\begin{ddefinition}[Finite trajectories and projected trace languages]
		\label{wn05-06:def:trace-language}
		Let $R_I\subseteq X_I^2$ be a local transition relation.
		For $S\subseteq I$, let $\pi_S^I:A^I\to A^S$ denote
		the coordinate projection.
		\begin{enumerate}[label=(\roman*),leftmargin=*,nosep]
			\item \textbf{Finite trajectories.}
			A trajectory of horizon $m\in\mathbb N_0$ is a
			sequence of $m+1$ states connected by $m$ permitted
			transitions. The finite trajectory set is
			\[
			\mathsf T_I(m):=
			\{(u_0,\ldots,u_m)\in X_I^{m+1}:
			(u_t,u_{t+1})\in R_I\ (0\leq t<m)\}.
			\]
			In particular, $\mathsf T_I(0)=X_I$ and
			$\mathsf T_I(1)=R_I$, under the natural identifications.
			
			\item \textbf{Trace projection and horizon trace language.}
			The trace projection onto $S$ maps each trajectory
			$(u_t)_{t=0}^m$ to $(\pi_S^I u_t)_{t=0}^m$.
			The horizon-$m$ trace language on $S$ is its image
			\[
			\mathsf L_I^S(m):=
			\{(\pi_S^I u_t)_{t=0}^m:
			(u_t)_{t=0}^m\in\mathsf T_I(m)\}
			\subseteq (A^S)^{m+1}.
			\]
			
			\item \textbf{Full finite-trace language.}
			The full finite-trace language on $S$ is
			\[
			\mathsf L_I^S:=
			\bigcup_{m\in\mathbb N_0}\mathsf L_I^S(m)
			\subseteq (A^S)^+,
			\]
			where $(A^S)^+$ denotes the set of nonempty finite
			words over $A^S$. For residual-language constructions,
			we adjoin the empty word $\varepsilon$.
		\end{enumerate}
		A finite trace represents a temporal assertion whose
		realizability requires a trajectory lift on $I$.
		Trace projection existentially eliminates the
		coordinates in $I\setminus S$.
	\end{ddefinition}

	These are trajectory assertions: recording the projected start and finish does not record the hidden intermediate states.
	In particular, relational composition and chart elimination do not generally commute.

	\begin{dlemma}[A temporal projection obstruction]\label{wn05-06:lem:projection-composition}
		For relations $R,S\subseteq X_I\times X_I$ and $U\subseteq I$,
		\[
		(S\circ R)|_U\ \subseteq\ (S|_U)\circ(R|_U),
		\]
		where composition is right to left: $(u,w)\in S\circ R$ means $\exists v:(u,v)\in R,(v,w)\in S$.
		The inclusion can be strict, even when $R=S$ is the graph of a total
		bijective update.
	\end{dlemma}

	\begin{proof}
		A hidden intermediate state $v$ witnessing a full-chart composition has an observable restriction $\pi_Uv$ witnessing the projected composition.
		For	strictness, take $X_I=\{0,1\}^2$, let $R=S$ be the graph of the swap $(a,b)\mapsto(b,a)$, and observe only coordinate $b$.
		Its one-step projected relation is full, so its square is full.
		But the projected two-step relation is the identity relation, since the swap squares to the identity.
	\end{proof}

	\begin{ddefinition}[Temporal exact globalization]
		\label{wn05-06:def:temporal-exact}
		For an atlas $\mathcal A=(X_I,R_I)_{I\in\mathcal H}$, let
		$R\subseteq M_{\mathcal H}^2$.
		It is an \textbf{exact globalization through horizon $m$} ($m\in\mathbb N_0$) if
		\[
		\pi_{I\times[0,h]}\mathsf T_R(h)=\mathsf T_I(h)
		\quad(I\in\mathcal H,\ 0\le h\le m).
		\]
		It is an \textbf{all-horizon exact globalization} if this holds for every $m\in\mathbb N_0$.
		The cases $h=0,1$ impose exact state and one-step realization, respectively.
	\end{ddefinition}

	\paragraph{Temporal trace defects and local recovery data.}	
	For a proposed join tree $T$ define, for an edge $I\!\! -\!\!J$ with separator $S=I\cap J$,
	\begin{equation}\label{wn05-06:eq:edge-trace-defect}
		\delta_{IJ}=\min\{m\ge0:\mathsf L_I^S(m)\ne\mathsf L_J^S(m)\}
		\ \in\mathbb N_0\cup\{\infty\},
	\end{equation}
	where the minimum of the empty set is $\infty$. Write
	$\delta_T(\mathcal A)=\min_{IJ\in E(T)}\delta_{IJ}$, with value $\infty$ for a one-chart atlas.
	This is a \emph{temporal trace-defect depth} of the
	specified transition atlas and the chosen join tree.
	It is not the intrinsic propagation--recovery rank $\operatorname{ipr}$ of \cite{MaPackaging2026}, nor the temporal observation depth of Section~\ref{sec:fixed}.
	A finite $\delta_{IJ}$ is accompanied by a distinguishing separator word,	thereby furnishing an explicit local temporal incompatibility certificate.

	\subsection{Existence on a prescribed state space}
	\label{subsec:v07-canonical}

	Fix a set $X$ of admissible assemblies.
	A missing local path can be certified by checking a single relation on $X$, even when many relations have the same one-step chart shadows.

	\begin{equation}\label{eq:canonical}
		\boxed{R_*=\{(x,y)\in X^2:(\pi_I x,\pi_I y)\in R_I\ \text{for every }I\in\mathcal H\}.}
	\end{equation}

	\begin{theorem}[Greatest relation and exact path realization]
		\label{thm:greatest-relation}
		For an atlas as in Definition~\ref{def:v07-atlas}, let $X\subseteq M$ be the stipulated subset of the compatible-state envelope.
		\begin{enumerate}[label=(\roman*),leftmargin=*,nosep]
			\item \emph{Greatest respecting law.} The relation $R_*$ in \eqref{eq:canonical} is the greatest relation on $X$ that respects every chart edge.
			
			\item \emph{Exactness criterion.} An exact realization of all finite trajectories on $X$ exists if and only if
			\[
			\pi_{I\times[0,m]}\mathsf T_{R_*}(m)=\mathsf T_I(m)
			\qquad(I\in\mathcal H,\ m\ge0).
			\]
			One-step exactness is the equality at $m=1$. State exactness is the separate equality at $m=0$.
			
			\item \emph{Finite decision.} For effectively given finite state spaces and transition tables, an automaton decides all-horizon exactness and gives a finite missing path on failure.
		\end{enumerate}
	\end{theorem}

	\begin{proof}
		\leavevmode
		\begin{enumerate}[label={(iii)},leftmargin=*,nosep]
			\item[(i)] Every respecting relation $R$ is contained in $R_*$ by the defining edge tests.
			
			\item[(ii)] On each chart and at each horizon one therefore has
			$\pi_I\mathsf T_R(m)\subseteq\pi_I\mathsf T_{R_*}(m)\subseteq\mathsf T_I(m)$.
			The second inclusion holds edge by edge, and at $m=0$ by $X\subseteq M$. If $R$ is exact, both inclusions are equalities. Conversely, the displayed equalities make $R_*$ itself a realization.
			
			\item[(iii)] In finite state spaces the two path languages are regular. Testing their difference for emptiness decides the condition and returns a witnessing word when it fails. \qedhere
		\end{enumerate}
	\end{proof}

	\begin{remark}[Universal extension]\label{rem:v07-UG-types}
		For a specified class $\mathscr D$, $\operatorname{UG}_{\mathscr D}(\mathcal H)$ quantifies over every compatible atlas in that class under the same global-state convention. Prescribing a smaller $X$ can exclude even a one-chart state.
	\end{remark}

	\subsection{Exact globalization of finite transition systems}
	\label{wn05-06:sec:exact-temporal-globalization}

	On a join tree, attaching a chart requires matching it only on its separator with its parent. We apply this extension first to states and edges, then to complete finite paths.

	\subsubsection{One-step gluing on acyclic covers}
	\label{wn05-06:subsec:one-step-relation}

	At a fixed step one can root a join tree at any prescribed chart transition. Equality on each edge of the tree extends that transition out to all charts.

	\begin{dtheorem}[Exact state and one-step globalization on a join tree]\label{wn05-06:thm:join-tree-one-step}
		Suppose $\mathcal H$ admits a join tree $T$ and $\mathcal A$ is state-compatible
		and one-step compatible on its edges, with nonempty local state sets and possibly empty transition relations. Then $M_{\mathcal H}$ is nonempty
		and satisfies $\pi_I M_{\mathcal H}=X_I$ for every $I$. Moreover,
		$R_{\mathcal A}|_I=R_I$ for every $I$, and
		$R_{\mathcal A}$ contains every other global relation with these same exact
		one-step restrictions. The global relation is nonempty if and only if every local transition relation is nonempty.
	\end{dtheorem}

	\begin{proof}
		Root the join tree. Fix any permitted state $u$ on the root chart.
		For a neighboring chart $J$, state compatibility on its separator with the parent guarantees a matching extension. Proceed recursively along the tree. The running-intersection property guarantees that a newly attached chart shares variables with previously attached charts only through its parent separator. The assembled state belongs to $M_{\mathcal H}$ and realizes the prescribed datum. Any chart can serve as root, proving exact state projection.
		For transition tuples the same argument applies to the doubled coordinates, using one-step compatibility. If one local transition relation is empty, its separator projection is empty (even on the empty separator), so emptiness propagates to every chart along the tree. In that case $R_{\mathcal A}=\varnothing$ and every one-step projection is still exact. Otherwise a prescribed root transition extends to a global transition as above, proving exactness and nonemptiness. Greatestness follows directly from Definition~\ref{wn05-06:def:join-transition}.
	\end{proof}

	This is the classical acyclic relational join theorem of Beeri--Fagin--Maier--Yannakakis~\cite{wn09:BeeriFaginMaierYannakakis1983}, applied to one-step transition data. By Definition~\ref{def:nu}, $\nu(\mathcal H)=\infty$ implies chordality and conformality, and hence
	a join tree (allowing redundant charts). It is precisely at longer horizons
	that a new dynamical compatibility question appears.

	\subsubsection{All-horizon trace globalization}
	\label{wn05-06:subsec:all-horizon}

	A whole path can be treated as one tuple whose entries are time strings.
	This changes the separator test from agreement on a pair of successive states to agreement on all complete words.

	\begin{dtheorem}[Acyclic temporal exactness criterion]
		\label{thm:acyclic-temporal-exactness-criterion}
		Let $\mathcal H$ admit a join tree $T$, and let $\mathcal A$ be a finite transition atlas with nonempty local state sets and possibly empty transition relations. The following are equivalent (the conditions at $m=0,1$ already impose state and one-step compatibility):
		\begin{enumerate}[label=(\roman*),leftmargin=2.3em]
			\item there exists $R\subseteq M_{\mathcal H}^{2}$ whose projected finite
			trajectory languages equal $\mathsf T_I(m)$ on every original chart
			$I$ and at every finite horizon $m$;
			
			\item for every edge $IJ$ of $T$ and every $m\ge0$, the separator path
			languages coincide:
			\begin{equation}\label{wn05-06:eq:edge-all-horizon}
				\mathsf L_I^{I\cap J}(m)=\mathsf L_J^{I\cap J}(m);
			\end{equation}
			
			\item $\delta_T(\mathcal A)=\infty$.
		\end{enumerate}
		When these conditions hold, the canonical greatest one-step extension
		$R_{\mathcal A}$ is an exact all-horizon globalization. Independently of these conditions, the canonical path relation satisfies
		\begin{equation}\label{wn05-06:eq:path-join}
			\mathsf T_{R_{\mathcal A}}(m)
			=\left\{(x_0,\ldots,x_m)\in(A^V)^{m+1}:
			(\pi_I x_0,\ldots,\pi_I x_m)\in\mathsf T_I(m)
			\ \forall I\in\mathcal H\right\}.
		\end{equation}
	\end{dtheorem}

	\begin{proof}
		\leavevmode
		\begin{enumerate}[label={},leftmargin=0pt,nosep]
			\item \textnormal{(i)$\Rightarrow$(ii).}\enspace For (i)$\Rightarrow$(ii), projection of any common global path set onto
			$S=I\cap J$ is independent of whether the projection is factored through
			$I$ or $J$. The two exact chart-path projections therefore have equal
			separator projections.
			
			\item \textnormal{(ii)$\Leftrightarrow$(iii).}\enspace Equivalence of (ii) and (iii) is Definition
			\eqref{wn05-06:eq:edge-trace-defect}.
			
			\item \textnormal{(ii)$\Rightarrow$(i).}\enspace For (ii)$\Rightarrow$(i), fix $m$. Regard a chart trajectory of length
			$m+1$ as a relation on the time-indexed coordinate set
			$I\times\{0,\ldots,m\}$; each such relation is $\mathsf T_I(m)$.
			The running-intersection property of $T$ persists after multiplying every
			chart by the same set of times. Equality~\eqref{wn05-06:eq:edge-all-horizon}
			expresses exact relational consistency on each tree separator. If any time-expanded relation is empty, its separator projection is empty, and equality propagates emptiness along the connected tree. Then every local relation and their join are empty, so exactness holds (including for a single chart). Otherwise all time-expanded relations are nonempty.
			The constructive join-tree argument in
			Theorem~\ref{wn05-06:thm:join-tree-one-step} therefore shows that the natural join of
			these time-expanded relations has exact projection onto every $\mathsf T_I(m)$. At $m=0$ the join is $M_{\mathcal H}$. At $m\ge1$ a global time word belongs to this join precisely when every consecutive global pair lies in $R_{\mathcal A}$, establishing~\eqref{wn05-06:eq:path-join} at all horizons.
			Thus $R_{\mathcal A}$ realizes every prescribed local trajectory at every
			horizon. The same global relation works for all $m$. \qedhere
		\end{enumerate}
	\end{proof}

	The join argument is the acyclic relational extension theorem applied to complete paths \cite{wn09:BeeriFaginMaierYannakakis1983}. Passing from paths to an effective separator test requires a finite state representation.

	\subsubsection{A finite temporal obstruction test}
	\label{wn05-06:subsec:finite-test}

	Finite local state sets give finite automata on the shared readings. A shortest distinguishing word provides both a counterexample and its first horizon.

	\begin{dtheorem}[Finite separator test and a bound for the first disagreement]\label{wn05-06:thm:finite-trace-test}
		Let $I,J$ have separator $S=I\cap J$ and nonempty finite local state sets of sizes $n_I,n_J$. One can decide whether their separator traces agree at every horizon. If they disagree, the first horizon at which they disagree satisfies
		\begin{equation}\label{wn05-06:eq:finite-depth-bound}
			\delta_{IJ}\le 2^{n_I}+2^{n_J}-2.
		\end{equation}
		On a finite join tree it suffices to apply the test on the tree edges.
	\end{dtheorem}

	\begin{proof}
		For $I$, use a fresh start state and one state for every subset of $X_I$. On the first separator symbol $s$, move to $\{u\in X_I:\pi_Su=s\}$. From a subset $U$, the transition on $s$ leads to $\{v:\pi_Sv=s,\ (u,v)\in R_I\text{ for some }u\in U\}$. Exactly the nonempty subsets accept. This machine recognizes the nonempty separator words of $I$. Construct the same machine for $J$.
		
		Take the disjoint union of the two machines and identify their rejecting empty-subset sinks. There are at most $D=1+2^{n_I}+2^{n_J}$ states. Partition the states first by acceptance (including the rejecting fresh starts), and then refine by the blocks of their successors for each symbol. Induction shows that the partition after $h$ refinements identifies exactly those states that agree on all words of length at most $h$. A stable partition never changes. Starting with two nonempty blocks allows at most $D-2$ strict refinements. Thus distinguishable start states admit a distinguishing word of length at most $D-2$. Subtracting one converts the word length to the horizon in \eqref{wn05-06:eq:finite-depth-bound}. A breadth-first search of pairs produces a shortest witness. The join-tree conclusion follows from Theorem~\ref{thm:acyclic-temporal-exactness-criterion}.
		
		The same refinement argument applies to minimal complete DFAs for prefix-closed languages containing the empty word. A nonuniversal such language has exactly one rejecting residual (the empty language): once a prefix fails, no continuation can be accepted. A universal language has one accepting state. Thus for differing languages with minimal state counts $u,v$, identify their rejecting sinks if both are nonuniversal. The combined automaton has
		\begin{equation}\label{rev:eq:residual-union-size}
			D(u,v)=u+v-\mathbf 1_{\{u>1,\,v>1\}}
		\end{equation}
		states. It has two initial acceptance blocks, so a separating word has length at most $D(u,v)-2$, and is nonempty. Its horizon is at most $D(u,v)-3$. This is the shared refinement estimate used in the intrinsic and canonical missing-path cutoffs below.
	\end{proof}
	This is the classical subset construction and deterministic state refinement \cite{wn09:MyhillNerode}, applied to the two nonempty path languages with a common rejecting sink; the horizon counts transitions, not letters.

	\subsubsection{A stronger, one-step-checkable lifting condition}
	\label{wn05-06:subsec:path-lifting}

	Trace equality permits a different hidden state for each whole word.
	To extend a specified current state step by step, we require a matching successor at every neighboring chart.

	\begin{ddefinition}[Forward separator-path lifting]
		\label{wn05-06:def:edge-lifting}
		Fix a join-tree edge $IJ$ and its separator $S=I\cap J$.
		The observation $\pi_S:X_I\to A^S$ has \textbf{forward separator-path lifting} if, for every $u\in X_I$ and $s'\in A^S$,
		\[
		(\pi_Su,s')\in R_I|_S\ \Longrightarrow\
		\exists v\in X_I:\ (u,v)\in R_I,\ \pi_Sv=s'.
		\]
		The lifting condition on the edge requires the analogous property on the $J$ side as well.
		A local relation $R_I$ is \textbf{serial} if
		$\forall u\in X_I\ \exists v\in X_I:(u,v)\in R_I$.
	\end{ddefinition}

	\begin{dtheorem}[Nonblocking exact globalization under edge lifting]\label{wn05-06:thm:fibration-gluing}
		Suppose $\mathcal H$ admits a join tree $T$, and the atlas is state-compatible and one-step compatible on its edges.
		Assume that every $R_I$ is serial and that forward separator-path lifting holds on both sides of every tree edge.
		Then:
		\begin{enumerate}[label=(\alph*),leftmargin=2.3em]
			\item the canonical global relation $R_{\mathcal A}$ is serial on all of	$M_{\mathcal H}$;
			
			\item for any chart $I$, any local finite or infinite trajectory from $\pi_I x$ and any compatible initial global state $x\in M_{\mathcal H}$, there exists a global trajectory from $x$ whose $I$-projection equals the prescribed local trajectory;
			
			\item consequently, the conditions of
			Theorem~\ref{thm:acyclic-temporal-exactness-criterion} hold, so $R_{\mathcal A}$ globalizes all finite horizons exactly.
		\end{enumerate}
	\end{dtheorem}

	\begin{proof}
		\leavevmode
		\begin{enumerate}[label={(iii)},leftmargin=*,nosep]
			\item[(a)] Fix a compatible global initial state $x$ and root the join tree at an arbitrary chart $I$. For (a), choose a successor of $\pi_Ix$, which exists
			by seriality.
			Move to each neighboring chart $J$, carrying the chosen parent transition's separator successor $s'$.
			The common separator pair $(\pi_Sx,s')$ belongs to $R_J|_S$ by one-step compatibility.
			Forward	lifting on $J$ supplies a successor from its specified current state $\pi_Jx$, with the required separator next state.
			Continue through the tree.
			Running intersection glues the chosen successor chart states into one global next state $y$, and $(x,y)\in R_{\mathcal A}$.
			
			\item[(b)] For (b), fix the entire local trajectory on root $I$.
			Repeat the same	construction at each successive time, keeping the already chosen local root trajectory fixed.
			For infinite time, this is an ordinary recursive
			construction (or, in finite state spaces, follows from finite-prefix extension and a finitely branching path tree).
			The running intersection argument guarantees consistency at every step and every chart.
			
			\item[(c)] Assertion (c) follows: the state-gluing part of Theorem~\ref{wn05-06:thm:join-tree-one-step} first extends any prescribed local initial state to a compatible global state.
			Then (b) extends its prescribed trajectory; the reverse inclusion of projected global paths into local paths is automatic. \qedhere
		\end{enumerate}
	\end{proof}

	The lifting test checks one state and one step at a time. It suffices for equality of the complete trace languages. Language equality alone allows different paths with the same observed prefix to use different hidden witnesses. For example, on binary charts $ab,bc$ with full state sets, take $(a',b')=(a,a)$ and $(b',c')=(c,c)$. Both $b$-languages consist of words constant from time one onward, with arbitrary initial bit, so the atlas globalizes exactly at every horizon. But a compatible state with $a\ne c$ has no global successor, although both local maps are total.

	\subsubsection{Explicit obstructions and independent global constraints}
	\label{wn05-06:subsec:temporal-counterexamples}

	A two-edge tree isolates a temporal obstruction.
	The example also explains why a separately prescribed global state constraint must appear in the join test.

	\begin{dexample}[A join tree with exact one-step gluing but no two-step path]
		\label{wn05-06:ex:one-not-two}
		Take binary coordinates $a,b,c$ and the acyclic cover
		$\mathcal H=\{\{a,b\},\{b,c\}\}$. Every local state is allowed. Define
		\[
		R_{ab}:\ (a',b')=(a,a),\qquad
		R_{bc}:\ (b',c')=(c,1-c).
		\]
		On the shared coordinate $b$ both local one-step projected relations are
		full, so one-step gluing is exact on the join tree. Yet a two-step path of
		$R_{ab}$ forces $b_2=b_1$, while a two-step path of $R_{bc}$ forces
		$b_2=1-b_1$. No global two-step path exists, although each local update is
		total and deterministic. The first edge trace defect is exactly $\delta_{ab,bc}=2$.
		The cover has $\nu(\mathcal H)=\infty$ by Definition~\ref{def:nu}, demonstrating that favorable static obstruction geometry alone does not control long-time exactness.
		Forward separator lifting fails in this example.
	\end{dexample}

	On a cyclic cover the defect can already occur at one step, as the incompatible triangle in Example~\ref{ex:pairwise-compatible-transition-relations} shows.

	\begin{dproposition}[One-step exactness with independent state guards]
		\label{wn05-06:prop:guarded-globalization}
		Suppose global admissible states are constrained by a finite family
		$C_J\subseteq A^J$, $J\in\mathcal G$:
		$X=\{x\in M_{\mathcal H}:\pi_Jx\in C_J\text{ for all }J\in\mathcal G\}$.
		On two-time coordinate sets $V^0\sqcup V^1$, lift each transition relation
		$R_I$ to scope $I^0\sqcup I^1$ and each state guard to the relation
		$C_J\times C_J$ on scope $J^0\sqcup J^1$. If this \emph{augmented}
		hypergraph admits a join tree and its lifted relations are pairwise
		compatible on the tree separators, then there is a global relation
		$R\subseteq X\times X$ whose restrictions to all original transition charts
		are exactly $R_I$. (The same conclusion holds for every guard relation on its doubled scope.) State exactness remains a separate condition.
	\end{dproposition}

	\begin{proof}
		Treat all lifted relations, including guards, as ordinary relations on the
		doubled attribute set. If any lifted relation is empty, separator consistency propagates emptiness along the tree, so every lifted relation and the join are empty. Otherwise the acyclic relational join argument in Theorem~\ref{wn05-06:thm:join-tree-one-step} gives a nonempty natural join with exact projections. Every joined current
		state and successor satisfies each guard $C_J$ and the original chart state requirements (because $R_I\subseteq X_I^2$). Hence the resulting relation lies in $X\times X$. Its projection onto each original transition scope is
		exactly $R_I$.
	\end{proof}

	The guard condition is a \emph{strong sufficient} test, not a necessary one.
	A diagonal admissible space $X=\{00,11\}$ with singleton chart laws
	$a'=1-a$ and $b'=b$ already shows that even one-step compatible local maps
	can fail on a constrained global state space (Section~\ref{sec:fixed}). For constraints
	whose augmented hypergraph is cyclic, new propagation and recovery
	hypotheses are needed; the static rank $\nu(\mathcal H)$ alone cannot test
	preservation of the independent admissibility relations.

	Example~\ref{wn05-06:ex:one-not-two} has a join tree and compatible one-step transitions, yet a two-step word fails to lift. This motivates keeping the size of a spatial obstruction separate from the first time at which a hidden witness fails.

	\section{Time-expanded covers and propagation of temporal assertions}
	\label{sec:covers-assertions}

	A time horizon replaces each coordinate by a string of observed cells. Theorem~\ref{thm:obstruction-rank-invariance-expansion} proves that this leaves the cover's obstruction rank unchanged. On the expanded charts, Theorem~\ref{wn08:thm:temporal-two-pass} computes the local paths that belong to a common global history. Theorem~\ref{wn08:thm:unbounded-temporal-delay} shows why the spatial rank alone cannot bound the first disagreement between two histories.

	At horizon $m$ the set $\mathcal M_m$ is the natural join of the complete chart path relations on the compatible-state envelope. A guard is an original chart containing the coordinates named by an assertion. The join supplies the semantics for testing that assertion.

	A finite temporal assertion specifies values or positive-event implications on the space-time cells $V\times\{0,\ldots,m\}$. For a local relation $T_I(m)$, a pattern $p$ is \emph{realizable} when some local path extends it, and is a \emph{nogood} otherwise. Valid implications and refutations follow by natural join, selection, and projection. Hidden intermediate states may introduce longer or higher-arity consequences when eliminated.
	
	For a causal persistent event system, a local closure $\gamma_I$ is empty preserving, causal, shift invariant and persistent; all semantic consequences are positive. The complete time-prefix closures are finite. For a monotone absorbing Horn update system, the eventual failed-event set is exactly the closure generated by its triggering implications from the initial seed set. The construction applies to monotone absorbing updates.

	\subsection{Temporal assertions and relational propagation}
	\label{wn08:sec:dynamic-assertions}

	At a fixed horizon an assertion specifies values on selected coordinate--time cells.
	Its validity depends on whether a complete path realizes those values.
	This interpretation allows hidden cells to be eliminated by projection.

	Let $V$ be a finite descriptor set, $A_v$ a nonempty finite alphabet for each $v\in V$, and let $\mathcal H$ be a finite chart cover.
	We retain the finite transition-atlas conventions from Section~\ref{sec:local-dynamical-atlases}:
	each chart $I$ has nonempty state set $X_I\subseteq A^I$ and local transition relation $R_I\subseteq X_I\times X_I$.
	The static packaged object supplies provenance; restriction to descriptor charts is additional data, not an axiom imposed on every intrinsic package.
	In particular, the temporal recovery invariants studied here are distinct from the static rank $\operatorname{ipr}$ of Paper~I~\cite{MaPackaging2026}.

	For a horizon $m\ge0$ put
	\[
	W_m=V\times\{0,1,\ldots,m\},\qquad
	W_{I,m}=I\times\{0,1,\ldots,m\}.
	\]
	For a set $W$ of space--time cells, write $A^W=\prod_{(v,t)\in W}A_v$.
	A \emph{space--time literal} is a typed expression $[v,t=a]$, with $v\in V$, $0\le t\le m$, and $a\in A_v$.
	A \emph{temporal pattern} $p$ is a finite partial assignment, that is, a set of literals containing at most one value for each pair $(v,t)$.
	Its \emph{space--time support} $\operatorname{st}(p)$ is the set of pairs $(v,t)$ appearing in $p$, and its spatial support is
	$\operatorname{sp}(p)=\{v:(v,t)\in\operatorname{st}(p)\}$.
	A pattern is \emph{chart-supported} when $\operatorname{sp}(p)\subseteq I$ for some original chart $I\in\mathcal H$.
	Literal and chart provenance remain part of every assertion; an unrecorded auxiliary chart is not an original chart.

	\begin{ddefinition}[Local and canonical global path relations]
		\label{wn08:def:temporal-model-relation}
		Let $\mathcal A$ be a transition atlas on $\mathcal H$
		and $m\in\mathbb N_0$ a temporal horizon.
		\begin{enumerate}[label=(\roman*),leftmargin=*,nosep]
			\item \textbf{Local path relation.}
			The local path relation on chart $I$ is
			\[
			Q_I^{(m)}:=\mathsf T_I(m)\subseteq A^{W_{I,m}},
			\]
			under the identification of trajectories with tuples
			on the time-indexed coordinates $W_{I,m}$.
			
			\item \textbf{Canonical global path relation.}
			The canonical global path relation of $\mathcal A$
			at horizon $m$ is
			\begin{equation}\label{wn08:eq:global-path-join}
				\begin{aligned}
					\mathcal M_m
					&:=\Join_{I\in\mathcal H}Q_I^{(m)}\\
					&=\{z\in A^{W_m}:z|_{W_{I,m}}\in Q_I^{(m)}
					\text{ for all }I\in\mathcal H\}.
				\end{aligned}
			\end{equation}
		\end{enumerate}
		Equivalently, $\mathcal M_m$ is the trajectory set
		$\mathsf T_{R_{\mathcal A}}(m)$ of the canonical global
		relation defined in Section~\ref{sec:local-dynamical-atlases}.
		The relation $\mathcal M_m$ may be empty.
	\end{ddefinition}

	For a prescribed admissible global state space $X\subseteq
	\mathcal M_0$, the corresponding global path relation is
	$\mathcal M_m\cap X^{m+1}$, under the time-coordinate
	identification.

	\begin{ddefinition}[Temporal implications, rooted circuits, and minimal nogoods]
		\label{wn08:def:temporal-assertions}
		For a finite trajectory relation $M\subseteq A^{W_m}$ and partial assignment $p$, define its witness set $M[p]:=\{w\in M:p\subseteq w\}$.
		A literal $q$ is a single-cell assignment.
		\begin{enumerate}[label=(\roman*),leftmargin=*,nosep]
			\item A \textbf{temporal implication}
			$p\Rightarrow_M q$ holds when $M[p]\subseteq M[q]$, for a literal $q$;
			$p\Rightarrow_M\bot$ is a \emph{refutation} or \emph{nogood} when $M[p]=\varnothing$.
			An implication is realizable when $M[p]\ne\varnothing$.
			This distinction rules out treating vacuous conclusions of an impossible premise as substantive dynamical recovery.
			
			\item A \textbf{rooted temporal circuit} is a realizable implication $p\Rightarrow_Mq$ with $q\notin p$ such that no proper subpattern
			$p'\subsetneq p$ entails $q$.
			Empty stems are permitted: a temporal literal
			may be forced without observations.
			
			\item A \textbf{minimal temporal nogood} is a pattern $p$ with $M[p]=\varnothing$ and $M[p']\ne\varnothing$ for every proper subpattern $p'\subsetneq p$.
		\end{enumerate}
		These are intrinsic to the specified semantic relation $M$.
		They do not depend on a chosen implicational presentation.
	\end{ddefinition}

	If $M=\varnothing$, even the empty premise refutes, and all material implications are formally true.
	We therefore keep realizability and refutation separate.  Unlike the empty-preserving static closure class in
	the original globalization theorem, temporal-assertion closures can have nontrivial empty-premise consequences.

	\begin{dproposition}[Semantic temporal closure and minimal cores]
		\label{wn08:prop:temporal-semantic-closure}
		Let $\mathcal L_m$ be the finite set of all typed literals over $W_m$.
		For any set $P\subseteq\mathcal L_m$, write $M[P]=\{w\in M:P\subseteq w\}$, including $M[P]=\varnothing$ when $P$ is inconsistent.
		Define
		\[
		c_M(P)=\{q\in\mathcal L_m:M[P]\subseteq M[q]\}.
		\]
		\begin{enumerate}[label=(\roman*),leftmargin=*,nosep]
			\item Then $c_M$ is an extensive, monotone, idempotent closure operator.
			
			\item Every	realizable implication $P\Rightarrow_Mq$ admits a rooted temporal circuit $P'\Rightarrow_Mq$ with $P'\subseteq P$, unless $q\in P$, in which case the conclusion is tautological.
			
			\item Every nonrealizable finite pattern contains an inclusion-minimal nogood.
		\end{enumerate}
		Neither statement bounds the spatial or temporal size of its core uniformly over all finite-state systems.
	\end{dproposition}

	\begin{proof}
		\leavevmode
		\begin{enumerate}[label={(iii)},leftmargin=*,nosep]
			\item[(i)] Extensivity and monotonicity follow from the definition: enlarging the
			premise shrinks the set of models.  Since all elements of $M[P]$ satisfy
			all literals of $c_M(P)$, the sets of models of $P$ and $c_M(P)$ coincide;
			this implies idempotence.
			\item[(ii)] To obtain the core of a realizable nontrivial
			implication, choose an inclusion-minimal subpattern of $P$ still entailing
			$q$.  It is realizable because the original premise is.
			\item[(iii)] The same finite
			minimality argument yields an inclusion-minimal nogood from any
			nonrealizable finite pattern. \qedhere
		\end{enumerate}
	\end{proof}

	Temporal consequences cannot in general be inferred by choosing the hidden intermediate state independently at each step.
	We instead use a typed finite relational calculus whose primitive operations are:
	\emph{time shift}, \emph{natural join}, \emph{selection} by a pattern, \emph{existential projection}, and \emph{relational composition}.
	Each operation preserves the actual coordinate/time typing.
	Time shift here transports a relation to the indicated time-labelled scope.
	It does not license arbitrary enlargement of its set of admissible histories. The cut rule
	\[
	\frac{p\Rightarrow q\qquad p\cup\{q\}\Rightarrow r}
	{p\Rightarrow r}
	\]
	and the refutation rule
	\[
	\frac{p\Rightarrow q\qquad p\cup\{q\}\Rightarrow\bot}
	{p\Rightarrow\bot}
	\]
	are sound for semantic implications; they are \emph{not} asserted to be a
	complete finite Horn-style calculus for arbitrary relational systems.

	\begin{dtheorem}[Exact finite-horizon propagation and witness-preserving elimination]
		\label{thm:exact-finite-horizon-propagation}
		For every finite transition atlas, every horizon $m$, every temporal pattern $p$, and every literal $q$, the following statements hold.
		\begin{enumerate}[label=(\roman*),leftmargin=*,nosep]
			\item The assertion tests are necessary and sufficient:
			\begin{align}
				\mathcal A\models_m p\Rightarrow\bot
				&\quad\Longleftrightarrow\quad
				\sigma_p\left(\Join_I Q_I^{(m)}\right)=\varnothing,
				\label{wn08:eq:exact-nogood}\\
				\mathcal A\models_m p\Rightarrow q
				&\quad\Longleftrightarrow\quad
				\sigma_p\left(\Join_I Q_I^{(m)}\right)
				\subseteq\sigma_q\left(\Join_I Q_I^{(m)}\right).
				\label{wn08:eq:exact-entailment}
			\end{align}
			Here $\sigma_p$ means restriction to tuples realizing $p$.
			To determine whether $p\Rightarrow q$ is nonvacuous, additionally test that the left-hand selection is nonempty.
			
			\item For an assertion involving only variables in $U\subseteq W_m$, the exact relation of feasible witnesses is
			\begin{equation}\label{wn08:eq:exact-elimination}
				\pi_U\left(\Join_IQ_I^{(m)}\right).
			\end{equation}
			
			\item In particular, for $R\subseteq A^{U\cup S}$ and $Q\subseteq A^{S\cup W}$ on disjoint $U,S,W$,
			\begin{equation}\label{wn11:eq:join-project}
				\pi_{U\cup W}(R\Join Q)=\{(u,w):\exists s\ [(u,s)\in R\ \text{and}\ (s,w)\in Q]\}.
			\end{equation}
			The same shared $S$-tuple must witness both constraints; independently projecting the endpoints can introduce false witnesses. Eliminate a hidden variable only after joining all constraints incident with it.
		\end{enumerate}
		Consequently, natural join, selection and projection suffice to decide every finite temporal implication or nogood in this model, with a
		finite global witness when the proposed assertion fails.
	\end{dtheorem}

	\begin{proof}
		\leavevmode
		\begin{enumerate}[label={(iii)},leftmargin=*,nosep]
			\item[(i)] The natural join in Definition~\ref{wn08:def:temporal-model-relation} is exactly
			the set of all compatible global finite trajectories.  Selection takes
			precisely those trajectories satisfying a premise, and the displayed
			inclusion is exactly the semantics of implication.  An empty selection
			witnesses refutation, whereas a tuple in the set difference of the two
			selections witnesses failure of implication.
			
			\item[(ii)] Projection existentially eliminates the unwanted space--time coordinates and preserves exactly the remaining feasible tuples.
			
			\item[(iii)] Equation~\eqref{wn11:eq:join-project} is its two-factor specialization: a common separator tuple is necessary and sufficient by the definition of join.
		\end{enumerate}
		All involved relations are finite, so this is a decision procedure.
		The result is completeness of the relational semantic calculus.
	\end{proof}

	\begin{dremark}[Temporal composition and hidden witnesses]
		\label{wn08:rem:temporal-composition}
		For a chart $I$ at horizon $m\ge1$, the path relation is computed by joining time-shifted copies of $R_I$ on consecutive pairs and projecting only after the full intermediate time coordinates are shared. At horizon zero the path relation is $X_I$.  If one first projects $R_I$ to a
		proper subchart and only then composes, intermediate hidden witnesses may
		be chosen inconsistently.  Thus
		$(R_I\circ R_I)|_S\subseteq(R_I|_S)\circ(R_I|_S)$,
		with possible strict inclusion as shown by the swap example in Section~\ref{sec:local-dynamical-atlases}.
		The order of existential elimination and temporal composition matters.
	\end{dremark}

	\subsection{Global-support recovery from original charts}
	\label{wn08:subsec:global-support}

	A local possibility survives the full system exactly when a global path projects to it.
	Projecting the history join back to each chart removes the unsupported tuples.

	\begin{ddefinition}[Global-support recovery and unsupported paths]
		\label{wn08:def:global-support-recovery}
		For an atlas $\mathcal A$, chart $I$ and horizon $m$, its global-support recovery operator is defined by
		\begin{equation}\label{wn08:eq:recovered-chart-paths}
			\operatorname{Rec}_I^{(m)}(\mathcal A)
			=\pi_{W_{I,m}}\mathcal M_m
			=\pi_{W_{I,m}}\left(\Join_JQ_J^{(m)}\right).
		\end{equation}
		It is the set of local paths that possess some compatible global realization.
		Its set of \textbf{unsupported local paths} is
		\[
		\mathcal D_I^{(m)}(\mathcal A)
		=Q_I^{(m)}\setminus\operatorname{Rec}_I^{(m)}(\mathcal A).
		\]
		The \textbf{global-support recovery profile} is $\operatorname{Rec}^{(m)}(\mathcal A):=(\operatorname{Rec}_I^{(m)}(\mathcal A))_{I\in\mathcal H}$.
		It is not a scalar intrinsic recovery rank.
		It is defined without choosing a presentation of any local rule set.
	\end{ddefinition}

	\begin{dproposition}[Maximality of recovered chart traces]
		\label{wn08:prop:maximal-supported}
		The family $\bigl(\operatorname{Rec}_I^{(m)}(\mathcal A)\bigr)_{I\in\mathcal H}$
		is the greatest, componentwise, family of subrelations $P_I\subseteq Q_I^{(m)}$
		for which there exists a \emph{single} global relation $G\subseteq A^{W_m}$
		with $\pi_{W_{I,m}}G=P_I$ for every original chart $I$.  In particular,
		exact temporal globalization through horizon $m$ holds if and only if
		$\mathcal D_I^{(h)}=\varnothing$ for all $I$ and all $0\le h\le m$.
	\end{dproposition}

	\begin{proof}
		If $G$ realizes the subrelations, every tuple of $G$ satisfies the
		original $Q_I^{(m)}$ restrictions; thus $G\subseteq\mathcal M_m$ and
		$P_I\subseteq\operatorname{Rec}_I^{(m)}$.  Conversely, $G=\mathcal M_m$
		realizes all recovered projections simultaneously, establishing
		maximality.  The final characterization follows from the definition of
		exact finite-horizon globalization and the fact that $\mathcal M_h$ is
		precisely the path relation of the canonical greatest one-step extension.
	\end{proof}

	\begin{dtheorem}[Two-pass temporal support recovery on a join tree]
		\label{wn08:thm:temporal-two-pass}
		Suppose the chart cover admits a join tree $T$.  At any fixed horizon $m$,
		start with the original path tables $Q_I^{(m)}$.  Root $T$ arbitrarily and
		perform two sweeps of semijoins:
		\begin{enumerate}[label=(\roman*),leftmargin=2.4em]
			\item in a bottom-up pass, for each child $J$ with parent $I$, replace
			$Q_I\leftarrow Q_I\ltimes Q_J$;
			\item in a top-down pass, replace each child
			$Q_J\leftarrow Q_J\ltimes Q_I$ using the already reduced parent.
		\end{enumerate}
		Here $Q\ltimes R$ retains precisely those tuples of $Q$ matching some
		tuple of $R$ on the separator, including all time-indexed separator
		coordinates.  The resulting table on every original chart $I$ is exactly
		\[
		Q_I^{\mathrm{final}}=\operatorname{Rec}_I^{(m)}(\mathcal A).
		\]
		Every surviving local trajectory extends to a global trajectory, and every
		removed tuple is unsupported.  If $\operatorname{Rec}_I^{(m)}=Q_I^{(m)}$ for all
		charts, the original atlas is exact at horizon $m$.
	\end{dtheorem}

	\begin{proof}
		\leavevmode
		\begin{enumerate}[label={(iii)},leftmargin=*,nosep]
			\item[(i)] For a rooted join tree, induction over a subtree shows that after the
			bottom-up pass, each rootward table retains exactly those tuples extendable
			to consistent tuples in every chart of its rooted subtree.  Running
			intersection ensures that subtrees interact with their complements only
			through the parent separator.
			\item[(ii)] The top-down pass imposes the extension
			requirements from the complementary parts of the tree.  After it, a tuple
			survives precisely when it extends across every chart, that is, when it is
			in $\pi_{W_{I,m}}\mathcal M_m$.  Conversely, a tuple in a genuine global
			trajectory matches all its neighbors and cannot be deleted by any semijoin. \qedhere
		\end{enumerate}
	\end{proof}

	Applying the full reducer of Yannakakis to complete history tuples preserves a common witness at every separator \cite{Yannakakis1981,wn09:BeeriFaginMaierYannakakis1983}. The first horizon at which it deletes a local tuple has a direct separator characterization.

	\begin{dcorollary}[Recovery defect equals first separator-language defect]
		\label{wn08:cor:first-recovery-defect}
		For a finite join-tree atlas, let
		\[
		\tau(\mathcal A)=\min\{m:\mathcal D_I^{(m)}\ne\varnothing
		\text{ for some original chart }I\},
		\]
		with value $\infty$ when no such horizon exists.  Then
		\begin{equation}\label{wn08:eq:tau-equals-delta}
			\quad \tau(\mathcal A)=\delta_T(\mathcal A),\quad
		\end{equation}
		where $\delta_T$ is the separator trace-defect depth of Section~\ref{sec:local-dynamical-atlases}.
		When finite, a shortest failure admits a certificate using the projected
		trajectory of one chart and the trace language of an adjacent chart.
	\end{dcorollary}

	\begin{proof}
		If every tree-edge separator language agrees through horizon $m$, the
		join-tree argument of Theorem~\ref{wn08:thm:temporal-two-pass} applied to the time-expanded relations realizes
		every original chart tuple, so all defects through that horizon vanish.
		Conversely, if $\mathsf L_I^S(m)\ne\mathsf L_J^S(m)$ on one edge $IJ$,
		choose a word present on one side but not the other.  Any local path
		realizing that word on its own side cannot extend to a global trajectory;
		hence its chart defect is nonempty.  Taking the first horizon gives the
		identity.  This local certificate involves the two adjacent original
		charts, even though the temporal word itself can be long.
	\end{proof}

	\begin{dproposition}[Horizon-prefix pruning and nonblocking]
		\label{wn08:prop:horizon-prefix}
		Let $\operatorname{pre}_m$ discard the final time slice.  Then
		\[
		\operatorname{pre}_m(\mathcal M_{m+1})\subseteq\mathcal M_m.
		\]
		Equality holds for every $m\ge0$ if and only if the canonical global
		transition relation $R_{\mathcal A}$ is serial on $M_{\mathcal H}$.
	\end{dproposition}

	\begin{proof}
		A global trajectory remains a trajectory after deleting its final state,
		proving inclusion.  If $R_{\mathcal A}$ is serial, extend the last state of
		any finite global path by one successor; this proves equality for all
		$m$.  Conversely, equality at $m=0$ requires every global state
		$x\in\mathcal M_0=M_{\mathcal H}$ to extend to a global one-step path,
		which is precisely seriality.
	\end{proof}

	\subsection{Obstruction-rank invariance and temporal obstructions}
	\label{wn08:subsec:temporal-cylinder}

	The time-expanded atlas has a precise relationship to the original static obstruction rank, but this relationship alone does not imply
	temporal exactness.

	\begin{dtheorem}[Obstruction-rank invariance under uniform time expansion]
		\label{thm:obstruction-rank-invariance-expansion}
		For any finite cover $\mathcal H$ of $V$ and any finite nonempty time index	set $T$, define the \emph{time-cylinder cover}
		\[
		\mathcal H^T=(I_j\times T)_{j\in J}
		\quad\text{of}\quad V\times T,\qquad \mathcal H=(I_j)_{j\in J}.
		\]
		For the obstruction rank $\nu$ of  Definition~\ref{def:nu},
		\begin{equation}\label{wn08:eq:cylinder-rank}
			\boxed{\nu(\mathcal H^T)=\nu(\mathcal H).}
		\end{equation}
	\end{dtheorem}

	\begin{proof}
		Let $G$ be the primal graph of $\mathcal H$.  The primal graph of
		$\mathcal H^T$ replaces each vertex $v\in V$ with a clique of
		$|T|$ \emph{true twins}: copies $(v,t)$ are adjacent, and two such
		vertex cliques are completely adjacent precisely when their original
		vertices were adjacent in $G$.  If $G$ has an induced cycle of length
		at least four, choosing one copy of each cycle vertex produces the same
		induced cycle in the expanded graph.  Conversely, if $G$ is chordal,
		take a perfect-elimination ordering of $G$~\cite{Dirac1961,FulkersonGross1965} and eliminate all time copies
		of each vertex consecutively in that order.  The later neighbors of
		every eliminated copy form a clique, so the expanded graph is chordal.
		Thus chordality is preserved in both directions.
		
		A clique $C'$ of the expanded primal graph projects to a clique
		$C=\pi_V(C')$ of $G$.  It is contained in an expanded chart $I\times T$
		if and only if $C\subseteq I$.  Hence an unguarded expanded clique
		projects to an unguarded original clique of no greater cardinality.
		Conversely, an unguarded clique in $G$ lifts, using one time index, to
		an unguarded expanded clique of the same cardinality.  The smallest
		unguarded-clique size is therefore unchanged.  Combining this with
		the chordality clause in the definition of $\nu$ proves
		\eqref{wn08:eq:cylinder-rank}, including the value $\infty$.
	\end{proof}

	This invariance isolates the issue: the original cover geometry is not
	deteriorated by copying every chart through time.  New temporal
	obstructions arise because a trajectory relation contains constraints
	connecting \emph{different time slices}, and local projection does not
	necessarily preserve their hidden intermediate witnesses.

	The preceding recovery operator detects an obstruction at a finite horizon
	when one exists.  Its first detection time, however, is not controlled by
	the original static cover geometry alone.

	\begin{dtheorem}[Arbitrarily delayed temporal obstruction on one acyclic cover]
		\label{wn08:thm:unbounded-temporal-delay}
		For every integer $N\ge3$, there is a time-homogeneous, finite, state-compatible and one-step-compatible transition atlas on the fixed cover
		\[
		\mathcal H=\{\{a,b\},\{b,c\}\},\qquad\nu(\mathcal H)=\infty,
		\]
		such that the overlap languages agree for every horizon $m<N-1$ but	disagree at $m=N-1$.
		Thus
		\[
		\delta_T(\mathcal A_N)=\tau(\mathcal A_N)=N-1.
		\]
		In particular no function of $\nu(\mathcal H)$, the number of charts, or the largest chart-cardinality alone can bound the first temporal obstruction depth in this class.
	\end{dtheorem}

	\begin{proof}
		Let the shared visible bit be $b\in\{0,1\}$.
		On chart $\{a,b\}$ use hidden run-count states
		\[
		X_{ab}=\{(0,0),(1,1),\ldots,(N-1,1)\}.
		\]
		From any state, the next state may be $(0,0)$, representing a zero.
		A next state $(h+1,1)$ is allowed if the current count is $h<N-1$.
		These rules	make the finite visible words precisely the binary words avoiding $N$ consecutive ones.
		On chart $\{b,c\}$ use the analogous counter up to
		$N$, permitting words avoiding $N+1$ consecutive ones.  Relabel its hidden-state coordinate as $c$; the order of coordinates is immaterial.
		Both charts have the full one-step projected relation on the visible bit	$b$ because $N\ge3$, and both project their state sets onto $\{0,1\}$.
		Every binary word of length less than $N$ is realizable in both charts,	whereas $1^N$ is forbidden in the first and realizable in the second.
		Therefore the earliest separator discrepancy occurs for a word of length $N$, i.e., at horizon $m=N-1$.
		The fixed two-edge-path cover is chordal and conformal, so its static obstruction rank is $\infty$.
		The	identity for $\tau$ follows from Corollary~\ref{wn08:cor:first-recovery-defect}.
	\end{proof}

	The construction uses hidden state spaces whose sizes grow with $N$;
	the subset-automaton cutoff of Section~\ref{sec:local-dynamical-atlases} depends on state complexity.
	Thus the static geometric rank and the temporal complexity of the local laws are independent inputs.

	The shortest temporal obstruction between two specified finite systems is bounded in terms of their state counts.
	That is different from bounding the lengths of all minimal forbidden words of one system.

	\begin{dproposition}[Three-state system with arbitrarily long minimal forbidden blocks]
		\label{wn08:prop:even-shift}
		There is a finite, three-state transition system with binary observation $b$ whose temporal observation language has forbidden blocks minimal under the contiguous-factor order, of unbounded length.
		In fact, for every $k\ge0$,
		\begin{equation}\label{wn08:eq:even-shift-forbidden}
			1\,0^{2k+1}\,1
		\end{equation}
		is forbidden and each proper contiguous factor is permitted.
	\end{dproposition}

	\begin{proof}
		Use the three labeled edges of a two-vertex graph:
		$E\xrightarrow{1}E$, $E\xrightarrow{0}O$, and
		$O\xrightarrow{0}E$.
		Treat these edges as the three states of a transition system; a state may follow another precisely when the target vertex of the first equals the source vertex of the second.
		The observed value is the edge label $0$ or $1$.
		A $1$-edge begins and ends at $E$;
		there must be an even number of $0$-edges before another $1$ can be read.
		Hence each word in~\eqref{wn08:eq:even-shift-forbidden} is impossible.
		Each proper contiguous factor contains at most one $1$ and can be read by starting in a suitable vertex; consequently it is permitted.
	\end{proof}

	This is the standard even shift in symbolic dynamics \cite{wn08:EvenShift}.
	Minimality here refers to deletion of a contiguous \emph{prefix or suffix} (the factor order), not necessarily to deletion of arbitrary time-indexed literals in Definition~\ref{wn08:def:temporal-assertions}.
	Consequently, defining a dynamical recovery rank merely as the maximum length of a minimal forbidden block would give value $\infty$ even for a three-state system.

	\subsection{Absorbing Horn propagation and original-chart recovery}
	\label{wn08:subsec:temporal-horn}

	Not every dynamical class needs the full relational apparatus.  The
	monotone absorbing fragment provides an exact connection to closure
	systems and to the package-collapse specialization of Section~\ref{subsec:monotone-cascades}.

	\begin{dtheorem}[Temporal Horn propagation and the absorbing closure]
		\label{wn08:thm:temporal-horn}
		Let $V$ be finite and let $\mathcal B$ be a finite family of rules
		$P\to i$ with $P\subseteq V\setminus\{i\}$, allowing $P=\varnothing$.
		Define the absorbing update
		\begin{equation}\label{wn08:eq:horn-dynamics}
			D_{t+1}=D_t\cup\{i:\exists(P\to i)\in\mathcal B,
			P\subseteq D_t\}.
		\end{equation}
		Then $D_t$ is exactly the set of events forced at time $t$ in the least model of the time-indexed Horn rules
		\[
		e_{i,t}\Rightarrow e_{i,t+1},\qquad
		\bigwedge_{j\in P}e_{j,t}\Rightarrow e_{i,t+1},
		\]
		with seed facts $e_{i,0}$ for $i\in D_0$.  The terminal set is the
		ordinary implicational closure $c_{\mathcal B}(D_0)$; it stabilizes in
		at most $|V|$ strict-update rounds.  For $D_0=\varnothing$, the generated
		closure is empty preserving precisely when the rule system has no empty-premise trigger (equivalently, no nonempty consequence of the empty seed).
	\end{dtheorem}

	\begin{proof}
		Induct on $t$.  At $t=0$ the least Horn model has exactly the seed
		facts; the rules for time $t+1$ produce exactly the old events and the
		heads of all rules whose premises already hold in $D_t$.  This is
		\eqref{wn08:eq:horn-dynamics}.  Iteration is the standard forward-chaining
		construction of the least closed set containing $D_0$, proving the
		closure identity.  Every strict round increases $|D_t|$, so there are
		at most $|V|$ strict rounds.  Finally, $c_{\mathcal B}(\varnothing)$ is
		empty precisely when no consequence arises from the empty seed; in the present finite forward system this is equivalent to the absence of every empty-premise trigger, since any such rule immediately adds its head.
	\end{proof}

	This is a standard property of Horn forward chaining~\cite{CaspardMonjardet2003,wn10:DatalogSurvey}.  It identifies a class where dynamical propagation	can be analyzed through the original kind of static implication
	closure, subject to the original class's empty-preserving hypotheses.
	For reversible and nonmonotone dynamics, a single absorbing-event
	closure generally loses essential information.

	\paragraph{Recovery from original charts.}

	The preceding results isolate four distinct ingredients of a general intrinsic dynamical propagation-recovery rank:
	\begin{enumerate}[label=(\roman*),leftmargin=2.6em]
		\item \emph{Semantic rooted temporal circuits and nogoods}, intrinsic to a
		finite trajectory language rather than a chosen syntactic rule basis;
		
		\item \emph{temporal composition and elimination} with typed time-indexed
		supports and a record of the original chart of each assertion;
		
		\item \emph{original-chart recovery}: a derived consequence must be
		realized in an original chart, not merely in an auxiliary amalgamation;
		
		\item \emph{guarding and class parameters}: a candidate classwise rank must
		specify the spatial guarding scale, temporal horizon, and admissible
		local-state complexity, and must distinguish individual and universal globalization, as in Paper~I~\cite{MaPackaging2026}.
	\end{enumerate}

	An atlas-specific shortest defect $\delta_T$, a temporal observation
	depth $\operatorname{tod}$, and the maximum duration of a minimal
	forbidden word are three different quantities, and none is declared
	here to be the intrinsic dynamical counterpart of
	$\operatorname{ipr}(\mathscr C)$.  In particular,
	Theorem~\ref{wn08:thm:unbounded-temporal-delay} shows no cover-only delay
	bound, while Proposition~\ref{wn08:prop:even-shift} shows a finite-state
	system may have infinitely many minimal forbidden blocks (under the
	contiguous-factor order).
	These facts motivate the rule-presentation-independent recovery parameters studied in Section~\ref{sec:width}, rather than a simple word-length rank.

	\begin{dremark}[Scope of the assertion calculus]
		The relational calculus, semijoin full reduction, finite-state trace
		recognition, even shift and Horn closure have well-established mathematical
		antecedents.  Section~\ref{sec:covers-assertions} provides a provenance-aware assertion language for the packaging specialization. Sharp classwise thresholds for specific recovery mechanisms are proved in Section~\ref{sec:global-theorems}; a single intrinsic dynamical rank covering all these mechanisms remains open (Section~\ref{sec:discussion}).
	\end{dremark}

	\section{Intrinsic temporal recovery and observation depth}
	\label{sec:width}

	How large a partial observation is needed to rule out a trajectory?
	We measure its size in named coordinate--time cells and its duration in time steps. Theorem~\ref{wn15:thm:effective-exact-rank} computes the exact cell bound of a finite-state observation language or detects that it is infinite.
	Residual languages give a separate bound on the duration of a first disagreement.

	\paragraph{Scope convention for intrinsic temporal recovery.}
	The primitive semantic input here is a previously specified observed trajectory language, whether from a realized system or an atlas-first chart.
	Temporal witness width and residual index are relative to the named observable coordinate atoms, not to an unobserved transition graph.
	The recognition results are all-horizon properties of a fixed finite-state language, not universal globalization over every atlas.
	Passage from local recovery certificates to globalization requires separately stated overlap, state-envelope, and cover-geometric hypotheses.
	These quantities need not distinguish marked interaction decompositions with identical aggregate dynamics.
	The essential read-set invariant of Proposition~\ref{dyn46:prop-native-read} is a different one-step property.

	Let $\Obs$ be a finite set of independently named observable coordinate atoms, each with a specified alphabet.
	In the finite-state recognition problems below all these alphabets are finite and nonempty.
	For each $m$, a partial observation pattern $p$ assigns values to some cells of $\Obs\times[0,m]$ and puts a wildcard elsewhere.
	Call $p$ \emph{minimally impossible} if it is unrealizable but erasing any single specified cell makes it realizable.

	\begin{definition}[Temporal witness width]
		\label{def:temporal-witness-width}
		Fix a finite set $\Obs$ of named observable atoms and a finite-horizon language family $(L(m))_{m\ge0}$.
		For a partial pattern $p$ at horizon $m$, let $\supp(p)$ be the set of specified cells.
		The \textbf{temporal witness width} is
		\[
		\tw(L;\Obs)=\max\{1,\sup_{m\ge0}\ \max_{p\text{ minimally impossible in }L(m)}|\supp(p)|\}\ \in\mathbb N\cup\{\infty\}.
		\]
		The maximum of an empty family is interpreted as zero. We also write $\operatorname{twit}(L;\Obs)=\tw(L;\Obs)$. The width is relative to the specified observable atoms.
	\end{definition}

	\begin{lemma}[Independent products and observation forgetting]
		\label{lem:product}
		For independent disjoint-coordinate trajectory systems whose path languages are nonempty at every finite horizon,
		$\tw(L_1\times L_2)=\max(\tw(L_1),\tw(L_2))$.
		The nonemptiness condition is needed for lifting a factor nogood while
		leaving the other factor unspecified. Without it, a terminating factor
		can remove the longer horizons at which the other factor has large
		witnesses. For any observed language and $U\subseteq\Obs$,
		$\tw(L|_U;U)\le\tw(L;\Obs)$.
	\end{lemma}

	\begin{proof}
		An impossible partial observation of an independent product is impossible in at least one factor. Minimality then forces all its specified cells to lie in that factor. Conversely every minimal impossible component pattern is still minimally impossible in the product. For coordinate forgetting, a minimal impossible pattern over $U$ remains minimally impossible when considered as a pattern over the larger observation set.
	\end{proof}
	Theorem~\ref{wn15:thm:effective-exact-rank} gives the all-horizon algorithm for this width, including a positive-cycle pumping certificate when it is infinite. The construction works with a specified initial subset and counts named cells, not whole observation letters. The fixed-bound recognizer of Theorem~\ref{wn14:thm:recognition} and the exact-width monitor of Theorem~\ref{wn15:thm:minimal-regular} serve different computational purposes.

	\begin{remark}
		The separate causal-Horn generating time-span invariant of Section~\ref{wn10:sec:dyn-recovery-global} is denoted $w(\gamma)=\genlen(\gamma)$, not $\tw$. A bounded primitive transition-rule span need not survive eliminating hidden variables. Likewise, a state-space bijection that merges observable coordinates can change $\tw$.
	\end{remark}

	\subsection{Intrinsic trace signatures and recovery depth}
	\label{wn09:sec:intrinsic-temporal-recovery}

	A partial history is useful for comparison only through its possible continuations.
	Residual languages record those continuations and bound the length of a first observable disagreement.

	\subsubsection{Observable trace languages and intrinsic residuals}
	\label{wn09:subsec:what-intrinsic-means}

	Paper~I~\cite{MaPackaging2026} defines its intrinsic propagation--recovery rank $\operatorname{ipr}(\mathscr C)$ for a restriction-stable class of finite empty-preserving closures:
	rooted circuits are the minimal-premise	implications of the closure itself, rather than those of a chosen rule basis.
	The rank measures how original-chart recovery after one-pivot composition interacts with spatial guarding.
	Its definition does not quantify over transition paths or observational automata.

	In a finite dynamical atlas $\mathcal A=(X_I,R_I)_{I\in\mathcal H}$, a semantically intrinsic temporal quantity should not depend on arbitrary
	labels or duplication of \emph{hidden states}, transition-rule syntax,
	or automata minimization choices. We construct an observable recovery signature rather than an unrestricted scalar analogue of
	$\operatorname{ipr}$. We distinguish three questions:
	(a) whether global paths realize all local paths,
	(b) how much history is required to expose a failure of realization, and
	(c) whether a composed consequence has a guard inside an original chart.
	The Section~\ref{sec:covers-assertions} support recovery operator	$\operatorname{Rec}_I^{(m)}=\pi_{I\times\{0,\ldots,m\}}\mathop{\bowtie}\limits_J\mathsf T_J(m)$ is retained throughout.

	Fix an original chart $I$ and any subset $U\subseteq I$.
	Let $\Sigma_U=\prod_{v\in U}A_v$, with $\Sigma_\varnothing=\{()\}$.
	The	coordinate alphabets $A_v$ are fixed as part of the chart typing.
	The	\emph{uninitialized observation language} of $I$ on $U$ is
	\begin{equation}\label{wn09:eq:trace-language}
		L_I^U=\{\varepsilon\}\cup
		\bigcup_{m\ge0}\bigl\{\,(\pi_Ux_0)\cdots(\pi_Ux_m):
		(x_0,\ldots,x_m)\in\mathsf T_I(m)\bigr\}
		\subseteq\Sigma_U^*.
	\end{equation}
	The empty word is included by convention. The language is prefix closed:
	for every $wv\in L_I^U$, its prefix $w$ also belongs to $L_I^U$. The
	observation is existential over hidden states; it need not define an
	autonomous system on $U$.

	\begin{ddefinition}[Intrinsic trace-residual index]
		\label{wn09:def:trace-residual-index}
		For any $u\in\Sigma_U^*$, define the right residual
		\[
		u^{-1}L_I^U=\{z\in\Sigma_U^*:uz\in L_I^U\}.
		\]
		The \textbf{trace-residual index} is
		\begin{equation}\label{wn09:eq:residual-index}
			\rho_I(U)=\bigl|\{u^{-1}L_I^U:u\in\Sigma_U^*\}\bigr|.
		\end{equation}
		This is the number of distinct Myhill--Nerode residuals of the observable language and equals the state count of its minimal complete deterministic acceptor \cite{Nerode1958,wn09:MyhillNerode}.
		It is determined by the full observation language, not by a particular hidden-state realization.
	\end{ddefinition}

	\begin{dlemma}[Finite memory and presentation invariance]
		\label{wn09:lem:finite-memory}
		For each finite nonempty chart state set $X_I$ and $U\subseteq I$,
		\begin{equation}\label{wn09:eq:subset-memory}
			1\le\rho_I(U)\le 1+2^{|X_I|}<\infty.
		\end{equation}
		If two chart systems have the same $U$-observation language, their trace-residual indices agree, even when their hidden state sets and	transition relations have different presentations or cardinalities.
		In particular, observation-preserving hidden-state splitting does not change $\rho_I(U)$.
	\end{dlemma}

	\begin{proof}
		To recognize~\eqref{wn09:eq:trace-language}, use an initial automaton
		position and a nondeterministic state for every $x\in X_I$. The first
		observed symbol $a$ chooses an arbitrary $x$ with $\pi_Ux=a$; each
		subsequent symbol chooses a successor along $R_I$ with that observation.
		The initial position and every reached chart state are accepting, while
		the empty set of hidden states is rejecting. Determinization has at most
		one initial position and $2^{|X_I|}$ hidden-state subsets, proving the
		bound. The residuals of a language are definitionally unchanged under
		any realization that preserves the language. The hidden-state splitting
		assertion is an instance of this invariance.
	\end{proof}

	\subsubsection{Spatial guarding and temporal defect signatures}
	\label{wn09:subsec:temporal-guard}

	Fix a join tree $T$ for the chart cover. For an oriented tree edge $e=(I,J)$ let $S_e=I\cap J$.
	We retain empty separators, since an empty-separator trace can detect mismatched permitted path lengths or deadlock.
	For every $U\subseteq S_e$, define
	\begin{equation}\label{wn09:eq:delta-U}
		\delta_e(U)=\min\{m\ge0:L_I^U\cap\Sigma_U^{m+1}
		\ne L_J^U\cap\Sigma_U^{m+1}\},
	\end{equation}
	with value $\infty$ if the languages agree at every horizon.

	\begin{ddefinition}[Intrinsic spatial--temporal recovery signature]
		\label{wn09:def:recovery-signature}
		Fix an atlas $\mathcal A$ on a cover with join tree $T$, with $S_e=I\cap J$ for each edge $e=IJ$.
		Using the separator defect horizon $\delta_e(U)$, define its \textbf{minimal observable guarding support} and its \textbf{full-separator defect depth} by
		\begin{align}
			g_e(\mathcal A)&=\min\{|U|:U\subseteq S_e,\ \delta_e(U)<\infty\},
			\label{wn09:eq:g-e}\\
			d_e(\mathcal A)&=\delta_e(S_e),\label{wn09:eq:d-e}
		\end{align}
		where $g_e=\infty$ if there is no disagreement on any $U$.
		The \textbf{spatial--temporal recovery signature} is \[ \mathbf{TR}_T(\mathcal A):=\bigl(
		\{(U,\delta_e(U),\rho_I(U),\rho_J(U)):
		U\subseteq S_e\}\bigr)_{e=IJ\in E(T)}.
		\]
		This is a signature indexed by original separators, not a scalar intrinsic dynamical propagation--recovery rank. In particular, $g_e=0$ is possible for an empty-separator deadlock mismatch, and $g_e=\infty$ signals equality of complete separator languages.
	\end{ddefinition}

	\begin{dproposition}[Semantic invariance and the shortest rooted refutation]
		\label{wn09:prop:semantic-invariance}
		\leavevmode
		\begin{enumerate}[label=(\roman*),leftmargin=*,nosep]
			\item \emph{Invariance.} The signature $\mathbf{TR}_T(\mathcal A)$ is invariant under
			observation-trace-equivalent replacements of local systems and under
			type-preserving renamings of chart and observable coordinates (with the
			corresponding renaming of the join tree).
			
			\item \emph{Shortest refutation.} If $d_e<\infty$, there is a
			nonempty separator word $w$ of length $d_e+1$ that is realizable on one
			side of $e$ and refuted on the other. Every proper prefix of $w$ is
			realizable on both sides. The resulting temporal nogood is
			prefix-minimal, but it need not be inclusion-minimal under arbitrary
			deletion of space--time literals.
		\end{enumerate}
	\end{dproposition}

	\begin{proof}
		\leavevmode
		\begin{enumerate}[label={(iii)},leftmargin=*,nosep]
			\item[(i)] Each entry of the signature is defined from observed trace languages,
			which are unchanged by the specified equivalences.
			
			\item[(ii)] For the second
			claim, choose a word of shortest length in the symmetric difference
			$L_I^{S_e}\triangle L_J^{S_e}$. Both languages contain $\varepsilon$,
			so the word is nonempty; its length is $d_e+1$. The language accepting
			$w$ also accepts its proper prefixes, since it is prefix closed. No
			shorter word separates the languages, so each proper prefix is accepted
			by the other side as well. This witnesses the promised temporal
			refutation with original-chart provenance. \qedhere
		\end{enumerate}
	\end{proof}

	\subsubsection{Finite observability and all-horizon globalization}
	\label{wn09:subsec:finite-observability-globalization}

	The trace languages above are regular.
	We now obtain a finite recovery cutoff which depends on their intrinsic residual indices, rather than the number of hidden states in a chosen representation.

	\begin{dtheorem}[Residual-index finite-witness theorem]
		\label{wn09:thm:finite-witness}
		Let $e=IJ$ be a tree edge and $S=I\cap J$, and write $u=\rho_I(S)$, $v=\rho_J(S)$ and $D=D(u,v)$ as in~\eqref{rev:eq:residual-union-size}.
		If $L_I^S\ne L_J^S$, the shortest separating observation word has length at most
		\begin{equation}\label{wn09:eq:bound-shortest}
			D-2\le u+v-2\le uv-1.
		\end{equation}
		Equivalently, whenever $d_e<\infty$,
		\[
		d_e \le D-3.
		\]
		Equality of the full-separator trace languages is therefore decidable by checking horizons $0\le m\le B_e:=\max\{0,D-3\}$.
	\end{dtheorem}

	\begin{proof}
		The two languages are prefix closed and both contain the empty word.
		Their minimal complete deterministic acceptors have $u,v$ states by Definition~\ref{wn09:def:trace-residual-index}.
		Apply the common rejecting-sink and acceptance-refinement estimate proved in Theorem~\ref{wn05-06:thm:finite-trace-test}.
		The length--horizon conversion subtracts one.
		Since $(u-1)(v-1)\ge0$, the linear bound is never weaker than the former product bound.
		If both languages are universal, there is no distinguishing word.
		The nonnegative cutoff convention covers this case as well.
	\end{proof}

	\begin{dtheorem}[Intrinsic finite-window recovery--globalization theorem]
		\label{wn09:thm:finite-window-globalization}
		Assume (G1) the cover $\mathcal H$ has $\nu(\mathcal H)=\infty$ and therefore admits a join tree $T$;
		(G2) each chart carries a finite, time-homogeneous relational transition system with nonempty state set (the transition relation may be empty), with no
		additional global state constraints beyond joining chart states; and
		(G3) every local transition relation is interpreted through its full trajectory language, including possible finite deadlocks.
		For $e=IJ\in E(T)$ let $S_e=I\cap J$ and set
		\[
		B(\mathcal A,T)=\max_{e=IJ\in E(T)}
		\max\{0,D(\rho_I(S_e),\rho_J(S_e))-3\},
		\]
		with $B=0$ if there are no tree edges.
		The following are equivalent:
		\begin{enumerate}[label=(\alph*),leftmargin=2.3em]
			\item The canonical natural-join global transition relation realizes every local chart path at every finite horizon.
			
			\item On every join-tree edge, full separator languages agree for every horizon $m$ with $0\le m\le B(\mathcal A,T)$.
			
			\item On every join-tree edge, full separator languages agree at every finite horizon.
		\end{enumerate}
		If $\rho_I(S_e)\le r$ on every directed separator, for some integer $r\ge2$, the uniform sufficient finite-window requirement
		\begin{equation}
			\label{wn09:eq:uniform-cutoff}
				\mathsf L_I^{S_e}(m)=\mathsf L_J^{S_e}(m)
				\quad(e=IJ\in E(T),\ 0\le m\le 2r-4)
		\end{equation}
		therefore guarantees exact globalization at every finite horizon.
		The conclusion is not asserted if additional independently prescribed global state constraints are present.
	\end{dtheorem}

	\begin{proof}
		By Definition~\ref{def:nu}, $\nu(\mathcal H)=\infty$ means chordality and conformality.
		For finite covers this is equivalent to the existence of a join tree~\cite{wn09:BeeriFaginMaierYannakakis1983}.
		A clique tree on maximal charts extends to the original indexed cover by attaching contained or duplicate charts as leaves (including empty supports).
		\begin{enumerate}[label={},leftmargin=0pt,nosep]
			\item \textnormal{(b)$\Leftrightarrow$(c).}
			For each tree edge, Theorem~\ref{wn09:thm:finite-witness}
			proves that equality through the stated finite window is equivalent to equality of the complete two projected trace languages.
			Thus (b) and (c) are equivalent.
			
			\item \textnormal{(c)$\Leftrightarrow$(a).}
			The all-horizon criterion of Theorem~\ref{thm:acyclic-temporal-exactness-criterion}, or the acyclic full-join argument on the time-expanded cover, identifies (c) with (a): a local path extends globally at each horizon
			exactly when all the tree-edge separator languages match.
			This proves	the equivalences.
		\end{enumerate}
		The $r$-uniform bound follows because two nonuniversal languages have $D-3=u+v-4\le2r-4$, while a universal and a nonuniversal language have $D-3\le r-2\le2r-4$.
		Two universal languages agree automatically.
		An extra global-state constraint need not be preserved under the joined transition map, so it is explicitly
		excluded.
	\end{proof}

	\begin{dcorollary}[Finite hidden-state and residual-index window bounds]
		\label{wn09:cor:hidden-complexity}
		Under Theorem~\ref{wn09:thm:finite-window-globalization}, if every original chart has at most $N$ hidden states, checking separator trace equality through horizon
		\[
		0\le m\le 2^{N+1}-2
		\]
		is sufficient.
		The stronger bound using the actual residual indices
		$\rho_I(S_e)$ is never worse, and can remain unchanged under arbitrary	finite hidden-state duplication.
	\end{dcorollary}

	\begin{proof}
		Apply Lemma~\ref{wn09:lem:finite-memory} to both sides of each tree edge and substitute $r=1+2^N$ into the uniform window in Theorem~\ref{wn09:thm:finite-window-globalization}.
		A hidden-state duplication that preserves the underlying observable trace language preserves its minimal residual indices.
		Arbitrarily changing the duplicated states' transitions need not preserve that language.
	\end{proof}

	\begin{dremark}[The residual index and rooted-circuit recovery]
		The finite-window condition is instance-sensitive:
		it bounds how long we must search for a separator trace defect, but does not itself force that no defect exists. The criterion also requires the finite trace-compatibility test.
		The result combines automata equivalence~\cite{wn09:MyhillNerode} with acyclic relational gluing~\cite{wn09:BeeriFaginMaierYannakakis1983}.
	\end{dremark}

	\subsubsection{Spatial guarding, temporal delay, and unguarded consequences}
	\label{wn09:subsec:spatial-temporal-examples}

	A fixed spatial cover can conceal an arbitrarily late trace disagreement.
	Conversely, a true assertion can involve distant coordinates even when every local path extends.

	\begin{dexample}[Arbitrarily late defect with one-dimensional guarding]
		\label{wn09:ex:delay-one}
		For each $N\ge3$, use the two local run-count systems of Section~\ref{sec:covers-assertions} on the
		fixed cover $\mathcal H=\{\{a,b\},\{b,c\}\}$. The chart $ab$ forbids
		$N$ consecutive observed $b=1$ values, while $bc$ forbids $N+1$.
		Both have full projected one-step $b$ transitions and agree on all
		shorter $b$ observation words. Their first separator defect is
		$d_e=N-1$, although $g_e=1$ and $\nu(\mathcal H)=\infty$.
		Their minimal trace DFAs have residual indices $N+1$ and $N+2$,
		respectively (including the rejecting sink). Hidden-state duplication
		changes the number of chart states but leaves these indices unchanged.
	\end{dexample}

	\begin{dexample}[A genuine two-coordinate temporal witness]
		\label{wn09:ex:two-coordinate}
		Use $\mathcal H=\{\{a,b,c\},\{b,c,d\}\}$ and let
		$S=\{b,c\}$.
		On the first chart, a hidden run counter forbids three consecutive observations $(b,c)=(0,0)$; on the second it forbids four.
		Both systems permit every length-one and length-two word over $\Sigma_S=\{0,1\}^2$, so they are state compatible and one-step compatible on $S$.
		All finite $b$-only and $c$-only trace words remain possible on both charts, by selecting the unobserved
		coordinate equal to one whenever necessary.
		But the length-three word $(00)(00)(00)$ is forbidden only in the first system.
		Therefore
		\[
		g_e=2,\qquad d_e=2,\qquad
			\delta_e(\{b\})=\delta_e(\{c\})=\infty.
		\]
		This distinguishes the required spatial witness size from the observation delay.
		The covered primal graph is chordal and conformal, so again $\nu(\mathcal H)=\infty$.
	\end{dexample}

	\begin{dproposition}[Exact global evolution need not guard all composed consequences]
		\label{wn09:prop:unguarded-consequence}
		There exists a finite relational atlas on the acyclic cover
		$\mathcal H=\{\{a,b\},\{b,c\}\}$ which globalizes exactly at every
		horizon, yet has a valid nontrivial temporal implication whose spatial
		support $\{a,c\}$ is not contained in either original chart.
	\end{dproposition}

	\begin{proof}
		Let all variables be binary and use the full local state spaces.
		On $\{a,b\}$ permit $((a,b),(a',b'))$ exactly when $b'=a$,
		leaving $a'$ unrestricted. On $\{b,c\}$ permit $((b,c),(b',c'))$
		exactly when $c'=b$, leaving $b'$ unrestricted. Either chart
		realizes \emph{every} finite observed $b$ word, so its separator trace
		language equals $\{0,1\}^*$. The acyclic temporal-globalization
		theorem gives exact global evolution at all horizons. But every global
		path obeys $b_1=a_0$ and $c_2=b_1$, hence $a_0=c_2$. The assertion
		$[a,0=1]\Rightarrow[c,2=1]$ is realizable and valid, but its
		spatial support $\{a,c\}$ is not guarded by either original chart.
	\end{proof}

	This example identifies the essential spatial qualifier in any proposed
	intrinsic recovery axiom. Even exact globalization does \emph{not}
	require every derived multi-chart consequence to reappear inside an
	original chart. The rank $\operatorname{ipr}$ of Paper~I~\cite{MaPackaging2026} accordingly restricts
	attention to composed rooted circuits whose combined support forms a
	primal clique. A future temporal analogue must impose a comparable
	spatial guard or an explicitly controlled enlargement of chart supports.

	The residual index of each original chart and the tree of its separator languages determine the signature $\mathbf{TR}_T(\mathcal A)$. The index bounds the length of a distinguishing word. Proposition~\ref{wn09:prop:unguarded-consequence} shows why that bound cannot by itself control the spatial support of a derived implication: two single-step updates can transmit information between coordinates that no original chart observes together. On a join tree the residual comparison and the acyclic join theorem give a finite equality test \cite{wn09:MyhillNerode,Yannakakis1981,wn09:BeeriFaginMaierYannakakis1983}.

	\subsection{Effective recognition of temporal semantic decomposability}
	\label{wn14:sec:effective-pbr}

	For a proposed width bound, failure means that an impossible partial word has no small impossible restriction. A finite automaton can recognize that failure over all word lengths.

	\subsubsection{The recognition problem}

	Section~\ref{wn13:sec:pb-r-wa} supplies a sharp globalization criterion for finite relational
	atlases whose \emph{complete} trajectory relations are projectively
	$r$-decomposable (\textsf{PB}$r$) for every finite horizon. Its hypothesis is
	semantic and independent of a preferred rule presentation; nevertheless,
	it quantifies over infinitely many horizons. We now give an exact decision
	procedure for finite time-homogeneous transition systems, including
	nondeterministic and reversible systems and non-Markovian projected
	observations. The construction uses ordinary regular-language operations.

	Fix a finite nonempty state set $X$, a transition relation
	$E\subseteq X\times X$, a finite coordinate set $I$, and a labeling
	$\lambda:X\to \prod_{v\in I}A_v$, where every $A_v$ is finite and
	nonempty. States need not be distinguished by the observation map. Fix the allowed initial subset $X_0\subseteq X$, with $X_0=X$ when no initialization is specified. Let $\mathsf T(m)$ be the relation of length-$(m+1)$ observation words $\lambda(x_0),\ldots,\lambda(x_m)$ for paths $x_0 E x_1 E\cdots E x_m$ with $x_0\in X_0$. This notation encompasses the special case
	$X\subseteq\prod_v A_v$ with coordinatewise observation.

	A \emph{partial observation} at one time is a letter in
	$\Gamma=\prod_{v\in I}(A_v\sqcup\{\star\})$; $\star$ means unspecified.
	A partial word $w\in\Gamma^{m+1}$ is \emph{feasible} if an actual
	trajectory agrees with its specified cells. We write
	$\operatorname{supp}(w)$ for its specified space--time cells and
	$w|_B$ for the word obtained by erasing every specified cell outside
	$B\subseteq\operatorname{supp}(w)$.

	\begin{ddefinition}[Temporal semantic defects and least defect horizon]
		\label{wn14:def:temporal-defect}
		Let $\Gamma$ be the partial-word alphabet of an observed
		system and fix an integer $r\geq1$.
		For $w\in\Gamma^+$, let $\operatorname{supp}(w)$ be its
		set of specified observable cells, and let $w|_B$
		retain only the cells in $B$ without changing the horizon.
		\begin{enumerate}[label=(\roman*),leftmargin=*,nosep]
			\item \textbf{Temporal semantic defect.}
			A partial word $w$ is a temporal semantic defect
			of order $r$ if it is infeasible but every restriction to at most $r$ specified cells is feasible.
			The corresponding defect language is
			\begin{equation}\label{wn14:eq:bad-r}
				\mathsf{Bad}_r:=
				\bigl\{w\in\Gamma^+:
				w\text{ is infeasible},\
				w|_B\text{ is feasible for every }
				B\subseteq\operatorname{supp}(w),\ |B|\leq r\bigr\}.
			\end{equation}
			
			\item \textbf{Least temporal defect horizon.}
			The least temporal defect horizon of order $r$ is
			\[
			\delta_r:=\min\{|w|-1:w\in\mathsf{Bad}_r\},
			\]
			where $|w|$ is the word length and
			$\min\varnothing=\infty$.
		\end{enumerate}
		For fixed observable coordinates, both $\mathsf{Bad}_r$
		and $\delta_r$ depend only on the observed trajectory
		language.
	\end{ddefinition}

	\begin{dlemma}[Partial-word characterization]
		\label{wn14:lem:partial-pbr}
		Assume $\mathsf T(m)\ne\varnothing$.
		Then $\mathsf T(m)$ is \textsf{PB}$r$ in the sense of Section~\ref{wn13:sec:pb-r-wa} (including every coordinate projection) if and only if no word of length $m+1$ belongs to $\mathsf{Bad}_r$.
	\end{dlemma}

	\begin{proof}
		A partial observation $w$ specifies a tuple on some subset $U$ of space--time coordinates.
		It is feasible exactly when this tuple belongs
		to the projection $\pi_U\mathsf T(m)$.
		By the definition of \textsf{PB}$r$, every tuple whose restrictions to subsets of size at most $r$ are all feasible must itself be feasible.
		Negating this assertion gives precisely~\eqref{wn14:eq:bad-r}.
	\end{proof}

	\subsubsection{An exact finite-state recognizer}

	For a partial letter $o\in\Gamma$, write
	$M(o)=\{x\in X:\lambda(x)\text{ agrees with }o\}$. The finite-state
	observer starts in $X_0$, filters at time zero by $M(o_0)$, and at every
	later time applies $\operatorname{Post}_E$ followed by filtering:
	\begin{equation}\label{wn14:eq:observer}
		Q_0=X_0\cap M(o_0),\qquad
		Q_{t+1}=\operatorname{Post}_E(Q_t)\cap M(o_{t+1}).
	\end{equation}
	The word is feasible if and only if its final observer subset is
	nonempty. This is a deterministic automaton with at most $2^{|X|}$
	observer subsets, plus a start convention.

	For an input partial word $w$, a second finite automaton guesses a
	\emph{mask}: at each position it may retain or erase each specified
	coordinate, spending one unit of a total budget of $r$ whenever a cell
	is retained. The guessed masked word is passed through the same
	subset observer. Acceptance occurs when a mask of cost at most $r$
	produces an empty final subset. Thus the second automaton recognizes
	\begin{equation}\label{wn14:eq:small-bad}
		\mathsf{SmallBad}_r=\{w:\exists B\subseteq\operatorname{supp}(w),
		|B|\le r,\ w|_B\text{ infeasible}\}.
	\end{equation}
	Both languages are regular. In particular,
	\begin{equation}\label{wn14:eq:lang-bad}
		\mathsf{Bad}_r=
		\mathsf{Infeasible}\setminus\mathsf{SmallBad}_r
	\end{equation}
	is effectively regular.

	\begin{dtheorem}[Effective all-horizon \textsf{PB}$r$ recognition]
		\label{wn14:thm:recognition}
		\leavevmode
		\begin{enumerate}[label=(\roman*),leftmargin=*,nosep]
			\item \emph{Recognition.}
			For each fixed integer $r\ge1$, there is a terminating algorithm
			that, given a finite transition relation, a finite observation map, and an initial set $X_0\subseteq X$ (default $X$), decides whether every nonempty complete finite-horizon observation relation is \textsf{PB}$r$.
			If not, it returns a finite-horizon partial observation which is impossible but all of whose $\le r$-cell restrictions are possible.
			
			\item \emph{Finite bound.}
			Let $n=|X|$. When a counterexample exists, one exists with fewer than
			\begin{equation}\label{wn14:eq:finite-window}
				N(n,r)=2^{\,1+n+(r+1)2^n}
			\end{equation}
			observation time slices.
			
			\item \emph{Presentation invariance.}
			The recognized semantic property, the least defect horizon, and the recognition verdict are invariant under finite-state presentations with the same named-coordinate observation language.
			Monitor size and running time may depend on the chosen presentation.
		\end{enumerate}
	\end{dtheorem}

	\begin{proof}
		\leavevmode
		\begin{enumerate}[label={(iii)},leftmargin=*,nosep]
			\item[(i)] For the decision statement, construct automata for $\mathsf{Infeasible}$ and $\mathsf{SmallBad}_r$ as above, complement the
			latter, intersect the two, and perform a finite-state emptiness test.
			This recognizes precisely $\mathsf{Bad}_r$. If it is nonempty, a
			shortest accepting word supplies the required witness. The preceding
			lemma makes emptiness equivalent to \textsf{PB}$r$ at every nonempty
			horizon.
			
			\item[(ii)] For the explicit size bound, determinize the mask automaton while
			tracking the full observer state $Q\subseteq X$. For each
			$j\in\{0,\ldots,r\}$, record the set $\mathcal M_j\subseteq 2^X$ of
			all final reachable-state subsets obtained from masks that retained
			exactly $j$ specified cells. The current deterministic monitor state
			consists of a start flag, $Q$, and the $r+1$ sets $\mathcal M_j$.
			There are at most $2\cdot 2^n(2^{2^n})^{r+1}=N(n,r)$ possible monitor
			states. After at least one input letter, a state is accepting exactly when
			$Q=\varnothing$ but $\varnothing\notin\mathcal M_j$ for every
			$j\le r$. The start flag excludes the empty word. Any accepting state reachable from the start is reachable
			by a simple path of length less than $N(n,r)$, proving the finite
			witness bound.
			
			\item[(iii)] Finally, feasibility of all partial observation words
			is determined by the observed trajectory language alone, so
			$\mathsf{Bad}_r$, $\delta_r$, and the recognition verdict are intrinsic
			to that language. \qedhere
		\end{enumerate}
	\end{proof}

	By Definition~\ref{def:temporal-witness-width}, the intrinsic temporal witness width is the maximum minimal-impossible support, with lower convention one. The equivalent fixed-bound formulation is
	\[
	\operatorname{twit}(X,E,\lambda)
	=\inf\{r\ge1:\mathsf{Bad}_r=\varnothing\}
	\in\mathbb N\cup\{\infty\},
	\]
	as proved in Lemma~\ref{wn15:lem:minimal-rank}. The formula also holds for empty path horizons under the support-zero convention; the nonemptiness gate for globalization is separate. It is distinct from the spatial guarding and fiber-Helly recovery ranks of Sections~\ref{wn12:sec:pb2-rel-global} and~\ref{wn13:sec:pb-r-wa}.

	The width is invariant under the named-coordinate language equivalence and nonincreasing under coordinate forgetting (Proposition~\ref{wn15:prop:invariance}). Theorem~\ref{wn14:thm:recognition} decides each fixed width bound; Theorem~\ref{wn15:thm:effective-exact-rank} decides finiteness and computes the least such bound.

	\begin{dremark}[Finite paths and the nonemptiness gate]
		\label{wn14:rem:nonempty-gate}
		Section~\ref{wn13:sec:pb-r-wa} defines \textsf{PB}$r$ for nonempty relations. For a specified initial set $X_0$, $\mathsf T(m)$ is nonempty for \emph{every} $m$ precisely when some directed cycle is reachable from $X_0$ in the finite transition digraph. For $X_0=X$ this is simply the existence of a directed cycle. This is effectively testable. A finite terminating
		system with an empty path relation at a later horizon is therefore
		reported as \emph{outside the Section~\ref{wn13:sec:pb-r-wa} class}, not as a positive
		certificate. The stronger condition that every local state has a
		successor is also decidable, but neither condition implies global
		nonblocking after chart gluing.
	\end{dremark}

	\paragraph{Examples of partial-word obstructions.}

	Independent bit flips and an XOR update have different minimal forbidden partial words. Their comparison shows what the recognition automaton measures.

	\begin{dexample}[Independent flips and XOR updates]
		\label{wn14:ex:flip-xor}
		On a single Boolean coordinate, the identity and bit-flip maps have
		complete trajectory relations of semantic width $2$: all time copies
		of the coordinate are determined by any one initial value, and any
		incompatible pair of observations is a two-cell nogood. Two
		coordinatewise independent bit flips are also temporally
		\textsf{PB2}. In contrast, the deterministic two-bit update
		\[
		(a_{t+1},b_{t+1})=(a_t\oplus b_t,b_t)
		\]
		already contains the ternary parity constraint
		$a_{t+1}\oplus a_t\oplus b_t=0$ at horizon one, so it fails the
		\textsf{PB2} recognition test. It is unnecessary to test later
		horizons to refute it.
	\end{dexample}

	\begin{dexample}[A finite-state system of unbounded temporal witness width]
		\label{wn14:ex:marker-parity}
		Consider the six-state labelled system (all states initially allowed)
		$X=\{s,e\}\cup\{d_{b,p}:b,p\in\{0,1\}\}$. The visible labels of
		$s,e,d_{b,p}$ are respectively $S,E,b$; the bit $p$ is a hidden
		parity register. Its edges are
		\[
		s\to d_{b,b},\quad
		d_{b,p}\to d_{c,p\oplus c},\quad
		d_{b,0}\to e,\quad e\to e
		\qquad(b,c,p\in\{0,1\}).
		\]
		This finite graph is nonblocking. For every $k\ge0$ the visible
		partial observation
		\begin{equation}\label{wn14:eq:unbounded-no-good}
			S\;1\;0^k\;E
		\end{equation}
		is impossible: the register has odd parity immediately before $E$.
		Erasing the leading $S$ permits a different starting parity;
		erasing the final $E$ permits continued data evolution; erasing any
		data symbol permits toggling that symbol to satisfy even parity.
		Hence~\eqref{wn14:eq:unbounded-no-good} is a minimal partial semantic
		nogood with support $k+3$. There is no finite $r$ for which all
		visible trajectory relations are \textsf{PB}$r$, even though the
		underlying system has only six states and is time homogeneous.
		This is stronger than simply observing that minimal forbidden
		\emph{contiguous} blocks may be long: here minimality concerns
		\emph{all} deletions of specified observation cells.
	\end{dexample}

	\begin{dexample}[A positive higher-arity system beyond one-step certification]
		Let $X\subseteq\{0,1\}^3$ consist of even-parity triples and let
		evolution be identity on $X$. Its temporal relations are
		\textsf{PB3}: the time-copy equalities are binary, and the three
		initial bits obey one ternary parity equation. They fail
		\textsf{PB2} already at horizon zero. The exact recognizer produces
		\emph{FAIL} at $r=2$ and a completed \emph{PASS} at $r=3$; this does
		not require guessing a near-unanimity polymorphism.
	\end{dexample}

	\subsubsection{A finite recognition pipeline for exact globalization}

	Let $\mathcal H$ be a finite package-chart cover. For each original
	chart $I$, suppose the full observed state is exactly the chart's coordinate tuple, and its complete trajectory relation consists of all paths of the declared first-order state and edge relations. Auxiliary hidden-state presentations may be used to compute these same relations and languages, but arbitrary hidden projections need not satisfy this first-order path condition. Let $q=\max_{v\in V}|A_v|$.

	\begin{dcorollary}[Effective globalization certificate]
		\label{wn14:cor:effective-global}
		Fix $r,q\ge2$, with $|A_v|\le q$ for every coordinate. Suppose the finite chart transition systems are nonblocking, their complete trajectory relations consist of all paths of the declared first-order relations on the full chart coordinate tuples, and there are no independent global-state constraints.
		The following conditions are jointly decidable:
		\begin{enumerate}[label=(\roman*),leftmargin=2.1em]
			\item every complete chart trajectory relation is \textsf{PB}$r$ at
			every finite horizon;
			
			\item on every chart intersection, the full projected trajectory
			languages agree at every horizon;
			
			\item $1+q(r-1)<\nu(\mathcal H)$.
		\end{enumerate}
		If all three hold, the canonical global transition relation realizes
		every prescribed local finite trajectory exactly.
	\end{dcorollary}

	\begin{proof}
		\leavevmode
		\begin{enumerate}[label={(iii)},leftmargin=*,nosep]
			\item[(i)] Condition (i) is decidable by Theorem~\ref{wn14:thm:recognition} and the
			nonemptiness gate.
			\item[(ii)] For condition (ii), use the deterministic
			subset-observers for each chart's trace language and decide their
			equality using a finite product automaton, as in Section~\ref{sec:width}.
			\item[(iii)] Condition (iii) is a finite combinatorial check on the cover's
			primal graph and guarded cliques.
		\end{enumerate}
		The sharp \textsf{PB}$r$ temporal globalization theorem (Theorem~\ref{wn13:thm:temporal-pbr}) then applies to the complete
		space--time trajectory atlases. Their canonical natural join is the
		path relation of the \emph{same} global one-step relation for every
		horizon, proving the claim.
	\end{proof}

	\begin{dremark}[Scope of the globalization certificate]
		The certificate is sufficient for an individual atlas; its strict cutoff is sharp for the universal class. The canonical one-step construction requires the first-order path hypothesis: an arbitrary projected regular language can contain every permitted consecutive pair of a word without containing the word itself (as the even-shift example, Proposition~\ref{wn08:prop:even-shift}, illustrates). It must instead be treated as history data or with a suitable declared state refinement. Nonblocking from every compatible state and exactness under additional independent global constraints require separate conditions. Arbitrary finite relational atlases can instead
		be decided by the canonical-global-relation procedure of Section~\ref{wn11:sec:nonhorn-relational}.
	\end{dremark}

	For a fixed $r$, the monitor tests all horizons in one finite reachability problem. The next construction computes the least successful width directly, by measuring the support of a minimal forbidden pattern.

	\subsection{Intrinsic temporal witness ranks and structural recognition}
	\label{wn15:sec:intrinsic-recognition}

	A monitor for minimal impossible words also computes the least width. Giving each transition the number of observed cells that it reads reduces the problem to detecting positive-weight cycles.

	\subsubsection{The rank-existence problem}

	Section~\ref{wn13:sec:pb-r-wa} will establish a sharp non-Horn globalization criterion for
	projectively $r$-decomposable (\textsf{PB}$r$) complete trajectory relations. Section~\ref{wn14:sec:effective-pbr} has recognized the all-horizon property for every \emph{fixed} $r$, but does not decide whether there exists \emph{any} finite uniform $r$. We solve the existence and exact-rank problems for arbitrary finite-state
	observed relational transition systems, without a fixed bound as input.
	
	Let $I$ be a finite observed-coordinate set, let $A_i$ be finite
	nonempty alphabets, and let $(X,E,\lambda)$ be a finite labelled
	transition system: $X$ is finite, $E\subseteq X\times X$, and
	$\lambda:X\to\prod_{i\in I}A_i$ is the observation map.
	Fix the initial subset $X_0\subseteq X$ as in Section~\ref{wn14:sec:effective-pbr}, with default $X_0=X$. The choice of $X$ and $\lambda$ may include internal, unobserved states.
	Write $d=|I|$ and $n=|X|$. To align with Sections~\ref{wn13:sec:pb-r-wa} and~\ref{wn14:sec:effective-pbr} we usually assume
	that paths of every finite length exist, which for a finite graph amounts to a directed cycle reachable from $X_0$.
	The rank algorithm itself also	handles systems without this property, but an empty horizon then fails the nonemptiness hypothesis of Theorem~\ref{wn13:thm:temporal-pbr}.
	
	Set
	\[
	\Gamma=\prod_{i\in I}(A_i\sqcup\{\star\}),\qquad
	\operatorname{wt}(a)=|\{i:a_i\ne\star\}|.
	\]
	For a nonempty partial observation word
	$p=a_0\cdots a_m\in\Gamma^{m+1}$, let
	$\operatorname{supp}(p)$ be its set of specified space--time cells;
	thus $|\operatorname{supp}(p)|=\sum_t\operatorname{wt}(a_t)$.
	It is \emph{feasible} when some path
	$x_0 E x_1 E\cdots E x_m$ with $x_0\in X_0$ agrees with every specified observation.
	Replacing one specified cell by $\star$ is a \emph{one-cell relaxation}
	$p^{-e}$, where $e\in\operatorname{supp}(p)$.
	Feasibility is downward closed under relaxation.

	\begin{ddefinition}[Minimal temporal nogoods]
		\label{wn15:def:minimal-nogood}
		A partial word $p$ of positive word length (possibly with zero specified cells) is a \textbf{minimal semantic temporal nogood} if $p$ is infeasible and $p^{-e}$ is feasible for every $e\in\operatorname{supp}(p)$.
		Write $\mathcal M(\mathbb D)$ for the language of all such words at all finite horizons.
		Under the nonemptiness hypotheses of the globalization applications, every all-$\star$ word is feasible and minimal nogoods have positive support.
		For the general recognition algorithm, an empty horizon instead	has an all-$\star$ minimal nogood of support zero.
		The rank convention still returns at least one, and a separate nonemptiness test is required.
	\end{ddefinition}

	\begin{dlemma}[Minimal-nogood characterization of the witness width]
		\label{wn15:lem:minimal-rank}
		For the intrinsic temporal witness width $\operatorname{twit}$ of Section~\ref{wn14:sec:effective-pbr},
		\begin{equation}\label{wn15:eq:rank-minimal}
			\operatorname{twit}(\mathbb D)=\max\left\{1,
			\sup_{p\in\mathcal M(\mathbb D)}|\operatorname{supp}(p)|\right\},
			\qquad\sup\varnothing:=0.
		\end{equation}
	\end{dlemma}

	\begin{proof}
		If an infeasible partial word $p$ has finite support, repeatedly erase specified cells while preserving infeasibility.
		This terminates at a minimal nogood $p'\preceq p$. Conversely, a minimal nogood with more
		than $r$ specified cells is infeasible while every subobservation with at most $r$ specified cells is feasible:
		any proper relaxation is feasible, since it relaxes at least one specified cell.
		Hence the $\mathsf{Bad}_r$ language from Section~\ref{wn14:sec:effective-pbr} is nonempty if and only if a minimal nogood has support exceeding $r$. Taking the least $r\ge1$ proves the formula.
	\end{proof}

	\subsubsection{The single-deletion observer}

	Let $\mathcal Q=2^X$. For an observation letter $a\in\Gamma$,
	the usual subset observer sends $Q\subseteq X$ to the set of states
	reachable in one transition and matching $a$; at the \emph{first}
	letter it only intersects the set of allowed initial states with the
	observation constraint. Denote this operation by $\delta_t(Q,a)$,
	where $t=0$ indicates the initial letter. Its extension to partial
	words sends $p$ to the final reachable subset $Q(p)$, and
	\[
	p\text{ is feasible}\quad\Longleftrightarrow\quad Q(p)\ne\varnothing.
	\]
	
	For each word $p$, define the \emph{single-deletion family}
	\[
	U(p)=\{Q(p^{-e}):e\in\operatorname{supp}(p)\}
	\subseteq2^X.
	\]
	At the start, $Q(\epsilon)=X_0$ and $U(\epsilon)=\varnothing$.
	If $a\in\Gamma$ is appended to $p$, then
	\begin{align}
		Q(pa)&=\delta_t(Q(p),a),\label{wn15:eq:single-full}\\
		U(pa)&=\{\delta_t(S,a):S\in U(p)\}\nonumber\\
		&\quad\cup
		\{\delta_t(Q(p),a^{-{i}}):i\in I,\ a_i\ne\star\},
		\label{wn15:eq:single-del}
	\end{align}
	where $a^{-{i}}$ replaces its $i$th entry by $\star$ and $t=0$
	only if $p=\epsilon$. The first term propagates a deletion at an
	earlier time; the second introduces the unique deletion at the current
	time. Thus the pair $(Q(p),U(p))$ is computable in one finite pass.

	\begin{dtheorem}[Regularity and state bound for minimal temporal nogoods]
		\label{wn15:thm:minimal-regular}
		The language $\mathcal M(\mathbb D)$ is effectively regular.
		A deterministic automaton recognizing it has at most
		\begin{equation}\label{wn15:eq:size}
			N(n)=2^{\,1+n+2^n}
		\end{equation}
		states (before removing unreachable states).
		The result holds with arbitrary hidden state and coordinate projection.
	\end{dtheorem}

	\begin{proof}
		Track the finite memory $(\textit{first},Q,U)$, where $Q\subseteq X$,
		$U\subseteq2^X$, and the Boolean flag marks whether the next letter is
		the first. Update by~\eqref{wn15:eq:single-full}--\eqref{wn15:eq:single-del}.
		Induction on word length proves the invariant that $Q$ represents the
		original word while $U$ collects the terminal state subsets for
		\emph{all} one-cell relaxations. Therefore the acceptance condition
		\begin{equation}\label{wn15:eq:accept}
			\textit{first}=0,\qquad Q=\varnothing,\qquad
			\varnothing\notin U
		\end{equation}
		holds precisely for minimal infeasible partial words.
		There are at most $2$ choices of the flag, $2^n$ choices of $Q$, and
		$2^{2^n}$ choices of $U$, yielding~\eqref{wn15:eq:size}. Observation
		compatibility with a hidden-state presentation changes only the
		observer transition function, not finiteness of its memory.
	\end{proof}

	\noindent\textbf{Convention for the empty word.}
	The first observation letter restricts the specified allowed initial-state set (possibly empty for rank recognition); subsequent letters impose actual transition steps.
	A wildcard removes one named coordinate observation, not an entire
	instantaneous symbol or time slice. The accepting-cycle criterion applies
	only to reachable, accepting-coaccessible monitor states. Exact rank
	invariance is asserted only for coordinatewise alphabet relabelings
	that preserve the set of named observable atoms; arbitrary state
	recodings may change the rank.

	\subsubsection{Exact recognition of the finite--infinite dichotomy}

	Let $\mathcal A_{\min}$ denote the reachable part of this automaton.
	Give each transition bearing letter $a$ the nonnegative integer
	weight $\operatorname{wt}(a)$. A state is \emph{coaccessible} if an
	accepting state is reachable from it.

	\begin{dtheorem}[Decidability and exact computation of temporal witness width]
		\label{wn15:thm:effective-exact-rank}
		For every finite-state initialized labelled transition system with a finite set of named observable coordinates and finite nonempty coordinate alphabets, as specified in this subsection, whether
		$\operatorname{twit}(\mathbb D)<\infty$ is decidable.
		More precisely, the following are equivalent:
		\begin{enumerate}[label=(\roman*),leftmargin=2em]
			\item $\operatorname{twit}(\mathbb D)=\infty$;
			
			\item $\mathcal A_{\min}$ has a positive-weight directed cycle
			that is reachable from its start and can reach an accepting state;
			
			\item there exist partial words $u,v,w$, with
			$|\operatorname{supp}(v)|>0$, such that
			$uv^kw\in\mathcal M(\mathbb D)$ for every $k\ge0$.
		\end{enumerate}
		If these conditions fail, then
		\begin{equation}\label{wn15:eq:finite-exact}
				\boxed{\operatorname{twit}(\mathbb D)=\max\left\{1,
				\max_{p\in\mathcal M(\mathbb D)}|\operatorname{supp}(p)|\right\}},
		\end{equation}
		where the maximum over an empty nogood language is zero. The maximum is computable by a weighted longest-path calculation on the accepting portion of the automaton after contracting strongly connected components; if no accepting state is reachable, the width is one.
	\end{dtheorem}

	\begin{proof}
		\leavevmode
		\begin{enumerate}[label={},leftmargin=0pt,nosep]
			\item \textnormal{(ii)$\Rightarrow$(iii)$\Rightarrow$(i).}\enspace A positive-weight cycle on a start-to-accept path yields words
			$uv^kw$ for all $k\ge0$ by repeatedly traversing the cycle. Every
			word is accepted by $\mathcal A_{\min}$, so is a minimal nogood,
			and its support equals
			$|\operatorname{supp}(u)|+k|\operatorname{supp}(v)|+
			|\operatorname{supp}(w)|$, which grows without bound.
			\item \textnormal{(i)$\Rightarrow$(ii).}\enspace Conversely, restrict the reachable automaton to vertices from which
			acceptance remains possible. If there is no positive-weight cycle,
			then every cycle in this subgraph has total weight zero, since all
			transition weights are nonnegative. Hence every edge internal to a
			strongly connected component has weight zero. After contracting these
			components, the remaining weighted graph is acyclic, so the weights
			of accepting paths admit an effectively computable maximum. By
			Lemma~\ref{wn15:lem:minimal-rank}, this maximum determines the exact
			$\operatorname{twit}$. If the maximum did not exist, finite-graph
			pigeonhole would force a positive-weight cycle on an accepting path.
			This proves all equivalences and the algorithm. \qedhere
		\end{enumerate}
	\end{proof}

	\begin{dcorollary}[Finite upper bound, independent of a guessed rank]
		\label{wn15:cor:rank-bound}
		If the intrinsic width is finite, then
		\begin{equation}\label{wn15:eq:bound}
			\operatorname{twit}(\mathbb D)
			\le \max\bigl\{1,\ d\bigl(N(n)-1\bigr)\bigr\}.
		\end{equation}
		In particular, a terminating exact algorithm needs no external search
		parameter $r$.
	\end{dcorollary}

	\begin{proof}
		If no accepting walk exists, the width is one by convention. Otherwise,
		a maximum-weight accepting walk can be shortened by removing its
		zero-weight cycles without reducing its weight. It therefore has a
		simple representative of at most $N(n)-1$ edges, each of weight at
		most $d$, proving the bound.
	\end{proof}

	\begin{dremark}[Computational complexity]
		The decision procedure explores the reachable single-deletion automaton, detects positive-weight cycles on accepting paths, and either produces a pumping witness or computes the maximum support of a minimal nogood. The stated state bound is doubly exponential in $|X|$, and the partial-observation alphabet has $\prod_{i\in I}(|A_i|+1)$ letters. A search stopped before completion gives no certified verdict.
	\end{dremark}

	\subsubsection{Semantic invariance and observation restriction}

	A witness is a set of observed cells. Relabeling values transports its support, while deleting coordinates restricts the possible witnesses.

	\begin{dproposition}[Semantic invariance and restriction monotonicity]
		\label{wn15:prop:invariance}
		The value of $\operatorname{twit}$, whether finite or infinite,
		only depends on the complete language of observable state sequences,
		not on any particular hidden-state presentation. It is unchanged
		by coordinatewise alphabet bijections preserving the named observable atoms,
		and cannot increase when some observed
		coordinates are hidden.
	\end{dproposition}

	\begin{proof}
		The minimal-nogood language is defined entirely through the
		feasibility of partial observations, which itself is determined by
		the observable trajectory language. A renaming bijection preserves
		feasibility and support. After discarding some coordinates, a partial
		observation is feasible exactly when the corresponding pattern with
		the discarded coordinates wildcarded is feasible in the original
		system. Any minimal forbidden projected pattern is therefore also a
		minimal forbidden original pattern, establishing monotonicity.
	\end{proof}

	\paragraph{Examples.}

	For the identity and flip on one bit, two incompatible time cells witness every failure. An unconstrained full shift has width one under the convention in Definition~\ref{def:temporal-witness-width}. The parity-marker observer in Section~\ref{wn14:sec:effective-pbr} instead admits the following family of minimal impossible words.
	
	\[
	S\,0\,0\,\bigl(0\bigr)^k\,1\,E
	\qquad(k\ge0).
	\]
	Each member is impossible, yet relaxing any one specified visible
	symbol makes it feasible. The support grows as $k+5$. This is a
	complete finite-state certificate of infinite width, not an
	extrapolation from a finite list of long examples.

	\subsubsection{Computing the guard threshold}

	For a finite dynamical chart atlas, compute the exact width
	$r_I=\operatorname{twit}(\mathbb D_I)$ for each chart. If every
	$r_I$ is finite, set $r=\max\{2,\max_I r_I\}$ and
	$q=\max(\{2\}\cup\{|A_v|:v\in V\})$. Assume all local trajectory relations are nonempty, chart-overlap trace languages agree at every horizon, and no independent global-state constraints are imposed. As in Corollary~\ref{wn14:cor:effective-global}, require that the complete local trajectory relations consist of all paths of the declared first-order state and transition relations on the chart coordinate tuples.
	Then the sharp PB$r$ theorem of Section~\ref{wn13:sec:pb-r-wa} gives the sufficient
	criterion
	\begin{equation}\label{wn15:eq:global-cert}
		1+q(r-1)<\nu(\mathcal H)
		\quad\Longrightarrow\quad
		\text{exact all-horizon relational globalization}.
	\end{equation}
	Exact temporal trace equivalence and the finite primal-graph bound are
	decidable by the earlier results, so the entire certificate is
	effective without prescribing $r$ in advance.

	This is not a universal globalization rank: a local chart may have
	infinite temporal witness width while a particular atlas still
	globalizes, and the displayed strict inequality is sufficient for an
	atlas rather than necessary for each individual atlas. The
	canonical-global-relation method from Section~\ref{wn11:sec:nonhorn-relational} remains the
	individual-atlas decision procedure for first-order chart data. For arbitrary projected regular histories the bound guarantees exact joins at each horizon, not realization by the greatest one-step relation on the visible tuples. For instance the run counter in Theorem~\ref{wn08:thm:unbounded-temporal-delay} with $N=3$, observed only through $b$, forbids $111$ but permits every one- and two-letter binary word. Its witness width is three (complete unspecified cells by zeros). Even on a single chart the visible one-step relation is full and therefore introduces the forbidden path $111$.

	The recognition arguments above use elementary regular-language closure and
	nonnegative weighted-graph boundedness arguments. Partial-word
	languages and decomposability of regular relations have substantial
	prior literature, including \cite{wn15:Dassow2014Partial,wn14:BarceloEtAl2019,wn15:BlanchetSadri2017Minimal}.

	\section{Recovery--globalization correspondences}
	\label{sec:global-theorems}

	Two extension problems lead to different sharp thresholds. For positive events, Theorem~\ref{thm:all-horizon-causal-Horn-globalization} constructs a least persistent causal closure on the nonnegative time axis from local semantic implication rules. For finite-valued relations, Theorem~\ref{thm:pbr} extends every original chart tuple using the natural join and geometric guarding of small obstructions. These are different semantic constructions; neither identifies implication determination with exclusion of a partial assignment. The later interface theorems examine witness width under assembly and feedback.

	\paragraph{Relation to local consistency.}
	Dechter obtained a domain-and-arity bound for extending locally consistent assignments in constraint networks. Dechter and van Beek developed the corresponding relational consistency methods \cite{Dechter1992,DechterVanBeek1997}. Here the bound on a chart applies after projection onto every subset of its coordinates. The conclusion concerns the projection of the join onto each original chart. In the Horn case the spatial threshold comes from the static closure theorem of Ma~\cite{MaPackaging2026}. Lemmas~\ref{wn10:lem:native-horn} and~\ref{wn10:lem:prefix-global} supply the temporal basis and its coherence under truncation.

	\paragraph{Projectively decomposable chart relations.}

	A partial assignment that cannot occur may already fail on a small subset of its cells.
	Requiring such a witness for every projection gives a condition that survives elimination of hidden coordinates.

	\begin{definition}[Projective $r$-decomposability]
		\label{def:projective-r-decomposability}
		Fix integers $q,r\ge1$.
		A nonempty finite chart relation $R\subseteq\prod_{v\in I}A_v$ is \textbf{projectively $r$-decomposable}, abbreviated $\mathsf{PB}r$, if every partial assignment on $U\subseteq I$ that fails to extend to $R$ has a restriction to at most $r$ coordinates that also fails. For every $U\subseteq I$, this is the projection identity
		\begin{equation}\label{wn13:eq:pb-r-def}
			\pi_U R=\left\{a\in\prod_{v\in U}A_v: a|_S\in\pi_S R\text{ for every }S\subseteq U,\ 1\le|S|\le r\right\}.
		\end{equation}
		For $U=\varnothing$ the identity uses nonemptiness of $R$.
		Write $\mathscr R_{q,r}$ for the class of these relations with $1\le |A_v|\le q$.
		In classwise atlas statements and recovery indices, the same notation denotes atlases with all local relations in this class.
	\end{definition}

	Theorem~\ref{thm:pbr} proves the exact classwise cutoff $1+q(r-1)<\nu(\mathcal H)$ when the coordinate alphabets have at most $q$ values, $q,r\ge2$.
	A single excluded value costs at most $r-1$ previously assigned coordinates.
	Guarding the union of at most $q$ such witnesses allows the next value to be chosen.
	Lemma~\ref{wn13:lem:ung-capacity} constructs the converse obstruction.
	The time-expanded form is Theorem~\ref{wn13:thm:temporal-pbr}.

	\paragraph{Composition, interfaces and temporal feedback.}

	The synchronized-admission construction of Definition~\ref{dyn46:def-admission-dynamics} provides the native structural interpretation when the separator and module states are declared operational states.
	The results below instead use explicitly specified relational interfaces and witness semantics.
	A union or sum of independently clocked mechanisms cannot generally replace synchronous composition (Example~\ref{dyn46:ex-clock-noncollapse}).
	
	The comparison below uses the declared synchronous interfaces and the corresponding witness families.

	\begin{example}[Feedback can destroy finite width]\label{ex:feedback}
		On binary coordinates, synchronizing a delayed equality $a_{t+1}=b_t$ with the instantaneous parity $b_t=a_t\oplus c_t$ yields an integrator $a_{t+1}=a_t\oplus c_t$ (addition modulo two). Finite component witness widths can thereby become infinite; the precise witness family and proof are given in Theorem~\ref{wn16:thm:feedback}. This concerns synchronization rather than an overlap-compatible atlas.
	\end{example}

	\subsection{Persistent causal closures and their least extension}
	\label{wn10:sec:dyn-recovery-global}

	A local rule records an event forced by observations at the same time or earlier. The semantic determination bound gives a canonical rule basis. This lets the finite-prefix extensions fit together into one least closure on infinite histories.

	\subsubsection{The persistent causal Horn class}

	Throughout, $V$ is a finite nonempty descriptor set, $\mathcal H$ is a finite cover of $V$, and $\mathbb N_0=\{0,1,\ldots\}$.
	Write
	\[
	W_{I,m}=I\times\{0,\ldots,m\},\qquad
	W_{I,\infty}=I\times\mathbb N_0,\qquad
	\sigma_h(i,t)=(i,t+h).
	\]
	The original packaged object retains its carriers, incidence and selected families.
	The descriptor cover and its rule-bearing charts are additional specialization data.
	All closure restrictions have the precise meaning
	$c|_S(A)=c(A)\cap S$ for $A\subseteq S$, as in the static closure theory \cite{MaPackaging2026}.

	\begin{ddefinition}[Stationary persistent causal closures, compatible atlases, and exact extensions]
		\label{wn10:def:causal-closure}
		Fix a finite indexed cover $\mathcal H$ of $V$ and
		an integer $k\geq1$.
		Write $W_{I,m}:=I\times\{0,\ldots,m\}$ and
		$W_{I,\infty}:=I\times\mathbb N_0$.
		\begin{enumerate}
			\item \textbf{Local causal closure.}
			A \textbf{stationary persistent causal closure} of uniform intrinsic determination order at most $k$ on chart $I$ is an extensive, monotone, and idempotent map $\gamma_I:2^{W_{I,\infty}}\to2^{W_{I,\infty}}$
			satisfying:
			\begin{enumerate}[label=(\roman*),leftmargin=*,nosep]
				\item \textbf{Empty preservation.}
				$\gamma_I(\varnothing)=\varnothing$.
				
				\item \textbf{Past causality.}
				For every $A\subseteq W_{I,\infty}$ and $m\geq0$,
				\[
				\gamma_I(A)\cap W_{I,m}
				=\gamma_I(A\cap W_{I,m})\cap W_{I,m}.
				\]
				
				\item \textbf{Stationarity.}
				For every $A\subseteq W_{I,\infty}$ and $h\geq0$,
				\[
				\gamma_I(\sigma_h A)=\sigma_h\gamma_I(A),
				\]
				where $\sigma_h(i,t):=(i,t+h)$.
				
				\item \textbf{Persistence.}
				For every $A\subseteq W_{I,\infty}$,
				\[
				(i,t)\in\gamma_I(A)
				\Longrightarrow(i,t+1)\in\gamma_I(A).
				\]
				
				\item \textbf{Uniform intrinsic determination.}
				For every $m\geq0$, $A\subseteq W_{I,m}$, and
				$q\in\gamma_I(A)\cap W_{I,m}$, there exists
				$P\subseteq A$ such that
				\[
				|P|\leq k,\qquad q\in\gamma_I(P).
				\]
			\end{enumerate}
			
			\item \textbf{Compatible closure atlas.}
			A causal closure atlas is an indexed family
			$(\gamma_I)_{I\in\mathcal H}$ of such closures.
			It is \emph{compatible} if, for every $I,J\in\mathcal H$
			and $A\subseteq W_{I\cap J,\infty}$,
			\begin{equation}\label{wn10:eq:temporal-compat}
				\gamma_I(A)\cap W_{I\cap J,\infty}
				=\gamma_J(A)\cap W_{I\cap J,\infty}.
			\end{equation}
			
			\item \textbf{Exact global extension.}
			A closure operator $\Gamma$ on $W_{V,\infty}$ is an
			\emph{exact global extension} of the atlas if
			\[
			\Gamma(A)\cap W_{I,\infty}=\gamma_I(A)
			\]
			for every $I\in\mathcal H$ and
			$A\subseteq W_{I,\infty}$.
			It is \emph{dynamically admissible} if it additionally satisfies conditions (i)--(iv).
		\end{enumerate}
	\end{ddefinition}

	Condition (v) bounds the premise size of every semantic
	consequence, rather than merely that of the rules in a
	chosen Horn basis.
	Since each time prefix is finite, past causality also implies that the closure is finitary.

	\begin{dremark}[Positive-event semantics]
		A closed set $F\subseteq W_{I,\infty}$ is interpreted as a permissible
		history of true or activated events.  Its complement need not
		mean that those events are impossible; the closure expresses positive
		implications.  Any premise is consistent (the all-true history is a
		model), so there are no nontrivial Horn nogoods in this restricted
		language.  This distinguishes the present native core from the
		positive/negative state literals and realizable refutations of
		Section~\ref{sec:covers-assertions}.  The distinction is essential, not a syntactic convenience.
	\end{dremark}

	\begin{dlemma}[Native causal Horn basis]\label{wn10:lem:native-horn}
		For $\gamma_I$ satisfying Definition~\ref{wn10:def:causal-closure}, let
		$\mathcal R_I$ be the collection of all nontrivial semantic rules
		\begin{equation}\label{wn10:eq:native-rules}
			P\longrightarrow q,\quad
			1\le|P|\le k,\quad q\in\gamma_I(P)\setminus P,\quad
			\max\{t:(i,t)\in P\}\le t(q).
		\end{equation}
		Then $\mathcal R_I$ generates exactly $\gamma_I$.  Every rule of
		$\mathcal R_I$ is a genuine local semantic consequence, and
		$\mathcal R_I$ is invariant under nonnegative time shifts and under
		backward shifts when all its atoms remain in $\mathbb N_0$.
	\end{dlemma}

	\begin{proof}
		Each rule is valid, so the closure it generates is contained in
		$\gamma_I$.  Conversely, suppose $q=(i,t)\in\gamma_I(A)\setminus A$.
		Past causality gives $q\in\gamma_I(A\cap W_{I,t})$, and the uniform
		Carath\'eodory bound provides $P\subseteq A\cap W_{I,t}$ with
		$1\le|P|\le k$ and $q\in\gamma_I(P)$.  This is a rule in
		$\mathcal R_I$, deriving $q$ directly from $A$.  Stationarity makes
		such valid rules equivariant under time translation in both directions
		whenever all shifted time labels remain nonnegative.  \qedhere
	\end{proof}

	\subsubsection{Canonical coherent globalization over infinite time}

	A causal rule whose conclusion lies in a finite prefix has all its premises there. This makes it possible to compare the least closures computed on successive prefixes.

	\begin{dlemma}[Prefix coherence of the least generated extension]
		\label{wn10:lem:prefix-global}
		For a compatible atlas let $\mathcal R=\bigcup_{I\in\mathcal H}\mathcal R_I$.
		Write $\Gamma$ for the closure on $W_{V,\infty}$ generated by
		$\mathcal R$, and let $C_m$ be the closure on $W_{V,m}$ generated by
		the rules in $\mathcal R$ whose heads have time at most $m$.
		Then for every $A\subseteq W_{V,\infty}$,
		\begin{equation}\label{wn10:eq:prefix-coherence}
			\Gamma(A)\cap W_{V,m}=C_m(A\cap W_{V,m}),
			\qquad C_{m+1}|_{W_{V,m}}=C_m.
		\end{equation}
		In addition, $C_m$ is the \emph{least generated extension} of the
		local finite-horizon closure atlas
		$(\gamma_I|_{W_{I,m}})_{I\in\mathcal H}$.
	\end{dlemma}

	\begin{proof}
		All rules in $\mathcal R$ are causal: no premise has time greater than its head.
		During forward chaining, a conclusion at time at most $m$
		therefore depends only on facts at times at most $m$; induction on the finite length of a Horn proof gives the first identity.
		The second follows by restricting the first identity to $W_{V,m}$.
		
		By Lemma~\ref{wn10:lem:native-horn}, the time-$m$ rules generate each $\gamma_I|_{W_{I,m}}$.
		Thus their union generates exactly the least closure containing all local finite-horizon implications, as in the closure-restriction and least-extension lemma of
		\cite{MaPackaging2026}. \qedhere
	\end{proof}

	\begin{dtheorem}[Sharp all-horizon causal Horn globalization]
		\label{thm:all-horizon-causal-Horn-globalization}
		Fix $k\geq1$ and a finite cover $\mathcal H$.
		The following conditions are equivalent:
		\begin{enumerate}[label=(\roman*),leftmargin=*,nosep]
			\item Every compatible atlas of stationary persistent
			causal closures of uniform intrinsic determination
			order at most $k$, as in
			Definition~\ref{wn10:def:causal-closure}, admits an
			exact dynamically admissible global extension.
			
			\item The primal graph of $\mathcal H$ is chordal,
			and every nonempty clique of size at most $2k+1$
			is guarded by an original chart.
			
			\item $\boxed{2k+1<\nu(\mathcal H)}$.
		\end{enumerate}
		When these conditions hold, the closure $\Gamma$ generated by the native local causal rules is the least exact global extension.
		It is dynamically admissible and finitary, with
		finite-horizon restrictions equal to the least exact
		extensions $C_m$.
		
		When the conditions fail, there exists a compatible atlas
		with no exact global extension, whose local closures are
		generated by finitely many stationary same-time Horn
		templates and one-step persistence templates.
	\end{dtheorem}

	\begin{proof}
		\leavevmode
		\begin{enumerate}[label={},leftmargin=0pt,nosep]
			\item \textnormal{(ii)$\Leftrightarrow$(iii).}\enspace
			Since $2k+1\geq3$, the equivalence follows directly
			from Definition~\ref{def:nu}.
			
			\item \textnormal{(iii)$\Rightarrow$(i).}\enspace
			For every finite horizon $m$,
			Theorem~\ref{thm:obstruction-rank-invariance-expansion} gives
			\[
			\nu(\mathcal H^{[0,m]})=\nu(\mathcal H)>2k+1.
			\]
			The restricted local closures are compatible,
			empty-preserving, and have Carath\'eodory number
			at most $k$.
			Hence Theorem~\ref{thm:source} ensures that their least generated extension $C_m$ is exact.
			
			By Lemma~\ref{wn10:lem:prefix-global}, these extensions are coherent and satisfy
			\[
			\Gamma(A)\cap W_{V,m}
			=C_m(A\cap W_{V,m}),
			\]
			where $\Gamma$ is generated by the union of the native local causal rules.
			Exactness of each $C_m$, together with local causality, gives
			\[
			\Gamma(A)\cap W_{I,m}
			=\gamma_I(A)\cap W_{I,m}
			\]
			for every $I$, $m$, and $A\subseteq W_{I,\infty}$.
			Taking the union over $m$ proves global exactness.
			
			Nonempty causal premises ensure empty preservation
			and causality.
			The native persistence rules ensure	persistence.
			By Lemma~\ref{wn10:lem:native-horn}, every derivation can be shifted forward in time.
			Conversely, a derivation from $\sigma_h A$ contains
			no events before time $h$, so it can be shifted
			backward by $h$.
			Thus $\Gamma(\sigma_h A)=\sigma_h\Gamma(A)$.
			Finitariness follows from the finite premises of
			the generating rules.
			Every exact global extension satisfies all native
			local rules and therefore contains $\Gamma$.
			
			\item \textnormal{(i)$\Rightarrow$(iii).}\enspace
			Suppose $2k+1\geq\nu(\mathcal H)$.
			By Theorem~\ref{thm:source}, there exists a compatible atlas $(c_I)_{I\in\mathcal H}$ of empty-preserving static closures of Carath\'eodory number at most $k$ with no exact global extension.
			
			For $A\subseteq W_{I,\infty}$, put
			$A_s:=\{i:(i,s)\in A\}$ and define
			\begin{equation}\label{wn10:eq:static-temporal-lift}
				\gamma_I(A):=
				\bigcup_{t\geq0}
				c_I\!\left(\bigcup_{s\leq t}A_s\right)
				\times\{t\}.
			\end{equation}
			Extensivity and monotonicity follow from those of $c_I$.
			The closed cumulative prefixes are nested, so
			idempotence follows from that of $c_I$.
			The formula also ensures empty preservation,
			causality, stationarity, and persistence.
			
			If $(i,t)\in\gamma_I(A)$, the static Carath\'eodory
			bound supplies at most $k$ labels in
			$\bigcup_{s\leq t}A_s$ determining $i$.
			Choosing one occurrence of each label gives at most
			$k$ events determining $(i,t)$.
			Thus the uniform intrinsic determination bound is
			preserved.
			Static compatibility implies temporal
			compatibility at every time.
			
			Restricting any exact global temporal extension
			to time zero would give an exact global extension
			of $(c_I)_{I\in\mathcal H}$, a contradiction.
			
			Finally, each $c_I$ admits a finite implicational
			basis with nonempty premises of size at most $k$
			by Lemma~4.13 of \cite{MaPackaging2026}.
			Repeating these rules at every time and adding
			the persistence templates $i_t\to i_{t+1}$
			generates exactly
			\eqref{wn10:eq:static-temporal-lift}:
			persistence transports all previously occurring
			labels to time $t$, and the same-time rules
			generate their static closure.
			This proves sharpness within the stated subclass. \qedhere
		\end{enumerate}
	\end{proof}

	\begin{dcorollary}[Exact finite and infinite propagation-history realization]
		\label{wn10:cor:history-exact}
		Under the hypotheses and inequality of
		Theorem~\ref{thm:all-horizon-causal-Horn-globalization}, let
		$\mathsf{Hist}_I(m)=\operatorname{Fix}(\gamma_I|_{W_{I,m}})$,
		and let $\mathsf{Hist}_V(m)=\operatorname{Fix}(C_m)$.
		\begin{enumerate}[label=(\roman*),leftmargin=*,nosep]
			\item \emph{Finite histories.} Then for every $m$ and every original chart $I$,
			\begin{equation}\label{wn10:eq:history-projection}
				\pi_{W_{I,m}}\mathsf{Hist}_V(m)=\mathsf{Hist}_I(m),
				\qquad
				\mathsf{Hist}_V(m)=\Join_{I\in\mathcal H}\mathsf{Hist}_I(m).
			\end{equation}
			
			\item \emph{Infinite extension.} Every closed finite-prefix history on an original chart extends to an \emph{infinite} global closed-event history with that prescribed prefix; every infinite locally closed-event history on an original chart likewise has an infinite global extension.
			
			\item \emph{Finite-state presentation.} If the local closures admit stationary generating templates of uniformly bounded temporal span $r$, the global history language has a finite-window presentation and can be encoded by a finite-state higher-block transition system.
		\end{enumerate}
	\end{dcorollary}

	\begin{proof}
		\leavevmode
		\begin{enumerate}[label={(iii)},leftmargin=*,nosep]
			\item[(i)] The fixed-set characterization of the least generated closure gives
			the second equality in~\eqref{wn10:eq:history-projection}.  Exact
			restriction $C_m|_{W_{I,m}}=\gamma_I|_{W_{I,m}}$ and the standard
			fixed-set restriction lemma give the first.  In fact, if $F_I$ is
			locally closed, then $C_m(F_I)$ extends it without adding facts on $I$.
			\item[(ii)] For an infinite local history $F_I$, the same argument with $\Gamma(F_I)$
			constructs a global closed extension; for a finite prefix use
			$\Gamma(F_I)$ with no initially asserted later events, together with
			prefix coherence.  Conversely a globally closed history restricts to
			locally closed histories.
			\item[(iii)] For stationary generating templates of span at most $r$, translate each earliest atom to time zero. Finiteness of $V\times[0,r]$ gives finitely many such normalized rules. Their shifts test windows of length at most $r+1$. Storing the last at most $r$ truth sets $F_t\subseteq V$ and checking newly completed rule windows gives a finite-state history presentation with the specified empty initial memory.  This last assertion concerns
			positive Horn history languages, not arbitrary local relational
			transition atlases. \qedhere
		\end{enumerate}
	\end{proof}

	\subsubsection{A native original-chart recovery rank}

	A rooted causal temporal core of $\gamma_I$ is a nontrivial implication $P\to x$ valid in $\gamma_I$ with $P$ inclusion-minimal and contained in times at most $t(x)$.  Empty preservation excludes empty-premise cores.
	Intrinsic determination implies $|P|\le k$.
	The spatial support of such a core is
	$\operatorname{sp}(P\cup\{x\})=\{i:(i,t)\in P\cup\{x\}\}$.
	A supported one-pivot composition consists of cores
	\[
	A\longrightarrow x\quad\text{in chart }I,
	\qquad B\cup\{x\}\longrightarrow y\quad\text{in chart }J,
	\quad y\notin A\cup B,
	\]
	where $x\notin B$, and the roots/premises are correctly time typed.
	The resulting semantic assertion is $A\cup B\to y$ in the generated global closure.
	Its \emph{composed spatial support} is
	$Q=\operatorname{sp}(A\cup B\cup\{x,y\})$.

	\begin{ddefinition}[Classwise two-cut temporal guarding rank]
		\label{wn10:def:two-cut-rank}
		For a class $\mathscr D$ of compatible causal Horn atlases with a declared uniform finite semantic determination bound, the \textbf{two-cut temporal guarding rank} $\operatorname{tprk}_2(\mathscr D)$ is the least integer $p\ge2$ with the following property: for every finite cover whose every primal clique of size at most $p$ is contained in an original chart, and
		for every atlas in $\mathscr D$ on that cover, every supported one-pivot composition with $Q$ a primal clique admits an original-chart rooted recovery $D\to y$ with
		$D\subseteq A\cup B$ and $|D|$ no larger than the class's
		semantic determination bound.
	\end{ddefinition}

	\begin{dtheorem}[Sharp original-chart two-cut recovery]
		\label{wn10:thm:sharp-temporal-rank}
		For the class $\mathscr T_k$ of all compatible temporal Horn closure atlases of Definition~\ref{wn10:def:causal-closure},
		\begin{equation}\label{wn10:eq:sharp-tprk}
			\operatorname{tprk}_2(\mathscr T_k)=2k+1.
		\end{equation}
		The upper bound follows by counting the composed support.
		The matching lower bound uses semantic rooted cores on the original charts and excludes every smaller guard threshold.
		Consequently,
		\begin{equation}\label{wn10:eq:tprk-vs-nu}
			\operatorname{UG}^{\mathrm{causal}}_{\mathscr T_k}(\mathcal H)
			\quad\Longleftrightarrow\quad
			\operatorname{tprk}_2(\mathscr T_k)<\nu(\mathcal H).
		\end{equation}
		Here $\operatorname{UG}^{\mathrm{causal}}$ denotes the classwise exact-globalization property of Theorem~\ref{thm:all-horizon-causal-Horn-globalization}.
		The equality concerns the class $\mathscr T_k$.
	\end{dtheorem}

	\begin{proof}
		Let $A\to x$ and $B\cup\{x\}\to y$ be as above.
		By intrinsic determination, $|A|\le k$ and $|B|\le k-1$, so
		\[
		|Q|\le |A|+|B|+2\le 2k+1.
		\]
		If all cliques through size $2k+1$ are guarded and $Q$ is a primal
		clique, an original chart $K$ contains $Q$.  Each parent assertion,
		including its time labels, is supported on the intersection of $K$
		with the chart in which it originated.  Exact local compatibility
		transfers both assertions to $\gamma_K$.  Its closure axioms then
		entail $A\cup B\to y$.  The local semantic Carath\'eodory bound
		supplies $D\subseteq A\cup B$ with $|D|\le k$ and
		$y\in\gamma_K(D)$; minimizing $D$ produces a rooted causal core.
		This proves $\operatorname{tprk}_2\le2k+1$.
		
		For the lower bound, let $Q=Z\sqcup W\sqcup\{x,y\}$ with
		$|Z|=k$, $|W|=k-1$, and take the simplex-boundary cover
		$\{Q\setminus\{z\}:z\in Q\}$.  Every clique of size at most $2k$
		is guarded, whereas the full primal clique $Q$ is unguarded.
		Assign to each facet exactly the supported members of the two
		implications $Z\to x$ and $W\cup\{x\}\to y$; since their supports
		jointly occupy $Q$, no facet receives both.  As in the proof of Theorem~4.22 of Paper~I~\cite{MaPackaging2026}, the local closures are compatible and
		have Carath\'eodory number at most $k$.
		Lift them by~\eqref{wn10:eq:static-temporal-lift}.  At time zero the two
		native roots compose to $Z\cup W\to y$, but the facet $Q\setminus\{x\}$
		has identity closure at time zero.  A recovery $D\to y$ with
		$D\subseteq Z\cup W$ on \emph{any} original facet would transfer,
		by compatibility, to that identity facet, a contradiction.
		Thus no $p\le2k$ suffices.  Together with the upper bound, this
		proves~\eqref{wn10:eq:sharp-tprk}.  Equation~\eqref{wn10:eq:tprk-vs-nu} is
		Theorem~\ref{thm:all-horizon-causal-Horn-globalization} with the just-computed
		rank substituted. \qedhere
	\end{proof}

	\begin{dremark}[Relation to the original intrinsic recovery rank]
		Paper~I establishes $\operatorname{ipr}(\mathscr C_k)=2k+1$
		for bounded-Carath\'eodory closures \cite{MaPackaging2026}.
		The time-zero embedding gives the same value in~\eqref{wn10:eq:sharp-tprk}. The temporal rank counts guarded compositions of two native causal implications through a common pivot.
	\end{dremark}

	\paragraph{Examples of the sharp causal cutoff.}

	Positive propagation can be exact on a cover with an unguarded clique. The size of that clique determines which locally bounded implications still have a common extension.
	
	\begin{dexample}[A positive non-join-tree example]
		\label{wn10:ex:tetra-positive}
		Let $V=\{1,2,3,4\}$ and let $\mathcal H$ consist of the four
		three-element facets $V\setminus\{i\}$.  Its primal graph is
		$K_4$, the four-clique is unguarded, and $\nu(\mathcal H)=4$;
		the cover has no join tree.  With $k=1$ the threshold $3<4$ holds.
		Take stationary templates $i_t\to4_{t+1}$ for $i=1,2,3$,
		plus all persistence rules.  Assign a rule to every original chart
		containing its full spatial support.  Each local semantic closure
		has Carath\'eodory number one, their overlap restrictions agree,
		and Theorem~\ref{thm:all-horizon-causal-Horn-globalization} gives an exact coherent
		global propagation-history system.  This example goes beyond the
		join-tree hypotheses of Sections~\ref{sec:local-dynamical-atlases}--\ref{sec:width}.
	\end{dexample}

	\begin{dexample}[Sharp failure at the boundary]
		For $k=1$ take $\mathcal H=\{\{a,b\},\{b,c\},\{a,c\}\}$.
		Assign the local rule $a_t\to b_t$ on $\{a,b\}$,
		$b_t\to c_t$ on $\{b,c\}$ and identity closure at each time on
		$\{a,c\}$; include persistence everywhere.  The local temporal
		closures agree on singleton overlaps.  Nevertheless global closure
		forces $a_0\to c_0$, contradicting the identity restriction on
		$\{a,c\}$.  Here $\nu(\mathcal H)=3=2k+1$, showing strictness.
		For $k=2$, the five-vertex simplex boundary with two disjoint
		binary-premise rules, as in the proof of
		Theorem~\ref{wn10:thm:sharp-temporal-rank}, gives the analogous failure
		at $\nu=5=2k+1$.
	\end{dexample}

	\begin{dremark}[Primitive arity versus semantic determination]
		With rules $\{a,b\}\to x$ and $\{x,c\}\to y$, the consequence
		$\{a,b,c\}\to y$ may require all three original premises.
		Thus a basis containing only two-premise rules need not have
		Carath\'eodory number two.  The \emph{local} uniform semantic
		Carath\'eodory bound, not the primitive rule arity, is indispensable.
		Nor does the theorem assert a uniform Carath\'eodory bound for the
		global extension: a global basis of small-premise rules can have
		composed consequences with larger minimal determining sets.
	\end{dremark}

	\begin{dremark}[Temporal extension of a static closure theorem]
		The geometric cutoff follows from the static closure theorem \cite{MaPackaging2026}. Horn generation and least-fixed-point semantics are classical \cite{CaspardMonjardet2003,wn10:DatalogSurvey}. The temporal argument gives exact restriction at every prefix and a stationary closure on the whole time axis.
	\end{dremark}

	\subsubsection{An independent intrinsic temporal-memory coordinate}

	A spatial recovery scale alone does not measure how far a consequence	propagates in time.
	At the same time, the maximal time span of a semantic rooted circuit is an unsuitable invariant:
	persistence implies $i_0\to i_t$ at all horizons $t$.

	\begin{ddefinition}[Intrinsic causal generation width and recovery profile]
		\label{wn10:def:temporal-width}
		Let $\gamma$ be a temporal closure on
		$W=V\times\mathbb N_0$, where $V$ is finite,
		and let $t(p)$ denote the time coordinate of an event $p$.
		\begin{enumerate}[label=(\roman*),leftmargin=*,nosep]
			\item \textbf{Causal rule and temporal span.}
			A nontrivial valid causal rule of $\gamma$ is an
			implication $P\to q$, where $P\subseteq W$ is finite
			and nonempty, $q\in\gamma(P)\setminus P$, and
			$t(p)\leq t(q)$ for every $p\in P$.
			Its temporal span is
			\[
			\operatorname{span}(P\to q):=
			t(q)-\min_{p\in P}t(p).
			\]
			For $r\in\mathbb N_0$, let
			$\mathcal R_{\leq r}(\gamma)$ be the family of all
			such rules with temporal span at most $r$.
			
			\item \textbf{Intrinsic causal generation width.}
			The intrinsic causal generation width of $\gamma$ is
			\begin{equation}\label{wn10:eq:intrinsic-width}
				w(\gamma)=\genlen(\gamma):=
				\min\bigl\{r\in\mathbb N_0:
				\operatorname{cl}\langle
				\mathcal R_{\leq r}(\gamma)\rangle=\gamma\bigr\},
			\end{equation}
			where $\operatorname{cl}\langle\mathcal R\rangle$
			denotes the least closure operator validating every
			rule in $\mathcal R$, and $\min\varnothing=\infty$.
			
			\item \textbf{Classwise generation width and recovery profile.}
			For a class $\mathscr D$ of temporal closures, define
			\[
			w(\mathscr D):=\sup_{\gamma\in\mathscr D}w(\gamma).
			\]
			For a class of atlases, the supremum is taken over
			all local chart closures of its members.
			The two-parameter recovery profile of $\mathscr D$ is
			\[
			\bigl(\operatorname{tprk}_2(\mathscr D),
			w(\mathscr D)\bigr),
			\]
			combining its two-cut temporal guarding rank and
			intrinsic causal generation width.
		\end{enumerate}
	\end{ddefinition}

	The recovery profile measures spatial guarding complexity and temporal generation width.

	\begin{dproposition}[Semantic invariance and hidden-state width growth]
		\label{wn10:prop:width-growth}
		The width $w(\gamma)$ is independent of the chosen Horn-rule
		presentation and is preserved under time-preserving, typed semantic
		isomorphisms.  For fixed $k=1$, however, $w$ is \emph{not monotone
			under spatial restriction}.  Indeed, on $V=\{a,b,c\}$ take the
		stationary persistent causal closure generated by
		\[
		a_t\longrightarrow b_{t+1},\qquad
		b_t\longrightarrow c_{t+1},\qquad
		i_t\longrightarrow i_{t+1}\quad(i=a,b,c).
		\]
		The full closure has intrinsic width $1$, while its restriction to
		$\{a,c\}\times\mathbb N_0$ has intrinsic width $2$.
	\end{dproposition}

	\begin{proof}
		The definition uses the \emph{complete semantic set} of valid
		bounded-span rules, so it is invariant under any different rule basis
		with the same closure.  The stated full basis has width one, and
		persistence precludes width zero.
		
		Under restriction to $\{a,c\}$, the closure is generated by
		$a_t\to c_{t+2}$ together with persistence for $a$ and $c$;
		this proves width at most two.  To see that width one fails, start
		from the sole premise $a_0$.  In the restricted closure, $c_2$ is
		entailed but $c_0,c_1$ are not.  Any derivation using only semantically
		valid rules of span at most one, starting from $a_0$, could derive
		$a_1,a_2,\ldots$, but the first possible consequence of $a_s$ on $c$ is $c_{s+2}$. Because the restriction has Carath\'eodory number one, any nonempty valid premise of only $a$-events forcing $c_t$ contains a singleton that already forces $c_t$, necessarily at time at most $t-2$. Thus no valid rule of span at most one, even with several premises, can introduce the first $c$-event.  Thus $w=2$. \qedhere
	\end{proof}

	\begin{dremark}[Classwise temporal uniformity]
		The broad class $\mathscr T_k$ has no finite uniform width bound:
		stationary unary implications $a_t\to b_{t+N}$ with persistence
		have $w=N$ for arbitrary $N$.  A subclass defined by the
		\emph{semantic} bound $w(\gamma)\le r$ has a finite temporal
		coordinate, but is generally \emph{not stable under arbitrary spatial
			restriction}; Proposition~\ref{wn10:prop:width-growth} explains why.
		The sharp globalization theorem does not require a finite width
		bound, since it controls temporal consequences through the stronger
		uniform semantic Carath\'eodory hypothesis at every horizon.
	\end{dremark}

	A positive closure always admits the all-true history. A transition relation can rule out a joint pattern or fail to compose two steps that use different hidden states. The relational calculus below keeps these witnesses during elimination. If an additional global admission relation is imposed, it enters the canonical test as part of the state space.

	\subsection{Witness-preserving relational globalization beyond Horn histories}
	\label{wn11:sec:nonhorn-relational}

	A relational path can constrain several values at once. We keep the witnesses used at intermediate times and compare them when two local path conditions are joined.

	\subsubsection{Relational constraints and time projection}

	Positive Horn propagation, studied in Section~\ref{sec:global-theorems}, has a least-closure
	semantics.  General finite transition relations instead admit forbidden
	configurations, mutually exclusive values, reversible dynamics, and
	existentially hidden intermediate states.  Consequently, composing
	implications in a closure algebra is insufficient.  This section develops
	an exact finite-state relational semantics and a structural non-Horn class
	in which original-chart recovery is controlled by cycle transports.
	The original packaged structure and its finite cover are retained;
	the local finite-state relations are additional specialization data.
	
	Let $V$ be a finite descriptor set, $\mathcal H$ a cover of $V$, and
	$A_v$ a finite nonempty alphabet for each $v$.  Write
	$A^I=\prod_{v\in I}A_v$.  For each $I\in\mathcal H$ choose nonempty
	$X_I\subseteq A^I$ and $R_I\subseteq X_I\times X_I$.  Independently,
	let $X\subseteq\{x\in A^V:x|_I\in X_I\text{ for all }I\}$ be a set of
	admissible global states.  The extra constraint $X$ is part of the
	problem data and need not follow from local-chart relations.
	Take the greatest respecting relation $R_*$ of~\eqref{eq:canonical}.
	Include the empty word in every finite-path language by convention.
	For a finite horizon $m\ge0$, denote the local path relation by
	$\mathsf T_I(m)$ and the canonical global path relation by
	$\mathsf T_*(m)$.  The latter consists of complete, jointly consistent
	sequences $(x_0,\dots,x_m)$, not independently projected one-step
	transitions.

	\begin{ddefinition}[Witness semantics of temporal assertions]
		\label{wn11:def:witness-semantics}
		Fix a path relation $\mathsf T_*(m)$ at horizon $m$.
		Let $p$ be a partial assignment on
		$D\subseteq V\times\{0,\ldots,m\}$.
		\begin{enumerate}[label=(\roman*),leftmargin=*,nosep]
			\item \textbf{Witness set.}
			The witness set of $p$ is
			\[
			\mathsf{Wit}_m(p):=
			\{w\in\mathsf T_*(m):w|_D=p\}.
			\]
			
			\item \textbf{Nogood.}
			A partial assignment $p$ is a nogood if
			$\mathsf{Wit}_m(p)=\varnothing$.
			
			\item \textbf{Realizable semantic consequence.}
			Let $q$ be a literal specifying a value at a single
			space--time cell $c$.
			A rooted assertion $p\Rightarrow q$ is a realizable
			semantic consequence if
			\[
			\mathsf{Wit}_m(p)\ne\varnothing,
			\qquad
			\mathsf{Wit}_m(p)\subseteq\mathsf{Wit}_m(q).
			\]
			Here $q$ is regarded as a singleton partial assignment.
			
			\item \textbf{Original-chart guarding.}
			A rooted assertion $p\Rightarrow q$ is original-chart
			guarded if there exists $I\in\mathcal H$ such that
			\[
			D\cup\{c\}\subseteq I\times\{0,\ldots,m\}.
			\]
		\end{enumerate}
		These semantic notions depend only on the time-indexed
		path relations and the declared original charts.
	\end{ddefinition}

	The witness-preserving elimination rule is the two-factor specialization of Theorem~\ref{thm:exact-finite-horizon-propagation}, equation~\eqref{wn11:eq:join-project}.
	The following example exhibits the false witnesses created by premature projection.

	\begin{dexample}[A missing intermediate witness]
		Let $U=\{a\}$, $S=\{b\}$, $W=\{c\}$ and all alphabets be binary.
		Take $R=\{(0,0),(1,1)\}$ on $(a,b)$ and
		$Q=\{(0,1),(1,0)\}$ on $(b,c)$.  Their natural join enforces
		$c=1-a$, giving the two endpoint pairs $(0,1),(1,0)$.
		But projecting $R$ to $a$ and $Q$ to $c$ separately permits all four
		endpoint pairs.  The additional two have incompatible hidden
		$b$-witnesses.  Exactly the same issue occurs when $b$ is a hidden
		intermediate state at successive time steps.
	\end{dexample}

	\subsubsection{A general finite-state exactness and decidability theorem}

	The canonical relation gives a decision test for one atlas. Minimizing the two path automata bounds the length of a missing trajectory, while the universal statement characterizes the covers on which no compatible instance can fail.

	\begin{dcorollary}[Length of a missing chart trajectory]\label{wn11:thm:greatest-relation-exact}
		On a finite atlas and its prescribed $X$, let $u$ and $v$ be the sizes of the minimal complete deterministic automata for a given chart's local path language and its path language projected from $R_*$. Include the empty word in both. If exactness fails on the chart, a missing path has horizon at most $D(u,v)-3$, where $D$ is given by~\eqref{rev:eq:residual-union-size}. In particular this is no larger than $\min\{uv-2,u+v-3\}$.
	\end{dcorollary}

	\begin{proof}
		Theorem~\ref{thm:greatest-relation} identifies failure with a local word missing from the projected language; the reverse difference is empty. Both languages are prefix closed and contain the empty word by convention. Apply the common-sink refinement estimate of Theorem~\ref{wn05-06:thm:finite-trace-test}, or its residual-index form in Theorem~\ref{wn09:thm:finite-witness}. The distinguishing word is therefore missing in the required direction. Since $u+v-3\le uv-2$, the claimed comparison follows.
	\end{proof}

	\begin{dtheorem}[Sharp universal criterion for unrestricted relational globalization]
		\label{wn11:thm:all-relational-classwise}
		Fix a finite indexed cover $\mathcal H$ of a finite set $V$.
		Allow arbitrary finite nonempty coordinate alphabets,
		nonempty local state sets $X_I$, and transition relations
		$R_I\subseteq X_I^2$.
		An atlas is \emph{all-horizon compatible} if
		\[
		\mathsf L_I^{I\cap J}(m)=\mathsf L_J^{I\cap J}(m)
		\]
		for every $I,J\in\mathcal H$ and $m\in\mathbb N_0$,
		including empty intersections.		
		In the absence of additional global-state constraints,
		\begin{equation}\label{wn11:eq:universal-nu-inf}
			\begin{gathered}
					\text{Every all-horizon compatible atlas admits}\\
					\text{an all-horizon exact globalization}
					\quad\Longleftrightarrow\quad
					\nu(\mathcal H)=\infty.
			\end{gathered}
		\end{equation}
	\end{dtheorem}

	\begin{proof}
		\leavevmode
		\begin{enumerate}[label={},leftmargin=0pt,nosep]
			\item \textnormal{Sufficiency.}\enspace
			Suppose $\nu(\mathcal H)=\infty$.
			By Definition~\ref{def:nu} and the classical acyclicity
			characterization~\cite{wn09:BeeriFaginMaierYannakakis1983}, the cover $\mathcal H$ admits a join tree.
			The same tree is a join tree for every time-expanded
			cover $\mathcal H^{[0,m]}$.
			
			All-horizon compatibility makes the local path
			relations $\mathsf T_I(m)$ compatible on their
			intersections.
			Acyclic relational gluing therefore gives
			\[
			\pi_{I\times[0,m]}
			\left(\Join_{J\in\mathcal H}\mathsf T_J(m)\right)
			=\mathsf T_I(m)
			\]
			for every $I$ and $m$.
			By Definition~\ref{wn08:def:temporal-model-relation},
			\[
			\Join_{J\in\mathcal H}\mathsf T_J(m)
			=\mathsf T_{R_{\mathcal A}}(m).
			\]
			Thus the single canonical global relation
			$R_{\mathcal A}$ realizes every local path relation
			at every finite horizon.
			
			\item \textnormal{Necessity.}\enspace
			Suppose $\nu(\mathcal H)<\infty$.
			The classical acyclicity characterization
			\cite{wn09:BeeriFaginMaierYannakakis1983}
			provides pairwise-compatible nonempty finite relations
			$S_I\subseteq A^I$ having no exact global relation.
			Set
			\[
			X_I:=S_I,\qquad
			R_I:=\{(u,u):u\in S_I\}.
			\]
			Every local trajectory is constant. Hence state
			compatibility implies trace compatibility at every
			finite horizon, including on empty intersections.
			An all-horizon exact globalization would, at horizon
			zero, give a global state set with projections $S_I$,
			contradicting the choice of these relations.
			\qedhere
		\end{enumerate}
	\end{proof}

	\begin{dremark}[Path lengths on a disjoint separator]
		On an empty separator there is one word of each admissible length. A chart with a self-loop has paths of all lengths, whereas a chart with a single directed edge between distinct states has none of horizon two. Even disjoint charts must agree on admissible lengths to have an exact joint path law. On a prescribed smaller global state set one must also check Theorem~\ref{thm:greatest-relation}.
	\end{dremark}

	\subsubsection{A non-Horn class: permutation-transport dynamical packages}

	A useful non-Horn specialization has native, perfectly functional
	relationships but global obstructions.  Let $G=(V,E)$ be a finite,
	connected simple graph with $|V|\ge2$.  All vertices have finite
	alphabets $A_v$ of the same size $q\ge2$.
	For every oriented edge $u\to v$ choose a bijection
	$g_{uv}:A_u\to A_v$, with $g_{vu}=g_{uv}^{-1}$.  At each vertex
	choose an arbitrary map $\tau_v:A_v\to A_v$ such that
	\begin{equation}\label{wn11:eq:time-square}
		g_{uv}\circ\tau_u=\tau_v\circ g_{uv}
		\qquad(uv\in E).
	\end{equation}
	For each edge $uv$, the original chart has admissible states
	$X_{uv}=\{(a,g_{uv}a):a\in A_u\}$ and deterministic update
	$(a,g_{uv}a)\mapsto(\tau_u a,\tau_v g_{uv}a)$.  This is an
	\emph{edge-permutation dynamical atlas}.  Its state predicates enforce
	negative exclusions (a vertex cannot have two different labels), so
	it is not interpreted as an empty-preserving positive Horn closure
	on activated events.

	\begin{dlemma}[Complete separator-trace compatibility]
		\label{wn11:lem:perm-separator}
		Every pair of edge charts sharing a vertex $v$ has identical
		singleton-$v$ observation languages at every finite horizon.  These
		languages consist of the deterministic $\tau_v$-trajectories from
		arbitrary initial values in $A_v$.
	\end{dlemma}

	\begin{proof}
		An arbitrary $v$-label extends uniquely across each incident edge
		using $g_{uv}$ or its inverse.  Commutation~\eqref{wn11:eq:time-square}
		ensures that the extension evolves as prescribed and remains on the
		edge relation.  Hence both incident edges realize the same complete
		$\tau_v$-trajectory language.
	\end{proof}

	Fix a spanning tree $T$ of $G$ and a root $r$.  For $v\in V$, let
	$p_v:A_r\to A_v$ be the bijection obtained by composing edge
	transports along the unique $T$-path from $r$ to $v$; in particular,
	$p_r=\mathrm{id}$.  For a non-tree edge $e=u\to v$ define its
	\emph{fundamental holonomy}
	\begin{equation}\label{wn11:eq:fundamental-hol}
		h_e=p_v^{-1}\circ g_{uv}\circ p_u:A_r\to A_r.
	\end{equation}
	Let $H_r\le\operatorname{Sym}(A_r)$ be the group generated by all
	$h_e$ and define its common fixed-point set
	\begin{equation}\label{wn11:eq:common-fix}
		F_r=\operatorname{Fix}(H_r)=\{a\in A_r:h(a)=a\ \forall h\in H_r\}.
	\end{equation}
	The generators give exactly the subgroup of all closed-walk transport permutations based at $r$: inserting the tree-path transports and their inverses at intermediate vertices expresses any closed walk as a product of the fundamental holonomies and their inverses. Hence this subgroup is independent of the chosen spanning tree.  Under a relabeling of
	vertex alphabets it changes by conjugation at the root; its triviality
	and the cardinality $|F_r|$ are therefore intrinsic.

	\begin{dtheorem}[Native guarded path recovery and exact holonomy globalization]
		\label{wn11:thm:holonomy}
		For the edge-permutation dynamical atlas, the following hold.
		\begin{enumerate}[label=(\roman*),leftmargin=2.35em]
			\item Every globally admissible state is uniquely of the form
			\begin{equation}\label{wn11:eq:global-parameter}
				x(a)=(p_v(a))_{v\in V},\qquad a\in F_r.
			\end{equation}
			The global transition sends $x(a)$ to $x(\tau_r a)$, and
			$\tau_r(F_r)\subseteq F_r$.
			\item The following conditions are equivalent:
			\begin{enumerate}[label=(\alph*),leftmargin=2.1em]
				\item exact state globalization;
				\item exact globalization of all finite trajectories;
				\item $F_r=A_r$;
				\item $h_e=\mathrm{id}$ for every non-tree edge;
				\item every tree-path assertion $x_v=p_vp_u^{-1}(x_u)$ whose
				endpoints $u,v$ belong to an original chord chart $uv$ is recovered
				\emph{exactly} by that chart's native relation.
			\end{enumerate}
			\item If $F_r\ne A_r$, the atlas is compatible on \emph{all}
			original overlap trajectory languages but is not exactly
			globalizable, even at horizon zero.  If $F_r=\varnothing$ there
			is no globally admissible state; if $\varnothing\ne F_r\subsetneq A_r$
			there are global states but some native local states have no global
			witness.
		\end{enumerate}
	\end{dtheorem}

	\begin{proof}
		\leavevmode
		\begin{enumerate}[label={(iii)},leftmargin=*,nosep]
			\item[(i)] A global state satisfying all tree-edge constraints is uniquely
			parametrized by $a=x_r\in A_r$, with $x_v=p_v(a)$.  For a non-tree
			edge $e=u\to v$, its constraint is $p_v(a)=g_{uv}p_u(a)$,
			equivalent to $h_e(a)=a$.  Thus the globally admissible states
			are exactly~\eqref{wn11:eq:global-parameter} with $a\in F_r$.
			Tree-path transports intertwine $\tau_r$ and $\tau_v$ by repeated
			applications of~\eqref{wn11:eq:time-square}.  The same identity on a
			non-tree edge implies $h_e\tau_r=\tau_rh_e$.  Hence $\tau_r$
			preserves the common fixed-point set $F_r$, and the stated update
			is the canonical global evolution.
			\item[(ii)] (a)--(e). For any original edge, its local states are parametrized by every
			$a\in A_r$ through its endpoint transport, whereas its globally
			realizable projections are parametrized only by $a\in F_r$.
			For (a)$\Leftrightarrow$(c), because all transports are bijective, exact state projection on
			any edge is equivalent to $F_r=A_r$.  For (c)$\Rightarrow$(b), in that case every local trajectory, from every local initial state, lifts via~\eqref{wn11:eq:global-parameter}; hence all-horizon exactness holds.  For (b)$\Rightarrow$(a), all-horizon exactness entails exact state projection.  For (c)$\Leftrightarrow$(d), the equality $F_r=A_r$ holds precisely when every generator $h_e$ is the identity.
			For (d)$\Leftrightarrow$(e), a non-tree edge $uv$ has native relation is
			$x_v=g_{uv}(x_u)$, whereas the composed tree-path assertion is
			$x_v=p_vp_u^{-1}(x_u)$.  These are the same rooted relation for
			\emph{every} initial value if and only if $h_e=\mathrm{id}$.
			Tree edges satisfy this recovery automatically, proving the final
			equivalence.
			\item[(iii)] Part~(iii) follows from the same parameterization and Lemma~\ref{wn11:lem:perm-separator}; disjoint charts also agree on all empty-separator words since every local map is total. \qedhere
		\end{enumerate}
	\end{proof}

	\paragraph{Transport fixed points versus temporal equilibria.}
	The parameter set $F_r=\operatorname{Fix}(H_r)$ tests spatial transport
	compatibility. It need not consist of equilibria of the temporal law.
	The parameterization in Theorem~\ref{wn11:thm:holonomy} gives
	\[
	\operatorname{Fix}(F_{\mathrm{glob}})
	=\{x(a):a\in F_r,\ \tau_r(a)=a\},\qquad
	\{x:F_{\mathrm{glob}}^n(x)=x\}
	=\{x(a):a\in F_r,\ \tau_r^n(a)=a\}.
	\]
	Indeed $F_{\mathrm{glob}}^n(x(a))=x(\tau_r^n(a))$ and $a\mapsto x(a)$
	is injective. Thus trivial holonomy guarantees exact state gluing but
	does not force temporal fixed points: a nontrivial cyclic update on
	the root alphabet may have none. See Definition~\ref{dyn45:fp-stationary}
	and Section~\ref{dyn46:subsec-weighted-equilibria} for the distinction between the operators
	whose fixed points are being considered.

	\begin{dcorollary}[Universal exact globalization of edge-permutation atlases]
		\label{wn11:cor:classwise-forest}
		For fixed $q\ge2$, every edge-permutation dynamical atlas on the
		connected two-uniform cover $\mathcal H_G=\{\{u,v\}:uv\in E\}$
		admits exact all-horizon globalization if and only if
		\begin{equation}\label{wn11:eq:forest-nu}
			G\text{ is a forest}
			\quad\Longleftrightarrow\quad
			\nu(\mathcal H_G)=\infty.
		\end{equation}
		When the graph has a cycle, there are compatible atlases whose
		native local trajectories do not globalize.
	\end{dcorollary}

	\begin{proof}
		A forest has no non-tree edge, so all holonomies are trivial and
		Theorem~\ref{wn11:thm:holonomy} applies.  If $G$ has a cycle, choose a
		spanning tree, give its edges identity transports, and choose a
		non-tree edge to carry a permutation with no fixed points (for
		example a nontrivial cyclic shift of a $q$-element alphabet).
		Let every vertex update be the identity.  Then the common fixed set
		is empty, so the atlas is not globally realizable despite complete
		overlap-trajectory compatibility.
		For a two-uniform edge cover, every forest is chordal and conformal.
		If $G$ contains a triangle, that three-clique is unguarded by any
		edge chart; otherwise a shortest cycle is induced and has length at
		least four, so the primal graph is nonchordal.  Definition~\ref{def:nu} therefore gives $\nu=3$ on every graph with
		a cycle, and $\nu=\infty$ on forests.
	\end{proof}

	\begin{dremark}[Interpretation of the classwise statement]
		This is a sharp, \emph{non-Horn} classwise result with the same
		universal-quantifier structure as the static globalization theorem of Paper~I~\cite{MaPackaging2026}.
		Unlike the Horn result $2k+1<\nu$, its graph-cover criterion is
		$\nu=\infty$, not a finite recovery threshold.  Flatness gives exact
		globalization even on cyclic graphs with $\nu=3$: geometry is a
		universal-class obstruction, while the actual holonomy is an
		atlas-specific invariant. Trivial holonomy is the balance condition for the associated gain graph with permutation gains~\cite{Zaslavsky1989}; it is also related to permutation constraint satisfaction and flat discrete connections.
	\end{dremark}

	\paragraph{Time-dependent transport and the temporal square condition.}

	A change of coordinates between consecutive times must transport both the spatial labels and their admissible updates. Commuting squares express precisely this compatibility.

	\begin{dproposition}[Holonomy transport under reconfiguration]
		\label{wn11:prop:time-dependent-hol}
		Allow the edge bijections to depend on time, $g_{uv,t}$, and suppose
		$\tau_{v,t}:A_v\to A_v$ are bijections satisfying each temporal
		square identity
		\[
		g_{uv,t+1}\tau_{u,t}=\tau_{v,t}g_{uv,t}.
		\]
		Then every root-based closed-walk holonomy satisfies
		\begin{equation}\label{wn11:eq:hol-conjugacy}
			h_C(t+1)=\tau_{r,t}\,h_C(t)\,\tau_{r,t}^{-1}.
		\end{equation}
		In particular, the full holonomy group is conjugated through time,
		its common fixed set moves by $\tau_{r,t}$, and trivial holonomy at
		one time is equivalent to trivial holonomy at every time.
	\end{dproposition}

	\begin{proof}
		Compose the temporal square identities around the directed cycle;
		all intermediate $\tau_{v,t}$ cancel, leaving
		$h_C(t+1)\tau_{r,t}=\tau_{r,t}h_C(t)$.  Invertibility gives
		\eqref{wn11:eq:hol-conjugacy}.  Conjugation preserves triviality and
		transports common fixed points bijectively.
	\end{proof}

	\subsubsection{Spatial and temporal obstructions to unrestricted globalization}
	\label{dyn45:sec:examples}

	The following examples separate failure of simultaneous spatial transport from disagreement of separator histories. They also show that large geometric obstruction rank alone cannot control the latter.
	\begin{dexample}[Full overlap-trace compatibility with cyclic failure]
		On the triangle $G=K_3$, let each $A_v=\{0,1\}$, all vertex
		updates be the identity, and set two edge transports equal to the
		identity and the third to the bit flip.  Every singleton overlap
		trace language is compatible at every horizon, but the fundamental
		holonomy is the bit flip, with no fixed point.  No global state
		exists.  On a three-symbol alphabet, replacing the flip by a
		transposition yields one global state, but not exact state
		projection.  With transports consisting of consistent cyclic
		translations on $\mathbb Z/3\mathbb Z$ and nontrivial vertex updates
		$a\mapsto a+1$, the same triangle has flat holonomy and exact
		all-horizon dynamics.
	\end{dexample}

	\begin{dexample}[Delayed temporal failure on an acyclic cover]
		\label{wn11:ex:delayed-join-tree}
		Example~\ref{wn05-06:ex:one-not-two} gives two binary charts on a join tree whose one-step overlaps agree, but whose two-step separator words require simultaneously $b_2=b_1$ and $b_2=1-b_1$. Here $\nu=\infty$: spatial cycle holonomy is absent but temporal separator consistency fails.
	\end{dexample}

	The frustrated triangle above has a spatial witness defect: its singleton
	overlaps agree, but its global state set is empty, with $\nu(\mathcal H)=3$.
	In contrast, Example~\ref{wn11:ex:delayed-join-tree} has a temporal
	witness defect at horizon two despite $\nu(\mathcal H)=\infty$.
	Neither defect is classified solely by primitive rule arity.

	\begin{dproposition}[Pairwise compatibility with arbitrarily large
		geometric obstruction rank]\label{wn11:prop:large-nu-parity}
		For every $N\ge4$ there is a finite binary, identity-evolution
		relational atlas on the simplex-boundary cover
		\[
		V=\{1,\dots,N\},\qquad
		\mathcal H_N=\{V\setminus\{i\}:i\in V\}
		\]
		with $\nu(\mathcal H_N)=N$, whose original-chart trajectory
		languages agree on every pairwise overlap at every finite horizon,
		but which has no exact global realization even at horizon zero.
	\end{dproposition}

	\begin{proof}
		On every $(N-1)$-variable facet take the binary tuples with
		\emph{even parity}, and use identity evolution on those tuples.
		A proper intersection of two distinct facets has $N-2$ variables.
		The projection of an even-parity relation onto these $N-2$ variables
		is full, since the omitted one variable can restore even parity.
		Thus all pairwise overlap state projections agree; under identity
		evolution, their full finite trajectory languages agree as well.
		A global assignment must satisfy
		$\sum_{j\ne i}x_j=0\pmod2$ for each $i$.  Subtracting any two
		equations gives $x_i=x_j$ for all $i,j$.  Hence the global state set
		has at most two elements (and exactly one when $N$ is even), whereas
		each local even-parity facet has $2^{N-2}>2$ states for every $N\ge4$.
		Therefore some local states admit no global realization.
		The primal graph is $K_N$, whose sole unguarded minimal clique
		is $V$; hence $\nu(\mathcal H_N)=N$.
	\end{proof}

	\begin{dremark}[Growth of the parity witness]
		The parity relations have premise scopes of size $N-1$, increasing
		with $N$.  The example shows that geometric rank alone is insufficient for unrestricted non-Horn relational atlases. The growing premise scopes motivate the uniform semantic recovery bounds below.
	\end{dremark}

	Theorem~\ref{thm:greatest-relation}, Corollary~\ref{wn11:thm:greatest-relation-exact} and Theorem~\ref{wn11:thm:holonomy} give a
	witness-sensitive non-Horn globalization theory for two clearly
	specified settings: completely general finite-state decision, and
	structural classification of invertible edge transports.  Their
	building blocks belong to established areas---relational database
	joins, finite automata, permutation constraints, and sheaf-style
	global-section obstructions
	\cite{Yannakakis1981,wn09:BeeriFaginMaierYannakakis1983,wn11:AbramskyBrandenburger2011,wn11:AtseriasKolaitis2020}.
	Theorem~\ref{wn11:thm:holonomy} yields a native recovery criterion for composed path assertions.  For restriction-stable projectively decomposable classes, the following
	sections obtain guard bounds from semantic witnesses. Theorems~\ref{wn12:thm:helly-classwise}
	and~\ref{wn13:thm:classwise-higher-suff} state the additional fiber-Helly
	hypotheses under which these bounds apply to the original charts.

	\subsection{Bounded semantic witnesses and sharp relational globalization}
	\label{wn12:sec:pb2-rel-global}

	A binary local obstruction rules out a pair of observations. To extend a partial assignment by one coordinate, the proof counts the candidate values excluded by those pairs.

	\subsubsection{The problem and a first obstruction}

	The canonical relational semantics developed in Section~\ref{wn11:sec:nonhorn-relational} allows
	negative constraints and actual intermediate-state witnesses.
	It does not by itself imply a threshold of the form
	$\operatorname{dprk}(\mathscr D)<\nu(\mathcal H)$.
	In particular, a bound on the arity of local semantic nogoods
	fails to control the simultaneous satisfiability of many compatible local witness demands.
	This section proves that failure constructively, then establishes sharp criteria for projectively decomposable relations.

	Let $V$ be a finite nonempty descriptor set, $\mathcal H$ a finite cover of $V$, and $A_v$ a nonempty finite alphabet for each $v\in V$.
	For $U\subseteq V$ write $A^U=\prod_{v\in U}A_v$, with the empty product a singleton.
	A local relational atlas is a family of nonempty $R_I\subseteq A^I$ indexed by $I\in\mathcal H$.
	It is \emph{projectively compatible} if
	\begin{equation}\label{wn12:eq:pb2-compatible}
		\pi_{I\cap J}R_I=\pi_{I\cap J}R_J\qquad(I,J\in\mathcal H).
	\end{equation}
	The canonical global relation is $M=\Join_{I\in\mathcal H}R_I$.
	We say the atlas \emph{globalizes exactly} if
	$\pi_I M=R_I$ for every original chart.
	No independent global-state restriction is imposed in this section: adding one can destroy exactness even on a join-tree cover.

	\begin{ddefinition}[Supported assignments, semantic nogoods, and rooted assertions]
		\label{wn12:def:pb2-nogood}
		Let $R\subseteq A^I$ be a relation and $p\in A^U$
		a partial assignment, where $U\subseteq I$.
		The \emph{witness set} of $p$ in $R$ is
		\[
		\mathsf{Wit}_R(p):=\{x\in R:x|_U=p\}.
		\]
		\begin{enumerate}[label=(\roman*),leftmargin=*,nosep]
			\item \textbf{Supported assignment.}
			The assignment $p$ is supported if
			$\mathsf{Wit}_R(p)\ne\varnothing$,
			equivalently, $p\in\pi_U R$.
			
			\item \textbf{Semantic and minimal nogoods.}
			The assignment $p$ is a semantic nogood if
			$\mathsf{Wit}_R(p)=\varnothing$.
			It is a minimal nogood if, additionally,
			\[
			\mathsf{Wit}_R(p|_B)\ne\varnothing
			\qquad\text{for every }B\subsetneq U.
			\]
			
			\item \textbf{Realizable rooted assertion.}
			For $v\in I$ and $a\in A_v$, a rooted assertion
			$p\Rightarrow x_v=a$ is realizable if
			\[
			\varnothing\ne\mathsf{Wit}_R(p)
			\subseteq\{x\in R:x_v=a\}.
			\]
		\end{enumerate}
		These notions depend only on the relation $R$,
		not on its rule presentation.
	\end{ddefinition}

	\begin{dproposition}[No-go: binary local nogoods and unbounded
		geometric rank]\label{wn12:prop:pb2-no-go}
		For every $n\ge3$ there is a compatible atlas on
		\[
		V_n=[n],\qquad\mathcal H_n=\{V_n\setminus\{i\}:i\in V_n\},
		\qquad\nu(\mathcal H_n)=n,
		\]
		whose chart relations and \emph{all their restrictions} have
		minimal nogoods of size at most two and admit full local
		extension of every supported partial assignment, but whose
		canonical global relation is empty.
		The obstruction persists under nonblocking, reversible identity
		transitions at every finite time horizon.
	\end{dproposition}

	\begin{proof}
		Let $A_v=[n-1]$ and take $R_I$ to be all injective assignments
		$I\to[n-1]$. Any assignment of distinct colors on $U\subseteq I$
		extends to an injective assignment on $I$, since $|I|=n-1$.
		Consequently, $\pi_U R_I$ consists of all injective maps
		$U\to[n-1]$; every minimal nogood repeats one color on a pair.
		The projections on $I\cap J$ are identical for adjacent charts, so
		\eqref{wn12:eq:pb2-compatible} holds. But an assignment in $M$ would
		inject $n$ variables into $n-1$ values, impossible. The primal graph
		is complete, and $V_n$ is its smallest unguarded clique, hence
		$\nu(\mathcal H_n)=n$. Identity transitions give precisely the constant paths with state in $R_I$; Lemma~\ref{wn12:lem:time-diagonal} proves that these remain PB2 at every horizon.
	\end{proof}

	\begin{dremark}[A nonempty but inexact global relation]
		Let $n=4$, let the alphabet be $\{0,1,2,3\}$, and on each
		three-variable facet permit either the all-zero tuple or an
		injective tuple using only $\{1,2,3\}$.
		Every projection is described by unary/binary restrictions:
		selected variables must be either all zero or nonzero and pairwise
		distinct. Pairwise atlas compatibility holds. The sole global
		assignment is all zero, although each chart permits six other tuples.
		Thus failure of exact globalization need not mean that no global
		state exists at all.
	\end{dremark}

	\subsubsection{Restriction-stable projective binary decomposability}

	For a relation $R\subseteq A^I$ define its binary reconstruction on
	$U\subseteq I$ by
	\begin{equation}\label{wn12:eq:binary-recon}
		B_2(R,U)=
		\left\{p\in A^U:
		p|_S\in\pi_S R\quad\text{for every }S\subseteq U,
		\ 1\le |S|\le2\right\}.
	\end{equation}

	\begin{ddefinition}[Projectively binary-decomposable relations]
		\label{wn12:def:pb2}
		With $B_2(R,U)$ the unary/binary consistency set defined above, a relation $R\subseteq A^I$ is \textbf{projectively binary-decomposable}, abbreviated \textsf{PB2}, if
		\begin{equation}\label{wn12:eq:pb2-definition}
			\pi_U R=B_2(R,U)\qquad\text{for every }U\subseteq I.
		\end{equation}
		The class $\mathscr R_q$ consists of all nonempty finite PB2 relations on finite coordinate sets, with each $|A_v|\le q$, where $q\ge2$.
		For atlas quantifiers it denotes the induced class of atlases.
		The relations may exclude tuples and need not be Horn closures.
	\end{ddefinition}

	\begin{dlemma}[Semantic characterization and restriction stability]
		\label{wn12:lem:pb2-nogood}
		A nonempty finite relation $R$ is PB2 if and only if every nogood in every
		semantic projection $\pi_U R$ has a nogood restriction supported
		on at most two coordinates (equivalently, every minimal nogood in
		every projection has size at most two). The PB2 property is stable
		under arbitrary coordinate projection, coordinate relabeling, and
		Cartesian products on disjoint coordinate sets.
	\end{dlemma}

	\begin{proof}
		A tuple fails membership in $B_2(R,U)$ exactly when one of its
		unary or binary restrictions fails membership in the corresponding
		projection of $R$. This proves the equivalence with
		\eqref{wn12:eq:pb2-definition}, including the empty-scope convention.
		If $S\subseteq I$, any further projection of $\pi_S R$ is
		$\pi_U R$ with $U\subseteq S$, so all the required equations persist.
		Relabeling merely transports the equations bijectively. For a
		product of relations on disjoint coordinate sets, all supports
		factor; forbidden patterns have a forbidden restriction in one
		factor, proving the last claim.
	\end{proof}

	The hereditary requirement in Definition~\ref{wn12:def:pb2} is essential.
	For example, over $A=\{0,1,2\}$ the four-variable relation
	\begin{equation}\label{wn12:eq:nonhereditary}
		R=\{(x,y,z,w):w\ne x,\ w\ne y,\ w\ne z\}
	\end{equation}
	is the conjunction of three binary inequalities and therefore
	binary-decomposable on its full scope. But its projection onto
	$(x,y,z)$ excludes exactly the assignments using all three values;
	all of that projection's binary projections are full. Hence the
	projection is not binary-decomposable. The missing condition
	is stability after hiding a common witness $w$.

	\subsubsection{An intrinsic fiber-Helly recovery index}

	For a compatible atlas, define the intrinsic unary and binary
	relations $D_v=\pi_{\{v\}}R_I$ and
	$P_{uv}=\pi_{\{u,v\}}R_I$ whenever a chart contains the indicated
	coordinates. Compatibility makes these independent of the chosen
	original chart. For $a\in D_u$ and $uv$ an edge of the primal graph,
	write
	\begin{equation}\label{wn12:eq:fiber}
		S_{v\leftarrow u}(a)=\{b\in D_v:(a,b)\in P_{uv}\}.
	\end{equation}
	This set is nonempty because $P_{uv}$ projects onto $D_u$.
	It is the complete set of values at $v$ compatible with one
	observed neighbor value $u=a$.

	\begin{ddefinition}[Fiber-Helly numbers, class index, and recovery rank]
		\label{wn12:def:fiber-rank}
		Let $\mathcal A$ be a compatible finite relational atlas
		on a coordinate set $V$. For each $v\in V$, define
		the fiber family
		\[
		\mathcal F_v(\mathcal A):=
		\{S_{v\leftarrow u}(a):
		uv\text{ is a primal edge},\ a\in D_u\},
		\]
		where the fibers are given by~\eqref{wn12:eq:fiber}.
		\begin{enumerate}[label=(\roman*),leftmargin=*,nosep]
			\item \textbf{Fiber-Helly number.}
			The fiber-Helly number $h_v(\mathcal A)$ at $v$
			is the least integer $h\geq1$ such that, for every
			finite subfamily $\mathcal F\subseteq\mathcal F_v(\mathcal A)$,
			\[
			\bigcap_{S\in\mathcal F}S=\varnothing
			\Longrightarrow
			\exists\,\mathcal G\subseteq\mathcal F:
			\quad |\mathcal G|\leq h,\quad
			\bigcap_{S\in\mathcal G}S=\varnothing.
			\]
			Thus $h_v(\mathcal A)=1$ if no subfamily has
			empty intersection.
			
			\item \textbf{Class fiber-Helly index.}
			For a restriction-stable class $\mathscr D$ of finite
			relational atlases, its fiber-Helly index is
			\[
			h(\mathscr D):=
			\sup_{\mathcal A\in\mathscr D}\max_{v\in V(\mathcal A)}
			h_v(\mathcal A),
			\]
			where $V(\mathcal A)$ is the coordinate set of $\mathcal A$.
			
			\item \textbf{Fiber recovery rank.}
			The specialized fiber recovery rank of $\mathscr D$ is
			\[
			\operatorname{dprk}_{\rm fib}(\mathscr D)
			:=\max\{3,\,1+h(\mathscr D)\}.
			\]
		\end{enumerate}
		These quantities depend only on the unary and binary
		projections of the relations in the class.
	\end{ddefinition}

	Any nonempty family of subsets of a set of at most $q$ elements has Helly number at most $q$:
	if its total intersection is empty, choose for each potential value one set excluding that value.
	Consequently $h(\mathscr R_q)\le q$.
	The injectivity	atlas of Proposition~\ref{wn12:prop:pb2-no-go} contains $q$ fibers
	$A\setminus\{1\},\ldots,A\setminus\{q\}$ with empty total intersection, but no smaller inconsistent subfamily.
	Hence
	\begin{equation}
		\label{wn12:eq:rank-value}
		h(\mathscr R_q)=q,
		\qquad
		\operatorname{dprk}_{\rm fib}(\mathscr R_q)=q+1.
	\end{equation}

	\begin{dlemma}[Guarded witness recovery]
		\label{wn12:lem:guarded-recovery}
		Let $\mathcal A$ be a compatible PB2 atlas, and suppose
		$a_u\in D_u$ has been assigned on a clique $N$ of the primal
		graph, with every pair among $N$ belonging to its permitted
		binary relation. Let $v$ be adjacent to every vertex of $N$.
		If the fiber-Helly number at $v$ is at most $h$ and every primal
		clique of size at most $h+1$ is guarded, then
		\begin{equation}\label{wn12:eq:common-fiber}
			\bigcap_{u\in N}S_{v\leftarrow u}(a_u)\ne\varnothing.
		\end{equation}
	\end{dlemma}

	\begin{proof}
		If the intersection were empty, select a subfamily indexed by
		$U\subseteq N$, $|U|\le h$, with empty intersection.
		The set $B=U\cup\{v\}$ is a primal clique and hence lies in an
		original chart $J$. The assigned $U$-tuple satisfies all unary
		and binary projections of $R_J$, because they are the canonical
		relations $D_w,P_{ww'}$. The PB2 property of $\pi_U R_J$ therefore
		lifts that tuple to some element $r\in R_J$. Its value $r_v$
		belongs to every fiber indexed by $U$, a contradiction.
		Thus~\eqref{wn12:eq:common-fiber} holds.
	\end{proof}

	\subsubsection{The sharp non-Horn recovery--globalization theorem}

	The obstruction rank $\nu(\mathcal H)$ of Definition~\ref{def:nu} equals three for
	nonchordal primal graphs and otherwise is the cardinality of the
	smallest unguarded clique, with value $\infty$ if none exists
	\cite{MaPackaging2026}. Hence $\nu(\mathcal H)>r\ge3$ implies
	primal chordality and guarding of every clique of size at most $r$.

	\begin{dtheorem}[Fiber-Helly recovery implies exact globalization]
		\label{wn12:thm:helly-classwise}
		Let $\mathscr D$ be a restriction-stable class of finite PB2
		relational atlases with uniformly finite intrinsic fiber-Helly
		number $h=h(\mathscr D)$.
		\begin{enumerate}[label=(\roman*),leftmargin=*,nosep]
			\item If $h=1$, every compatible atlas in $\mathscr D$ globalizes exactly on every finite cover.
			\item For $h\ge2$, every compatible atlas in $\mathscr D$ on a finite cover $\mathcal H$ globalizes exactly if
			\begin{equation}\label{wn12:eq:helly-suff}
				\operatorname{dprk}_{\rm fib}(\mathscr D)
					<\nu(\mathcal H).
			\end{equation}
		\end{enumerate}
		No independent global-state constraint is allowed.
	\end{dtheorem}

	\begin{proof}
		\leavevmode
		\begin{enumerate}[label={(iii)},leftmargin=*,nosep]
			\item[(i)] First suppose $h=1$. Every individual fiber is nonempty, so the defining Helly property implies that no finite family of fibers can have empty intersection. Begin with any specified original-chart tuple and assign the remaining coordinates in any order. At each new vertex intersect the fibers prescribed by its already assigned neighbors (or take the whole attained domain if none); the intersection is nonempty without a clique hypothesis on those neighbors. All guarded unary and binary constraints are then satisfied, hence every original-chart tuple constraint is satisfied by PB2. This proves exact globalization for every cover.
			\item[(ii)] For $h\ge2$, the inequality implies chordality of the primal graph and guards
			all cliques of size at most $h+1$. Fix an original chart $I$ and
			$a\in R_I$. A chordal graph admits a perfect elimination ordering~\cite{Dirac1961,FulkersonGross1965} which leaves the chosen clique $I$ to the end: root a clique tree at a maximal clique containing $I$ (different graph components can be connected along empty separators), eliminate vertices from leaf cliques outside the root, and eliminate the remaining root vertices outside $I$ before eliminating $I$.
			Reverse this ordering to assign variables, beginning with $a$ on
			$I$. The previously assigned neighbors $N$ of each new vertex
			$v$ form a clique. Maintain the unary and binary compatibility
			relations on all assigned coordinates. By
			Lemma~\ref{wn12:lem:guarded-recovery}, some value at $v$ satisfies every
			neighbor constraint. If $N=\varnothing$, choose any value in $D_v$.
			Induction produces a tuple $x\in A^V$ extending $a$ and
			satisfying all canonical unary and binary relations on the primal
			graph. For each original chart $J$, its restriction $x|_J$
			satisfies all unary/binary projections of $R_J$. Since $R_J$ is
			PB2, $x|_J\in R_J$. Thus $x\in\Join_J R_J$ extends $a$.
			Every original chart tuple has a global extension, which is exact
			projective globalization. \qedhere
		\end{enumerate}
	\end{proof}

	\begin{dtheorem}[Sharp bounded-domain projective-binary classification]
		\label{wn12:thm:pb2-sharp}
		For each integer $q\ge2$, the class $\mathscr R_q$ satisfies
		\begin{equation}\label{wn12:eq:sharp}
			\operatorname{UG}_{\mathscr R_q}(\mathcal H)
				\quad\Longleftrightarrow\quad
				q+1<\nu(\mathcal H).
		\end{equation}
		The threshold is classwise and sharp on \emph{every} cover,
		not merely on a family of extremal examples.
	\end{dtheorem}

	\begin{proof}
		Sufficiency follows from Theorem~\ref{wn12:thm:helly-classwise}
		and~\eqref{wn12:eq:rank-value}. Conversely, assume
		$\nu(\mathcal H)\le q+1$.
		If the primal graph is nonchordal, choose an induced chordless
		cycle $C$. Assign binary alphabets to its vertices and singleton
		alphabets outside $C$. Since each original chart is a clique, it
		contains at most two adjacent vertices of $C$. Impose equality
		on every cycle edge except one, on which impose inequality;
		all other local coordinates are unconstrained. These relations
		are nonempty, PB2 under every restriction, and agree on overlaps;
		but no global assignment can satisfy the twisted cycle.
		
		Otherwise the primal graph is chordal and has an unguarded clique
		$K$ of size $r\le q+1$. Give the vertices of $K$ an alphabet of
		$r-1\le q$ colors and all other vertices singleton alphabets.
		On chart $I$ impose injectivity on $I\cap K$, leaving the singleton
		coordinates fixed. Since $K$ is unguarded, $|I\cap K|\le r-1$;
		every locally supported partial injection extends on that chart.
		Therefore all local relations are PB2 and pairwise compatible.
		No assignment can inject all $r$ vertices of $K$ into $r-1$
		colors. Thus exact globalization fails in either case, proving
		necessity.
	\end{proof}

	\begin{dremark}[Classwise versus atlas-specific ranks]
		The number $q+1$ is the sharp \emph{universal} cutoff for
		$\mathscr R_q$, not necessarily the minimal cutoff for any
		proper subclass. In the PB2 setting, the fiber-Helly index gives a semantic local witness-recovery bound and an associated classwise guarding criterion based on fiber intersections, in contrast to the original-chart closure recovery underlying $\operatorname{ipr}$ in Paper~I~\cite{MaPackaging2026}. The proof uses the relation between consistency and domain size developed by Dechter and van Beek \cite{Dechter1992,DechterVanBeek1997}.
	\end{dremark}

	\subsubsection{All-horizon dynamical globalization}

	Let $X_I\subseteq A^I$ be a finite local state set and $R_I\subseteq X_I\times X_I$ its time-homogeneous transition relation.
	Define the local path relation $\mathsf T_I(m)$ on $I\times\{0,\ldots,m\}$.
	Its members are complete local trajectories, not separately chosen one-step transitions.
	Call the system \emph{temporally PB2} if $\mathsf T_I(m)$ is PB2 on its full space--time coordinate set for every $m\ge0$.
	The class is restriction-stable as a class of \emph{observed trace families}:
	spatial or temporal projection preserves PB2, although an observed projection need not itself possess a first-order Markov presentation.

	\begin{dlemma}[Time-diagonal stability of PB2]
		\label{wn12:lem:time-diagonal}
		Let $R$ be PB2 on $I$. For each horizon $m$, the relation of constant paths $\Delta_m(R)$ on $I\times[0,m]$ is PB2.
	\end{dlemma}

	\begin{proof}
		This is the $r=2$ identity-update case of the diagonal stability argument in Lemma~\ref{wn13:lem:pbr-diagonal}. The proof there uses only projective decomposability and equality of repeated time copies, not the globalization theorem.
	\end{proof}

	\begin{dtheorem}[Sharp all-horizon temporal PB2 globalization]
		\label{wn12:thm:temporal-pb2}
		Fix $q\ge2$. Consider the class of finite, time-homogeneous relational transition atlases with $|A_v|\le q$ that are temporally PB2 (in particular, each chart has paths at every finite horizon) and for which
		\begin{equation}\label{wn12:eq:temporal-compat}
			\pi_{(I\cap J)\times[0,m]}\mathsf T_I(m)
			=\pi_{(I\cap J)\times[0,m]}\mathsf T_J(m)
			\quad\forall I,J,\ m\ge0.
		\end{equation}
		No global admissible-state constraints are imposed beyond the natural join of local states.
		Universal exact globalization at every finite horizon holds for this class if and only if
		\begin{equation}\label{wn12:eq:temporal-sharp}
			q+1<\nu(\mathcal H).
		\end{equation}
		When this holds, the canonical global transition relation is an exact extension of every prescribed finite local trajectory.
		The same universal equivalence holds if the class is restricted to nonblocking local laws.
	\end{dtheorem}

	\begin{proof}
		For every $m$, apply Theorem~\ref{wn12:thm:pb2-sharp} to the
		finite time-cylinder cover
		$\mathcal H^{[0,m]}=\{I\times[0,m]:I\in\mathcal H\}$.
		Its obstruction rank equals $\nu(\mathcal H)$, since duplicating
		each primal vertex into a clique of time copies preserves
		chordality and the smallest unguarded clique size (as proved in
		Section~\ref{sec:covers-assertions}). The alphabet bound $q$ is unchanged because each
		space--time coordinate still has alphabet $A_v$. The PB2 and
		compatibility assumptions therefore make the join of local path
		relations exact at every horizon. This join is precisely the
		trajectory relation of the canonical global one-step relation,
		which demands simultaneous satisfaction of every chart update.
		Thus one and the same global transition system works for all $m$.
		Exactness of observed paths does not in general assert that every
		compatible global initial state has a successor; that stronger
		property requires a separate global nonblocking condition.
		
		For sharpness, use either of the static counterexamples from the
		preceding theorem and put identity dynamics on every local
		state space. These systems are nonblocking and reversible.
		Each temporal path relation is the diagonal time-copy relation
		of the original PB2 state relation, which remains PB2 by
		Lemma~\ref{wn12:lem:time-diagonal}.
		Overlap compatibility holds at every horizon, while failure of
		exact state globalization at horizon zero persists. Hence
		$\nu(\mathcal H)\le q+1$ cannot guarantee exact all-horizon
		globalization.
	\end{proof}

	\begin{dlemma}[Majority-operation certificate for temporal PB2 globalization]
		\label{wn12:lem:majority}
		Suppose each alphabet $A_v$ carries a ternary majority operation
		$m_v$ with
		$m_v(a,a,b)=m_v(a,b,a)=m_v(b,a,a)=a$.
		If each local state relation $X_I$ and its transition relation
		$R_I$ are closed under the coordinatewise product majority
		operation, then every nonempty local trajectory relation
		$\mathsf T_I(m)$ is PB2. If, in addition, the local state sets are
		nonempty and the local laws are nonblocking, the positive implication of
		Theorem~\ref{wn12:thm:temporal-pb2} applies whenever
		$\nu(\mathcal H)>q+1$ and the full trace-compatibility condition
		holds. Closure under an operation alone does not exclude empty
		trajectory relations.
	\end{dlemma}

	\begin{proof}
		Intersection of majority-closed cylinder relations remains
		majority-closed; existential projection also preserves
		majority closure. Hence each $\mathsf T_I(m)$ and all of its
		projections are majority-closed.
		To see that a majority-closed relation $S$ is binary-decomposable,
		induct on its arity $d$. For $d\le2$ the claim is immediate.
		If a tuple $x$ has all its pairs in the corresponding projections,
		then by induction each deletion $x|_{[d]\setminus\{i\}}$ belongs to
		that projection. Choose three distinct indices $i,j,k$ and full
		witnesses in $S$ agreeing with $x$ except possibly at those
		respective coordinates. Their coordinatewise majority equals $x$,
		so $x\in S$. Apply the argument to every projection.
	\end{proof}

	\begin{dexample}[Nontrivial reversible evolution]
		For $q=2$ and the four-facet cover
		$\mathcal H=\{[4]\setminus\{i\}:i\in[4]\}$, take on every chart
		the binary states with all coordinates equal, and update every
		coordinate by bit flip. The complete local trajectory relations
		are PB2, the atlas is compatible at all horizons, and
		$\nu(\mathcal H)=4>3$. The canonical global system has exactly
		two alternating global trajectories, and their projections
		recover all prescribed local paths. This is reversible dynamics,
		not a positive Horn closure.
	\end{dexample}

	The binary criterion combines guarded elimination with bounded-domain consistency \cite{Dechter1992,DechterVanBeek1997}. The majority operation supplies one algebraic route to the needed decomposability \cite{BakerPixley1975}. For larger witnesses, one excluded value can require several previously observed coordinates.

	\subsection{Higher-arity semantic witness amalgamation and sharp globalization}
	\label{wn13:sec:pb-r-wa}

	For an obstruction of size at most $r$, excluding one candidate value may use $r-1$ other coordinates. Guarding their union leads to the sharp domain-and-arity threshold.

	\subsubsection{The intrinsic problem}

	The preceding PB2 theorem concerns finite relational charts whose semantic
	projections are completely determined by unary and binary marginals. We now
	allow arbitrary higher-order exclusions, including ternary parity, and require
	that the relevant width remain bounded \emph{after every coordinate
		projection}. This is a semantic property rather than a bound on the arity of a
	chosen rule presentation. We work with the cover obstruction rank
	$\nu(\mathcal H)$ of~\cite{MaPackaging2026}: it is $3$ when the primal graph
	has an induced cycle of length at least four, and otherwise the size of the
	smallest unguarded primal clique, with $\infty$ if no such clique exists.
	
	Fix integers $q,r\ge2$. Let $V$ be a finite variable set with
	$1\le |A_v|\le q$, and let $\mathcal H$ be a finite cover of $V$.
	A \emph{relational atlas} has nonempty relations
	$R_I\subseteq\prod_{v\in I}A_v$, one for each $I\in\mathcal H$.
	It is \emph{compatible} if $\pi_{I\cap J}R_I=\pi_{I\cap J}R_J$
	for all $I,J$. Put $M=\Join_{I\in\mathcal H}R_I$.
	An atlas \emph{globalizes exactly} when $\pi_I M=R_I$ for every
	original chart $I$. As in Section~\ref{wn12:sec:pb2-rel-global}, no independent constraint on
	admissible global states is imposed: such a constraint can separately destroy
	exactness even when $\nu(\mathcal H)=\infty$.
	
	Write $\mathscr R_{q,r}$ for the resulting class in Definition~\ref{def:projective-r-decomposability}. The projection formula~\eqref{wn13:eq:pb-r-def} is the property needed at each elimination step.

	\begin{dlemma}[Presentation independence and hereditary projection]
		\label{wn13:lem:hereditary-pbr}
		The \textsf{PB}$r$ property is invariant under coordinate relabeling and
		is inherited by all coordinate projections. Its validity is determined by
		the semantic relation, not by the rules used to present it. For $r=2$ it is
		exactly the \textsf{PB2} property of Section~\ref{wn12:sec:pb2-rel-global}.
	\end{dlemma}

	\begin{proof}
		The condition in~\eqref{wn13:eq:pb-r-def} for a projection $\pi_U R$
		only concerns projections $\pi_S R$ with $S\subseteq U$;
		these are unchanged when first projecting to $U$. A tuple outside
		$\pi_U R$ satisfies all its $\le r$ restrictions exactly when it has no
		nogood restriction of size at most $r$, proving the equivalent nogood
		characterization. Coordinate renaming transports all the equations.
	\end{proof}

	\subsubsection{A semantic higher-arity fiber-Helly invariant}

	When $S\subseteq V$ is a primal clique of size at most $r$, assume it
	is guarded by a chart and denote its canonical relation by
	$P_S=\pi_SR_I$, for any $I\supseteq S$.
	Compatibility makes this independent of the chart; $P_\varnothing$ consists of the unique empty tuple. For
	$v\notin B$, $|B|\le r-1$, and a supported tuple
	$a_B\in P_B$ on a guarded clique $B\cup\{v\}$, define the
	\emph{semantic witness fiber}
	\begin{equation}\label{wn13:eq:higher-fiber}
		\mathcal F_{v\mid B}(a_B)
		=\bigl\{b\in A_v:(a_B,b)\in P_{B\cup\{v\}}\bigr\}.
	\end{equation}
	The fiber is nonempty because $P_{B\cup\{v\}}$ projects onto $P_B$.

	\begin{ddefinition}[Witness-Helly index and witness-amalgamation rank]
		\label{wn13:def:witness-index}
		For a chart atlas $\mathcal A$ and a vertex $v$, its \textbf{witness-Helly number} $h_v(\mathcal A)$ is the least $h\ge1$ such that any subfamily of the nonempty semantic fibers~\eqref{wn13:eq:higher-fiber} with empty intersection already has a subfamily of at most $h$ members with empty intersection.
		If no inconsistent family exists, put $h_v=1$.
		For a restriction-stable class $\mathscr D$ of \textsf{PB}$r$ atlases define its \textbf{witness-Helly index} and \textbf{witness-amalgamation rank} by
		\[
		h(\mathscr D)=\sup_{\mathcal A\in\mathscr D}\max_v h_v(\mathcal A),
		\qquad
		\operatorname{dprk}_{\mathrm{wa}}(\mathscr D)
		=\max\{3,\,1+(r-1)h(\mathscr D)\}.
		\]
		In the supremum all admitted coordinate restrictions are included.
		This \emph{witness-amalgamation rank} is semantic and presentation independent.
		It supplies a sufficient recovery bound for the stated class.
	\end{ddefinition}

	A family of subsets of a nonempty set with at most $q$ values has Helly
	number at most $q$: select, for each value, a set that excludes it.
	Thus $h(\mathscr R_{q,r})\le q$. The reverse inequality is attained:
	choose $q$ disjoint $(r-1)$-element groups $B_c$, one for each color $c$,
	and a pivot $v$. On the chart $B_c\cup\{v\}$ use assignments with at
	most $r-1$ occurrences of the color $c$; all other colors are
	unrestricted. These charts have compatible full unary $v$-projections,
	are \textsf{PB}$r$, and the fiber conditioned on $B_c$ being entirely
	colored $c$ is $A_v\setminus\{c\}$. All $q$ fibers are required to
	witness empty intersection. Consequently,
	\begin{equation}\label{wn13:eq:rank-qr}
		h(\mathscr R_{q,r})=q,\qquad
			\operatorname{dprk}_{\mathrm{wa}}(\mathscr R_{q,r})
			=1+q(r-1).
	\end{equation}

	\begin{dlemma}[Guarded higher-arity witness amalgamation]
		\label{wn13:lem:high-amalgamation}
		Let $\mathcal A$ be compatible and \textsf{PB}$r$ on a cover for which
		every primal clique of size at most $d=1+h(r-1)$ is guarded.
		Suppose a tuple has been assigned on a clique $N$ and every restriction
		to a clique of size at most $r$ is semantically supported.
		Let $v$ be adjacent to all of $N$, and suppose the witness-fiber
		Helly index at $v$ is at most $h$. Then some $b\in A_v$ makes every
		restriction to at most $r$ coordinates of $N\cup\{v\}$ semantically
		supported.
	\end{dlemma}

	\begin{proof}
		For every $B\subseteq N$ of size at most $r-1$, form the
		nonempty fiber $\mathcal F_{v\mid B}(a_B)$; each scope $B\cup\{v\}$
		is a guarded clique since $d\ge r$. Suppose their total intersection
		were empty. Pick at most $h$ such fibers with empty intersection
		and put $U=\bigcup B$. Then $|U\cup\{v\}|\le1+h(r-1)=d$ and
		$U\cup\{v\}$ is a guarded clique, say in chart $J$.
		The assigned tuple on $U$ has all its $\le r$ restrictions supported
		in $\pi_S R_J$, by the invariant and overlap compatibility. The
		\textsf{PB}$r$ property of $\pi_U R_J$ therefore supplies a witness
		in $R_J$ extending the $U$-tuple. Its $v$-value belongs to each of the
		selected fibers, a contradiction. Any value in the total intersection
		is the required extension.
	\end{proof}

	\subsubsection{Exact globalization beyond binary decomposability}

	To extend a chart tuple, eliminate vertices in an order that leaves its chart last. The witness-fiber bound gives a simultaneous value at each extension step.

	\begin{dtheorem}[Semantic witness-amalgamation sufficiency]
		\label{wn13:thm:classwise-higher-suff}
		Let $\mathscr D$ be a restriction-stable class of compatible \textsf{PB}$r$ relational atlases with finite classwise witness-Helly
		index $h$.
		If
		\begin{equation}\label{wn13:eq:general-rank-suff}
			\operatorname{dprk}_{\mathrm{wa}}(\mathscr D)
				<\nu(\mathcal H),
		\end{equation}
		then every atlas in $\mathscr D$ on $\mathcal H$ globalizes exactly.
	\end{dtheorem}

	\begin{proof}
		The strict inequality implies that the primal graph is chordal and
		that every clique of size at most $1+h(r-1)$ is guarded. Fix any
		original chart $I$ and any $a\in R_I$. A perfect-elimination order can end with $I$: while vertices outside $I$ remain, remove a simplicial vertex outside $I$. In a complete graph any such vertex works. A noncomplete chordal graph has two nonadjacent simplicial vertices \cite{Dirac1961}, and at most one can lie in the clique $I$. Induced subgraphs remain chordal. Read the resulting order backwards,
		starting with $a$ on $I$. At each step, the previously assigned
		neighbors $N$ of the new vertex $v$ form a clique. Maintain that all
		assigned subsets which form a clique and have at most $r$ vertices
		are supported by their canonical projections. This is initially true
		inside $I$. Lemma~\ref{wn13:lem:high-amalgamation} supplies a value at
		$v$ preserving the invariant; when $N=\varnothing$ choose any value
		from its nonempty unary projection. At the end, every $\le r$ subset
		of each original chart $J$ is supported. By~\eqref{wn13:eq:pb-r-def},
		the constructed global tuple restricts to a member of $R_J$. Thus
		it belongs to $M=\Join_JR_J$ and extends $a$. Since $I$ and $a$
		were arbitrary, $\pi_I M=R_I$ for every original chart.
	\end{proof}

	\begin{dcorollary}[Bounded-alphabet witness-amalgamation criterion]
		\label{wn13:cor:pb-r-upper}
		For every $q,r\ge2$, compatible \textsf{PB}$r$ atlases with
		$|A_v|\le q$ globalize exactly whenever
		\begin{equation}\label{wn13:eq:pbr-upper}
			1+q(r-1)<\nu(\mathcal H).
		\end{equation}
	\end{dcorollary}

	\begin{proof}
		Every family of nonempty fibers at a coordinate with at most $q$ values has Helly number at most $q$: if the intersection is empty, choose one excluding fiber for each value. Hence $h\le q$, and Theorem~\ref{wn13:thm:classwise-higher-suff} applies. This is its bounded-alphabet specialization.
	\end{proof}

	\subsubsection{Sharpness on every finite cover}

	The preceding bound is not merely sufficient.
	An exact counterexample uses finite cardinality capacities rather than parity constraints.

	\begin{dproposition}[Balanced capacity obstruction]
		\label{wn13:prop:balanced-capacity}
		For any $q,r\ge2$, put $s=1+q(r-1)$. On the simplex-boundary cover
		$\mathcal H_s=\{[s]\setminus\{i\}:i\in[s]\}$, give each coordinate
		$q$ possible colors. On every facet permit the assignments with
		exactly $r-1$ occurrences of each color. The resulting atlas is
		nonempty, compatible, and \textsf{PB}$r$, but its global join is empty.
		Moreover $\nu(\mathcal H_s)=s$.
	\end{dproposition}

	\begin{proof}
		Every facet has $q(r-1)=s-1$ coordinates. Its balanced relation is
		nonempty. For any subset $U$ of a facet, a partial coloring extends
		to a balanced coloring if and only if each color occurs at most
		$r-1$ times on $U$. Thus every forbidden partial coloring contains
		$r$ occurrences of one color, making every semantic projection
		$r$-decomposable. The two adjacent facets have the same projection
		onto their common $(s-2)$-coordinate intersection: the color counts
		there are bounded by $r-1$, with precisely one color occurring
		$r-2$ times. Hence the atlas is compatible.
		If a full assignment satisfied every facet, the color multiset of
		the facet omitting $i$ would be the same for every $i$. Comparing
		facets omitting $i$ and $j$ forces the values at $i$ and $j$ to be
		equal. Hence all coordinates carry one color, contradicting the
		requirement that every facet uses all $q\ge2$ colors. The primal
		graph is complete, and the full $s$-clique is the smallest unguarded
		clique. Thus $\nu=s$.
	\end{proof}

	\begin{dlemma}[Capacity obstruction on every unguarded clique]
		\label{wn13:lem:ung-capacity}
		Let $q,r\ge2$ be integers. Suppose the primal graph contains an unguarded clique $K$ of size
		$3\le s\le1+q(r-1)$. Then there is a nonempty, compatible,
		\textsf{PB}$r$ atlas on the original cover with alphabets of size at most $q$ whose global join is empty.
	\end{dlemma}

	\begin{proof}
		Replace $K$ by an inclusion-minimal unguarded subset $K'$ of it and write $s'=|K'|$. Every proper subset of $K'$ is guarded. Singletons are covered, and every edge of a primal clique is guarded by definition, so $3\le s'\le s$. Choose integer capacities $0\le c_j\le r-1$, $1\le j\le q$, of sum $s'-1$, at least two of them positive. This is possible by allocating one unit to each of two colors and distributing the remainder within the available capacities.
		
		Use the $q$ colors on $K'$ and singleton alphabets elsewhere. On $I$ allow exactly the assignments with at most $c_j$ occurrences of each color $j$ on $I\cap K'$. Since $|I\cap K'|\le s'-1$, a partial assignment that respects all capacities extends to $I$ by filling unused capacity. Thus the same inequalities describe every overlap projection, the charts are nonempty, and every locally forbidden assignment has a witness on at most $c_j+1\le r$ cells. In a full coloring of $K'$ some color $j$ occurs more than $c_j$ times. Since at least two capacities are positive, $c_j+1\le s'-1$. Select $c_j+1$ such cells. They form a proper subset of $K'$, hence lie together in an original chart. That chart rejects the coloring. This proves that the join is empty.
	\end{proof}

	\begin{dtheorem}[Sharp classwise higher-order relational globalization]
		\label{thm:pbr}
		For all integers $q,r\ge2$ and every finite cover $\mathcal H$,
		\begin{equation}\label{wn13:eq:sharp-pbr}
			\operatorname{UG}_{\mathscr R_{q,r}}(\mathcal H)
				\quad\Longleftrightarrow\quad
				1+q(r-1)<\nu(\mathcal H).
		\end{equation}
		Thus $\operatorname{dprk}_{\mathrm{wa}}(\mathscr R_{q,r})
		=1+q(r-1)$ is the exact universal geometric cutoff for this class.
	\end{dtheorem}

	\begin{proof}
		Sufficiency follows from
		Corollary~\ref{wn13:cor:pb-r-upper}. Conversely assume
		$\nu(\mathcal H)\le1+q(r-1)$. If the primal graph is chordal,
		its smallest unguarded clique has size
		$3\le s=\nu(\mathcal H)\le1+q(r-1)$; apply
		Lemma~\ref{wn13:lem:ung-capacity}. If the primal graph is nonchordal,
		choose an induced chordless cycle. On its vertices use two colors and outside it singleton alphabets, with equality on every cycle edge except one, which enforces inequality. Each chart meets that induced cycle in at most one edge,
		so its relation is nonempty, projectively binary-decomposable
		(hence \textsf{PB}$r$), and the local overlap projections agree.
		The cycle constraints cannot all hold globally. This supplies a
		counterexample in either case and proves sharpness on every cover.
	\end{proof}

	\begin{dexample}[The first genuinely ternary sharp boundary]
		\label{wn13:ex:ternary-capacity}
		For $q=2,r=3$, take five Boolean variables and every four-variable
		facet as a chart. Require each facet to contain exactly two zeros and
		two ones. Every local relation has six tuples; a partial assignment
		is locally supported if and only if it contains no three equal bits,
		so every local minimal nogood has size three. The three-coordinate
		overlap projections are the common \emph{not-all-equal} relation.
		Nevertheless, a global assignment cannot make every four-coordinate
		facet balanced. Here
		\[
		\nu(\mathcal H)=5=1+2(3-1),\qquad M=\varnothing.
		\]
		The upper bound is thus sharp already for ternary semantic recovery
		on Boolean domains, not merely for $r=2$.
	\end{dexample}

	\subsubsection{Time-indexed relational globalization}

	Let each original chart $I$ carry a finite, time-homogeneous
	relational transition system with state set $X_I\subseteq A^I$ and
	relation $E_I\subseteq X_I^2$. Write $\mathsf T_I(m)$ for its complete
	local path relation on $I\times[0,m]$. Say the atlas is
	\emph{temporally PB$r$} if \emph{every} $\mathsf T_I(m)$ is
	\textsf{PB}$r$ on its individual space--time coordinates, for every
	$m\ge0$. Require exact equality of the projected local path
	relations on every chart intersection, for every horizon.

	\begin{dlemma}[Diagonal and coordinatewise bijective trajectory stability]
		\label{wn13:lem:pbr-diagonal}
		For any \textsf{PB}$r$ relation $R$, its constant-path diagonal
		relation on $I\times[0,m]$ is \textsf{PB}$r$ for all $m\ge0$.
		More generally, this holds for the trajectory relation obtained
		from $R$ by coordinatewise \emph{bijective} state updates that
		preserve $R$.
	\end{dlemma}

	\begin{proof}
		A partial space--time tuple consistent with every $\le r$ marginal
		must assign compatible values to all observed time copies of the
		same spatial variable, since $r\ge2$. In the bijective-update case,
		transport each observed value back to time zero. The resulting
		assignment on the represented spatial variables is supported on
		every subset of size at most $r$; the \textsf{PB}$r$ property of
		$R$ supplies a common global state witness. Transport it through
		the coordinatewise updates. All other implications follow by
		projection. Constant paths correspond to identity updates.
	\end{proof}

	\begin{dtheorem}[Sharp all-horizon temporal PB$r$ globalization]
		\label{wn13:thm:temporal-pbr}
		Fix $q,r\ge2$. For finite, time-homogeneous temporally \textsf{PB}$r$ relational transition atlases with at most $q$ values per spatial coordinate, nonempty local path relations at every horizon, complete overlap-trace compatibility at all horizons, and no independent global-state restrictions, universal exact finite-horizon globalization holds if and only if
		\begin{equation}\label{wn13:eq:sharp-temporal}
			\boxed{1+q(r-1)<\nu(\mathcal H).}
		\end{equation}
		Under the strict inequality the canonical global one-step relation realizes all prescribed finite local paths exactly. The same universal equivalence holds in the subclass of nonblocking local laws, since necessity is witnessed by identity dynamics. Exactness does not automatically imply nonblocking from \emph{every} compatible global state.
	\end{dtheorem}

	\begin{proof}
		For a fixed horizon $m$, the time-expanded atlas lives on
		$\mathcal H^{[0,m]}=\{I\times[0,m]:I\in\mathcal H\}$.
		By the temporal-cylinder identity~\eqref{wn08:eq:cylinder-rank} its obstruction rank is
		$\nu(\mathcal H)$; each time-coordinate alphabet retains size
		at most $q$. Apply Theorem~\ref{thm:pbr} to obtain an exact
		natural join of full local path relations for every $m$. That joint
		relation is precisely the horizon-$m$ path relation (with $m+1$ time slices) of the canonical
		one-step global transition relation: a global path satisfies every
		chart's transition constraints at every step exactly when each of
		its chart restrictions is a permitted local path. Thus the same
		global relation works for all $m$.
		For necessity, use the static counterexample from the preceding
		theorem and give each chart reversible identity dynamics.
		The corresponding constant-path relations remain \textsf{PB}$r$
		by Lemma~\ref{wn13:lem:pbr-diagonal} and satisfy complete trace
		compatibility, yet they fail at horizon zero.
	\end{proof}

	\begin{dcorollary}[A one-step near-unanimity certificate]
		\label{wn13:cor:nu-cert}
		Suppose each alphabet $A_v$ carries an $(r+1)$-ary near-unanimity
		operation $\mu_v$ and, for every original chart $I$, both its state
		relation $X_I$ and its one-step transition relation $E_I$ are closed
		under the corresponding coordinatewise operations. Then every nonempty
		complete local trajectory relation $\mathsf T_I(m)$ is \textsf{PB}$r$.
		With nonempty local state sets and nonblocking local laws, exact
		all-horizon globalization follows from full separator-trace
		compatibility and $1+q(r-1)<\nu(\mathcal H)$. The nonemptiness
		condition is separate from closure under the operations.
	\end{dcorollary}

	\begin{proof}
		The complete path relation is an intersection of cylinder relations
		expressing the initial-state and successive transition constraints.
		Each cylinder is closed under the product near-unanimity operation;
		so is the intersection, as are all its coordinate projections.
		The classical near-unanimity decomposability argument applies also to these coordinate-dependent operations~\cite{wn13:BartoKrokhinWillard2017,BakerPixley1975}: induct on the arity of a nonempty projected relation. For arity at most $r$ the claim is immediate. At larger arity, a tuple whose restrictions to at most $r$ cells are supported has, by induction, a witness agreeing off each single deleted coordinate. Choose $r+1$ distinct deleted coordinates and combine their witnesses using the product near-unanimity operation. At every coordinate at least $r$ witnesses agree with the tuple, so the resulting tuple is in the relation. This proves $r$-decomposability of every projection.
		Apply Theorem~\ref{wn13:thm:temporal-pbr}. This is a sufficient checkable
		certificate, not a characterization of all temporally
		\textsf{PB}$r$ systems. In particular, all-horizon separator-trace
		compatibility remains a separate hypothesis.
	\end{proof}

	\begin{dremark}[Why one-step arity is insufficient]
		\label{wn13:rem:arity-blow-up}
		Let $b=x_1\oplus x_2$ and $b=x_3\oplus x_4$ be two ternary
		Boolean constraints sharing a hidden witness $b$.
		Joining them and then eliminating $b$ produces
		$x_1\oplus x_2\oplus x_3\oplus x_4=0$, which has full
		three-coordinate projections but is not \textsf{PB3}.
		Hence it is essential to require the bound for the complete
		semantic trajectory relations at all horizons; a bound on the
		arity of primitive one-step clauses does not suffice.
	\end{dremark}

	Guarded elimination uses the same domain-size principle as classical relational consistency \cite{Dechter1992,DechterVanBeek1997}. Near-unanimity yields sufficient decomposability hypotheses \cite{wn13:BartoKrokhinWillard2017}. Under temporal projection one must verify the property for each complete history relation.

	\subsection{Compositional recovery and structural rank inequalities}
	\label{wn16:sec:composition-rank}
	
	Independent products preserve the larger of two witness widths. A shared coordinate changes the problem because the local witnesses must now agree on a common interface history.

	Theorem~\ref{wn15:thm:effective-exact-rank} established that the intrinsic
	temporal witness width $\operatorname{twit}(\mathbb D)$ of a finite-state
	observed transition system is computable in $\mathbb N_{\ge1}\cup\{\infty\}$: it is the maximum support
	size of a minimal infeasible partial temporal observation, with the convention
	that a system without nogoods has width one. We now study its behaviour under
	composition, synchronization, observation restriction, and feedback.
	
	We retain the exact observation granularity. Let
	$\mathbb D=(X,E,(A_v)_{v\in V},\lambda)$ be finite, with
	$\lambda:X\to\prod_{v\in V}A_v$, and all initial states allowed.
	An observation word $p=(p_0,\ldots,p_m)$ has entries in
	$\prod_v(A_v\sqcup\{\star\})$; it is \emph{feasible} if a path
	$x_0Ex_1E\cdots Ex_m$ agrees with its specified cells.
	Write $\operatorname{Feas}(\mathbb D)$ for the set of all feasible partial
	words, $\mathcal M(\mathbb D)$ for the minimal infeasible partial words,
	and $\operatorname{wt}(p)$ for its number of specified space--time cells.
	For applications of Theorem~\ref{wn13:thm:temporal-pbr}, assume paths of every
	finite length exist; stronger nonblocking assumptions are stated when needed.
	The concept of minimality refers to removal of \emph{one specified cell},
	not to deletion of an entire time slice. This distinction is critical.

	\subsubsection{Independent products and observation restriction}

	For two systems on disjoint named coordinate sets, their independent
	clock-synchronous product has product states and simultaneous component
	transitions; its global trajectory language is the Cartesian product of the
	component trajectory languages at every horizon.

	\begin{dtheorem}[Independent-product maximum law]\label{wn16:thm:independent-product}
		If $V_1\cap V_2=\varnothing$ and both factors admit paths at every
		horizon, then
		\begin{equation}\label{wn16:eq:product}
			\ \operatorname{twit}(\mathbb D_1\times\mathbb D_2)
				=\max\{\operatorname{twit}(\mathbb D_1),
				\operatorname{twit}(\mathbb D_2)\}.\
		\end{equation}
		The equality includes infinite widths. The same conclusion holds for
		finitely many mutually independent factors.
	\end{dtheorem}

	\begin{proof}
		The two-factor assertion, including infinite width, is Lemma~\ref{lem:product} applied to the product trajectory languages. Their nonemptiness at every horizon allows induction over any finite number of independent factors.
	\end{proof}

	\begin{dproposition}[Witness-width monotonicity under restriction and subsampling]
		\label{wn16:prop:projection}
		\leavevmode
		\begin{enumerate}[label=(\roman*),leftmargin=*,nosep]
			\item Let $U\subseteq V$ be a subset of the \emph{named observed coordinates}.
			Then
			\begin{equation}\label{wn16:eq:projection}
				\operatorname{twit}(\mathbb D|_U)
					\le\operatorname{twit}(\mathbb D).
			\end{equation}
			\item For each $k\ge1$, let $\mathbb D^{[k]}$ have the same underlying state set $X$ and observation map $\lambda$, with transition relation $E^k$ (relational composition of $k$ steps). Then
			\begin{equation}\label{wn16:eq:sampling}
				\operatorname{twit}(\mathbb D^{[k]})\le
				\operatorname{twit}(\mathbb D).
			\end{equation}
		\end{enumerate}
	\end{dproposition}

	\begin{proof}
		\leavevmode
		\begin{enumerate}[label={(iii)},leftmargin=*,nosep]
			\item[(i)] For the first assertion, extend a partial $U$-word by $\star$ on the
			other coordinates at every time. Its feasibility in $\mathbb D$ equals
			that of the original word in $\mathbb D|_U$. Thus every minimal
			impossible $U$-word extends to a minimal impossible full word with the same
			support.
			\item[(ii)] For the second assertion, embed a word of length $m+1$ for
			$E^k$ into a word of length $km+1$ for $E$ by recording its letters at
			times $0,k,\ldots,km$ and putting $\star$ at intervening times.
			Relational path composition gives equivalence of feasibility and preserves
			one-cell relaxation minimality. Both inequalities follow. \qedhere
		\end{enumerate}
	\end{proof}

	\begin{dremark}[Two different meanings of restriction]
		Coordinate hiding or temporal subsampling only removes observations.
		Imposing new state or transition constraints is \emph{not} an observation
		restriction; the latter operation can increase the witness width, as
		proved below. Likewise, an arbitrary many-to-one \emph{symbol recoding}
		need not be a coordinate projection and is not covered by
		Proposition~\ref{wn16:prop:projection}.
	\end{dremark}

	\subsubsection{Synchronization and the interface-amalgamation defect}

	Let $\mathbb D_A$ and $\mathbb D_B$ have coordinate sets $V_A,V_B$
	with identical alphabets on overlaps. Their \emph{state-synchronized product} $\mathbb C=\mathbb D_A\bowtie\mathbb D_B$ consists of pairs of states whose observations agree on the shared coordinates, with simultaneous component transitions between such pairs.
	All its observable coordinates are those in $V=V_A\cup V_B$.
	For a partial word $p$ on $V$, let $p_A$ and $p_B$ denote its coordinate
	restrictions. Always
	\begin{equation}\label{wn16:eq:feas-one-way}
		p\in\operatorname{Feas}(\mathbb C)
		\Longrightarrow
		p_A\in\operatorname{Feas}(\mathbb D_A)
		\ \text{and}\ 
		p_B\in\operatorname{Feas}(\mathbb D_B).
	\end{equation}
	The converse can fail: the two component realizations may use different values for the same currently unspecified interface cells.

	\begin{ddefinition}[Interface-amalgamation nogoods and intrinsic width]
		\label{wn16:def:amw}
		Let $\mathbb C=\mathbb D_A\bowtie\mathbb D_B$ be a
		synchronized product. For a partial assignment $p$,
		let $p_A$ and $p_B$ denote its restrictions to the
		respective component coordinates at the same horizon.
		\begin{enumerate}[label=(\roman*),leftmargin=*,nosep]
			\item \textbf{Interface-amalgamation nogood.}
			A minimal nogood $p\in\mathcal M(\mathbb C)$ is an
			interface-amalgamation nogood if both component
			restrictions are feasible. The family of such nogoods is
			\[
			\mathcal A(A,B):=
			\{p\in\mathcal M(\mathbb C):
			p_A\in\operatorname{Feas}(\mathbb D_A),\
			p_B\in\operatorname{Feas}(\mathbb D_B)\}.
			\]
			
			\item \textbf{Intrinsic interface-amalgamation width.}
			The interface-amalgamation width of the two components
			under their specified interface is
			\begin{equation}\label{wn16:eq:amw}
				\operatorname{amw}(A,B):=
				\sup_{p\in\mathcal A(A,B)}\operatorname{wt}(p),
				\qquad \sup\varnothing:=0,
			\end{equation}
			where $\operatorname{wt}(p)$ counts the specified
			observable cells of $p$.
		\end{enumerate}
		The width depends only on the component observable languages
		and their common-witness semantics under the specified
		interface. It is invariant under coordinatewise alphabet
		bijections and equivalent hidden-state presentations,
		but may change under regrouping of observable coordinates.
	\end{ddefinition}

	\begin{dtheorem}[Compositional rank inequality and exact witness condition]
		\label{wn16:thm:amw-bound}
		Suppose the synchronized product has at least one path at every horizon.
		\begin{enumerate}[label=(\roman*),leftmargin=*,nosep]
			\item \emph{Width bound.}
			Then
			\begin{equation}\label{wn16:eq:main-bound}
				\boxed{\operatorname{twit}(\mathbb C)
					\le\max\{\operatorname{twit}(\mathbb D_A),
					\operatorname{twit}(\mathbb D_B),
					\operatorname{amw}(A,B)\}.}
			\end{equation}
			
			\item \emph{Exact witness condition.}
			Moreover, $\operatorname{amw}(A,B)=0$ if and only if the	\emph{partial-observation  witness-amalgamation property}
			\begin{equation}\label{wn16:eq:amalg-property}
				\ p_A\in\operatorname{Feas}(\mathbb D_A),\quad
					p_B\in\operatorname{Feas}(\mathbb D_B)
					\ \Longrightarrow\ p\in\operatorname{Feas}(\mathbb C)\
			\end{equation}
			holds for every finite partial temporal observation $p$.
			
			\item \emph{Finiteness.} In particular, when both component widths are finite,
			\begin{equation}
				\label{wn16:eq:finiteness-equivalence}
				\operatorname{twit}(\mathbb C)<\infty
				\quad\Longleftrightarrow\quad
				\operatorname{amw}(A,B)<\infty.
			\end{equation}
		\end{enumerate}
	\end{dtheorem}

	\begin{proof}
		\leavevmode
		\begin{enumerate}[label={(iii)},leftmargin=*,nosep]
			\item[(i)] Let $p$ be a minimal impossible word of $\mathbb C$. If $p_A$ is
			impossible in $A$, minimality of $p$ forces its entire specified support
			to lie in $V_A$ and forces $p_A$ to be minimal impossible in $A$:
			otherwise erasing a cell outside $V_A$, or erasing a cell of a
			nonminimal impossible $p_A$, leaves the composite infeasible. Hence
			$\operatorname{wt}(p)\le\operatorname{twit}(A)$.
			The same argument applies if $p_B$ is impossible in $B$. In the only
			remaining case both projected partial words are feasible separately;
			then $p\in\mathcal A(A,B)$ and
			$\operatorname{wt}(p)\le\operatorname{amw}(A,B)$. Take the supremum.
			
			\item[(ii)] If~\eqref{wn16:eq:amalg-property} fails, take a partial word feasible in
			both components but impossible in the synchronized product, and erase
			specified cells until it becomes minimal impossible for the product.
			Every restriction of a feasible component word remains feasible, so the
			result belongs to $\mathcal A(A,B)$, which is therefore nonempty.
			Conversely, any member of $\mathcal A(A,B)$ violates
			\eqref{wn16:eq:amalg-property}. Because the product has a path at every horizon, all minimal nogoods have positive support. Thus $\operatorname{amw}(A,B)=0$ precisely when $\mathcal A(A,B)$ is empty, which is precisely the witness property.
			
			\item[(iii)] Finally $\operatorname{amw}(A,B)\le
			\operatorname{twit}(\mathbb C)$ by definition; combine this with
			\eqref{wn16:eq:main-bound} when both factor ranks are finite to obtain
			\eqref{wn16:eq:finiteness-equivalence}. \qedhere
		\end{enumerate}
	\end{proof}

	\begin{dremark}[Why gluing partial witness languages matters]
		For arbitrary temporal languages, feasibility does not distribute through synchronous conjunction:
		$\operatorname{Feas}(L_A\cap L_B)$ need not equal
		$\operatorname{Feas}(L_A)\cap\operatorname{Feas}(L_B)$.
		This is an existential-witness issue, not a problem of primitive constraint arity.
		The condition~\eqref{wn16:eq:amalg-property} is strong;
		when it holds,~\eqref{wn16:eq:main-bound} reduces to the maximum bound.
		These observations are closely related to established database join dependencies and compositional automata theory.
	\end{dremark}

	\subsubsection{A finite-state feedback catastrophe}

	Let $a,b,c$ be Boolean coordinates.
	Module $A$ observes $(a,b)$, allows all four states, and imposes the one-step constraint
	\begin{equation}
		\label{wn16:eq:delay}
		a_{t+1}=b_t,
	\end{equation}
	with no restriction on the next value of $b$.
	Module $B$ observes $(a,b,c)$, allows states satisfying
	\begin{equation}
		\label{wn16:eq:xor}
		b_t=a_t\oplus c_t,
	\end{equation}
	and permits all transitions between its four states.
	Both systems are finite-state and nonblocking. Their exact widths are
	\begin{equation}
		\label{wn16:eq:factor-ranks}
		\operatorname{twit}(A)=2,
		\qquad \operatorname{twit}(B)=3.
	\end{equation}
	Indeed each failed partial observation of $A$ contains a violated two-cell equality across consecutive times; each failed partial observation of $B$ contains a violated three-cell parity equation at one time.
	Their synchronized product has only four states and
	obeys
	\begin{equation}\label{wn16:eq:integrator}
		a_{t+1}=a_t\oplus c_t.
	\end{equation}

	\begin{dtheorem}[Finite-rank components with infinite-rank feedback]
		\label{wn16:thm:feedback}
		For the two systems in~\eqref{wn16:eq:delay}--\eqref{wn16:eq:xor},
		\begin{equation}\label{wn16:eq:feedback-infinite}
			\operatorname{twit}(A\bowtie B)
				=\operatorname{amw}(A,B)=\infty,
		\end{equation}
		even though the factors have widths $2$ and $3$, and the
		composite has only four states.
	\end{dtheorem}

	\begin{proof}
		For each $m\ge1$, consider a partial observation that specifies
		\[
		a_0=0,\qquad a_m=1,\qquad
		c_0=\cdots=c_{m-1}=0,
		\]
		and leaves every other cell, including every $b$-cell, free.
		Equation~\eqref{wn16:eq:integrator} makes this impossible: it forces
		$a_m=a_0$. But erasing either endpoint permits a different endpoint,
		and erasing any one $c_j$ permits $c_j=1$ and hence $a_m=1$.
		Therefore it is a minimal impossible word of support $m+2$. This family witnesses unbounded width.
		It is separately feasible in $A$ (choose suitable unconstrained
		$b$-values) and in $B$ (choose the $a$-values independently over time).
		Thus it belongs to $\mathcal A(A,B)$ for every $m$.
		The widths in~\eqref{wn16:eq:factor-ranks} follow from the independent
		pairwise and ternary constraints already noted.
	\end{proof}

	\begin{dremark}[No contradiction with earlier atlas globalization]
		The two factors in this feedback example are \emph{not} an exact,
		overlap-compatible atlas as defined in Section~\ref{sec:local-dynamical-atlases}: their complete
		$\{a,b\}$-trajectory languages differ. They are independently
		specified modules composed by synchronizing shared signals.
		The example disproves a blanket \emph{compositional rank bound},
		not the earlier exact-atlas theorems.
	\end{dremark}

	\subsubsection{Effective recognition of the amalgamation defect}

	The Section~\ref{wn15:sec:intrinsic-recognition} observer tracks $(\textit{first},Q,U)$, where $Q$ is
	the set of reachable states for a partial word and $U$ is the family
	of reachable-state sets for its one-cell relaxations. For the
	synchronized product, augment this memory with two more subsets
	$Q_A,Q_B$, tracking feasibility of the same partial word projected
	onto the factors.

	\begin{dtheorem}[Effective interface-defect rank]
		\label{wn16:thm:effective-amw}
		For finite-state modules $A,B$, the language $\mathcal A(A,B)$ is
		effectively regular. Therefore finiteness of $\operatorname{amw}(A,B)$
		is decidable. When finite, the exact value is computable; when
		infinite, a positive-weight cycle yields an explicit pumping family
		of minimal composite nogoods that are feasible separately.
	\end{dtheorem}

	\begin{proof}
		Run the Section~\ref{wn15:sec:intrinsic-recognition} single-deletion monitor for the synchronized
		product in parallel with the ordinary subset observers for
		$A$ and $B$. At the end accept if the composite partial word is
		minimal impossible and both factor state subsets are nonempty.
		The finite product monitor recognizes exactly $\mathcal A(A,B)$.
		Give each transition the number of specified observed coordinates as
		weight. The same positive-weight-cycle / longest-accepting-path
		criterion as in Section~\ref{wn15:sec:intrinsic-recognition} proves the dichotomy and computes the
		exact bound. A search stopped before completion gives no certified verdict.
	\end{proof}

	\subsubsection{Changing observation granularity can destroy invariance}

	The composite system~\eqref{wn16:eq:integrator} has an instantaneous
	state represented by two bits $(a,c)$; the value of $b$ is then
	uniquely determined. Its temporal witness width on these two
	\emph{named} coordinates is infinite. Now regard the very same four
	states as four atomic labels in \emph{one} coordinate $s$.

	\begin{dproposition}[Witness width depends on observable coordinate atomization]
		\label{wn16:prop:atomization}
		The two presentations above describe isomorphic unlabelled
		transition graphs under a bijection of the four instantaneous states,
		yet
		\begin{equation}\label{wn16:eq:atomization}
			\operatorname{twit}(\mathbb C;(a,c))=\infty,
				\quad \operatorname{twit}(\mathbb C;s)=2.
		\end{equation}
	\end{dproposition}

	\begin{proof}
		The first claim follows from Theorem~\ref{wn16:thm:feedback} because the
		minimal nogoods in its proof never specify $b$.
		For the atomic four-symbol presentation, the allowed adjacent state
		pairs are the edges of the state graph. From any state $(a,c)$, every
		state $(a'',c'')$ is reachable in two steps: select the intermediate
		$c'$ so that $a''=a\oplus c\oplus c'$.
		Every state has a predecessor and a successor, and the two-step relation is full; hence the relation of paths of any length at least two is full on endpoint pairs. With all initial states allowed, a partial word of atomic states is therefore feasible exactly when no two specified states at \emph{adjacent time slices} violate the one-step relation. Longer gaps, leading blanks and trailing blanks can be filled. Such a violation is a minimal two-cell nogood.
		There are no larger minimal nogoods. The symbol map
		$s\leftrightarrow(a,c)$ is bijective as a map of instantaneous states,
		but it changes the decomposition into independently erasure-capable
		observable coordinates.
	\end{proof}

	\begin{dremark}[Meaning of intrinsic]
		The temporal witness width is invariant under presentation changes
		\emph{preserving its named coordinate observation semantics}, such
		as coordinatewise alphabet bijections and state-space duplications
		with the same observable trace language. It is \emph{not} invariant
		under arbitrary regrouping or splitting of observable coordinates.
		A general dynamical packaging theory must therefore include its
		coordinate/support structure in the object or explicitly restrict
		the morphisms under which a rank is asserted to be intrinsic.
	\end{dremark}

	\subsubsection{Compositional globalization certificate}

	Let $\mathscr D$ be a finite relational atlas whose original chart
	dynamics each arise by synchronizing bounded-rank modules. Fix a binary synchronization order in each chart. If all component ranks and the interface-amalgamation widths at each successive binary composition (including intermediate products) are bounded by $r$, and the resulting charts admit paths of every finite length, Theorem~\ref{wn16:thm:amw-bound} inductively bounds the complete chart trajectory relations by \textsf{PB}$r$ at every horizon. Under the remaining assumptions of Section~\ref{wn13:sec:pb-r-wa} (compatible chart traces at all horizons, alphabet bound $|A_v|\le q$ with $q\ge2$, complete chart histories given by first-order relations on the full chart tuples, and no additional global-state constraints), one gets
	\begin{equation}\label{wn16:eq:global-corollary}
		\ 1+q(\max\{2,r\}-1)<\nu(\mathcal H)
			\Longrightarrow
			\text{the specified atlas $\mathscr D$ is all-horizon exact}.\
	\end{equation}
	The rank $r$ may be computed by Section~\ref{wn15:sec:intrinsic-recognition}'s exact width algorithm
	and Theorem~\ref{wn16:thm:effective-amw}.

	\begin{dremark}[Partial words and synchronized products]
		The positive product, automata-intersection, and weighted-cycle
		steps use classical regular-language and automata constructions.
		Shared-witness consistency is closely related to database joins,
		constraint satisfaction, and synchronous transition-system products.
		The interface-amalgamation width records those minimal forbidden patterns
		that become impossible only after the shared histories are synchronized.
	\end{dremark}

	\paragraph{Examples of interface-amalgamation width.}

	The bounds give the following widths for fixed observable coordinates. In the feedback example a minimal inconsistent pattern grows by one cell at each time step.

	\begin{center}
		\small
		\begin{tabular}{@{}p{7.2cm}c@{}}
			\toprule
			Finite-state dynamical package & $\operatorname{twit}$ \\
			\midrule
			Delayed-copy factor $A$: $a_{t+1}=b_t$ & $2$ \\
			Static parity factor $B$: $b_t=a_t\oplus c_t$ & $3$ \\
			Synchronized feedback $A\bowtie B$ on $(a,b,c)$ & $\infty$ \\
			Feedback observed only through $(a,c)$ & $\infty$ \\
			Identical four-state graph encoded as one atomic coordinate & $2$ \\
			Feedback sampled every second step and observed via $(a,c)$ & $1$ \\
			Independent identity/flip product & $2$ \\
			\bottomrule
		\end{tabular}
	\end{center}
	For independent modules, width is the maximum of the component widths. In a synchronized product the additional obstruction is the existence of a common separator history, as measured by the interface-amalgamation width.

	\section{Exact temporal formation descent and reconfiguration}
	\label{sec:temporal-formation}

	This section studies temporal descent through a partial or multivalued formation relation. Theorem~\ref{dyn47:thm:rel-descent} characterizes existence, the greatest output law, forced edges and finite minimum-cost realizations. Theorem~\ref{dyn48:thm:reconfiguration} extends the characterization to changing state types and transition masks; the cancellation criterion then determines when exact descent is coherent along parallel paths.

	\subsection{Exact temporal descent of relational formation}
	\label{dyn47:subsec:relational-descent}

	A native realization is generally a relation, not a function. Retain a
	declared interpretation $J\subseteq S\times Y$ of source formation states
	and realized output states. For instance this can be a finite family of
	actual formation judgments with their native marks, or an admitted
	operational relation when used as a relation between state sorts. We do
	not identify those two native predicates. Given a source law $E\subseteq
	S^2$, define a \emph{relational descent} to be an output law $D\subseteq
	Y^2$ satisfying the exact square in the category of relations
	\begin{equation}\label{dyn47:eq:descent}
		J\circ E=D\circ J.
	\end{equation}
	Composition is in right-to-left order. This compares reachable output
	sets from a formation state. It does not require every source step to
	lift from every chosen current output representative; that forward-lift
	condition alone also does not imply the square.

	\begin{theorem}[Exact relational descent: existence, solution family, rigidity, and cost]
		\label{dyn47:thm:rel-descent}
		Write $J(s)=\{y:(s,y)\in J\}$ and
		$B_s=\bigcup_{t\in E(s)}J(t)$. For $y\in Y$ put
		\[
		H_y=\bigcap_{s:y\in J(s)}B_s,
		\qquad D^*(y)=H_y,
		\]
		with the empty intersection taken to be $Y$.
		\begin{enumerate}[label=(\roman*),leftmargin=*,nosep]
			\item \emph{Existence.} A relational descent
			exists if and only if
			\begin{equation}\label{dyn47:eq:desc-test}
				\boxed{\forall s\in S,\quad B_s=\bigcup_{y\in J(s)} H_y.}
			\end{equation}
		
			\item \emph{Greatest law and all solutions.} If it exists, $D^*$ is the unique greatest descent. All descents are
			exactly the subrelations of $D^*$ meeting every demand
			\[
			C_{s,z}=\{(y,z):y\in J(s),\ z\in H_y\}
			\qquad(s\in S,\ z\in B_s).
			\]
			
			\item \emph{Rigidity.} When the square is solvable, an edge in $D^*$ is forced in every descent
			if and only if some demand is its singleton. In particular, in the
			solvable case the descent is unique exactly when every edge of $D^*$
			is forced.
		
			\item \emph{Minimum cost.} If $S$ and $Y$ are finite and nonnegative edge costs are declared, the minimum total cost of a descent is the corresponding finite hitting-set optimum; its value is zero when there are no demands.
		\end{enumerate}
	\end{theorem}

	\begin{proof}
		\leavevmode
		\begin{enumerate}[label={(iii)},leftmargin=*,nosep]
			\item[(i)] The square says $\bigcup_{y\in J(s)}D(y)=B_s$. Its inclusion from left
			to right forces each $D(y)\subseteq H_y$, and that inclusion is sufficient
			for soundness. Thus $D^*$ is the greatest candidate, and its coverage
			of all demands is both necessary and sufficient.
			
			\item[(ii)] For a subrelation of
			$D^*$ coverage is precisely the listed hitting condition.
			
			\item[(iii)] An edge is
			forced if a singleton demand requires it; otherwise removing just that
			edge from $D^*$ leaves every demand met. This proves rigidity.
			
			\item[(iv)] The minimum cost is the hitting-set optimum by (ii); with no demands the empty law has zero cost.
		\end{enumerate}
		The proof
		permits empty rows and outputs outside $J(S)=\bigcup_s J(s)$; rows of
		the latter are unconstrained unless a further marking is imposed.
	\end{proof}

	\begin{corollary}[Functional descent, composition, iterates, and path lifting]
		\label{dyn47:cor:desc-compose}
		\leavevmode
		\begin{enumerate}[label=(\roman*),leftmargin=*,nosep]
			\item \emph{Functional case.} For the graph $J$ of a surjection $q:S\to Y$, a descent exists exactly
			when $q(E(s))$ depends only on $q(s)$; then it is unique.
			
			\item \emph{Composition.} For composable
			relations $J:S\rightsquigarrow Y$ and $K:Y\rightsquigarrow Z$ with exact
			squares for $E,D$ and $D,F$, the composite $K\circ J$ has an exact square
			for $E,F$.
			
			\item \emph{Iterates.} Any exact square satisfies $J\circ E^n=D^n\circ J$ for every
			$n\ge0$.
			
			\item \emph{Path lifting.} Moreover an output path starting at any $y\in J(s)$ has a
			source lift starting at $s$ with all successive pairs in $J$.
		\end{enumerate}
	\end{corollary}

	\begin{proof}
		\leavevmode
		\begin{enumerate}[label={(iii)},leftmargin=*,nosep]
			\item[(i)] For a functional graph the union on the right of the square has only
			one row $D(q(s))$, and surjectivity forces every row.
			
			\item[(ii)] For composition,
			associativity gives $KJE=KDJ=FKJ$.
			
			\item[(iii)] The iterated square follows by
			induction, including the identity relation at $n=0$.
			
			\item[(iv)] For a chosen
			compatible pair $(s,y)$ and step $yDz$, the square gives a source
			successor $t$ with $z\in J(t)$. Iterate this implication along the
			specified output path. \qedhere
		\end{enumerate}
	\end{proof}

	\begin{example}[Exact descent is not representative-wise forward lifting]
		\label{dyn47:ex:desc-boundary}
		Let $S=\{s,t\}$, $Y=\{a,b,c\}$, $J(s)=\{a,b\}$, $J(t)=\{c\}$,
		$E=\{(s,t)\}$ and $D=\{(b,c)\}$. The square holds: both sides give
		$\{c\}$ at $s$ and the empty set at $t$. But from the compatible pair
		$(s,a)$ the source step to $t$ has no matching output step from $a$.
		The greatest descent adds $(a,c)$; it is not unique, as the demand at
		$(s,c)$ can be met by either output edge. Thus exact descent and the
		representative-wise lifting condition in the functional-realization discussion
		preceding Theorem~\ref{thm:exact-deterministic-chart-factor-criterion} must remain separate.
	\end{example}

	\begin{proposition}[Representative-wise lifting, joint descent, rigidity, and cost]
		\label{dyn47:prop:representative-descent}
		In the notation of Theorem~\ref{dyn47:thm:rel-descent}, require for each
		$(s,y)\in J$ both of the following properties:
		\begin{enumerate}[label=(\roman*),leftmargin=2em]
			\item for every $t\in E(s)$ there is $z\in J(t)$ with $z\in D(y)$;
			
			\item for every $z\in D(y)$ there is $t\in E(s)$ with $z\in J(t)$.
		\end{enumerate}
		
		\begin{enumerate}[label=(\alph*),leftmargin=*,nosep]
			\item \emph{Two-way lifting.} Such a representative-wise two-way lifting law exists if and only if
			$H_y\cap J(t)\ne\varnothing$ for every $s,t,y$ with $t\in E(s)$ and
			$y\in J(s)$. In that case its greatest law is $D^*$, and its solutions
			are precisely the subrelations of $D^*$ hitting each set
			\[
			L_{s,t,y}=\{(y,z):z\in H_y\cap J(t)\}.
			\]
			
			\item \emph{Joint exactness.} If both representative-wise lifting and the exact descent square
			are required, use the union of these demands and the sets $C_{s,z}$
			in Theorem~\ref{dyn47:thm:rel-descent}. Existence is equivalent to
			nonemptiness of every demand in the combined library.
		
			\item \emph{Rigidity and cost.} In either solvable
			problem rigidity again has the singleton test; for finite $S,Y$ with
			declared nonnegative edge costs the minimum has the finite hitting-set
			formulation.
		\end{enumerate}
	\end{proposition}

	\begin{proof}
		\leavevmode
		\begin{enumerate}[label={(iii)},leftmargin=*,nosep]
			\item[(a)] Condition (ii) is equivalent to $D(y)\subseteq B_s$ for every
			compatible $(s,y)$, or to $D\subseteq D^*$. Under this bound, (i)
			is precisely the requirement to meet each $L_{s,t,y}$. Thus $D^*$ solves the two-way problem exactly when every such demand is nonempty; empty source transitions create no forward demands.
			
			\item[(b)] For joint exactness both libraries must be met. Their demands are subsets of $D^*$, so $D^*$ solves the combined problem precisely when each is nonempty.
			
			\item[(c)] For either library the singleton and finite hitting-set arguments of Theorem~\ref{dyn47:thm:rel-descent} apply.
		\end{enumerate}
		The conditions are not equivalent to the exact
		square alone: Example~\ref{dyn47:ex:desc-boundary} gives one failure.
		Conversely for $J(s)=\{a\}, J(t)=\{b,c\}$, $E=\{(s,t)\}$, and
		$D=\{(a,b)\}$ both representative-wise conditions hold, but
		$D(J(s))=\{b\}\ne\{b,c\}=J(E(s))$.
	\end{proof}

	These exact squares compose along the declared formation relations. Structural signatures with several input sorts
	can use their product as the source, with any additional admission
	condition retained in the source states and relation $J$.

	\subsection{Relational reconfiguration, rigidity and path coherence}
	\label{dyn48:subsec:reconfiguration}

	For a finite directed reconfiguration graph with vertices $i$, let $S_i$ and $Y_i$ be the source and packaged-output state sets of type $i$, and let $J_i\subseteq S_i\times Y_i$ be its possibly partial, multivalued native realization relation. For each directed edge $e:i\to j$ prescribe a source transition relation $E_e\subseteq S_i\times S_j$; seek an output relation $D_e\subseteq Y_i\times Y_j$. The required square is $J_j\circ E_e=D_e\circ J_i$; any additional admissible output edges may be specified by a fixed mask $M_e\subseteq Y_i\times Y_j$, with $M_e=Y_i\times Y_j$ when no restriction is imposed. Composition is in temporal order from right to left. This is a typed extension of the descent law in Section~\ref{dyn47:subsec:relational-descent}, not a map requiring unchanged carriers across time.

	\begin{theorem}[Masked reconfiguration descent: existence, rigidity, cost, and composition]
		\label{dyn48:thm:reconfiguration}
		For each edge $e:i\to j$ put
		\[
		B_e(s)=\bigcup_{t\in E_e(s)}J_j(t),\qquad
		H_e(y)=M_e(y)\cap\!\bigcap_{s:\,y\in J_i(s)}B_e(s),
		\]
		with empty intersection $Y_j$. Here $M_e(y)$ is the allowed output fiber.
		\begin{enumerate}[label=(\roman*),leftmargin=*,nosep]
			\item \emph{Existence.} The exact masked square is solvable if and only if
			\begin{equation}\label{dyn48:eq:reconfiguration}
				B_e(s)=\bigcup_{y\in J_i(s)} H_e(y)\qquad(s\in S_i).
			\end{equation}
			
			\item \emph{All solutions.} If solvable, $H_e$ is the greatest output law; all solutions are precisely the subrelations $D_e\subseteq H_e$ meeting every demand
			\[
			\{(y,z): y\in J_i(s),\ z\in H_e(y)\}\qquad(s\in S_i,\ z\in B_e(s)).
			\]
			
			\item \emph{Rigidity.} For a solvable square, a greatest-law edge $(y,z)$ is forced exactly when some demand consists of it alone. Such a square is rigid exactly when every edge of $H_e$ is forced.
			
			\item \emph{Minimum cost.} For finite $S_i,S_j,Y_i,Y_j$ with declared nonnegative edge costs, the minimum is the finite hitting-set minimum for this demand family.
			
			\item \emph{Independent edges.} If the reconfiguration graph imposes no further identifications among its edge laws, all solutions on the graph form the product of the individual solution families.
			
			\item \emph{Path composition.} For every path of edges their composed output relation is an exact descent of the composed source relation between its endpoint types.
		\end{enumerate}
	\end{theorem}

	\begin{proof}
		\leavevmode
		\begin{enumerate}[label={(iii)},leftmargin=*,nosep]
			\item[(i)] Every output edge $y\to z$ must be allowed by the mask and must have $z\in B_e(s)$ for each source representative containing $y$, so every solution lies in $H_e$. That maximal relation is sound, and the reverse inclusion in the square is exactly the coverage test \eqref{dyn48:eq:reconfiguration}.
			
			\item[(ii)] For a subrelation the coverage condition is equivalent to meeting each displayed demand.
			
			\item[(iii)] A singleton demand forces its edge. Conversely if an edge is never a singleton demand, removing it from the greatest solution leaves all demands met, proving rigidity.
			
			\item[(iv)] The finite cost problem is the hitting-set problem by (ii).
			
			\item[(v)] Independent edge constraints give the stated product.
			
			\item[(vi)] Along two composable edges, substitute the first square into the second, obtaining $J_k E_f E_e=D_f D_e J_i$; iteration gives the path assertion. A path of length zero has the identity square. No equality of different output paths was assumed. \qedhere
		\end{enumerate}
	\end{proof}

	\begin{proposition}[Exact criterion for cancellation and path coherence]
		\label{dyn48:prop:path-coherence}
		\leavevmode
		\begin{enumerate}[label=(\roman*),leftmargin=*,nosep]
			\item \emph{Cancellation.} For any $J\subseteq S\times Y$, the right-composition comparison $D\mapsto D\circ J$ on relations $D\subseteq Y\times Z$ is injective for every nonempty $Z$ if and only if each $y\in Y$ has a singleton source witness $s$ with $J(s)=\{y\}$.
			
			\item \emph{Path coherence.} Consequently, if two parallel reconfiguration paths with the same endpoint realization data have the same source composite and their common initial realization has this witness property, their exact output composites must coincide.
			
			\item \emph{Masked uniqueness.} More generally, if that endpoint square has a unique exact descent in a declared mask containing both output composites, those composites coincide.
		\end{enumerate}
		Without such a condition, exact squares along every edge do not imply path coherence.
	\end{proposition}

	\begin{proof}
		\leavevmode
		\begin{enumerate}[label={(iii)},leftmargin=*,nosep]
			\item[(i)] A singleton witness reads off $D(y)$ as $(D\circ J)(s)$, giving injectivity. If a witness for some $y$ is missing, fix $z\in Z$. Compare the relation having the single output $z$ from every $y'$ with that relation after deleting $(y,z)$. A nonempty source fiber containing $y$ also contains a different output, so the two composites agree everywhere, including empty fibers. The relations differ.
			
			\item[(ii)] For the parallel paths the endpoint squares of Theorem~\ref{dyn48:thm:reconfiguration} have the same left side $J_j\circ E_{\rm path}$ and the same initial $J_i$ on the right. Apply injectivity to these squares.
			
			\item[(iii)] If the endpoint square has a unique exact descent within the declared mask, both output composites belong to the same singleton solution family.
		\end{enumerate}
		For failure take $S_i=\{s\}$, $Y_i=\{a,b\}$, $J_i(s)=\{a,b\}$ and $S_j=\{t\}$, $Y_j=\{z\}$, $J_j(t)=\{z\}$. Two parallel labeled edges with identical source relation $s\to t$ can have respective exact descents $a\to z$ and $b\to z$. The squares are exact but the output relations differ. The same issue arises with diamond paths by inserting intermediate identity steps. This proves the necessity of additional path comparisons in a general reconfiguration theory.
	\end{proof}

	These results admit empty realization fibers, changing source and output types and masks. They are relative to the specified reconfiguration graph and realization data. Simultaneously imposing arbitrary path equations on all unknown edge laws is not reduced to independent minima; the cancellation test is a sufficient way of ensuring a prescribed equality of source paths descends. Stronger representative-wise path lifting is still a separate condition as in Section~\ref{dyn47:subsec:relational-descent}.

	\section{Gluing probability laws on states and histories}
	\label{sec:gluing-probability-laws}

	This section develops globalization for probability laws, from single-time marginals to complete histories and finite-state generators. Theorem~\ref{dyn47:thm:weighted-join} identifies the canonical gluing and the dimension of the full marginal fiber on a join tree. Theorem~\ref{dyn48:thm:sharp-wug} gives the exact guard threshold on general covers for small independent history blocks. Theorem~\ref{dyn48:thm:projective} constructs the canonical complete-history law; the final results treat continuous-time globalization and finite-system reduction.

	Throughout Sections~\ref{sec:gluing-probability-laws}--\ref{sec:structural-admission-obstructions}, use nonempty finite indexed covers and nonempty finite coordinate alphabets, unless expressly stated otherwise. An empty chart in the unmarked word-law contract has the unique law on its one-letter alphabet and imposes no further probability constraint; marked chart data are not discarded by this convention. Fix the observed alphabets, restrictions and initial laws. A single probability
	row, a family of rows for every starting state, and a complete process law
	give three different gluing problems. A join tree controls each under its
	overlap hypotheses. Away from join trees the independent-block hypothesis
	gives an additional classwise theorem. For the canonical history law,
	conditioning on an entire separator history can change the distribution of
	an earlier prefix.

	\subsection{Weighted structural composition and exact rigidity}
	\label{dyn47:subsec:weighted-compose}

	Let a nonempty finite site family $\mathcal H$ on $V=\bigcup_I I$ have
	join tree $\mathcal U$. Sites are occurrences; their alphabets agree on
	common coordinates. A probability law $p_I$ on $A^I$ may be supported
	on a specified admissible state set. Write $Z_I=\supp p_I$. The family
	is separator consistent if for every edge $I-J$ its marginals on
	$B=I\cap J$ coincide, denoted $p_B$ with support $Z_B$. A global law
	must have all these specified marginals, including their zeros. Its
	possible support is
	$Z=\Join_I Z_I=\{x:x|_I\in Z_I\ \forall I\}$.

	\begin{theorem}[Weighted join: existence, entropy characterization, and dimension]
		\label{dyn47:thm:weighted-join}
		\leavevmode
		\begin{enumerate}[label=(\roman*),leftmargin=*,nosep]
			\item Separator consistency is necessary and sufficient for a common global probability law~\cite{Vorobev1962}. The canonical join $p^*$ is zero outside $Z$ and on
			$Z$ is
			\begin{equation}\label{dyn47:eq:weighted-join}
				p^*(x)=\frac{\prod_{I\in\mathcal H}p_I(x_I)}
				{\prod_{I-J\in\mathcal U}p_{I\cap J}(x_{I\cap J})}.
			\end{equation}
			It is strictly positive on $Z$, has the prescribed marginals, and is
			independent of the choice of root (indeed, of the join tree).
			
			\item For any other gluing $p$ let $H$
			denote Shannon entropy using natural logarithms, with $0\log0=0$.
			Then
			\begin{equation}\label{dyn47:eq:entropy-defect}
				D_{\rm KL}(p\Vert p^*)=\sum_I H(p_I)-\sum_{I-J}H(p_{I\cap J})-H(p).
			\end{equation}
			In particular $p^*$ is the unique entropy-maximizing gluing and the
			unique gluing satisfying conditional independence of each newly added
			private block from the preceding variables given its separator in a
			rooted parent-before-child ordering.
			
			\item The affine dimension of the
			entire gluing polytope is exactly
			\begin{equation}\label{dyn47:eq:glue-dimension}
				d=|Z|-\sum_I|Z_I|+\sum_{I-J\in\mathcal U}|Z_{I\cap J}|.
			\end{equation}
			The marginal data therefore determine the global law without any
			conditional-independence assumption if and only if $d=0$.
		\end{enumerate}
	\end{theorem}

	\begin{proof}
		\leavevmode
		\begin{enumerate}[label={(iii)},leftmargin=*,nosep]
			\item[(i)] Necessity follows by taking marginals. Root the tree and add sites in
			parent-before-child order. The intersection of a new site $I$ with the
			preceding union is its separator $B$ with its parent by the running
			intersection property. On every separator value of positive mass,
			attach the new private variables using conditional probability
			$p_I(x_I)/p_B(x_B)$. All positive partial configurations then extend,
			all prescribed marginals are preserved, and multiplying the factors
			gives~\eqref{dyn47:eq:weighted-join}. The expression is symmetric in the
			unrooted tree and has no zero denominators on $Z$. It also shows that
			every compatible partial configuration on a connected subtree extends
			to $Z$.
			
			\item[(ii)] Every other gluing is supported on $Z$. Integrating $\log p^*$
			against it yields~\eqref{dyn47:eq:entropy-defect}. Nonnegativity and the
			zero case of relative entropy follow from $\log t\le t-1$ (or the
			finite log-sum inequality). The conditional-independence characterization
			is the conditional chain rule for the same construction; it only
			requires conditions on positive-mass separator fibers. The unique entropy maximizer is specified by the marginal constraints alone, so it is independent of the choice of join tree.
		
			\item[(iii)] For the dimension, let $M:\mathbb R^Z\to\bigoplus_I\mathbb R^{Z_I}$
			be the linear marginal map. Its rank is the dimension of the sum of
			spaces of functions of individual site values on $Z$. Each site space
			has dimension $|Z_I|$ because the construction projects $Z$ onto $Z_I$.
			In the rooted ordering the intersection of the next site function space
			with the preceding sum consists exactly of functions of the separator.
			Indeed, a function in that intersection is both a function of the new
			site and a function of the preceding variables. For each separator
			value, all compatible new and preceding private values can be combined
			and extended to $Z$. Hence that function must be constant across the
			new private values at fixed separator. Conversely every separator
			function belongs to the parent space as well as the new site space.
			This proves $\operatorname{rank} M=\sum_I|Z_I|-\sum_{I-J}|Z_{I\cap J}|$.
			The equality constraints define $p^*+\ker M$; since $p^*$ is strictly
			positive on $Z$, small perturbations in every kernel direction remain
			probability laws. Thus nonnegativity does not lower the affine
			dimension. The polytope is a singleton exactly when that dimension
			vanishes. \qedhere
		\end{enumerate}
	\end{proof}

	\begin{corollary}[Two-site rigidity]
		\label{dyn47:cor:two-site-rigid}
		For two sites $I,J$ with separator $B$, let $m_s,n_s$ be the sizes of
		the two positive support fibers above $s\in Z_B$. Then
		\[
		d=\sum_{s\in Z_B}(m_s-1)(n_s-1).
		\]
		The gluing is rigid exactly when at each positive separator value at
		least one side has just one supported private configuration.
	\end{corollary}

	\begin{proof}
		The joined support has size $\sum_s m_sn_s$ and the two site support
		sizes are $\sum_s m_s,\sum_s n_s$. Apply
		\eqref{dyn47:eq:glue-dimension} and nonnegativity of each product.
	\end{proof}
	The product-over-separators probability construction and entropy identity are classical facts about decomposable models~\cite{Lauritzen1996}. The dimension formula
	makes the \emph{remaining weighted coupling freedom} explicit in native
	site support sizes. All of these quantities transport under marked
	site and alphabet isomorphisms. The formula concerns an affine fiber,
	not a quotient by possible relabeling automorphisms.

	\begin{theorem}[Exact Markov globalization: existence, row dimension, and rigidity]
		\label{dyn47:thm:markov-join}
		On the same join tree prescribe nonempty local state sets $X_I\subseteq
		A^I$ and local stochastic kernels $P_I$ on $X_I$. Suppose their natural
		join $X=\Join_I X_I$ projects onto each $X_I$.
		\begin{enumerate}[label=(\roman*),leftmargin=*,nosep]
			\item \emph{Existence.} An exact global Markov
			kernel on $X$, having all $P_I$ as autonomous factors for every
			initialization, exists if and only if for each $x\in X$ and each tree
			edge $I-J$ the two local next-state laws $P_I(x_I,\cdot)$ and
			$P_J(x_J,\cdot)$ have equal separator marginals.
			
			\item \emph{Canonical rows.} In that case let
			$P^*(x,\cdot)$ be their canonical join. The row gluing dimensions $d_x$
			are given by~\eqref{dyn47:eq:glue-dimension} using the positive support
			of each next-state row and the corresponding separator supports.
			
			\item \emph{Dimension and rigidity.} The set of exact global kernels is the product of these row polytopes
			and has dimension $\sum_{x\in X}d_x$. In particular it is rigid if and
			only if every row is rigid.
			
			\item \emph{Entropy maximization.} The canonical kernel uniquely maximizes
			$\sum_x H(P(x,\cdot))$ over exact global kernels.
		\end{enumerate}
	\end{theorem}

	\begin{proof}
		\leavevmode
		\begin{enumerate}[label={(iii)},leftmargin=*,nosep]
			\item[(i)] The strong factor criterion of Proposition~\ref{dyn46:prop-weighted-factor}
			says exactly that every global row has the specified local row
			marginals. Necessity follows by taking separator marginals; sufficiency
			is Theorem~\ref{dyn47:thm:weighted-join} row by row. Their joined supports
			lie in $X$ by definition; the kernel rows are indexed by all of $X$, including states not reached from a specified initialization.
			
			\item[(ii)] The canonical row and $d_x$ are the canonical law and dimension supplied by Theorem~\ref{dyn47:thm:weighted-join}.
			
			\item[(iii)] Row constraints are independent, so the feasible kernel set is the product of row polytopes, with dimension $\sum_x d_x$. It is a singleton exactly when every row polytope is a singleton.
			
			\item[(iv)] The entropy statement follows from the unique row entropy maximizers in the same theorem and the independence of the row choices.
		\end{enumerate}
		Iterating the factor identity yields exact finite
		and all-time discrete path projections for each initial law. When
		initial local laws are prescribed, they must themselves be separator
		consistent; the same join theorem constructs a compatible initial law.
	\end{proof}
	This is a stochastic structural composition operation with exact
	projection, not independent-clock summation of generators. It is
	specific to the joined state space. Imposing further global admission
	constraints may destroy existence even if all row separator tests pass.

	\subsection{Weighted trajectory recovery: path laws and temporal exactness}
	\label{dyn47:subsec:path-recovery}

	At a fixed horizon a word at each coordinate can be treated as a single finite-valued observation. Passing from such finite tables to a process additionally requires consistency under truncation.

	\begin{theorem}[All-horizon probability globalization and uniqueness]
		\label{dyn47:thm:path-glue}
		Suppose on a finite join tree there are prescribed local probability
		laws $p_{I,n}$ of histories through times $0,\ldots,n$, for all
		$n\ge0$. They are supported on histories of the local state spaces.
		Assume prefix consistency at each site.
		\begin{enumerate}[label=(\roman*),leftmargin=*,nosep]
			\item \emph{Existence.} There is a global path law
			on $X^{\mathbb N_0}$, $X=\Join_I X_I$, whose local path marginals are
			all the prescribed laws if and only if for every $n$ the local
			history laws have equal edge-separator history marginals.
			
			\item \emph{Finite-horizon joins.} For a fixed
			horizon, the unique canonical history join and its rigidity dimension
			are given by Theorem~\ref{dyn47:thm:weighted-join} with each coordinate
			alphabet replaced by its $(n+1)$-history alphabet.
			
			\item \emph{Uniqueness.} If every horizon
			is rigid, the global path law is unique.
		\end{enumerate}
		The global law need not be
		Markov and the canonical joins at separate horizons need not be
		prefix consistent.
	\end{theorem}

	\begin{proof}
		\leavevmode
		\begin{enumerate}[label={(iii)},leftmargin=*,nosep]
			\item[(i)] Necessity is marginalization. At each horizon apply the join theorem
			to the coordinate histories; its support consists of histories in
			$X$. For each $N$ form such a gluing at horizon $N$ and take all its
			shorter prefix marginals. For any fixed horizon these probability
			vectors lie in a compact finite simplex. Successive subsequences and
			a diagonal extraction give limits at every fixed horizon, with the
			prescribed marginals. Prefix restriction is linear and continuous,
			so the limiting sequence is prefix consistent. It defines a path law:
			use its initial distribution and the successive conditional laws given
			positive-probability prefixes (arbitrary on null prefixes), or the
			Kolmogorov extension theorem for finite alphabets~\cite{Kallenberg}.
			
			\item[(ii)] At a fixed horizon the join theorem gives the canonical law and its dimension, with complete coordinate histories as letters.
			
			\item[(iii)] If every gluing polytope is a singleton the prefix marginals of any global law are
			forced, hence the law is unique. The following example shows why the
			canonical entropy maximizer cannot simply be iterated over horizons. \qedhere
		\end{enumerate}
	\end{proof}

	\begin{example}[Canonical weighted history joins need not commute with truncation]
		\label{dyn47:ex:prefix-entropy}
		Take the two sites $\{a,b\}$ and $\{b,c\}$, all binary, and initialize
		locally with $b_0=0$, $a_0$ and $c_0$ each uniform. The first site keeps
		$a$ fixed and sets $b'=a$; the second keeps $c$ fixed and sets $b'=c$.
		The separator history law is the same at every horizon: at positive
		horizons it is $(0,U,U,\ldots,U)$ with $U$ uniform. At horizon zero
		the canonical join makes $a_0,c_0$ independent. At horizon one every
		gluing requires $a_0=c_0$ almost surely, because the two values of
		$b_1$ must agree. Hence even the canonical join at horizon one has a
		correlated horizon-zero prefix, different from the horizon-zero
		canonical join. The global path law with $a_0=c_0$ exists. A kernel on
		the full joined state space with both exact local factors for every
		initialization does not exist: from a state with $a\ne c$, the next
		separator rows are incompatible. This separates law-specific history
		globalization from all-initial-state Markov globalization and from
		prefix compatibility of entropy-selected representatives.
	\end{example}

	\begin{example}[Weighted obstruction invisible to relational support]
		\label{dyn47:ex:cycle-weight}
		On the three binary pair sites $ab,bc,ca$, prescribe each pair law with
		probability $9/20$ on each unequal pair and $1/20$ on each equal pair.
		All single-coordinate marginals are uniform and all pair supports
		are full. Thus relational state globalization is unproblematic. If a
		joint probability law existed, the expected number of unequal pairs
		would be $27/10$. But on any binary triple the number is either zero
		or two. This contradiction shows that even strict positivity and
		consistent overlaps do not suffice on a general cyclic cover.
	\end{example}

	\subsection{A native weighted rank and a sharp globalization theorem}
	\label{dyn48:subsec:rank}

	Let $I$ be a nonempty finite native coordinate set with nonempty finite alphabets $A_i$. An initialized observed weighted object $\mathbf D$ specifies prefix-consistent laws $p_n$ on $\prod_{i\in I} A_i^{n+1}$, $n\ge0$ (e.g. the laws induced by an initialized finite observed Markov system). Here each spatial coordinate carries its \emph{whole history} at the given horizon. Its probability marginal on $J\subseteq I$ is denoted $p_{J,n}$. Empty-site laws can be ignored in the block-bounded universality statements, since their unique marginal has no dependence block. We use the laws themselves as the comparison data, not a chosen state presentation; no all-initializations quantifier is implicit. A block partition $\pi$ of $I$ is independent at horizon $n$ if $p_n=\bigotimes_{B\in\pi}p_{B,n}$, and independent for the process if this holds for all $n$. The empty-coordinate law has one value of mass one.

	\begin{definition}[Intrinsic weighted temporal dependence and recovery rank]
		\label{dyn48:def:wdrk}
		Let $\mathbf D$ be an initialized probability process on
		a finite nonempty coordinate set $I$, with finite nonempty
		alphabets $A_i$ and prefix-consistent laws $(p_n)_{n\geq0}$
		on $\prod_{i\in I}A_i^{n+1}$.
		\begin{enumerate}[label=(\roman*),leftmargin=*,nosep]
			\item \textbf{Mixed cumulant.}
			For $\varnothing\ne T\subseteq I$ and functions
			$g_i:A_i^{n+1}\to\mathbb R$, $i\in T$, define
			the mixed cumulant by
			\[
			K_{p_n}(T;g):=
			[t_T]\log\mathbb E_{p_n}
			\prod_{i\in T}(1+t_i g_i),
			\qquad t_T:=\prod_{i\in T}t_i,
			\]
			where $[t_T]$ extracts the coefficient of $t_T$
			and the logarithm is expanded in the commutative
			algebra with $t_i^2=0$.
			For singletons, this gives ordinary expectations
			\cite{Speed1983}.
			
			\item \textbf{Temporal dependence hypergraph.}
			The dependence hypergraph at horizon $n$ is
			$H_n=(I,\mathcal E_n)$, where
			\[
			\mathcal E_n:=
			\{T\subseteq I:|T|\geq2,\
			K_{p_n}(T;g)\ne0\text{ for some }g\}.
			\]
			The full temporal dependence hypergraph is
			\[
			H_\infty:=\left(I,\bigcup_{n\geq0}\mathcal E_n\right).
			\]
			Let $\pi_n$ and $\pi_\infty$ denote the connected-component
			partitions of $H_n$ and $H_\infty$, respectively,
			including isolated vertices.
			
			\item \textbf{History-block size and recovery rank.}
			The history-block size $b(\mathbf D)$ and the
			weighted dynamical recovery rank
			$\operatorname{wdrk}(\mathbf D)$ are defined by
			\begin{equation}\label{dyn48:eq:wdrk}
				b(\mathbf D):=\max_{B\in\pi_\infty}|B|,
				\qquad
				\operatorname{wdrk}(\mathbf D)
				:=\max\{2,2b(\mathbf D)-1\}.
			\end{equation}
			
			\item \textbf{Class recovery rank and bounded-block class.}
			For a nonempty class $\mathscr C$ of such processes,
			define its recovery rank by
			\[
			\operatorname{wdrk}(\mathscr C)
			:=\sup_{\mathbf D\in\mathscr C}
			\operatorname{wdrk}(\mathbf D).
			\]
			For $k\geq1$, the bounded-block class $\mathscr B_k$
			consists of all such processes satisfying
			$b(\mathbf D)\leq k$, with binary alphabets permitted.
		\end{enumerate}
		The rank counts spatial history-block complexity,
		rather than individual space--time cells.
	\end{definition}

	\begin{lemma}[Intrinsic finest partition and recovery from block marginals]
		\label{dyn48:lem:block}
		\leavevmode
		\begin{enumerate}[label=(\roman*),leftmargin=*,nosep]
			\item \emph{Finest independent partitions.} The partition $\pi_n$ is the unique finest independent partition at horizon $n$. The partition $\pi_\infty$ is the unique finest partition independent for the whole process. The partitions $\pi_n$ coarsen with increasing $n$, and eventually stabilize at $\pi_\infty$.
			
			\item \emph{Invariance and restriction.} The rank is invariant under coordinate and coordinatewise alphabet bijections of the observed laws, and does not increase on restriction to a coordinate subset.
			
			\item \emph{Recovery from marginals.} Among laws whose independent blocks have size at most $k$, marginals on all subsets of size at most $k$ determine the full law uniquely.
		\end{enumerate}
	\end{lemma}

	\begin{proof}
		\leavevmode
		\begin{enumerate}[label={(iii)},leftmargin=*,nosep]
			\item[(i)] For functions on a fixed finite coordinate set, put $M=\mathbb E\prod_i(1+t_i g_i)$ in the nilpotent algebra. Independence across a partition gives a product of the block moment polynomials, hence a sum of their logarithms; no mixed coefficient crossing blocks survives. Conversely, if all such mixed coefficients vanish for every choice of functions, $\log M$ is the sum of its block restrictions. Exponentiating recovers the product of block moment polynomials. Taking coordinate functions to be atom indicators shows that every joint atom probability factors across those blocks. Thus a partition is independent exactly when each dependence edge lies inside one block. This proves the finest-partition assertions, including the all-horizon assertion. The cumulants at shorter horizons are cumulants of prefix functions at longer horizons, so every earlier edge persists. There are finitely many coordinate subsets, hence their increasing union stabilizes.
			
			\item[(ii)] Restrictions preserve precisely the cumulant tests on retained coordinates; their components are contained in old components. Relabeling transports the test functions.
			
			\item[(iii)] Finally, in a law with blocks of size at most $k$ all cumulants of support larger than $k$ vanish. Its smaller cumulants are determined by its small marginals. Equal such marginals for two laws give equal coefficients of $\log M$ for every family of functions, hence equal atom probabilities. This last uniqueness is within the block-bounded class, not among arbitrary probability measures. \qedhere
		\end{enumerate}
	\end{proof}

	For example, an initialized process with independent coordinate histories has rank two. This includes independent time-dependent or non-Markov histories, not just independent one-time states. Conversely let two binary coordinates evolve by $(a,b)\mapsto(U,a)$ with a fresh fair independent bit $U$ and initially independent fair bits. At every time the two coordinates are independent, but the next $b$ equals the preceding $a$. Their joint histories form a block of size two and the rank is three. Even constant histories with correlated strictly positive initial distributions can have a larger rank than independent constant histories of exactly the same relational support. Thus the rank uses weights and complete histories, not just a static support or a snapshot independence test.

	\begin{proposition}[Finite-state recognition of the weighted rank]
		\label{dyn48:prop:rank-recognition}
		For an initialized $d$-state finite Markov presentation with a deterministic finite observation tuple, $d\ge1$, independence of a proposed bipartition of all coordinate histories can be decided by checking equality of joint and product-marginal word probabilities through word length $d+d^2-1$. Consequently $\pi_\infty=\pi_{d+d^2-2}$, and the rank is effectively computable for rational transition and initialization data by finitely many exact linear tests. The same bound holds for a $d$-state edge-emitting presentation, with its letter matrices in place of the matrices $D_aP$.
	\end{proposition}

	\begin{proof}
		With row initialization $\lambda$, transition matrix $P$, and diagonal observation indicators $D_a$, a word $a_1\cdots a_l$ has probability $\lambda D_{a_1}P\cdots D_{a_l}P\mathbf 1$. This uses $P\mathbf 1=\mathbf 1$. For edge emissions the same expression uses the given letter matrices instead of $D_aP$. Marginalization in either convention sums the full letter matrices over letters with the prescribed restriction. The product of two marginal word probabilities is represented by the tensor product of two $d$-dimensional representations. Their difference from the full word law has a linear representation of dimension $D=d+d^2$, formed by a block direct sum with opposite signs on one initialization. For any such representation the reachable row span is generated by words of length at most $D-1$: the spans obtained by adding one letter at a time increase until they stabilize, and equality at one stage makes the span invariant under every letter matrix. A nonzero initial row gives dimension at least one at length zero; the zero case is immediate. Thus vanishing against the terminal vector through length $D-1$ implies vanishing for all words. At the horizon in the statement every independent bipartition is therefore independent forever. The finest partition is the common refinement of independent bipartitions, by Lemma~\ref{dyn48:lem:block} (a component cannot be separated by any such split, while distinct components can). This proves the claim.
	\end{proof}

	For a nonempty finite cover $\mathcal H$ of a finite nonempty $V$, define $\operatorname{WUG}(\mathcal H,\mathscr C)$ to mean that \emph{every} $\mathscr C$-atlas of initialized local process laws on the sites, with agreement of all overlap history marginals at every horizon, has a common global probability process on the coordinate histories. The global law need not belong to $\mathscr C$, be Markov, or meet an additional global constraint not declared by the local laws. It is supported on the local joins at each horizon. Use the inherited $\nu$ of Section~\ref{sec:covers-assertions}: $\nu>m$ for $m\ge3$ means chordality and guarding of all primal cliques of size at most $m$. In particular $\nu=\infty$ is precisely the finite join-tree condition, by the classical chordal/conformal characterization of hypergraph acyclicity \cite{wn09:BeeriFaginMaierYannakakis1983}. This convention permits repeated and contained sites and disconnected components (empty separators).

	\begin{theorem}[Sharp weighted block-recovery globalization]
		\label{dyn48:thm:sharp-wug}
		\leavevmode
		\begin{enumerate}[label=(\roman*),leftmargin=*,nosep]
			\item \emph{Sharp block cutoff.} For every finite cover $\mathcal H$ and integer $k\ge1$,
			\begin{equation}\label{dyn48:eq:sharp-wug}
				\boxed{\operatorname{WUG}(\mathcal H,\mathscr B_k)
					\ \Longleftrightarrow\ \max\{2,2k-1\}<\nu(\mathcal H).}
			\end{equation}
			Equivalently, $k=1$ always globalizes, and for $k\ge2$ the criterion is chordality and guarding of all cliques of size at most $2k-1$. The maximum in the displayed rank convention retains the normalization in Definition~\ref{dyn48:def:wdrk}; it does not impose a restriction for $k=1$, since $\nu\ge3$. Here universality allows binary alphabets at every vertex; for arbitrary fixed alphabets the sufficient direction still holds, and the converse holds if every alphabet admits at least two symbols.
			
			\item \emph{Classwise sufficiency.} More generally, every nonempty subclass $\mathscr C$ satisfies
			\begin{equation}\label{dyn48:eq:class-wug}
				\operatorname{wdrk}(\mathscr C)<\nu(\mathcal H)
				\ \Longrightarrow\ \operatorname{WUG}(\mathcal H,\mathscr C).
			\end{equation}
			
			\item \emph{Unrestricted acyclic gluing.} If $\nu=\infty$, all compatible process atlases globalize, even when the classwise supremum is infinite.
		\end{enumerate}
	\end{theorem}

	\begin{proof}
		\leavevmode
		\begin{enumerate}[label={(iii)},leftmargin=*,nosep]
			\item[(i)] \emph{Sufficiency.} For $k=1$ every local history law is a product of coordinate marginals. Their consistent products on $V$ give a process, on every cover. For $k\ge2$ assume $2k-1<\nu(\mathcal H)$. Fix a finite horizon, treating histories as single alphabet symbols. Every primal clique of size at most $2k-1$ is guarded. For a maximal primal clique $C$, define a hypergraph on $C$ by including $T\subseteq C$, $2\le |T|\le k$, whenever its marginal, taken from any guarding chart, has a nonzero mixed cumulant. This definition is independent of the guard by overlap agreement. If one component has more than $k$ vertices, choose a connected sequence of its edges, beginning with one edge and each subsequent edge meeting the accumulated union, until the union first exceeds $k$. Such a sequence exists by connectivity. The first exceeding union $U$ has $k<|U|\le2k-1$: the preceding union has at most $k$ vertices, and the new edge adds at most $k-1$. A chart guards $U$. Every selected edge is a dependence edge of the chart law, by compatibility and marginal locality. They are connected there, contradicting its bound $b\le k$.
			
			Therefore each component $B$ of the clique hypergraph has $|B|\le k$. Take its prescribed marginal $p_B$ from any guard and form the clique law $q_C=\bigotimes_B p_B$. For a chart $I$, its law restricted to $I\cap C$ has all dependence edges of size at most $k$, since it has block size at most $k$. Every such edge lies within a constructed component of $C$. It consequently factors across the intersections of these components with $I$ by Lemma~\ref{dyn48:lem:block}. Each factor equals the corresponding restriction of $p_B$, by overlap agreement of its guard with $I$. Thus $q_C$ has exactly the prescribed restriction on $I\cap C$. In particular every original site, being a clique, is extended by any maximal clique containing it. All $q_C$ have blocks of size at most $k$. On an overlap of two maximal cliques their marginals on every subset of size at most $k$ agree, because each such subset is guarded and both clique laws extend its prescribed marginal. The uniqueness clause of Lemma~\ref{dyn48:lem:block} gives full agreement on the overlap. The same uniqueness gives uniqueness of each block-bounded clique completion among laws extending the prescribed chart restrictions, since all subsets of size at most $k$ are guarded. It also shows that its successive history laws are prefix consistent: a prefix of a block-bounded law is block-bounded, and all prescribed small marginals restrict correctly.
			
			The primal graph is chordal. Its maximal cliques admit a join tree (with empty separators between components); apply Theorem~\ref{dyn47:thm:weighted-join} to the completed clique laws. There is a gluing at each horizon. The compact prefix argument in Theorem~\ref{dyn47:thm:path-glue}, or the canonical process construction below, gives a single global path law.
			\emph{Necessity: induced cycle.} If the primal graph has an induced cycle on $r\ge4$ vertices, every site meets that cycle in at most an adjacent pair. Prescribe binary fair singletons, and on every cycle edge the uniform law satisfying a prescribed parity of its two bits, choosing parities whose sum modulo two is one. On each site take the indicated marginal on its intersection with the cycle, and independent fair bits off the cycle. These local laws agree on all overlaps. A global sample would have to satisfy every edge parity, which is impossible since summing all edge parities gives zero. The local laws have blocks of size at most two. Make every coordinate constant in time; this is a compatible process atlas and no global process exists already at horizon zero.
			
			\emph{Necessity: missing clique.} Otherwise failure of the threshold for $k\ge2$ gives an unguarded clique $M$ of size at most $2k-1$. Pass to an inclusion-minimal unguarded subclique, still denoted $M$, of size $r\ge3$. All its proper subsets are guarded. Choose proper subsets $U,W$ with $U\cup W=M$, $|U\cap W|=1$ and $2\le |U|,|W|\le k$. Such sizes exist precisely for $3\le r\le2k-1$. On every proper subset $J$ of $M$ use uniform bits subject to even parity on $U$ if $U\subseteq J$, or on $W$ if $W\subseteq J$, and use independent uniform bits otherwise. A proper subset cannot contain both. Marginalizing a uniform parity law to any strict subset of its parity support gives independent fair bits there; hence these prescriptions are consistent under restriction. Each site meets $M$ in a proper subset; extend the prescribed marginal there by independent fair bits off $M$. All sites are compatible and have blocks of size at most $k$. Both $U$ and $W$ have prescribed even parity, so a global law must have even parity on $U\triangle W$. But $U\triangle W$ is a guarded proper subset containing neither $U$ nor $W$, and has the uniform law. This contradicts its odd-parity probability $1/2$. Constant-in-time processes again transport this obstruction to every horizon. Binary values may be embedded into any larger alphabets. This proves necessity and the equivalence.
			
			\item[(ii)] For a class with finite rank supremum strictly below $\nu$, take the corresponding uniform bound on $b$.
			
			\item[(iii)] The separate assertion for an unrestricted class on a join tree follows directly from the all-horizon join theorem. \qedhere
		\end{enumerate}
	\end{proof}

	\begin{corollary}[Sharpness, full-support obstructions and unrestricted class]
		\label{dyn48:cor:rank-sharpness}
		\leavevmode
		\begin{enumerate}[label=(\roman*),leftmargin=*,nosep]
			\item For $k\ge2$ the least integer $m\ge3$ for which $\nu(\mathcal H)>m$ guarantees $\operatorname{WUG}(\mathcal H,\mathscr B_k)$ on all finite covers is $2k-1$.
			
			\item Every counterexample above may be made strictly positive on each local instantaneous binary product without changing its block bound.
			
			\item For the class of \emph{all} initialized finite-alphabet processes, universal weighted globalization holds exactly on join-tree covers.
		\end{enumerate}
	\end{corollary}

	\begin{proof}
		\leavevmode
		\begin{enumerate}[label={(iii)},leftmargin=*,nosep]
			\item[(i)] The simplex boundary of size $2k-1$ has obstruction rank $2k-1$, and the missing-clique witness applies. This excludes every smaller eligible guard threshold; for $k=2$ the permitted minimum is already three.
			\item[(ii)] For the cycle example replace each hard edge parity by a fair conditional parity mixture with probability $\eta>0$ of violation, with $r\eta<1$. The sum of violation indicators is at least one on every global assignment, contradicting its prescribed expectation. For the clique example give each of the two parity constraints odd-parity probability $0<\eta<1/4$. The odd-parity indicator on $U\triangle W$ is bounded above by the sum of the two odd-parity indicators, yet its prescribed expectation is $1/2>2\eta$. Proper restrictions of these parity mixtures are still uniform and the block bound persists; all initial local atoms now have positive probability. Thus the constant-history supports are exactly those of independent constant-history laws.
			\item[(iii)] For an arbitrary non-join-tree cover, the cycle or minimal unguarded clique construction uses a finite $k$. Join trees globalize the unrestricted class by Theorem~\ref{dyn47:thm:path-glue}. \qedhere
		\end{enumerate}
	\end{proof}

	\begin{remark}[Rank interpretation and limits]
		\label{dyn48:rem:rank-limits}
		The guard calibration in \eqref{dyn48:eq:wdrk} is sharp for the \emph{full} block-bounded class. The rank is determined by connected dependence components of the initialized observed process, rather than by the largest individual cumulant or a support-only witness width. The same local bound can force a larger global component. On the chain of sites $ab,bc$ with perfect equalities and fair bits, the local block bound is two, yet every globalization has a three-coordinate dependence block. Even a two-variable marginal-order rank of two would incorrectly predict globalization on the triangle with frustrated parity, whose $\nu$ is three. The weighted proof uses clique completion and parity transport in place of the static native-recovery signature and obstruction portability of Paper I. A class allowing only special subclasses of these laws may have a better threshold. Corollary~\ref{dyn53:cor:finite-wug} treats finite autonomous realization, and Theorem~\ref{dyn53:thm:block-preserving} characterizes preservation of the block bound itself.
	\end{remark}

	\subsection{Canonical projective globalization of weighted histories}
	\label{dyn48:subsec:projective}

	The entropy-maximizing finite-horizon joins of Section~\ref{dyn47:subsec:path-recovery} need not be prefix consistent. We resolve the canonical-selection problem for complete compatible path laws on a finite join tree. This selection is allowed to use the entire supplied temporal law, not only the observations up to the prefix being selected.

	\begin{theorem}[Canonical full-history gluing: convergence, uniqueness, and naturality]
		\label{dyn48:thm:projective}
		Assume the finite-alphabet compatible all-horizon data of Theorem~\ref{dyn47:thm:path-glue}, with join tree $\mathcal U$. Let $q_N$ be their canonical finite-history join through horizon $N$.
		\begin{enumerate}[label=(\roman*),leftmargin=*,nosep]
			\item \emph{Limit and existence.} For every fixed $n$, the prefix marginals $q_N|_{0:n}$ converge as $N\to\infty$ in total variation to a law $q^\infty_n$. The resulting family is projectively consistent, has all prescribed local marginals, and defines a canonical global process $q^\infty$. This construction is independent of the choice of join tree.
			
			\item \emph{Characterization.} For every tree edge, let the two components of the cut have coordinate unions $V_1,V_2$ and shared separator $B=V_1\cap V_2$. The law $q^\infty$ is the unique gluing for which the two full histories on $V_1,V_2$ are conditionally independent given the \emph{entire} history on $B$, at every edge of a fixed join tree.
			
			\item \emph{Naturality.} The canonical law is natural under native coordinate and alphabet bijections; it need not be Markov or causal with respect to the shorter separator prefixes.
		\end{enumerate}
	\end{theorem}

	\begin{proof}
		Write $\mu_I$ for the full history law on $I$, which exists from the prefix data by the finite-alphabet extension construction in Theorem~\ref{dyn47:thm:path-glue}. These spaces are countable products of finite discrete spaces, hence standard Borel; conditional laws given full separator histories exist. For completeness, the conditional expectations needed below can be obtained by taking the limit of conditional expectations on finite separator prefixes. For disintegration on standard Borel spaces, see \cite{Kallenberg}. For a bounded random variable, L\'evy's upward martingale convergence theorem (see \cite{Kallenberg}) gives convergence in $L^2$ to its conditional expectation on the full history sigma algebra. Equivalently, the conditional expectations are orthogonal projections onto increasing closed subspaces whose union is dense in the space measurable with respect to the full history. Density follows from cylinder simple functions and the fact that the cylinders generate the sigma algebra. For conditional atom probabilities on the countable cylinder algebra take a common almost-sure version of the finite probability identities; these consistent finite distributions give a conditional law by the same finite-alphabet extension theorem.
		
		Root the tree. Sample the root full history with its given law; for each child, condition its history on the full history of its separator with the parent, using its prescribed conditional law, and attach the remaining coordinates. All such conditional laws are used only over separator histories distributed according to their common marginal. The running-intersection property ensures that a new site shares exactly its parent separator with the previously attached union. Thus the construction is defined almost surely and has the correct site laws.
		\begin{enumerate}[label={(iii)},leftmargin=*,nosep]
			\item[(i)] For the prefix-limit assertion, we compare finite canonical joins with this full-history construction. Fix bounded cylinder functions $g_I$ of the site histories, all depending on times at most $n$. Their product under the canonical finite join $q_N$, $N\ge n$, has expectation computed by the rooted tree recursion. At a nonroot site $I$ with parent separator $B$, the message is
			\[
			m_{I,N}=\mathbb E_{\mu_I}\!\left[g_I\prod_{J\text{ child of }I}m_{J,N}
			\ \middle|\ \mathcal F_{B,N}\right],
			\]
			where $\mathcal F_{B,N}$ is the sigma algebra of the separator prefix, and messages on child separators are pulled back to the site. The root integrates the same product without a parent condition. This recursion is just the finite conditional-product construction in \eqref{dyn47:eq:weighted-join}. The messages for the full-history construction satisfy the identical recursion using full separator sigma algebras. All messages are uniformly bounded by the product of the relevant bounds for the $g_I$. At a leaf the increasing-subspace projection convergence just proved gives $L^1$ convergence to the full-history message. Inductively the product in the expectation converges in $L^1(\mu_I)$: use the finite telescoping product inequality and the common separator marginal to transfer the $L^1$ bounds from each child. Conditional expectation is an $L^1$ contraction, and conditioning the fixed limiting bounded product on growing separator prefixes converges in $L^1$ by the same projection argument. This proves convergence of every message and then of the root integral. Global finite cylinder atom indicators are products of site cylinder indicators (assigning each coordinate to a containing site, or using consistent copies). Hence all atom probabilities of the fixed prefix converge to those of the full-history construction, and its finite alphabet implies total-variation convergence.
			
			Each $q_N$ is the unique finite entropy maximizer among laws with the prescribed site marginals by Theorem~\ref{dyn47:thm:weighted-join}; this property is independent of the choice of join tree. Its prefix limit therefore is independent of that choice. Its prefixes are projectively consistent because they come from a single probability process (also because restriction is continuous in these finite spaces).
			
			\item[(ii)] For a fixed rooted ordering, successive conditional attachments are independent given the already sampled histories except for their stated separators. Changing the root across an adjacent edge does not alter the initial joint law of those two sites: it is their relatively independent join over the common full separator, symmetric in the two factors. The remaining branches can be attached in the same order and use the same conditional laws; attachments of different branches commute by conditional independence. Repeating adjacent root changes proves root independence. Rooting on either side of any edge now shows that the two full subtrees are conditionally independent given the full separator. Conversely a gluing with all these edge independences must attach every leaf to the remainder by precisely its conditional law on the full separator. Removing leaves successively gives uniqueness: the edge independence statements for the remaining tree follow by marginalizing the deleted leaf from the corresponding cut side. This describes a unique full-history gluing for the specified tree.
			
			\item[(iii)] Relabeling commutes with the finite canonical joins and their limits. This proves all assertions. \qedhere
		\end{enumerate}
	\end{proof}

	\begin{corollary}[Canonical globalization at the weighted rank threshold]
		\label{dyn48:cor:rank-canonical}
		For $\mathscr B_k$ with $\max\{2,2k-1\}<\nu(\mathcal H)$ there is a canonical global probability process constructed by unique block-bounded completion on each maximal clique, then the full-separator gluing of Theorem~\ref{dyn48:thm:projective}. For $k=1$ the construction is simply the product process, even on a nonchordal cover. For $k\ge2$ the result is independent of the clique-tree choice. It retains all local history weights; it is not in general a unique globalization without the declared canonical rule.
	\end{corollary}

	\begin{proof}
		In the sufficiency proof of Theorem~\ref{dyn48:thm:sharp-wug}, the block-bounded completions at each horizon are unique, prefix consistent and agree on all clique overlaps. Chordality gives the join tree; apply the preceding theorem. Independent coordinate products give the $k=1$ case.
	\end{proof}

	\begin{example}[No universal selector using only the current prefixes]
		\label{dyn48:ex:no-prefix-selector}
		Use Example~\ref{dyn47:ex:prefix-entropy}: the sites $ab$ and $bc$ have $b_0=0$ and independent-in-their-respective-site fair private bits $a_0,c_0$. In one continuation keep the private bits fixed and let $b'=a$ at the first site, $b'=c$ at the second. In a second continuation keep them fixed and let the first site again set $b'=a$ but the second set $b'=1-c$. Both continuations have identical initial local data and compatible separator histories at every horizon. The first forces $a_0=c_0$ in every global process, the second forces $a_0=1-c_0$. There is no choice of initial global law, based only on the initial local laws, that works in both cases. Thus the full-law selector above does not solve a causal online selection problem for arbitrary compatible continuations; that stronger problem is impossible without additional conditions. Even local deterministic Markov laws with prescribed initializations exhibit this obstruction.
	\end{example}

	\subsection{Continuous-time weighted globalization and finite-system reduction}
	\label{dyn48:subsec:continuous}

	Let $X=\Join_I X_I$ be the finite nonempty joined state space, with each local state projection onto $X_I$. Prescribe conservative finite-state generators $Q_I$ on the local states. Exact autonomous globalization means a generator $Q$ on $X$ satisfying the strong generator factor criterion for each site, for every initialization. For $x\in X$ and $z_I\ne x_I$, the criterion is
	\begin{equation}\label{dyn48:eq:rate-marginals}
		\sum_{y\in X:\ y_I=z_I}Q(x,y)=Q_I(x_I,z_I).
	\end{equation}
	The omitted diagonal condition follows from zero row sums. For a separator, interpret the projected rates from a local row as the sums into separator values distinct from the present value, not as a normalized probability distribution.

	\begin{theorem}[Continuous-time globalization and uniformized generator polytopes]
		\label{dyn48:thm:continuous-join}
		\leavevmode
		\begin{enumerate}[label=(\roman*),leftmargin=*,nosep]
			\item On a finite join tree, the local generators admit an exact global generator if and only if for every $x\in X$ the projected rates into every distinct separator output agree across each tree edge.
			
			\item More generally on an arbitrary cover, choose a constant
			\[
			\Lambda>0,\qquad \Lambda\ge \max_{x\in X}\sum_I[-Q_I(x_I,x_I)].
			\]
			Set $P_I=\mathrm{Id}+Q_I/\Lambda$. The map $Q\mapsto \mathrm{Id}+Q/\Lambda$ is an affine bijection between exact global generators and exact discrete-time global kernels for the $P_I$. Thus the row-polytope description and dimension of Theorem~\ref{dyn47:thm:markov-join} apply on join trees, including the rigidity criterion.
			
			\item On any cover, if the uniformized local next-state laws for each $x$ agree on overlaps and have dependence blocks of size at most $k$ in the sense of Lemma~\ref{dyn48:lem:block}, the threshold $\max\{2,2k-1\}<\nu(\mathcal H)$ also suffices. It applies separately to any prescribed compatible initial laws satisfying the same block bound. The resulting global continuous-time process has the prescribed local generators as factors at every time.
		\end{enumerate}
	\end{theorem}

	\begin{proof}
		\leavevmode
		\begin{enumerate}[label={},leftmargin=0pt,nosep]
			\item \textnormal{(i)--(ii).}\enspace Necessity of separator-rate consistency follows by summing the global off-diagonal rates into distinct separator outputs. Every global jump changes at least one coordinate in a covering site. Nonnegativity and \eqref{dyn48:eq:rate-marginals} therefore give, for every exact generator row,
			\[
			\max_I[-Q_I(x_I,x_I)]\le -Q(x,x)\le\sum_I[-Q_I(x_I,x_I)].
			\]
			The first inequality counts only jumps visible at one site and the second counts each global jump at least once. Thus uniformization at the stated $\Lambda$ gives a stochastic kernel. Its exact factor identity is equivalent to the generator identity by subtracting the identity matrices and multiplying by $\Lambda$. Conversely any such stochastic kernel yields a conservative generator by that inverse affine transformation. For two local uniformized rows, separator-rate consistency on distinct outputs is exactly equality of their separator probability marginals (the probability at the current separator value is determined by normalization). Apply Theorem~\ref{dyn47:thm:markov-join} on join trees. The affine bijection also gives the dimension and rigidity results.
			
			\item \textnormal{(iii).}\enspace For the additional bounded-block assertion apply the fixed-horizon clique-completion proof of Theorem~\ref{dyn48:thm:sharp-wug} to each next-state row separately. Since these rows are supported on their local state sets, all global rows are supported on the joined state space. Choose any such rows and the corresponding compatible initial distribution if required. The factor identities for generators intertwine their powers and hence their matrix exponentials $e^{tQ}$, $e^{tQ_I}$ for all $t\ge0$. These semigroup identities give the exact site process laws by the Markov property, for any compatible initial law.  All statements concern the actual joined state space; extra global constraints require separate checks. \qedhere
		\end{enumerate}
	\end{proof}

	The choice of uniformization rate is part of the canonical row construction, not an intrinsic continuous-time coupling law. For example, on two disjoint binary singleton sites let each jump from zero to one at rate $a$ and $b$, respectively, and absorb at one. At state $(0,0)$ the canonical product of uniformized rows gives a simultaneous jump to $(1,1)$ of rate $ab/\Lambda$, which varies with $\Lambda$. Independent-clock gluing would give rate zero for that simultaneous jump. Both have the exact local factors. The generator-polytope classification itself does not depend on $\Lambda$.

	\begin{theorem}[Continuous-time reduction: error, exactness, and path-law transfer]
		\label{dyn48:thm:continuous-error}
		Let $Q$ be a finite fine generator on $X$, $m:X\twoheadrightarrow Y$ a declared finite reduction, and $\overline Q$ a proposed coarse generator on $Y$. Define the aggregate signed row $q_x(y)=\sum_{z:m(z)=y}Q(x,z)$ (including the diagonal contribution), and let
		\begin{equation}\label{dyn48:eq:rate-defect}
			d=\tfrac12\sup_x\sum_{y\in Y}|q_x(y)-\overline Q(m(x),y)|.
		\end{equation}
		Let $\delta$ be the total-variation distance of the projected fine and coarse initial laws. The two path laws on $Y$ through time $T$ have total-variation distance at most
		\begin{equation}\label{dyn48:eq:continuous-error}
			1-(1-\delta)e^{-Td}\ \le\ \min\{1,\delta+Td\}.
		\end{equation}
		They coincide for all initializations with matched projected laws if and only if $d=0$. For a size-indexed family the sufficient condition $\delta_N+T_N d_N\to0$ transfers any common measurable path statistic's limiting distribution from one construction to the other on the stated time windows. The bound is valid without Markov autonomy of the fine projection; it does not assert that the projected law has a generator on $Y$.
	\end{theorem}

	\begin{proof}
		Uniformize both generators with a common rate $\Lambda>0$ at least both maximum exit rates. The fine and coarse chains at the Poisson clock have transition kernels $P=\mathrm{Id}+Q/\Lambda$ and $\overline P=\mathrm{Id}+\overline Q/\Lambda$. If the current states match under $m$, the total-variation distance of their next projected rows is at most $d/\Lambda$. A maximal coupling of those rows, followed by sampling a fine next state conditional on its projected value, gives conditional probability at least $1-d/\Lambda$ of continued equality. This is a legal coupling even when the projected fine process alone is not Markov. Couple the initial projected states maximally and lift to the fine initial law. Use one common Poisson clock, independent of the embedded coupled chain. Conditional on $l$ clock events the probability that the projected paths match throughout is at least $(1-\delta)(1-d/\Lambda)^l$. Averaging over the Poisson count of mean $\Lambda T$ gives $(1-\delta)e^{-Td}$. The coupling inequality bounds the path total variation, and $1-e^{-u}\le u$ gives the second bound. Exactness for $d=0$ follows; conversely equality of projected laws from every point initialization implies the row identities by differentiating at time zero. Pushforward by any common statistic contracts total variation, proving the limiting assertion. Sharpness of the first bound as a general bound already occurs when the fine binary process jumps from zero to absorbing one at rate $d$ while the coarse process starting at zero stays there ($\delta=0$): the path laws differ by $1-e^{-Td}$.
	\end{proof}

	\section{Finite autonomous globalization of observed word laws}
	\label{sec:Finite-utonomous-globalization}

	Theorem~\ref{dyn53:thm:universal} constructs a finite global hidden Markov model from finite local models on a join tree. By \emph{unmarked finite autonomous globalization} we mean one finite homogeneous stochastic presentation whose projected initialized word laws are exactly the prescribed ones, with no requirement to preserve the hidden presentations or to meet an additional admission mask. Its compatibility test on join trees is equality of separator word probabilities at every length. This differs from a Markov coupling on the original pairs of states in \cite{CraneHolroydRussell2026}, which under the Markov-image hypothesis also permits conditional independence given the whole common trajectory. The construction passes through a common positive separator presentation and then recovers the prescribed initialized local word laws. Theorem~\ref{dyn53:thm:block-preserving} subsequently characterizes when the independent history-block bound can be retained in the global law. Strong lifting of given hidden presentations imposes further constraints; the intermediate common-factor and driver constructions below retain their separate hypotheses.

	\subsection{Predictive word signatures and finite autonomous realization}
	\label{dyn49:subsec:finite}

	For a finite joint observation alphabet $\Sigma$ write $p(w)$ for the
	probability of a finite prefix word $w$, with $p(\epsilon)=1$ and
	$\sum_{a\in\Sigma}p(wa)=p(w)$. This is the same initialized law as the
	history family above, now indexed by successive joint observations.
	For every native subset $J$, the analogous word function is obtained by
	summing over the omitted coordinates at each letter. Its full predictive
	signature consists of the Hankel array $H_J(u,v)=p_J(uv)$, with rows and
	columns indexed by finite words. Its rank means the supremum of the ranks
	of finite subarrays. This gives an intrinsic profile of linear predictive
	resources on the native subsets.

	An edge-emitting autonomous realization with $d$ hidden states consists
	of a nonnegative probability row $\lambda$ and nonnegative $d\times d$
	matrices $M_a$, with $P=\sum_a M_a$ row stochastic, such that
	$p(w)=\lambda M_w\mathbf1$. The letter is emitted during a state
	transition. A deterministic-observation state presentation can be
	converted to this form without increasing its state count, using
	$M_a=D_aP$; conversely an edge-emitting presentation can be converted to
	at most $d|\Sigma|$ deterministic-observation states $(a,j)$ with
	initial law $(\lambda M_a)_j$ and transition $(a,j)\to(b,k)$ of
	probability $(M_b)_{jk}$. These are different state-budget conventions.
	The following results use the edge-emitting budget unless noted.

	\begin{proposition}[Canonical linear realization and finite positive realization criterion]
		\label{dyn49:prop:residual-space}
		Set $r_u(v)=p(uv)$ and $\mathcal R=\operatorname{span}\{r_u\}$ in the
		space of real word functions.
		\begin{enumerate}[label=(\roman*),leftmargin=*,nosep]
			\item \emph{Linear realization.} Then $\dim\mathcal R=\operatorname{rank}H_I$,
			and $L_a r(v)=r(av)$ is well-defined on $\mathcal R$, with
			$L_a r_u=r_{ua}$. If $\dim\mathcal R=d<\infty$, this is a canonical
			$d$-dimensional linear realization, and every linear realization has
			at least this dimension. In particular every $d'$-state autonomous
			probability realization satisfies $\operatorname{rank}H_I\le d'$.
			Finite rank by itself is not a positivity assertion.
			
			\item \emph{Positive realization.} More precisely, a $d'$-state edge-emitting realization exists if and only
			if there are at most $d'$ normalized prefix-consistent nonnegative word
			functions $g_i$, a probability vector $\lambda$, and nonnegative
			coefficients $m_{a,ij}$ such that
			\[
			p=\sum_i\lambda_i g_i,\qquad L_ag_i=\sum_jm_{a,ij}g_j.
			\]
			The $g_i$ may lie outside $\mathcal R$. This formulation gives the exact
			minimum number of generators of such a finite invariant predictive
			family. This is the classical positive realization criterion in word-function form (cf. Heller's invariant-cone theorem \cite{dyn49:Heller,dyn49:Schoenhuth}).
		\end{enumerate}
	\end{proposition}

	\begin{proof}
		\leavevmode
		\begin{enumerate}[label={(iii)},leftmargin=*,nosep]
			\item[(i)] Rows of the Hankel array are the residual functions, so their span has
			the stated dimension (linear independence of finitely many functions is
			witnessed by finitely many evaluation columns). Prepending a letter to
			the argument is a linear operation on functions and takes each row to the
			next residual. Starting from $r_\epsilon$ and applying successive $L_a$
			in temporal order gives the word residual; evaluation at $\epsilon$
			gives its mass. A linear representation of dimension $d'$ factors the
			Hankel array through that space, yielding the lower bound.
			
			\item[(ii)] For a nonnegative realization take $g_i(v)=e_i^\top M_v\mathbf1$; these
			are normalized laws of the hidden-state initializations, and the displayed
			relations hold. Conversely evaluate the relations at $\epsilon$. Since
			$\sum_a g_i(a)=1$ and $g_j(\epsilon)=1$, the matrices
			$M_a=(m_{a,ij})$ have row-stochastic sum. Induction on word length gives
			$g_i(w)=e_i^\top M_w\mathbf1$, and then $p(w)=\lambda M_w\mathbf1$.
			Taking minima proves the exact characterization, including infinite
			minimum if no finite family exists. The cone generated by the $g_i$ is invariant under all $L_a$; the generators are normalized by evaluation at the empty word, not by a choice of coordinates in a minimal linear realization. Marginalization replaces letter
			matrices by sums, so all native marginal Hankel ranks are lower bounds
			for the joint presentation budget. \qedhere
		\end{enumerate}
	\end{proof}

	Finite word presentations can be compared effectively by finite-dimensional linear tests; for probabilistic automata, see~\cite{Tzeng1992}. The next theorem uses this fact for partial native data rather than assuming a full global law in advance.

	\begin{theorem}[Finite autonomous state-budget test and horizon compactness]
		\label{dyn49:thm:budget}
		Let each site $I$ of a finite cover have an initialized local word law
		with a given linear representation of dimension $d_I$, and assume these
		are probability laws. Fix $d\ge1$ and a nonempty allowed joint observation
		alphabet $A\subseteq\prod_{i\in V} A_i$ (the unrestricted product is
		permitted).
		\begin{enumerate}[label=(\roman*),leftmargin=*,nosep]
			\item \emph{Finite test.} There exists an edge-emitting $d$-state global probability
			realization with all prescribed local laws and letters in $A$ if and only
			if the following finite feasibility problem has a solution:
			\begin{align*}
				&\lambda\ge0,\quad\lambda\mathbf1=1,\quad M_a\ge0\ (a\in A),
				\quad \sum_{a\in A}M_a\mathbf1=\mathbf1,\\
				&\lambda\widehat M_{I,w}\mathbf1=p_I(w)
				\qquad (I\in\mathcal H,\ |w|\le d+d_I-1),
			\end{align*}
			where $\widehat M_{I,b}=\sum_{a\in A:a|_I=b}M_a$, and products follow
			the order of letters.
			
			\item \emph{Minimum budget and decidability.} Thus the minimum autonomous edge-state budget is
			the least feasible $d$ (or infinity). For rational input data, each fixed
			budget is decidable by real polynomial feasibility; no finite upper
			bound on the required budget or general decision procedure for infinite
			minimum is asserted.
			
			\item \emph{Horizon compactness.} Moreover a fixed $d$ suffices exactly when every
			finite horizon separately admits some $d$-state autonomous realization
			of the local prefixes, with the same alphabet mask, even if these
			finite-horizon realizations are not coherently chosen.
		\end{enumerate}
	\end{theorem}

	\begin{proof}
		\leavevmode
		\begin{enumerate}[label={(iii)},leftmargin=*,nosep]
			\item[(i)] Necessity is immediate by summing over the omitted letters. For
			sufficiency subtract the candidate marginal series from the prescribed
			site series using their block-direct-sum linear representations. Its
			dimension is at most $D=d+d_I$. The span of the reachable rows from the
			initial row is generated by words of length at most $D-1$: once the
			increasing row spans stop growing they are invariant under all letter
			matrices, while at length zero the nonzero row spans at least one
			dimension. Hence zero evaluation against the terminal vector at lengths
			through $D-1$ forces zero evaluation at all lengths. This is the same
			linear test mechanism as Proposition~\ref{dyn48:prop:rank-recognition}.
			The site marginals are therefore exact.
			
			\item[(ii)] The finite constraints are
			polynomials with rational coefficients when the prescribed data are
			rational, so the decision procedure for real closed fields applies~\cite{BasuPollackRoy2006}. A realization with fewer
			states can be padded by inaccessible states whose rows are stochastic (choose any allowed letter for their self-loops).
			
			\item[(iii)] For the final assertion,
			choose a horizon at least the largest required word length minus one;
			its feasible presentation passes the finite test. Alternatively,
			the parameters form a compact set of stochastic matrices and initial
			probability rows, and the nested prefix-equality constraints have the
			finite intersection property. Either argument establishes a common
			homogeneous finite presentation. \qedhere
		\end{enumerate}
	\end{proof}

	\begin{theorem}[Finite hidden-state gluing through common separator factors]
		\label{dyn49:thm:hidden-glue}
		Let a cover have a join tree. At site $I$ take a finite initialized
		Markov presentation $(S_I,P_I,\lambda_I,r_I)$ with deterministic native
		observation $r_I:S_I\to\prod_{i\in I}A_i$. For each tree edge $e=I-J$
		assume there is a nonempty finite shared separator state set $Z_e$,
		surjections $q_{I,e}:S_I\to Z_e$, $q_{J,e}:S_J\to Z_e$, and a common
		Markov kernel $P_e$ of which these maps are strong factors. Assume their
		initial pushforwards agree and that the overlap observations on both
		sites are the same function of $Z_e$. Then an initialized finite Markov
		globalization exists on
		\[
		S=\{(s_I)_I:q_{I,e}(s_I)=q_{J,e}(s_J)\text{ for every }e=I-J\}.
		\]
		It realizes the native local word laws exactly, with at most $|S|$ hidden
		states in the deterministic-observation convention. A global letter
		admission mask is satisfied if every tuple of $S$ has its joined
		observation in that mask; weaker reachable-mask conditions may suffice.
	\end{theorem}

	\begin{proof}
		Use root sampling of the local initial distributions, followed by
		conditional sampling in each child given its parent's separator value;
		matching separator marginals ensures the resulting law has the specified
		site marginals and is supported on $S$. At each $s\in S$ the next-state
		site distributions $P_I(s_I,\cdot)$ have equal separator pushforwards
		by the common strong factor and the seam equalities at $s$. The same
		root construction gives a next-state distribution on $S$. This defines
		a Markov kernel whose projection to each $S_I$ is a strong factor.
		Inductively its entire local state paths, hence native observations, have
		the prescribed laws. Join-tree running intersection and agreement on
		edge overlaps make the joined native observation well-defined. Each site-state projection is onto: starting with any prescribed state at that site, root the tree there and extend recursively by surjectivity of the separator maps. The
		mask claim and state bound follow directly. This is a hidden-separator
		version of the rowwise weighted gluing criterion, not a consequence of
		mere equality of separator word laws.
	\end{proof}

	\begin{theorem}[Canonical history globalization need not be finite-state]
		\label{dyn49:thm:greatest-relation-not-finite}
		There are two binary-coordinate sites $AS$ and $SB$ with identical
		separator history laws and finite initialized Markov presentations for
		which the canonical complete-separator join of
		Theorem~\ref{dyn48:thm:projective} has infinite Hankel rank. More precisely, for every $d>0$ a Hankel
		submatrix of size $d$ using word lengths through $2d-1$ has full rank. The same
		atlas nevertheless admits a finite autonomous (noncanonical) global
		realization.
	\end{theorem}

	\begin{proof}
		Choose $0<s<r<1$. At the first site, a fair bit $A$ is constant through
		time. Given $A=1$, let the separator history be $1^n0^\infty$ with
		probability $(1-r)r^n$ for $n\ge0$; given $A=0$ use $(1-s)s^n$.
		The geometric run length is realized by a four-state Markov chain with
		state $(A,S)$, $A$ frozen, $S=0$ absorbing and the probability of
		$1\to1$ equal to $r$ or $s$ according to $A$. Initialize the first bit of
		the separator accordingly. At the second site use the identical law with
		$B$ replacing $A$. Both project to the same separator law. In the canonical
		join $A$ and $B$ are conditionally independent given the full separator
		history. The probability of the global prefix in which $A=B=1$ and the
		separator word is $1^n0$ equals
		\[
		f_n=\frac{\bigl((1-r)r^n/2\bigr)^2}
		{(1-r)r^n/2+(1-s)s^n/2}
		=\frac{c r^n}{1+C q^n},
		\quad c=(1-r)/2,\ C=(1-s)/(1-r),\ q=s/r\in(0,1).
		\]
		Observing the first zero determines the entire separator history almost
		surely, justifying this computation from complete-history conditionals.
		For a finite linear presentation this sequence has the form
		$\lambda M_{111}^{n} M_{101}\mathbf1$ (letters ordered $A,S,B$), so it
		satisfies a nontrivial constant-coefficient linear recurrence by
		Cayley--Hamilton. Suppose $\sum_{j=0}^k c_j f_{n+j}=0$ for all sufficiently
		large $n$, with not all $c_j$ zero. Divide by $c r^n$; the rational
		function
		\[
		R(z)=\sum_{j=0}^k\frac{c_j r^j}{1+Cq^j z}
		\]
		vanishes at infinitely many distinct $z=q^n$, hence is identically zero.
		Its terms have distinct simple poles, so multiplying by the corresponding
		linear denominator and evaluating at that pole gives $c_j=0$ for each
		$j$, a contradiction. Infinite Hankel rank follows from the linear
		residual construction in Proposition~\ref{dyn49:prop:residual-space}.
		For the finite certificates, the submatrix with row words $(111)^i$
		(where the notation means $i$ repetitions of the joint letter $111$)
		and column words consisting of $j$ repetitions of $111$ followed by
		$101$, for $0\le i,j<d$, has entries $f_{i+j}$. If a linear combination
		of its rows is zero, dividing column $j$ by $c r^j$ gives a rational
		function of $z=q^j$ with at most $d$ distinct simple poles and numerator
		of degree at most $d-1$. Its $d$ distinct zeros force it identically
		zero, and the pole argument forces all coefficients zero. Thus this
		submatrix has rank $d$. Its longest tested concatenation has length
		$2d-1$. In particular matching this finite portion of the canonical
		law already requires at least $d$ hidden edge states.
		But choosing $B=A$ and using the original four-state chain glues the
		same local laws, without conditional independence of the private bits.
		This distinguishes failure of a canonical finite realization from
		nonexistence of every finite realization of an atlas.
	\end{proof}

	\subsection{Exact finite-state feasibility by shared matrices}
	\label{dyn51:subsec:matrix}

	The finite budget test of Theorem~\ref{dyn49:thm:budget} uses matrices
	for full joint letters. On a join tree there is a smaller exact
	formulation in terms of local letters and matrix-valued marginal
	compatibility. Here the local presentations are variables, not marks
	that have to be inherited from input presentations.

	\begin{theorem}[Shared-state matrix globalization on join trees and arbitrary covers]
		\label{dyn51:thm:matrix}
		Let $p_I$ be the specified initialized local word laws on a finite
		join tree, given with linear dimensions $d_I$. Fix $d\ge1$.
		\begin{enumerate}[label=(\roman*),leftmargin=*,nosep]
			\item \emph{Join-tree test.} A
			$d$-state unmarked autonomous globalization exists if and only if there
			are a common probability row $\lambda$ of size $d$ and nonnegative
			$d\times d$ matrices $A^I_b$ for all local letters satisfying:
			\begin{align}
				&\sum_b A^I_b\mathbf1=\mathbf1,\qquad
				\lambda A^I_w\mathbf1=p_I(w)\quad(|w|\le d+d_I-1),
				\label{dyn51:eq:matrix-native}\\
				&\sum_{b:b|_S=s}A^I_b=\sum_{c:c|_S=s}A^J_c
				\quad(S=I\cap J,\ I-J\text{ a tree edge},\ s\in\Sigma_S).
				\label{dyn51:eq:matrix-compat}
			\end{align}
			\item \emph{Minimum budget.} In particular the exact unmarked state minimum is the least feasible
			$d$, or infinity. The unknown letter matrices number
			$\sum_I|\Sigma_I|$, rather than $|\prod_i A_i|$; the condition is
			finite real polynomial feasibility for rational input.
			
			\item \emph{Arbitrary covers and masks.} For arbitrary
			covers with an allowed joint-letter set $A$, replace the tree equalities
			by the following linear marginal feasibility for each pair of hidden
			states $h,j$: the vector $((A^I_b)_{hj})_{I,b}$ belongs to the cone
			generated by the local letter-indicator vectors of $a\in A$.
			The latter condition is exact, including on cyclic covers and with
			structural letter admission.
		\end{enumerate}
	\end{theorem}

	\begin{proof}
		\leavevmode
		\begin{enumerate}[label={(iii)},leftmargin=*,nosep]
			\item[(i)] A global realization gives the stated matrices by summing its letter
			matrices, so the conditions are necessary. The direct-sum word-span
			test in Theorem~\ref{dyn49:thm:budget} shows that the native constraints
			hold at all word lengths. On the join tree, summing the matrix overlap
			equalities over $s$ identifies a common total transition matrix $P$;
			this includes empty separators. For a hidden edge $h,j$ with
			$P_{hj}>0$, the local numbers $(A^I_b)_{hj}/P_{hj}$ are probability
			laws on letters with consistent separator marginals. The weighted
			join-tree construction of Theorem~\ref{dyn47:thm:path-glue}, or its
			one-letter restriction, supplies a global letter distribution $\mu_{hj}$
			with these marginals. Set $(M_a)_{hj}=P_{hj}\mu_{hj}(a)$. If $P_{hj}=0$,
			all corresponding local entries vanish by nonnegativity and all global
			entries are set to zero. Then $\sum_a M_a=P$ is stochastic and the
			local projected letter matrices are precisely the chosen $A^I_b$.
			This gives the exact local word laws.
			
			\item[(ii)] The least feasible $d$ is the minimum by (i). The displayed system uses one matrix per local letter and polynomial constraints with rational coefficients for rational input.
			
			\item[(iii)] For a general letter mask, a
			nonnegative coefficient vector of the indicated cone supplies the
			entries $(M_a)_{hj}$ directly. Conversely any such entries give the
			cone representation. Normalization follows from the local row sums. \qedhere
		\end{enumerate}
	\end{proof}

	The theorem also admits a simultaneous Hankel-factor formulation:
	for each site use the finite prefix and suffix index sets of length
	at most $L_I=d+d_I-1$. Take nonnegative row factors $R_I(u)$ and
	column factors $C_I(v)$ with shared $R_I(\epsilon)=\lambda$,
	$C_I(\epsilon)=\mathbf1$, and impose
	\[
	R_I(u)C_I(v)=p_I(uv),\quad
	R_I(ub)=R_I(u)A^I_b,\quad C_I(bv)=A^I_b C_I(v)
	\]
	whenever the respective indices belong to the sets, along with the
	matrix gluing constraints and stochastic normalization. These equations
	are necessary by taking the actual reachable rows and future columns,
	and sufficient since $v=\epsilon$ and the shift equations give
	\eqref{dyn51:eq:matrix-native}. Thus a simultaneous nonnegative
	factorization must retain the shift operators and their matrix-valued
	compatibility. Ordinary or nonnegative ranks of individual or mixed
	Hankel matrices alone do not assert those compatibility conditions.

	\subsection{Shared stochastic drivers}
	\label{dyn52:subsec:driver}

	Write a local initialized word law as $p_I(w)=\lambda_I M^I_w\mathbf1$,
	with $d_I$ hidden states, nonnegative letter matrices and stochastic total
	matrix, as in Section~\ref{dyn51:subsec:matrix}. Let $W$ be the set of
	coordinates belonging to at least two sites, and $S_I=I\cap W$. The
	coordinates outside $W$ belong to exactly one site. Let $(\eta,B_c)$ be
	a finite hidden driver on the shared-coordinate alphabet $\Sigma_W$,
	with $r$ hidden states. Write $\widehat M^I_s$ and $\widehat B^I_s$ for
	the sums of letter matrices projecting to $s\in\Sigma_{S_I}$.
	If $W$ is empty, the shared alphabet consists of one empty letter. The stochastic links below are letterwise intertwinings, related to the classical Markov-function and intertwining constructions of~\cite{RogersPitman1981,DiaconisFill1990}.

	\begin{theorem}[Shared stochastic driver globalization]
		\label{dyn52:thm:driver}
		Suppose for each site there is a row-stochastic $d_I\times r$ matrix $C_I$
		satisfying
		\begin{equation}\label{dyn52:eq:intertwine}
			\lambda_I C_I=\eta,\qquad
			\widehat M^I_s C_I=C_I\widehat B^I_s
			\quad(s\in\Sigma_{S_I}).
		\end{equation}
		\begin{enumerate}[label=(\roman*),leftmargin=*,nosep]
			\item \emph{Existence and bound.} Then the prescribed site word laws have a finite autonomous globalization
			whose state budget is at most
			\begin{equation}\label{dyn52:eq:driver-budget}
				\sum_{h=1}^r\prod_I |\{i:C_I(i,h)>0\}|\ \le r\prod_I d_I.
			\end{equation}
			
			\item \emph{Stationarity.} If all local laws are stationary, the joining can be made stationary without
			increasing the budget.
			
			\item \emph{Effective links.} For fixed rational presentations and a fixed rational
			driver, the link test is a finite linear feasibility problem.
		\end{enumerate}
		This theorem
		allows arbitrary finite covers with a supplied common shared-coordinate
		driver, not only join trees. It does not assert causal conditional independence
		of the private histories or exact pointwise lifting of the original hidden
		states.
	\end{theorem}

	\begin{proof}
		\leavevmode
		\begin{enumerate}[label={(iii)},leftmargin=*,nosep]
			\item[(i)] Use the hidden tuples $(h,(i_I)_I)$ with $C_I(i_I,h)>0$ for all $I$.
			Choose $h$ initially with law $\eta$, and, given $h$ with $\eta(h)>0$,
			choose the local hidden states independently with probabilities
			$\lambda_I(i)C_I(i,h)/\eta(h)$. The initial joint law of $(i_I,h)$ is
			$\lambda_I(i_I)C_I(i_I,h)$, including at zero-mass $h$.
			At a tuple first sample the shared letter $c$ and next driver state $h'$
			with weights $B_c(h,h')$. Given this choice, write $s=c|_{S_I}$ and
			choose local letter $b$ and next hidden state $j$ at site $I$ with weights
			\begin{equation}\label{dyn52:eq:disintegration}
				K_I(b,j\mid i,s,h')=
				\frac{M^I_b(i,j)C_I(j,h')}{(\widehat M^I_s C_I)(i,h')}
				\quad(b|_{S_I}=s).
			\end{equation}
			The denominator is positive whenever the driver choice has positive weight:
			by~\eqref{dyn52:eq:intertwine} it is at least
			$C_I(i,h) B_c(h,h')>0$. The weights sum to one and every positive next
			choice has $C_I(j,h')>0$. Conditional choices at distinct sites may be
			independent since the remaining letter coordinates are disjoint. Together
			with $c$, they specify a unique global letter.
			
			For the marginal on site $I$ and driver hidden state, the letter transition
			from $(i,h)$ to $(j,h')$ is the product of
			$\widehat B^I_s(h,h')$ and~\eqref{dyn52:eq:disintegration}. Averaging the
			transition over $h$ with weights $C_I(i,h)$ gives exactly
			$M_b^I(i,j) C_I(j,h')$: if the denominator vanishes, the nonnegative numerator is zero and every driver weight from a state with $C_I(i,h)>0$ for that separator letter and destination is zero. Such conditional choices are unused. Otherwise the denominator cancels by the intertwining equation.
			Induction therefore gives the joint weight of local word $w$ and terminal
			$(j,h')$ as $(\lambda_I M^I_w)_j C_I(j,h')$. Summing proves all prescribed
			word laws. Counting the tuple set gives~\eqref{dyn52:eq:driver-budget}.
			
			\item[(ii)] For stationary marginals, each shifted initialization of the constructed global hidden matrix gives the same local word laws. Ces\`aro average these shifted initializations and take a convergent subsequence in the finite probability simplex. Its limit is stationary by telescoping and preserves every local word law by linearity. For rational matrices the stationary choice can instead be made by rational linear feasibility, imposing stationarity and preservation of each local word law through the finite linear comparison bound.
			
			\item[(iii)] Finally all equations in $C_I$ are linear when the two presentations are given. The normalization and nonnegativity constraints are linear as well; a feasible rational link system has a rational solution and the sampling formulas then use only rational operations. \qedhere
		\end{enumerate}
	\end{proof}

	\subsection{A common finite positive refinement without a span hypothesis}
	\label{dyn53:subsec:refinement}

	We use the edge-emitting convention $p_I(w)=\lambda_I M^I_w\mathbf1$ of Section~\ref{dyn52:subsec:driver}. All alphabets and state sets are nonempty and finite. A star atlas here has common separator $S$ and disjoint remaining coordinate sets. Write $L^I_s=\widehat M^I_s$ for the separator letter matrices. The following elementary construction is a stochastic version of a common positive simulation in weighted automata; related span-shaped simulation results over nonnegative reals and rationals are treated in \cite{dyn53:SokolovaWoracek}. We give the construction and its budget because its direction and normalization are important for globalization.

	\begin{theorem}[Common positive separator refinement]
		\label{dyn53:thm:refinement}
		Let $k$ finite stochastic presentations $(\lambda_I,L^I_s)$ have identical separator word law. Define the compatible belief polytope
		\begin{equation}\label{dyn53:eq:polytope}
			D=\{(x_I)_I:x_I\ge0,\ x_I\mathbf1=1,\
			x_I L^I_w\mathbf1=x_J L^J_w\mathbf1\quad\forall I,J,w\}.
		\end{equation}
		\begin{enumerate}[label=(\roman*),leftmargin=*,nosep]
			\item \emph{Finite polytope.} It is nonempty and defined by finitely many linear equations and inequalities.
			
			\item \emph{Common presentation.} There exist a common stochastic separator presentation $(\eta,B_s)$ with $r=|\operatorname{Vert}D|$ states and stochastic matrices $C_I$ of size $r\times d_I$ such that
			\begin{equation}\label{dyn53:eq:outgoing}
				\eta C_I=\lambda_I,\qquad B_s C_I=C_I L^I_s.
			\end{equation}
			
			\item \emph{State bounds.} In particular
			\begin{equation}\label{dyn53:eq:refinement-bound}
				r\le \prod_I(2^{d_I}-1).
			\end{equation}
			If $D=\{x\ge0:Ex=b\}$ in concatenated coordinates of dimension $D_0=\sum_I d_I$ and $\rho=\operatorname{rank}E$, also $r\le\sum_{j=1}^{\rho}\binom{D_0}{j}$.
			
			\item \emph{Converse.} Conversely the existence of~\eqref{dyn53:eq:outgoing} implies equality of all separator word laws.
			
			\item \emph{Stationarity and rationality.} If all $\lambda_I$ are stationary for their total hidden matrices, $\eta$ may be chosen stationary without changing the bound. Rational data give a rational construction.
		\end{enumerate}
	\end{theorem}

	\begin{proof}
		\leavevmode
		\begin{enumerate}[label={(iii)},leftmargin=*,nosep]
			\item[(i)] The equal-future equations need only compare each site with one fixed reference site, through word length $d_I+d_{\rm ref}-1$: use the finite linear word test from Section~\ref{dyn49:subsec:finite}. Equivalently compute the span of the column vectors $(L^I_w\mathbf1,L^{\rm ref}_w\mathbf1)$ by successive multiplication until it stabilizes. The polytope contains the tuple of initializations.
			
			\item[(ii)] For $x\in D$ the mass $m_s=x_I L^I_s\mathbf1$ is independent of $I$. If $m_s>0$, the tuple $(x_I L^I_s/m_s)_I$ belongs to $D$, by applying its defining equations to the word $sw$. If $m_s=0$, each nonnegative row $x_I L^I_s$ vanishes. For each vertex $v$ of $D$ choose nonnegative barycentric weights for the normalized updated tuple and multiply them by $m_s$; these are the row $B_s(v,\cdot)$. Since $\sum_s m_s=1$, the total matrix is stochastic. Let $C_I(v,\cdot)=v_I$ and express the tuple $(\lambda_I)_I$ as the mixture $\eta$ of vertices. The two equations follow by construction, including the zero-mass case. Induction then proves equality of every word law for the common presentation and the local presentations.
			
			\item[(iii)] Here $Ex=b$ includes each site's normalization and the equal-future constraints. A vertex of $\{x\ge0:Ex=b\}$ has linearly independent columns of $E$ on its positive support: otherwise a nonzero kernel vector supported there permits a two-sided small perturbation. In particular the support determines the vertex uniquely, and its size is at most $\rho$. Each site has nonempty support, giving both bounds.
			
			\item[(iv)] The converse follows by multiplying the intertwining identities and using $C_I\mathbf1=\mathbf1$.
			
			\item[(v)] For stationary local hidden rows, $(\eta B^n) C_I=\lambda_I$ for $B=\sum_s B_s$. A Ces\`aro limit of these rows is stationary and keeps the required images. For rational data the vertex enumeration, convex decompositions, and stationary choice are rational linear feasibility problems; real feasibility of rational linear constraints implies rational feasibility. This proves the claims. \qedhere
		\end{enumerate}
	\end{proof}
	
	The arrow in~\eqref{dyn53:eq:outgoing} is opposite to~\eqref{dyn52:eq:intertwine}. This alone is not the latter driver construction. In particular, private histories can introduce future functions outside the separator's intrinsic reachable predictive span. The next step converts the comparison direction exactly, rather than assuming that this extra predictive information disappears.

	\begin{theorem}[Stationary finite autonomous globalization on stars and join trees]
		\label{dyn53:thm:stationary}
		\leavevmode
		\begin{enumerate}[label=(\roman*),leftmargin=*,nosep]
			\item Every compatible atlas of stationary finite hidden Markov word laws on a finite join-tree cover has a stationary finite hidden Markov globalization. Here stationarity of a word law means invariance under a one-letter shift, $\sum_a p(aw)=p(w)$ for every $w$.
			
			\item For a star in the common-separator sense of Section~\ref{dyn53:subsec:refinement}, with presentation sizes $d_I$, a constructive bound is
			\begin{equation}\label{dyn53:eq:star-budget}
				d_{\rm global}\le r\prod_I d_I\le\prod_I d_I(2^{d_I}-1),
			\end{equation}
			where $r$ is the vertex count of the compatible belief polytope in Theorem~\ref{dyn53:thm:refinement}, computed using the local stationary supports. The driver may have fewer than $r$ states after deletion of its zero-stationary-mass states.
			
			\item For a general join tree, order sites root-first, with parent preceding child; put $D_1=d_1$ and
			\begin{equation}\label{dyn53:eq:recursive-budget}
				D_{j+1}=D_jd_{j+1}(2^{D_j}-1)(2^{d_{j+1}}-1).
			\end{equation}
			Then $D_k$ suffices.
		\end{enumerate}
		No predictive-rank, finite-residual, ergodicity or common-span assumption is imposed. Rational local presentations admit a rational stationary globalization.
	\end{theorem}

	\begin{proof}
		\leavevmode
		\begin{enumerate}[label={(iii)},leftmargin=*,nosep]
			\item[(i)] For existence, the star construction in (ii) followed by the tree induction in (iii) gives a stationary finite globalization.
			\item[(ii)] First make every supplied hidden initialization stationary by Ces\`aro averaging its shifts; because the observable word law is stationary, each shift has the same observable law. The stationary initialization can equivalently be found by linear feasibility, imposing nonnegativity, normalization, stationarity and equality of the word laws through a finite linear recognition bound. For rational input use an exact rational feasible row; all support deletions below refer to exactly zero entries, not a numerical tolerance. Ces\`aro averaging establishes existence and is not a prescribed approximate algorithm. Discard zero stationary hidden states, which cannot be reached from positive stationary states. Apply Theorem~\ref{dyn53:thm:refinement} to the common separator of a star, with stationary hidden rows $\lambda_I$ and $\eta$. Discard zero stationary states of $B$ as well. Each discarded state has no incoming transition from a positive stationary state, since the total matrix is nonnegative and its stationary mass is zero. The retained letter matrices are therefore closed on the positive support, and restricting the rows of $C_I$ to this support preserves both equations in~\eqref{dyn53:eq:outgoing}. Let $D_\lambda$ denote the diagonal matrix with the entries of a strictly positive probability row $\lambda$. For any stationary letter presentation define its reversed letter matrices by
			\[
			M_b^{I,\leftarrow}=D_{\lambda_I}^{-1}(M_b^I)^{\mathsf T}D_{\lambda_I},
			\qquad B_s^{\leftarrow}=D_\eta^{-1}B_s^{\mathsf T}D_\eta.
			\]
			They are nonnegative and their total matrices are stochastic by stationarity. The word probability of the reversed presentation is the probability of the word in reversed order, by telescoping the diagonal factors. Define the Bayes transpose of $C_I$ by
			\[
			C_I^\dagger=D_{\lambda_I}^{-1}C_I^{\mathsf T}D_\eta.
			\]
			The normalization $\eta C_I=\lambda_I$ gives $C_I^\dagger\mathbf1=\mathbf1$, and $C_I\mathbf1=\mathbf1$ gives $\lambda_I C_I^\dagger=\eta$. Transposing $B_s C_I=C_I L^I_s$ and multiplying by the two diagonal factors gives exactly
			\begin{equation}\label{dyn53:eq:reverse-links}
				(L^I_s)^{\leftarrow}C_I^\dagger=C_I^\dagger B_s^{\leftarrow}.
			\end{equation}
			Thus the reversed full local letter laws satisfy the hypotheses of Theorem~\ref{dyn52:thm:driver}. Its construction gives a finite globalization of all reversed word laws, which can be made stationary without increasing the budget. Reverse this stationary global presentation (again discard any zero stationary states). The result has the required forward word laws. This proves~\eqref{dyn53:eq:star-budget}. Empty shared coordinate sets use the one-letter separator, which also satisfies the argument.
			\item[(iii)] For a general join tree, attach one child at a time. By running intersection the new site's intersection with the previously attached union is exactly its parent separator. The already constructed global law on that union has the correct separator word law. Apply the two-site instance of the star construction, keeping all previously achieved marginals. Induction yields the bound~\eqref{dyn53:eq:recursive-budget}. All matrices used in the rational case are constructed by rational linear feasibility, finite multiplication, positive division or rational diagonal reversal; the stationary choices are rational linear feasibility as in the preceding theorem. \qedhere
		\end{enumerate}
	\end{proof}

	\subsection{Removing stationarity by a common reset and harmonic reweighting}
	\label{dyn53:subsec:reset}

	A stationary extension of a reset-augmented process need not have a globally constant reset hazard conditional on all private observations. Simply deleting reset letters and renormalizing their matrices row by row could change the prescribed marginals. The following argument avoids that difficulty.

	\begin{lemma}[Harmonic extraction of an initialized word law]
		\label{dyn53:lem:harmonic}
		Let $\alpha$ be a probability row, $N_a$ finitely many nonnegative matrices, $N=\sum_a N_a$, and suppose $\alpha N^n\mathbf1=1$ for all $n\ge0$. Suppose for every prescribed site projection and local word $w$,
		\begin{equation}\label{dyn53:eq:wildcard}
			\alpha\widehat N^I_w N^n\mathbf1=p_I(w)\qquad(n\ge0).
		\end{equation}
		There is a finite stochastic presentation of a global law with exactly these local laws, using no more states than the matrices $N_a$. It can be constructed rationally if the matrices, $\alpha$, and finite presentations of the $p_I$ are rational.
	\end{lemma}

	\begin{proof}
		Restrict to states reachable from $\alpha$ under the total matrix $N$. For each such state $i$, choose $t$ with $(\alpha N^t)_i>0$. Since all entries are nonnegative,
		\[
		(N^n\mathbf1)_i\le \frac{\alpha N^{t+n}\mathbf1}{(\alpha N^t)_i}=\frac1{(\alpha N^t)_i}.
		\]
		Hence $N^n\mathbf1$ is uniformly bounded on the retained state set. Take a convergent subsequence of the Ces\`aro means of these vectors. Its limit $h\ge0$ satisfies $Nh=h$, by telescoping the means, and $\alpha h=1$. By~\eqref{dyn53:eq:wildcard}, $\alpha\widehat N_w^I h=p_I(w)$ for every local word. Put $Z=\{i:h_i>0\}$. If $h_i=0$, no positive $N$ transition leads from $i$ to $Z$, by $Nh=h$ and nonnegativity. On $Z$ use the matrices
		\begin{equation}\label{dyn53:eq:doob}
			M_a=D_{h|_Z}^{-1}(N_a)_{Z,Z}D_{h|_Z},\qquad \lambda=(\alpha_i h_i)_{i\in Z}.
		\end{equation}
		They have stochastic total matrix and probability initialization, and their projected word probabilities are $\alpha\widehat N_w^I h$. All summands with an intermediate zero-$h$ state and a positive-$h$ endpoint vanish, so restriction does not lose a term in this identity. For rational construction, seek $h\ge0$ satisfying $Nh=h$, $\alpha h=1$, and the word equalities through length $d+d_I-1$ for each local representation size $d_I$ and retained matrix size $d$. These finite linear constraints are feasible by the proof; the word test for linear representations guarantees all remaining equalities. A rational solution exists and~\eqref{dyn53:eq:doob} is rational. In the effective rational case solve these equations and inequalities exactly and use the strictly positive support of that exact solution; an approximate Ces\`aro iterate need not satisfy either the harmonic equation or the word equalities. Harmonic normalization supplies prefix consistency of the resulting joint word law.
	\end{proof}

	\begin{theorem}[Unrestricted finite autonomous join-tree globalization]
		\label{dyn53:thm:universal}
		Let $\mathcal H$ be a finite join-tree cover of finite native coordinates with finite alphabets. At every site prescribe an initialized finite hidden Markov word law, with no assumption of stationarity.
		\begin{enumerate}[label=(\roman*),leftmargin=*,nosep]
			\item \emph{Existence.}
			\[
			\boxed{\begin{gathered}
				\text{A common finite hidden Markov realization exists}\\
				\Longleftrightarrow\quad
				\text{The prescribed separator word laws agree.}
			\end{gathered}}
			\]
			Here the realization has exactly the prescribed initialized
			local word laws.
			
			\item \emph{State bounds.} The bounds~\eqref{dyn53:eq:star-budget} (stars) and~\eqref{dyn53:eq:recursive-budget} (general join trees) remain valid using the original local state sizes.
			
			\item \emph{Rational solutions.} Rational inputs admit a rational solution. In particular all ergodic infinite-residual separator cases in the unrestricted unmarked join-tree question are included.
		\end{enumerate}
	\end{theorem}

	\begin{proof}
		\leavevmode
		\begin{enumerate}[label={(iii)},leftmargin=*,nosep]
			\item[(i)] Necessity is projection. For sufficiency, fix $0<\varepsilon<1$ and write $c=1-\varepsilon$. Adjoin one common binary reset marker to every site. When the marker is zero, a local letter $b$ is emitted with matrix $c M_b^I$. When the marker is one, all native coordinates are at a fixed default letter and the single local reset matrix is $R_I=\varepsilon\mathbf1\lambda_I$. All other marker-one letters have zero matrices. Defaults are chosen consistently per coordinate. The new total matrix is $cP_I+R_I$, $P_I=\sum_b M_b^I$. Its stationary row is
			\begin{equation}\label{dyn53:eq:reset}
				\pi_I=\varepsilon\sum_{n\ge0}c^n\lambda_I P_I^n
				=\varepsilon\lambda_I(\mathrm{Id}-cP_I)^{-1}.
			\end{equation}
			It is a probability row, and $\pi_I(cP_I+R_I)=\pi_I$. The reset-augmented separator word laws agree. Indeed for a separator block $w$ with no reset, its stationary probability is
			\[
			c^{|w|}\varepsilon\sum_{n\ge0}c^n
			\sum_{|v|=n}q_S(vw),
			\]
			where $q_S$ is the common original initialized separator word law (the sum is over all separator words of length $n$, including the empty word when $n=0$). For words with resets, the segment before the first reset has the preceding stationary probability, each reset contributes $\varepsilon$ and resets the hidden row to $\lambda_I$, and each subsequent reset-free segment contributes its initialized separator word probability times the corresponding power of $c$. The common marker also makes this argument valid for an originally empty separator. Hence every augmented separator is compatible. Adding the same coordinate to every site preserves the running-intersection property.
			
			Apply Theorem~\ref{dyn53:thm:stationary} to the augmented atlas. Its global presentation has the same finite state bound as for the original local state sizes. Let $\widetilde\lambda$ be its initial row, $G_*$ the global reset matrix and $G_a$ the matrices of marker-zero global letters $a$ (letters with a reset and nondefault native coordinates have probability zero and may be discarded on reachable states). The initial reset event has probability $\varepsilon$; put $\alpha=\widetilde\lambda G_*/\varepsilon$ and $N_a=G_a/c$. For a local word $w$ of length $\ell$ and $n$ subsequent marker-zero letters with arbitrary native values, the probability of this event preceded by the reset is $\varepsilon c^{\ell+n}p_I(w)$, by the exact augmented local law. Since the reset marker is common, summing over all unobserved native letters gives precisely~\eqref{dyn53:eq:wildcard}. Taking the empty word gives $\alpha N^n\mathbf1=1$ for every $n$. Lemma~\ref{dyn53:lem:harmonic} yields the desired finite stochastic law on the original global alphabet, without increasing the state budget.
			
			\item[(ii)] The alphabet augmentation changes letter matrices but introduces no extra local hidden states; the support restrictions and harmonic extraction never increase their respective state counts.
			
			\item[(iii)] Choosing rational $\varepsilon$ and using the rational construction of the preceding theorem and lemma proves the last claim. \qedhere
		\end{enumerate}
	\end{proof}

	\paragraph{Comparison with Markov couplings over a Markov image.}
	The qualitative existence in Theorem~\ref{dyn53:thm:universal} can alternatively be obtained from the common positive refinement and Crane--Holroyd--Russell's Theorem~1 \cite{CraneHolroydRussell2026}. For a common separator, use the outgoing links $B_s C_I=C_I L_s^I$, $\eta C_I=\lambda_I$ from Theorem~\ref{dyn53:thm:refinement}. For each individual site lift to pairs $(h,x)$ with $C_I(h,x)>0$, initialized by $\eta(h)C_I(h,x)$. For a local letter $b$ restricting to $s$, define the lifted letter matrix by
	\[
	K_b^I((h,x),(h',y))=
	\frac{B_s(h,h')M_b^I(x,y)C_I(h',y)}{(B_s C_I)(h,y)},
	\]
	with zero terms assigned zero. Whenever $C_I(h,x)M_b^I(x,y)>0$, the denominator is positive by the outgoing link identity. Summing over $h'$ recovers $M_b^I(x,y)$, so the total lifted matrix is stochastic and the original site letter law is exact. Conversely, summing over $x$ with weights $C_I(h,x)$ and over letters restricting to $s$ gives $B_s(h,h')C_I(h',y)$. By induction over a driver path, the conditional local hidden row after its last edge is $C_I(h',\cdot)$; the driver edges therefore have precisely the common law initialized by $\eta$ and the matrices $B_s$. Regard the positive labeled transitions of each lifted process as the states of a finite Markov chain, with next-transition probabilities determined by the current endpoint. The image of this chain given by $(h,s,h')$ is the common Markov driver-edge chain. The theorem of \cite{CraneHolroydRussell2026} joins the two enlarged chains over this Markov image; projecting their emitted letters gives the required finite joint word law. Induction along the join tree completes this alternative existence argument. It is a consequence of the positive-refinement and Markov-image results together, not of the Markov-image theorem applied directly to the prescribed separator. The preceding proof does not rely on that coupling theorem and supplies its stated quantitative and rational conclusions.

	\paragraph{A six-state joining of the non-Markov-image example.}
	In Example~1 of \cite{CraneHolroydRussell2026}, the two chains have states $1,2,3,4$, begin at $2$ and $1$, and their binary images are zero exactly at $2$ and $1$, respectively. No Markov coupling on the original pairs exists. Nevertheless, with the six hidden states having ordered output pairs
	\[
	((4,3),(3,2),(3,4),(2,1),(1,2),(3,2)),
	\]
	take initial state $4$ and transition matrix
	\[
	Q=\begin{pmatrix}
		0&1/2&1/2&0&0&0\\
		1/2&0&0&1/2&0&0\\
		1&0&0&0&0&0\\
		0&0&0&0&1/3&2/3\\
		0&0&0&1&0&0\\
		3/4&0&0&1/4&0&0
	\end{pmatrix}.
	\]
	The first coordinate follows the original first chain: on observing $3$, the hidden row is concentrated at state $6$ if the preceding first-coordinate value was $2$, or uniform on states $2,3$ if it was $4$. In both cases the next first-coordinate value is $4$ with probability $3/4$ and $2$ with probability $1/4$. Each of its other values has a single hidden representative, with the required transition row. The second coordinate on observing $2$ has hidden row $(1/3,2/3)$ on states $5,6$ if preceded by $1$, or is concentrated at state $2$ if preceded by $3$. Both rows give the next second-coordinate values $1,3$ with probabilities $1/2,1/2$. Each other second-coordinate value has a single representative. These conditional computations verify both full initialized Markov laws by induction, and all output pairs have equal binary images. In particular the repeated pair $(3,2)$ has two distinct hidden representatives. This explicit coupling illustrates why finite hidden realization and a Markov coupling of the original pairs are different contracts.

	\begin{corollary}[Effective unmarked existence and exact minimum on join trees]
		\label{dyn53:cor:decision}
		\leavevmode
		\begin{enumerate}[label=(\roman*),leftmargin=*,nosep]
			\item On a given finite join tree with rational finite local presentations, existence of an unmarked finite autonomous globalization is decidable by finitely many separator word comparisons, through length $d_I+d_J-1$ on each edge.
			
			\item If incompatible, a word of at most that length is a certificate of nonexistence even of an arbitrary probability globalization.
			
			\item If compatible, the finite budget bounds above and the exact bounded matrix test of Theorem~\ref{dyn51:thm:matrix} decide the minimum unmarked state budget. The feasibility search over $1,\ldots,D_k$ therefore terminates, giving a computable bound in terms of the local presentation sizes.
		\end{enumerate}
		The exact test uses real polynomial feasibility; rational feasibility at the minimum budget is a separate question. In the unconstrained join-tree category there are no compatible instances requiring infinitely many states in every realization.
	\end{corollary}

	\begin{proof}
		\leavevmode
		\begin{enumerate}[label={(iii)},leftmargin=*,nosep]
			\item[(i)] Apply the finite linear equivalence bound to the projected local presentations. Necessity and sufficiency are Theorem~\ref{dyn53:thm:universal}.
			
			\item[(ii)] A distinguishing separator word contradicts the equality of marginals of any probability extension.
			
			\item[(iii)] The fixed-budget real feasibility test is exact, and the bound makes the search finite; for a stationary minimum one may add the stationarity equation or stationarize a witness as above. \qedhere
		\end{enumerate}
	\end{proof}

	\begin{corollary}[Sharp finite autonomous weighted recovery correspondence]
		\label{dyn53:cor:finite-wug}
		Let $\mathscr B_k^{\rm fin}$ be the finite hidden Markov subclass of $\mathscr B_k$ from Definition~\ref{dyn48:def:wdrk}. Let $\operatorname{FWUG}(\mathcal H,\mathscr C)$ mean that every compatible atlas of initialized local laws in $\mathscr C$ has an unmarked finite autonomous globalization, with no requirement that the global law belongs to $\mathscr C$. Allow binary alphabets as in Theorem~\ref{dyn48:thm:sharp-wug}. Then
		\begin{equation}\label{dyn53:eq:finite-wug}
			\operatorname{FWUG}(\mathcal H,\mathscr B_k^{\rm fin})
			\quad\Longleftrightarrow\quad \max\{2,2k-1\}<\nu(\mathcal H).
		\end{equation}
		For the entire finite hidden Markov class, the universal property holds exactly for join-tree covers.		
	\end{corollary}

	\begin{proof}
		For $k=1$ take the product of the finite single-coordinate processes. For $k\ge2$ at the asserted threshold, the sufficiency proof of Theorem~\ref{dyn48:thm:sharp-wug} constructs a block-bounded completion on each maximal primal clique, with compatible prefix laws. By Lemma~\ref{dyn48:lem:block}, its independent partitions stabilize and the process is the product of the marginal processes on the resulting blocks of size at most $k$. Each block is guarded by an original site, and its process is exactly a restriction of that site's law. It is therefore finite hidden Markov, and their finite product gives a finite presentation of the completed clique law. Use Theorem~\ref{dyn53:thm:universal} on the clique tree. The original local laws are its exact marginals. The converse is given by the frozen-in-time binary parity witnesses of Theorem~\ref{dyn48:thm:sharp-wug}; these have finite hidden presentations and do not extend even at time zero. The unrestricted case follows by choosing a finite block bound large enough for the witness on any non-join-tree cover, and by the universal theorem on join trees. No global block-rank bound is implied.		
	\end{proof}

	\subsection{Intrinsic nonfinite certificates and the limits of the unmarked theorem}
	\label{dyn53:subsec:certificates}

	Theorem~\ref{dyn53:thm:universal} settles the existence problem for unmarked finite join-tree input. The situation for additional admission, causality, exact canonical joining, cyclic covers and minimum native-carrier realization is different. The next statement gives a precise finite obstruction test for probability extension and distinguishes it from nonfinite realization of a selected process.

	\begin{proposition}[Finite marginal infeasibility certificates and nonfinite realization]
		\label{dyn53:prop:certificates}
		\leavevmode
		\begin{enumerate}[label=(\roman*),leftmargin=*,nosep]
			\item \emph{Finite certificates.} On any finite cover, with arbitrary compatible or incompatible prefix-consistent finite-alphabet local laws and prefix-closed admission masks, nonexistence of an admitted global probability process is equivalent to failure of the admitted marginal polytope test of Theorem~\ref{dyn49:thm:admission} at some finite horizon.
			\item \emph{Effective semidecision.} With rational local finite presentations and decidable masks this is a semidecision by enumeration and rational linear programming.
			\item \emph{Unmarked join trees.} On a join tree without additional masks, equality of all separator words through the bound of Corollary~\ref{dyn53:cor:decision} is sufficient to exclude all such incompatibility and all unmarked nonfinite-realization obstructions.
			\item \emph{Nonfinite selected laws.} For a separately specified global word law, nonsingular Hankel minors of arbitrarily large size certify that no finite autonomous presentation realizes that particular law. If its Hankel rank is finite, the rank alone is only a linear realization bound; a finite positive realization still requires an invariant finitely generated positive cone, in the classical sense of \cite{dyn49:Heller}.
		\end{enumerate}
		No equivalence between nonexistence of an admitted finite presentation and a finite marginal infeasibility certificate is asserted on arbitrary covers.
	\end{proposition}

	\begin{proof}
		\leavevmode
		\begin{enumerate}[label={(iii)},leftmargin=*,nosep]
			\item[(i)] The admitted polytope criterion and prefix compactness give the first equivalence.
			\item[(ii)] A failed rational finite polytope test has a rational separating certificate.
			\item[(iii)] The join-tree assertion is Corollary~\ref{dyn53:cor:decision}.
			\item[(iv)] If a word law has $d$ hidden states, its Hankel matrix factors into prefix rows and suffix columns of dimension $d$, so every finite minor of larger order vanishes. The invariant positive cone distinction is the classical positive realization criterion reviewed in Section~\ref{dyn49:subsec:finite}; a finite representation implies finite linear dimension, but positivity is not entailed by linear dimension in general. For example, the canonical infinite-rank join of Theorem~\ref{dyn49:thm:greatest-relation-not-finite} coexists with a finite noncanonical join; it would be invalid to infer a failure of atlas globalization from its minors. This verifies the separation of certificate types. \qedhere
		\end{enumerate}
	\end{proof}

	\subsection{Exact preservation of the independent history-block bound}

	The local bound in Theorem~\ref{dyn48:thm:sharp-wug} controls extension, but an extension may merge independent blocks. We now characterize when a globalization can retain the same bound. This requires both the guard threshold and a condition on the connected components of the primal graph.

	\begin{theorem}[Block-preserving weighted globalization]
		\label{dyn53:thm:block-preserving}
		Let $k\ge1$ and let $\mathcal H$ be a finite cover of a nonempty finite coordinate set $V$. Consider a compatible atlas of initialized local laws in $\mathscr B_k$. Form the hypergraph $\mathcal D$ on $V$ by taking the union of the local dependence hypergraphs over all sites and all finite horizons, as in Definition~\ref{dyn48:def:wdrk}. Include isolated vertices as components.
		\begin{enumerate}[label=(\roman*),leftmargin=2em]
			\item Suppose the atlas has a probability extension. It has an extension in $\mathscr B_k$ if and only if every connected component of $\mathcal D$ has at most $k$ vertices. If the atlas has a finite hidden Markov extension, the same condition is equivalent to existence of an extension in $\mathscr B_k^{\rm fin}$. A rational finite extension, when given, yields a rational block-preserving extension.
			\item Universality ranges over finite alphabets allowing two symbols at each coordinate, as in Theorem~\ref{dyn48:thm:sharp-wug}. Every compatible $\mathscr B_k$-atlas has an extension in $\mathscr B_k$ precisely under the following conditions: for $k=1$ there is no restriction; for $k\ge2$,
			\begin{equation}\label{dyn53:eq:block-preserving}
				\boxed{\begin{gathered}2k-1<\nu(\mathcal H),\\\text{every primal connected component of size greater than $k$ is complete.}\end{gathered}}
			\end{equation}
			The identical geometric criterion holds for extension of every compatible $\mathscr B_k^{\rm fin}$-atlas within $\mathscr B_k^{\rm fin}$. For rational finite local input satisfying the criterion, a rational finite block-preserving extension exists.
		\end{enumerate}
	\end{theorem}

	\begin{proof}
		\leavevmode
		\begin{enumerate}[label={(iii)},leftmargin=*,nosep]
			\item[(i)] A dependence edge in any local law is a dependence edge of every global extension, because mixed cumulants on the specified coordinates depend only on their joint marginal. Hence every component of $\mathcal D$ must lie in a single independent block of any global law. This proves necessity. Conversely let $p$ be a probability extension, let $\mathcal P$ be the component partition of $\mathcal D$, and define the global process $q=\bigotimes_{B\in\mathcal P}p_B$, where $p_B$ is the complete marginal process of $p$ on $B$. By Lemma~\ref{dyn48:lem:block}, each local law factors over the nonempty intersections of its site with $\mathcal P$: all its dependence edges are contained in these blocks. Consequently for each site $I$,
			\[
			q_I=\bigotimes_{B\in\mathcal P,\,I\cap B\ne\varnothing}p_{I\cap B}=p_I.
			\]
			If every $B$ has size at most $k$, this is an extension in $\mathscr B_k$. Marginalization and finite independent products preserve finite hidden Markov presentations, by summing letter matrices and taking tensor products, respectively. They also preserve rationality. This proves all assertions in (i).
			\item[(ii)] For $k=1$, take the product of the common single-coordinate processes; the same construction works for finite presentations. For $k\ge2$, assume~\eqref{dyn53:eq:block-preserving}. Theorem~\ref{dyn48:thm:sharp-wug} supplies a probability extension; Corollary~\ref{dyn53:cor:finite-wug} supplies a finite extension in the finite-input case. We show that every component of $\mathcal D$ has size at most $k$. Every local dependence edge has size at most $k$ and is contained in a chart, so every component lies in a primal connected component. If a component of $\mathcal D$ contains more than $k$ vertices, the containing primal component is complete by hypothesis. Select a connected sequence of dependence edges, adding each edge so that it meets the accumulated union, until the union $U$ first exceeds size $k$. Then $k<|U|\le2k-1$. Because $U$ is a clique, the guard inequality places it in an original chart $J$. For each selected edge, overlap agreement identifies its entire history marginal in its original chart with that in $J$. Thus all these edges belong to the all-horizon dependence hypergraph of $J$. They are connected and their union has more than $k$ vertices, contradicting $b(p_J)\le k$. Part (i) now gives block-preserving extension. For rational finite inputs the clique-completion proof of Corollary~\ref{dyn53:cor:finite-wug} uses only rational restrictions, finite products and the rational join-tree construction of Theorem~\ref{dyn53:thm:universal}, so part (i) also gives the rational assertion.
			
			For necessity, universal block-preserving extension implies ordinary universal extension; the frozen binary obstructions in Theorem~\ref{dyn48:thm:sharp-wug} give the guard inequality, also for finite inputs. Suppose a primal connected component has more than $k$ vertices but is not complete. It contains an induced path of length two (take the first three vertices on a shortest path between nonadjacent vertices). Starting with those three vertices, add adjacent vertices in the component until a connected set $U$ of size $k+1$ is obtained. This is possible since $k\ge2$. The nonadjacent endpoints of the induced path ensure that no chart contains $U$. Let all coordinates of $U$ be copies of one fair bit, and let all coordinates outside $U$ be independent fair bits, also independent of that bit. Make each bit constant in time. Prescribe the local restrictions of this global process. They are compatible, finite and rational. Each chart meets $U$ in at most $k$ vertices, so every local law belongs to $\mathscr B_k^{\rm fin}$. But every global extension must give identical nondegenerate bits on every primal edge within $U$. Connectedness forces all bits of $U$ to be equal almost surely. Each such edge is a dependence edge, so the global dependence component contains $U$ and has size at least $k+1$. No extension can belong to $\mathscr B_k$. This proves necessity for both universality statements. \qedhere
		\end{enumerate}
	\end{proof}

	Whenever an extension exists, the construction in (i) attains the smallest possible largest independent block size, namely $\max_{B\in\mathcal P}|B|$. The partition of $\mathcal D$ is the finest common coarsening of the local finest independent partitions extended by singleton coordinates outside their respective sites. For finite rational input, Proposition~\ref{dyn48:prop:rank-recognition} computes those partitions, so the component-size test in (i) is effective when existence of an extension is given. For example, the two-site cover $\{ab,bc\}$ with $k=2$ admits a probability and a finite autonomous extension for every compatible block-bounded atlas. It fails the stronger condition because its three-vertex primal component is not complete. Conversely, a chordal conformal cover whose connected components are complete always preserves the block bound. The additional condition in~\eqref{dyn53:eq:block-preserving} measures the possibility of merging dependence blocks across distinct charts.

	\section{Structural admission and obstructions to finite realization}
	\label{sec:structural-admission-obstructions}

	Probability extension and finite autonomous realization are different requirements. Theorem~\ref{dyn49:thm:admission} first characterizes admitted probability extension by finite-prefix feasibility. Theorems~\ref{dyn54:thm:triangle-infinite} and~\ref{dyn54:thm:conditional-spectrum} then establish the infinite-Hankel-rank obstruction and its exact cover classification; Corollary~\ref{dyn54:cor:binary-stationary} strengthens the obstruction to binary stationary ergodic local laws. Finally, Corollary~\ref{dyn54:cor:regular-mask} applies the cyclic construction to operational admission on a join-tree observation cover.

	\subsection{Weighted globalization with general structural admission}
	\label{dyn49:subsec:admission}

	We now allow a global admission condition in addition to the local
	laws. At horizon $n$ let $A_n\subseteq\prod_{i\in V}A_i^{n+1}$ be a finite
	allowed history set, not assumed to be a join of local supports. Assume
	prefix truncation carries $A_{n+1}$ into $A_n$. The constraints may come
	from operational admission, forbidden transitions, or other specified
	history conditions. These allowed sets are data, not conclusions of
	weighted compatibility. Write $p_{I,n}$ for the required site marginals.
	For a history $x\in A_n$ let $\ell(x)$ be its vector of atom indicators
	on the local sites (one coordinate for each local atom). Let $b_n$ be
	the vector of the local probability entries.

	\begin{theorem}[Exact weighted admission certificates at all horizons]
		\label{dyn49:thm:admission}
		For each $n$ the following are equivalent:
		\begin{enumerate}[label=(\roman*),leftmargin=*]
			\item a probability law supported on $A_n$ has all the prescribed marginals;
			\item $b_n\in\operatorname{conv}\{\ell(x):x\in A_n\}$;
			\item for every family of real local history functions $f_I$ with
			$\sum_I f_I(x_I)\ge0$ for all $x\in A_n$, one has
			$\sum_I\mathbb E_{p_{I,n}} f_I\ge0$.
		\end{enumerate}
		The empty set has empty convex hull. If the local laws are prefix
		consistent, a global probability process satisfying admission at every
		horizon exists if and only if these tests hold for every horizon.
		In particular every failure has a finite-horizon linear certificate
		on native site histories. For rational marginal data the functions in
		a strict certificate can be taken rational. The certificate width can
		be measured by the number of sites with nonconstant contributions,
		but is relative to this atlas and admission library, not an intrinsic
		rank of an individual local process.
	\end{theorem}

	\begin{proof}
		\leavevmode
		\begin{enumerate}[label={},leftmargin=0pt,nosep]
			\item \textnormal{(i)$\Leftrightarrow$(ii).}\enspace A probability law is a convex combination of the allowed atom vectors,
			proving (i)--(ii).
			\item \textnormal{(i)$\Rightarrow$(iii).}\enspace Taking expectations proves (iii).
			\item \textnormal{(iii)$\Rightarrow$(ii).}\enspace If the vector lies
			outside the nonempty compact convex hull, strict separation gives a
			linear combination of site functions with its value on $b_n$ strictly
			below the minimum on allowed atoms. Subtract this minimum as a constant
			on one site; the sum becomes nonnegative on $A_n$ with negative local
			expectation, contradicting (iii). If $A_n$ is empty the constant $-1$
			is already a certificate. Strict separation permits rational approximation
			of the coefficients followed by subtraction of a rational minimum over
			the finitely many atom vectors; strict negativity persists.
		\end{enumerate}
		For the
		all-horizon claim necessity follows by projection. Conversely take
		feasible laws at arbitrarily long horizons and restrict them to shorter
		ones. The prefix condition ensures admissibility, and local prefix
		consistency ensures the required marginal values. A diagonal subsequence
		in the compact finite probability simplices yields consistent limiting
		prefix laws; their supports stay in the finite allowed sets. These give
		a probability measure on the full history space, with admission almost
		sure at each horizon, hence at all horizons by countability.
	\end{proof}

	\begin{corollary}[Exact admitted strong Markov globalization and dimension]
		\label{dyn49:cor:admitted-markov}
		Let $X\subseteq\prod_i A_i$ be a nonempty finite admitted state set,
		and let $R\subseteq X^2$ be an additional allowed transition mask.
		\begin{enumerate}[label=(\roman*),leftmargin=*,nosep]
			\item \emph{Existence.} For
			site kernels $P_I$ on the attained site sets $X_I=\pi_I X$ and local initial distributions $\lambda_I$, a global
			Markov kernel on $X$ supported in $R$ with all site projections strong
			factors (including at unreachable states), and a matching initial law,
			exists exactly when the marginal test of Theorem~\ref{dyn49:thm:admission}
			holds for the initial distributions on $X$ and, for each $x\in X$, for
			the row distributions $P_I(x_I,\cdot)$ on $R(x)$.
			\item \emph{Larger nominal sites.} If a site kernel is originally given on a larger state set containing $X_I$, the same tests (using its full marginal alphabet) characterize row intertwining and matching initialization: feasibility forces its attained image to be closed under the site kernel and its initial law supported there. The restrictions to those images are strong factors; the maps to the larger sets need not be surjective.
			\item \emph{Dimension and rigidity.} For any nonempty such
			polytope let $Z$ be the allowed atoms positive in at least one feasible
			law, and $L_Z$ the marginal map restricted to those columns. Its affine
			dimension is $|Z|-\operatorname{rank}L_Z$. The total dimension for the
			initial law and all marked kernel rows is the sum of these dimensions;
			rigidity is equivalent to the sum being zero.
		\end{enumerate}
	\end{corollary}

	\begin{proof}
		\leavevmode
		\begin{enumerate}[label={(iii)},leftmargin=*,nosep]
			\item[(i)] Each strong factor equation for a kernel is exactly one row marginal
			constraint; initialization adds an independent set of marginal constraints.
			Feasibility row by row therefore defines a total Markov kernel and an
			initial law, and every such globalization occurs in their Cartesian
			product.
			\item[(ii)] If a larger target is used, every supported global row and the global initial law project into the attained image, so the corresponding local rows and initial law assign zero mass to its complement. This proves the stated restriction and distinguishes the full-target intertwining from a surjective factor.
			\item[(iii)] In any feasible polytope average one feasible law per atom in
			$Z$ to obtain a law strictly positive on $Z$. All other atom masses
			are forced zero. Thus the affine hull is exactly that law plus
			$\ker L_Z$, because sufficiently small variations along every kernel
			vector preserve positivity. Normalization is included in the site
			marginals. Dimensions add over independent rows and initialization. \qedhere
		\end{enumerate}
	\end{proof}

	For example let $X=\{00,01,10\}$ with the two singleton sites, and allow
	all transitions within $X$. Both sites start at zero. Both local kernels
	have the same next-one probability $3/4$ at every row. On the unrestricted
	product the sites admit a strong independent join; on $X$ no row is
	feasible, since $x_1+x_2\le1$ on every admitted state but the required
	expectation is $3/2$. Even disjoint sites with trivial overlap compatibility
	can therefore have an admission obstruction. The affine dimension formula
	also applies when the structural admission set is not a tree join, unlike
	the explicit tree cardinality formula of Section~\ref{dyn47:subsec:weighted-compose}.

	\subsection{Joint Hankel completion ranks and obstructions to finite realization}
	\label{dyn54:subsec:factorization}

	\paragraph{Joint Hankel completion ranks.}

	Write $\Sigma=\prod_{i\in V}A_i$ for the joint alphabet. An atlas specifies cylinder probabilities of local words but not those of all joint words. It is important not to mistake a nonnegative factorization of an individual or mixed native Hankel array for a common stochastic shift action. For later use, the following intrinsic hierarchy bounds state resources without imposing a common stochastic shift action.

	\begin{definition}[Joint Hankel completion ranks and state-budget lower bounds]
		\label{dyn54:def:completion-rank}
		\leavevmode
		\begin{enumerate}[label=(\roman*),leftmargin=*,nosep]
			\item \textbf{Joint Hankel completion rank.} For $h\ge1$ let $\mathcal P_{2h}$ be the compact polytope of prefix-consistent joint probability tables through length $2h$ that realize every prescribed local word probability through that length. If structural admission is part of the comparison, set the forbidden joint word probabilities equal to zero. The joint Hankel completion rank is \begin{equation}\label{dyn54:eq:completion-rank}
				\mathfrak r_h=\min_{p\in\mathcal P_{2h}}
				\operatorname{rank}_+\bigl(p(uv)\bigr)_{|u|,|v|\le h},
			\end{equation}
			where $\operatorname{rank}_+$ is nonnegative matrix rank and an empty completion set has value $+\infty$. These ranks are functions only of the declared native laws and admission mask (and finite horizon), not of their positive presentations.
			\item \textbf{Exact state-resource invariant.} The invariant $d_*$ is the least size of an unmarked stochastic presentation of any admitted globalization, with value $+\infty$ if none exists.
			\item \textbf{State-budget lower bounds.} Prefix restriction and the submatrix inequality for nonnegative rank give $\mathfrak r_h\le\mathfrak r_{h+1}$; factoring every finite Hankel block of a $d_*$-state realization into nonnegative prefix and suffix factors gives $\mathfrak r_h\le d_*$.
		\end{enumerate}
		The definition does not impose a common shift action.
	\end{definition}

	\paragraph{Probability-extendible cyclic atlases.}

	The following obstruction is different from incompatible marginals or from a nonfinite \emph{choice} of canonical globalization. Every probability extension will be forced to have infinite Hankel rank. The mechanism is saturation of a prescribed overlap probability; no external admission constraint is needed.

	\begin{theorem}[Cyclic infinite-rank obstruction and finite state-budget certificates]
		\label{dyn54:thm:triangle-infinite}
		\leavevmode
		\begin{enumerate}[label=(\roman*),leftmargin=*,nosep]
			\item \emph{Infinite-rank atlas.} There is a rational finite hidden Markov atlas on the three sites $SA,SB,AB$, with $S$ ternary and $A,B$ binary, which admits a unique global probability law but no finite linear, and hence no finite autonomous positive, realization of that law. Every site admits at most four hidden edge states.
			\item \emph{Finite certificates.} For each integer $d\ge1$ there is an effective finite horizon $h$ at which $\mathfrak r_h\ge d$; thus arbitrarily large state demands are already witnessed by finite native probability data.
		\end{enumerate}
	\end{theorem}

	\begin{proof}
		\leavevmode
		\begin{enumerate}[label={(iii)},leftmargin=*,nosep]
			\item[(i)] Use the separator alphabet $\{0,1,\#\}$. Put $a=1/2$, $b=3/8$, $t=1-a-b=1/8$ and $c=a+b=7/8$. A finite binary word $w$ is emitted, followed by $\#^\infty$. For the two laws $\mu,\nu$ of the terminal binary word, put
			\begin{equation}\label{dyn54:eq:word-laws}
				\mu(w)=t a^{\#0(w)}b^{\#1(w)},\qquad
				\nu(w)=t b^{\#0(w)}a^{\#1(w)}.
			\end{equation}
			Each is normalized by the geometric sum of $t c^n$. At site $SA$, let $A$ be a constant fair bit, with $w$ of law $\mu$ when $A=1$ and $\nu$ when $A=0$. At site $SB$ let $B$ be constant fair, with $w$ of law $\nu$ when $B=1$ and $\mu$ when $B=0$. Each is an edge presentation with states consisting of a frozen regime and a two-state running/stopped flag: the running state emits the respective probabilities of $0,1$ without changing the flag, and emits $\#$ with probability $t$ changing it to stopped; stopped emits $\#$ forever. The shared law is $(\mu+\nu)/2$.
			
			For any probability globalization, write $r_{ij}(w)=\Pr(A=i,B=j,S=w\#^\infty)$. Nonnegativity and the two local word laws imply
			\[
			r_{11}(w)\le m(w):=\tfrac12\min(\mu(w),\nu(w)).
			\]
			We compute the sum of this upper bound exactly. Under $\mu$, the distribution $g_k$ of the signed count $K=\#0-\#1$ satisfies $g_k=t\mathbf1_{k=0}+a g_{k-1}+b g_{k+1}$. On $\ell^1(\mathbb Z)$ the recursion has a unique solution since the norm of the convolution part is at most $c<1$. Direct substitution and summation give
			\[
			g_k=\tfrac14(2/3)^k\ (k\ge0),\qquad
			g_{-k}=\tfrac14(1/2)^k\ (k\ge0).
			\]
			In particular $\Pr_\mu(K=0)=1/4$, $\Pr_\mu(K<0)=1/4$. Swapping zeros and ones exchanges $\mu$ and $\nu$; on $K>0$ the smaller word mass is $\nu$ and on $K<0$ it is $\mu$. Hence $\sum_w\min(\mu(w),\nu(w))=3/4$ and $\sum_w m(w)=3/8$.
			
			Prescribe the $AB$ site law to be constant bits with probabilities $p_{11}=p_{00}=3/8$, $p_{10}=p_{01}=1/8$. This is a four-state rational presentation, and the overlaps of all three sites agree. For any globalization the upper bound on each $r_{11}(w)$ must be saturated: the sum of their nonnegative deficits is $3/8-p_{11}=0$. Consequently all four conditional word masses are forced:
			\begin{equation}\label{dyn54:eq:forced-coupling}
				r_{11}(w)=r_{00}(w)=m(w),\quad
				r_{10}(w)=\tfrac12(\mu(w)-\nu(w))_+,\quad
				r_{01}(w)=\tfrac12(\nu(w)-\mu(w))_+.
			\end{equation}
			Here $z_+=\max(z,0)$. Conversely these nonnegative masses sum for each $w$ to $(\mu(w)+\nu(w))/2$ and have exactly the required marginals. They construct the unique extension, since all trajectories have constant private bits and a finite separator word before absorption almost surely. In particular the support of the $A=1,B=0$ terminal binary words is exactly $\{w:\#0(w)>\#1(w)\}$. Intersecting with $0^*1^*$ shows this language is nonregular (the prefixes $0^i$ have distinct continuation languages). The positive support of a finite nonnegative word representation is regular: retain positive initial states and positive letter transitions, and accept paths with positive terminal weight. Taking the indicated terminal-marker language and deleting the fixed marker and private labels preserves regularity. Hence the forced law already has no finite positive realization.
			
			There is a stronger linear obstruction. Order joint letters as $(S,A,B)$. Let the $i$th row prefix be $(0,1,0)^{i+1}$ and the $j$th column suffix be $(1,1,0)^j(\#,1,0)$, $0\le i,j<d$. Their Hankel submatrix has entries $r_{10}(0^{i+1}1^j)$, strictly positive exactly when $i+1>j$. This is lower triangular with strictly positive diagonal. Its rank is $d$, for every $d$. Hence no finite linear realization exists.
			\item[(ii)] To obtain finite certificates from the prescribed atlas, let $L\ge2d$ and take any prefix-consistent joint table matching local words through length $L$. In this paragraph $r_{ij}(w)$ denotes a terminal-event mass in this finite table, rather than a full-process mass. The terminal events with $|w|<L$ and the event of no $\#$ through length $L$ partition its probability, with the private flags constant on its support by the local laws. The latter event has mass $c^L$. On terminal events the same upper bounds $m(w)$ hold. Since the $AB$ mass is already fixed by a one-letter word, the sum of the deficits $m(w)-r_{11}(w)$ over $|w|<L$ is at most $c^L$: it equals the mass with both flags one on the no-stop event minus the tail sum of $m(w)$. Each deficit is nonnegative. On each terminal word $w$, the prescribed $SA$ marginal gives $r_{10}(w)=\mu(w)/2-r_{11}(w)$, so the increase of $r_{10}(w)$ over its forced value is exactly that word's deficit. All words used in the displayed $d\times d$ block terminate by length $2d\le L$. Thus every entry differs from its forced value by at most $c^L$. If $H_d$ denotes the forced rational triangular matrix, choose $L$ so large that $d c^L\|H_d^{-1}\|_\infty<1$. Then every possible completion has an invertible block by the Neumann-series criterion. Taking $h\ge L$ (so all row/column words occur in~\eqref{dyn54:eq:completion-rank} and its constraints include length $L$) gives $\mathfrak r_h\ge d$, since nonnegative rank is at least ordinary rank. This is an effective finite-horizon certificate for each proposed finite budget; no marginal-infeasibility certificate exists for this atlas. \qedhere
		\end{enumerate}
	\end{proof}

	\begin{theorem}[Universal finite realizability of probability-extendible atlases]
		\label{dyn54:thm:conditional-spectrum}
		For a finite cover $\mathcal H$ of native coordinates, consider the property that every atlas of finite hidden Markov word laws admitting a global probability process admits a finite autonomous globalization. Quantify over choices of nonempty finite coordinate alphabets.
		\begin{enumerate}[label=(\roman*),leftmargin=*,nosep]
			\item \emph{Geometric classification.} This property holds if and only if
			\begin{equation}\label{dyn54:eq:conditional-spectrum}
				\boxed{\nu(\mathcal H)=\infty,}
			\end{equation}
			that is, the cover admits a join tree.
			\item \emph{Sharp obstructions.} On each other cover a probability-extendible counterexample can be chosen with coordinate alphabets of size at most three (unused coordinates can be singleton), rational local presentations of at most four edge-emitting states per site, and every global law of infinite Hankel rank. The all-binary stationary ergodic strengthening is Corollary~\ref{dyn54:cor:binary-stationary}, without the four-state budget assertion.
		\end{enumerate}
	\end{theorem}

	\begin{proof}
		\leavevmode
		\begin{enumerate}[label={(iii)},leftmargin=*,nosep]
			\item[(i)] Sufficiency is Theorem~\ref{dyn53:thm:universal}: existence of a probability extension ensures separator compatibility. Necessity follows from the counterexamples in (ii).
			\item[(ii)] For necessity use the finite chordal/conformal characterization of join trees in Section~\ref{sec:covers-assertions}. Start from the unique probability extension of Theorem~\ref{dyn54:thm:triangle-infinite}. Coordinates not used below have singleton alphabets.
			
			If the primal graph has an induced chordless cycle of length at least four, partition its cyclically ordered vertices into three nonempty contiguous blocks and assign a copy of $S,A,B$ respectively to the three blocks. A chart meets this cycle in at most two consecutive vertices. Thus its restriction of the copied process is either a single base coordinate, identical copies of one coordinate, or one of the three finite base pair laws. All local laws are rational finite-state, by letter projection or duplication. The copied process supplies a probability extension. In any other extension, consecutive copies within a block must agree almost surely at all times, since the corresponding edge is guarded. The three boundary edges prescribe the three base pair laws on the representative variables. Projection to the three representatives would therefore give the infinite-Hankel law of Theorem~\ref{dyn54:thm:triangle-infinite}, which cannot be the letter projection of a finite linear presentation. This proves the claim for nonchordal covers.
			
			Otherwise there is a minimally unguarded clique $M$ of size at least three: all proper subsets are guarded and no site contains $M$. Partition $M$ into three nonempty disjoint groups $M_S,M_A,M_B$. Use the additive group $\mathbb Z/3\mathbb Z$ for letters of $S$ (identify $\#$ with $2$) and $\mathbb Z/2\mathbb Z$ for $A,B$. For each time and each group, given the corresponding base letter, sample shares uniformly over group assignments summing to that letter. Make these choices conditionally independent across groups and times. If a group has one vertex the share is just the base letter. If a site sees a proper subset of a group, its shares there have the full independent uniform law, independent of the base letter and of the shares of the other groups: each specified proper assignment has the same number of completions to any group sum. This remains true for every finite time word by the independent sampling. A chart cannot contain all three complete groups. Hence its process is obtained from a base law on at most two coordinates by independent finite letter-randomization and uniform extra shares. A letterwise stochastic output kernel applied to a finite nonnegative matrix law remains finite: replace each letter matrix by its weighted sums according to that kernel. Thus every chart law is rational finite-state with at most four edge states: letter randomization does not increase the dimension of the chosen base marginal presentation, including the trivial marginal when no group is fully seen. The randomized share process gives one probability extension.
			
			In any globalization, decode each base letter by summing the shares of its entire group. For each pair of groups, their union is a proper subset of $M$, hence guarded. Its chart prescribes exactly the base pair law after decoding. The decoded process must be the unique base extension. A deterministic letterwise map preserves finiteness of linear realizations by summing letter matrices over preimages. By Proposition~\ref{dyn49:prop:residual-space}, finite Hankel rank would give just such a finite linear realization, so every global law has infinite Hankel rank. The proof permits alphabets of size two or three on the chosen obstruction coordinates; it asserts no obstruction for an arbitrarily fixed singleton or binary alphabet on every coordinate. The two cases exhaust non-join-tree covers. \qedhere
		\end{enumerate}
	\end{proof}

	The conditional property differs from universal globalization of overlap-compatible atlases: here the failure witnesses have passed \emph{every} finite admitted marginal polytope test with the unconstrained mask. In particular a sequence of successful finite-window tests with state counts that grow cannot prove finite autonomous existence. This conditional classification is also distinct from the block-bounded spectrum of Corollary~\ref{dyn53:cor:finite-wug}; the counterexample class does not impose its dependence-block bound.

	\begin{corollary}[Binary stationary ergodic local obstructions]\label{dyn54:cor:binary-stationary}
		A finite cover admits a join tree if and only if every probability-extendible atlas of stationary ergodic finite hidden Markov laws on that cover with binary coordinate alphabets admits a finite hidden Markov globalization. On each cover without a join tree the local presentations in a counterexample can be rational. The obstruction also rules out nonstationary finite hidden Markov globalizations.
	\end{corollary}

	\begin{proof}
		On a join tree apply Theorem~\ref{dyn53:thm:universal}. For the converse first produce stationary ergodic local laws over finite alphabets on a connected component carrying an obstruction.
		
		Take the initialized rational finite local laws of Theorem~\ref{dyn54:thm:conditional-spectrum} and their probability extension, restricted to the component containing the constructed obstruction. The cycle or minimally unguarded clique lies entirely in that component, so the decoding contradiction of that theorem still applies there. Adjoin a reset bit to each coordinate in this component. At each reset, every coordinate emits its fixed default letter with marker one; all other times the marker is zero. Globally use iid reset events of probability $\varepsilon\in(0,1)$, with independent fresh samples of the original probability extension for the successive reset-free segments. A stationary version exists by using a bi-infinite iid sequence of reset events and independent segment samples: almost surely the last reset before any given time is at finite distance in the past. The restriction to each site is precisely the stationary reset-augmented presentation of Theorem~\ref{dyn53:thm:universal}, with matrices $c M_b^I$ and $\varepsilon\mathbf1\lambda_I$ ($c=1-\varepsilon$) and stationary initialization~\eqref{dyn53:eq:reset}. Take rational $\varepsilon$. On states reachable from $\lambda_I$, the augmented hidden total matrix is irreducible (reset to $\lambda_I$ followed by a positive path reaches every such state) and aperiodic (a state of positive $\lambda_I$ has a self-loop at reset). Thus each stationary local observable law is ergodic. Other components may use a stationary deterministic law.

		In any globalization on the obstruction component, all the reset bits must coincide almost surely at all times: each guarded edge enforces equality of its endpoint bits, and the primal component is connected. Therefore a putative finite globalization can be regarded as a finite presentation with one common reset marker. Condition on an initial reset, discard subsequent reset symbols and divide all remaining letter matrices by $c$. Write $G_*$ for the common reset matrix, $G_a$ for the marker-zero global letter matrices, and $\widetilde\lambda$ for the putative global initial row, which need not be stationary. The event of an initial reset has probability $\varepsilon$, so $\alpha=\widetilde\lambda G_*/\varepsilon$ is a probability row. Put $N_a=G_a/c$ and $N=\sum_a N_a$. The event consisting of an initial reset, a prescribed local marker-zero word $w$, and $n$ arbitrary marker-zero letters has probability $\varepsilon c^{|w|+n}p_I(w)$ by the stationary local reset law. Common markers allow all other native letters to be summed out. Thus $\alpha\widehat N_w^I N^n\mathbf1=p_I(w)$ for every $w,n$; taking $w$ empty gives $\alpha N^n\mathbf1=1$. Lemma~\ref{dyn53:lem:harmonic} would extract a finite initialized globalization of the original atlas, contrary to its forced infinite Hankel rank. 
		It remains to encode these stationary laws in binary without losing the observable block boundaries. Choose an integer $k\ge1$ large enough to label every letter of each coordinate alphabet injectively by $k$ bits $b_1,\ldots,b_k$. Encode a letter by the word
		\[
		11\,0\,b_1\,0\,b_2\,0\,\cdots\,b_k\,0,\qquad L=3+2k .
		\]
		In any concatenation the adjacent symbols $11$ occur exactly at block starts. Encode the stationary joint process above in blocks of length $L$, and use a common uniform phase in $\{0,\ldots,L-1\}$. Equivalently, use a two-sided stationary version and then restrict to nonnegative times. Each local law is stationary and ergodic. Indeed on the space of paths with a marked phase the shift advances the phase cyclically, and its $L$th power on phase zero is the original ergodic shift. An invariant event thus has probability zero or one on phase zero and the same probability on every phase. The binary observation is a factor of this system. A finite rational local presentation samples one original edge at each boundary and stores its output code, position and destination until the next boundary. The hidden transition between boundary states is the original stochastic total matrix; intermediate code positions advance deterministically. At a uniform phase its law is a rational mixture of the $L$ phase laws obtained by advancing the rational stationary boundary law. The encoded joint process supplies a probability extension.
		
		In any extension, block starts coincide almost surely throughout the obstruction component. Each coordinate displays its own starts, every guarded pair has simultaneous starts under its prescribed law, and the component is connected. Fix one coordinate there and condition on its first two bits being $11$, an event of probability $1/L$. This is the phase-zero event of every site in the component. Skip the first $L$ bits. The subsequent decoded site words have precisely the original stationary reset-augmented laws: the phase-zero event is determined within each site by its first two bits on any one coordinate, and the site marginal conditional on this event is the law of its original stationary block process. Skipping one entire block shifts that block process by one step and hence leaves its law unchanged. Equality of the phase events almost surely lets the same conditioning be used for all sites, even when the globalization itself is not stationary. If the binary extension had a finite hidden presentation, write $Q_a$ for its joint binary letter matrices and $Q=\sum_a Q_a$. For the fixed reference coordinate let $Q_{[1]}$ be the sum of $Q_a$ over letters having bit one there. Its hidden row after the conditioned prefix and the remaining $L-2$ unrestricted letters is $L\lambda Q_{[1]}^2Q^{L-2}$, a probability row because the conditioning event has probability $1/L$ (here $\lambda$ is the initialized global row). Group subsequent transitions in blocks of $L$ and sum their matrix products by decoded letter. Extend the decoder arbitrarily to invalid blocks, which have conditional probability zero on every site and hence on the entire finite component. Summation over all length-$L$ binary blocks gives the total matrix $Q^L$. These matrices have stochastic total and would give a finite globalization of the reset-augmented atlas. The preceding harmonic extraction rules this out. Coordinates on other components may be constant. The binary coding therefore gives the claimed counterexample on the original cover.
	\end{proof}

	\paragraph{Structural admission and finite realization.}

	A rule excluding certain global words can change which hidden-state realizations are possible.

	\begin{corollary}[Regular structural admission can force nonfinite realization on a tree]
		\label{dyn54:cor:regular-mask}
		There is an atlas of rational finite hidden Markov laws on three disjoint sites, with a regular prefix-closed letterwise structural admission mask, that admits a probability globalization but admits no finite autonomous globalization. The mask consists solely of equality of designated interface copies at each time.
	\end{corollary}

	\begin{proof}
		Use disjoint copies of the three sites of Theorem~\ref{dyn54:thm:triangle-infinite}. Require the two copies of each of $S,A,B$ to agree at every emitted letter. This is a letterwise diagonal admission constraint recognized by an automaton with an accepting and a rejecting state. Copying the global probability law of that theorem supplies an admitted extension. Every admitted finite extension would project to the three base coordinates and give a finite extension of the cyclic atlas, contradiction. The disjoint sites have a join tree with empty separators; the added operational equalities, not the declared observation geometry, impose the cyclic comparison. This example concerns admitted unmarked probability laws; it does not assert that the extension is strongly causal.
	\end{proof}

	\section{Compositional invariants, conservation, and structural resilience}
	\label{sec:compositional-invariants-conservation}

	Theorem~\ref{dinv:thm:complete} establishes complete temporal records for subgroup-matching continuations. Theorem~\ref{dinv:thm:quotient} constructs the conserved quotient of a finite group evolution. Theorem~\ref{dinv:thm:emergence} identifies the extra conserved observations created by linear synchronization. A system invariant compares systems under the specified equivalence; a conserved observable is constant along trajectories; a forward-invariant set contains every successor of its members (Sections~\ref{dyn45:eq-invariant-kinds} and~\ref{dyn46:subsec-coupling-invariants}). The finite and linear specializations retain declared algebraic structure and synchronization constraints. The final exposure model uses the state-set fixed-point theory to characterize survival and intervention costs (Theorems~\ref{thm:shield} and~\ref{thm:rescue}).

	\subsection{Finite group dynamics: invariant cores and complete temporal records}
	\label{dinv:sec:records}

	This specialization adds autonomous endomorphisms to the finite-group matching data of Structural Packaging (see \hyperref[comp:structural]{Paper~II}).
	Unlike a single-time record, a temporal record retains joint boundary trajectories.
	The distinction	between all-time feasibility, exact projection and conservation is essential.
	No update is inferred from the primitive packaging axioms.

	\begin{definition}[Initialized group fragments, temporal assemblies, and continuation tests]
		\label{dinv:def:fragment}
		\leavevmode
		\begin{enumerate}[label=(\roman*),leftmargin=*,nosep]
			\item \textbf{Initialized group fragment.}
			An initialized group fragment is a tuple
			\[
			\mathcal G=(G,F,J,p,(f_d)_{d\in D}),
			\]
			where $G$ is a finite group, $F:G\to G$ is an
			endomorphism, $J\leq G$ is an initial subgroup,
			$p:G\to P$ is a boundary homomorphism into a finite
			product of interface groups, and $f_d:G\to D_d$
			are marked homomorphic readouts indexed by a finite
			set $D$.
			The subgroup $J$ need not be $F$-invariant.
			An isomorphism of fragments is a group isomorphism
			preserving the marked boundary and readouts,
			intertwining the endomorphisms, and carrying the
			initial subgroup onto its counterpart.
			
			\item \textbf{Temporal readouts and premises.}
			Define the boundary-history map and time-indexed
			readouts by
			\[
			p^\infty(g):=(p(F^t g))_{t\geq0},
			\qquad f_{d,t}:=f_d\circ F^t.
			\]
			A temporal premise is a finite set
			$T\subseteq D\times\mathbb N_0$.
			Its joint readout is $f_T:=(f_{d,t})_{(d,t)\in T}$,
			and its agreement kernel on initial states is
			\[
			K(T):=J\cap\bigcap_{(d,t)\in T}\ker f_{d,t},
			\qquad K(\varnothing):=J.
			\]
			
			\item \textbf{Temporal assembly.}
			Given finitely many fragments $\mathcal G_i$,
			a matching subgroup $C\leq\prod_i P_i$, and an
			output boundary homomorphism $\rho:C\to P'$,
			put
			\[
			A:=\prod_i G_i,\qquad
			F_A:=\prod_i F_i,\qquad p_A:=(p_i)_i.
			\]
			The state group and initial subgroup of the
			assembly are
			\begin{equation}\label{dinv:eq:core}
				W:=\bigcap_{t\geq0}F_A^{-t}(p_A^{-1}(C)),
				\qquad
				J_W:=W\cap\prod_i J_i.
			\end{equation}
			The assembled fragment has state group $W$,
			initial subgroup $J_W$, temporal law $F_A|_W$,
			and boundary map $\rho\circ p_A|_W$.
			Its readouts are retained component readouts
			$f_{i,d}\circ\operatorname{pr}_i|_W$ and any
			specified additional readouts $h\circ p_A|_W$,
			where $h:C\to D_h$ is a homomorphism.
			Private readouts retain their occurrence labels.
			Initialization is imposed only at time zero.
			
			\item \textbf{Continuation language.}
			A continuation is a finite expression formed from
			these assemblies with one distinguished input
			fragment. The identity expression is permitted.
			The unrestricted algebraic subgroup-matching
			language permits all such fragments, matching
			subgroups, and attachments, including
			non-surjective attachments and subgroup probes.
			
			\item \textbf{Temporal tests and equivalence.}
			For a temporal premise $T$ and a requested
			time-indexed readout $u=(d,t)$, the recovery test
			asks whether
			\[
			\forall g,g'\in J,\quad
			f_T(g)=f_T(g')\Longrightarrow f_u(g)=f_u(g'),
			\]
			equivalently, $K(T)\subseteq\ker f_u$.
			For a finite set $U$ of requested readouts,
			a numerical test returns $|J|$ or
			\[
			\left|\{f_U(g):g\in J,\ f_T(g)=a\}\right|,
			\qquad a\in f_T(J).
			\]
			The latter count is independent of the attained
			premise value $a$, since the fibers are subgroup cosets.
			Two fragments are \emph{recovery equivalent}
			(respectively, \emph{numerically equivalent}) if
			their recovery (respectively, numerical) test
			outcomes agree in every common continuation.
		\end{enumerate}
		All comparisons retain the marked boundary and readout
		targets and their labels. Tests are invariant under
		renaming private readout values and do not compare
		literal values across different times. These equivalences
		are instances of Definition~\ref{dyn45:eq-contract}.
	\end{definition}

	\paragraph{Native interpretation.}
	Every finite group $G$ admits a table presentation with
	generators $\dot g$ ($g\in G$) and relations
	$\dot g\dot h=\dot{gh}$ and $\dot 1=1$.
	Every word reduces to its value-name, and homomorphisms
	induce maps between these presentations
	\cite{MaStructural2026}.
	Thus the algebraic fragments and subgroup probes admit
	native interpretations in contexts permitting these
	presentations and attachments.
	The subgroup-matching language includes finite graphs of
	equality constraints and always admits the identity match.
	Additional native markings, admission filters, or rejection
	tests define different comparison contracts.

	\begin{theorem}[Greatest invariant core and sharp finite safety bound]
		\label{dinv:thm:core}
		Put $W_h=\bigcap_{t=0}^h F_A^{-t}(p_A^{-1}(C))$ for $h\geq0$.
		The groups $W_h$ descend. Their first equality $W_h=W_{h+1}$ implies
		$W_h=W$; $W$ is the largest forward-invariant subset of $p_A^{-1}(C)$ and is a
		subgroup.
		If $\Omega(n)$ counts prime factors of the positive integer $n$ with
		multiplicity, with $\Omega(1)=0$, then
		\[
		W=W_{\Omega(|W_0|)}.
		\]
		The same conclusion holds for the sequence of kernels of full labelled
		observation prefixes on an ambient group, by taking $C=\{1\}$ for
		the joint observation map. These bounds do not require invertible laws.
		The safety bound is sharp, even for elementary abelian groups. An
		initialized restriction $J_W$ is imposed after the ambient iteration;
		an early equality of $J\cap W_h$ alone does not certify stabilization.
	\end{theorem}

	\begin{proof}
		Preimages of subgroups under homomorphisms are subgroups, and
		$W_{h+1}=W_0\cap F_A^{-1}(W_h)$. Equality of two successive groups
		makes the earlier group forward invariant and makes every subsequent
		iterate equal. Conversely every forward-invariant subset of $W_0$
		is contained in each $W_h$, by following its forward orbits. A strict inclusion has integer index at
		least two and consumes at least one prime factor of the original order.
		There are at most $\Omega(|W_0|)$ strict steps. For sharpness take
		$A=\mathbb F_p^{n+1}$, $F(x_0,\ldots,x_n)=(x_1,\ldots,x_n,0)$,
		$p_A(x)=x_0$, and $C=0$; then $W_h$ imposes exactly the first
		$\min(h+1,n+1)$ coordinates to be zero. Stabilization first occurs at
		$h=n=\Omega(|W_0|)$. In this same example with $n=2$, the initialized
		line spanned by the last basis vector has equal intersections for
		$h=0,1$ but not for $h=2$.
	\end{proof}

	Define the joint boundary-history records (subgroups are actual marked
	subgroups, not just abstract isomorphism classes)
	\begin{equation}\label{dinv:eq:record}
		\begin{aligned}
			B(T)&=p^\infty(K(T)),& N(T)&=K(T)\cap\ker p^\infty,\\
			V(T)&=\{u\in D\times\mathbb N_0:N(T)\leq\ker f_u\},&
			m(T)&=|N(T)|.
		\end{aligned}
	\end{equation}
	The qualitative record is $(B(T),V(T))_T$; the numerical record is
	$(B(T),m(T))_T$. Since the vertical kernels are nested finite groups,
	the latter determines the former by
	$u\in V(T)\Longleftrightarrow m(T)=m(T\cup\{u\})$. Under a fixed native isomorphism with marked
	readouts these records are unchanged. Every trajectory of a finite
	endomorphism is eventually periodic, with a common preperiod and period
	for all its states. The records are finite data in the following
	precise sense. Let $s$ be the first time the images $F^tG$
	stabilize, and let $a$ be the order of the resulting automorphism on
	$F^sG$ (including order one for the trivial group). Every map $F^t$
	is identical to one with $0\leq t<s+a$. Thus the time labels can be
	reduced to this interval, the boundary words represented by their
	prefixes of length $s+a$ with periodic extension, and the tables
	indexed by subsets of the resulting finite descriptor set. For comparing
	two fragments, the larger preperiod and the least common multiple of
	periods give one common coding convention. It gives an
	effective representation when the finite group tables and readout values
	have effective encodings with decidable equality. Effective evaluation
	of assemblies also requires effective constraint and readout data on
	the attained finite images. No bound for every temporal recovery
	invariant is asserted.

	\begin{theorem}[Composition and completeness of temporal recovery and numerical records]
		\label{dinv:thm:complete}
		For Definition~\ref{dinv:def:fragment}, both records in
		\eqref{dinv:eq:record} compose exactly and are complete for their
		respective continuation equivalences.
		\begin{enumerate}[label=(\roman*),leftmargin=*,nosep]
			\item \emph{Composition of records.} For a temporal premise in an assembly, write $T_i$ for its component
			premises and $T_*$ for the pairs $(h,t)$ naming new readouts. Regard
			$C^{\mathbb N_0}$ as a subgroup of $\prod_i P_i^{\mathbb N_0}$ by
			reordering factors. Write $\rho^\infty$ for the pointwise map on
			$C^{\mathbb N_0}$. Set
			\[
			M_T=C^{\mathbb N_0}\cap\prod_i B_i(T_i)
			\cap\bigcap_{(h,t)\in T_*}\ker(h\operatorname{ev}_t),
			\qquad L_T=M_T\cap\ker\rho^\infty.
			\]
			Then
			\begin{equation}\label{dinv:eq:compose}
				\boxed{B_W(T)=\rho^\infty(M_T),\qquad
				m_W(T)=|L_T|\prod_i m_i(T_i),\qquad
				|K_W(T)|=|M_T|\prod_i m_i(T_i).}
			\end{equation}
			\item \emph{Recovery tests.} A retained component time-readout $u=(d,t)$ from $i$ belongs to
			$V_W(T)$ exactly when
			\begin{equation}\label{dinv:eq:qualitative}
				u\in V_i(T_i)\quad\text{and}\quad
				\operatorname{pr}_i L_T\leq B_i(T_i\cup\{u\}).
			\end{equation}
			A new time-readout $(h,t)$ belongs to $V_W(T)$ exactly when
			$L_T\leq\ker(h\operatorname{ev}_t)$. Substituting $M_T$ for $L_T$
			gives the unconditioned recovery tests (without fixing the new boundary).
			\item \emph{Numerical tests.} For every attained premise value the conditional range of a requested
			finite time-readout family $U$ is
			\[
			R_W(U\mid T)=\frac{|K_W(T)|}{|K_W(T\cup U)|}.
			\]
			\item \emph{Completeness and minimality.} Forgetting readouts restricts the record. Equality of records is
			necessary and sufficient for equality of all continuation tests in the
			specified unrestricted language. Sufficiency for its tests persists
			under restriction of the language, but necessity requires enough
			separating probes. In the unrestricted language necessity can be
			witnessed by finite group-native initialized probes:
			the trivial boundary history and cyclic groups generated by individual
			prime-power-order boundary histories. Every exact summary of the tests
			therefore factors onto the attained record classes by a unique map,
			in the sense of Proposition~\ref{dyn45:eq-classification}. Completeness
			here refers to these test outcomes; it does not identify all-data
			isomorphism classes of native initialized group systems.
		\end{enumerate}
	\end{theorem}

	\begin{proof}
		\leavevmode
		\begin{enumerate}[label={(iii)},leftmargin=*,nosep,beginpenalty=10000]
			\item[(i)] The joint boundary-history map sends $K_W(T)$ onto $M_T$: any tuple
			of histories in $M_T$ has independent initialized lifts in the respective
			$K_i(T_i)$; membership in $C$ at every time makes the tuple of lifts an
			admissible initial state in \eqref{dinv:eq:core}. Conversely each lift
			in $K_W(T)$ gives exactly those constraints. Its kernel is
			$\prod_i N_i(T_i)$. Hence we have a short exact sequence
			\[
			1\longrightarrow\prod_i N_i(T_i)\longrightarrow K_W(T)
			\longrightarrow M_T\longrightarrow1.
			\]
			All groups in this sequence are finite. Restricting the last map to
			the preimage of $L_T$ gives the vertical kernel of $W$; counting fibers
			proves \eqref{dinv:eq:compose}.
			\item[(ii)] If $u\notin V_i(T_i)$, take a local
			vertical witness that does not kill $u$ and identity states elsewhere.
			This refutes vertical recovery in $W$. If $u\in V_i(T_i)$, its value
			on each fiber of the local boundary-history map is constant, because
			any two lifts differ by an element of $N_i(T_i)$. It is the identity
			exactly on histories in $B_i(T_i\cup\{u\})$. Every history in $L_T$
			has a lift, proving \eqref{dinv:eq:qualitative}. The new readouts are
			functions of the boundary history alone. The unconditioned assertion
			uses the same argument with $M_T$.
			\item[(iii)] The premise fibers in $J_W$ are
			cosets of $K_W(T)$, and their images under the requested homomorphism
			have size $[K_W(T):K_W(T\cup U)]$. This proves the range formula.
			\item[(iv)] These operations and forgetting preserve record equality, so induction
			on finite continuations proves sufficiency.
			
			For necessity, first attach the one-element fragment by equality of the
			entire boundary history. Recovery tests on this continuation recover
			$V(T)$; numerical tests recover $m(T)$ from its order and the joint
			range on $T$ (the order divided by that range). It remains to recover
			$B(T)$. Suppose the two candidate boundary-history subgroups differ
			at a premise. One of them contains an element $e$ of prime-power order
			not in the other: a finite group is generated by its prime-power-order
			elements, since each element is a product of its prime-power components.
			Write $E=\langle e\rangle$, of order $\ell^r$ with $r\geq1$.
			Let $\sigma$ be the left shift on boundary histories. The subgroup
			$H$ generated by all $\sigma^j E$ ($j\geq0$) is finite: $e$ belongs to
			one fragment's boundary image and all its entries and their powers
			lie in that fragment's finite subgroup $p(G)\leq P$. Its shifts share
			a uniform eventual period and preperiod bound. Words over this finite
			group with those bounds form a finite group, and $H$ is a shift-invariant
			subgroup. This argument does not require the whole target $P$ to be finite. Form a finite group-native probe on $H\times E$
			with law $(b,a)\mapsto(\sigma b,a)$, boundary $\operatorname{ev}_0(b)$,
			initial subgroup $\{(a,a):a\in E\}$ and readout
			$z(b,a)=aE^\ell$. Its initialized boundary histories are exactly $E$.
			Attach it by equality of all boundary values at every time. The premise
			$T$ fails to determine $z$ at time zero exactly when
			$B(T)\cap E\nsubseteq E^\ell$, equivalently $E\leq B(T)$ because
			$E^\ell$ is the unique maximal subgroup of $E$. This distinguishes the
			two records. No equivariant surjectivity of the probe or invariance of its initial
			subgroup is required. The probe and attachment belong to the unrestricted algebraic language by definition. For its native interpretation, the table construction following Definition~\ref{dinv:def:fragment} realizes the finite group and the maps when the context admits them. Numerical tests also distinguish
			recovery (range one). Thus both
			necessity assertions hold. If an exact summary identifies two fragments,
			all tests agree, hence so do the records. This defines the unique
			surjection from attained summary codes onto attained record classes. \qedhere
		\end{enumerate}
	\end{proof}

	\begin{example}[Static equivalence need not imply temporal equivalence]
		\label{dinv:ex:static}
		Use $G=\mathbb F_p^2$, $J=G$, $p(x,y)=x$ and no private descriptors.
		The identity law and the swap $F(x,y)=(y,x)$ have identical static
		boundary records, but their boundary-history groups are respectively
		constant words and alternating words. The latter contains a nonconstant
		word of order $p$; the preceding probe separates them. Even agreement
		kernels at each separate time do not, in general, test conservation:
		on $\mathbb F_p$ for odd $p$, the laws $x\mapsto x$ and $x\mapsto-x$
		with private identity readout and trivial boundary have identical
		temporal recovery and numerical records. Only the former conserves the
		actual readout value. Thus literal-value conservation is outside the equality-of-fibers test language.
	\end{example}

	\paragraph{Exact trajectory coverage and conservation tests.}
	For initialized matching, the image of $J_W$ in component $i$ is the
	preimage in $J_i$ of $\operatorname{pr}_i M_\varnothing$ under
	$p_i^\infty$. All initialized native trajectories of component $i$
	(including their initial states) are realized exactly when
	$\operatorname{pr}_i M_\varnothing=B_i(\varnothing)$:
	the independent-lifting argument in the proof supplies each initialized
	state, not just a representative of its boundary history. Failure is an
	exact obstruction to coverage by any compatible initialized source.
	For designated temporal recovery targets, Theorem~\ref{dinv:thm:complete}
	gives the exact tests on the maximal initialized matching source; a
	missing counterexample cannot be restored by deleting initial states.
	Additional targets (in particular recovery theories whose positive rules
	are not forced on that source) need separate tests. These statements do
	not imply full-cut common information, uniqueness of arbitrary native
	sources or preservation of generator markings under replacement.
	To test literal conservation of a homomorphic readout $f$, expose it
	as an extra marked boundary coordinate before computing the record.
	Then $f(F^t g)=f(g)$ for every initialized state and all $t\geq0$ if
	and only if that coordinate of $B(\varnothing)$ consists of constant
	words. Thus value-sensitive boundary tests, not private agreement kernels
	alone, detect conservation.

	\subsection{Observable dynamics, conserved quotients and asymptotic structure}
	\label{dinv:sec:asymptotic}

	An endomorphism induces an update on observable cosets precisely when it preserves the relevant kernel.
	Quotienting further by the change produced at each step yields a conserved observation.

	\begin{theorem}[Observed and conserved quotients; stable image and entry time]
		\label{dinv:thm:quotient}
		\leavevmode
		\begin{enumerate}[label=(\roman*),leftmargin=*,nosep]
			\item \emph{Observed quotient.} For an autonomous finite group endomorphism $F:G\to G$ and a finite
			joint readout $r:G\to R$ that is a group homomorphism, put
			$N=\bigcap_{t\geq0}\ker(rF^t)$. The quotient $Q=G/N$, with induced
			endomorphism $\bar F$ and readout, is the coarsest surjective
			group-dynamical factor retaining $r$, unique up to the unique isomorphism
			commuting with its quotient map and marked readout. It is canonically isomorphic to
			the marked group $r^\infty(G)$ with shift. Its observation prefixes
			stabilize by $\Omega(|\ker r|)$ in the sense of
			Theorem~\ref{dinv:thm:core}. Every exact surjective homomorphic factor
			with a descended evolution and readout has a unique surjective
			comparison onto $Q$. This is the group-homomorphic instance of
			Proposition~\ref{dyn45:eq-observation-quotient}: the comparison is between
			readout-retaining factors of the fixed source. These factors classify
			observable behavior in that comparison, not unobserved native state
			structure.
			
			\item \emph{Conserved quotient.} For a group endomorphism $L:Q\to Q$, let $D_L$ be the normal closure
			in $Q$ of $\{x^{-1}L(x):x\in Q\}$. Then $Q/D_L$ is the universal
			conserved homomorphic quotient of $(Q,L)$: every homomorphism $a:Q\to H$
			into a group $H$ with $aL=a$ factors uniquely through it. For the
			observed quotient above this describes conserved homomorphisms factoring
			through $Q$, not necessarily all conserved homomorphisms on $G$.
			
			\item \emph{Stable image and entry time.} The stable image $S=L^sQ$ at the first
			stabilization time $s$ is the unique largest invariant subset on which
			$L$ is bijective; every state enters $S$ within $s$ steps, and
			$s\leq\Omega(|Q|)$. The function
			\[
			h(x)=\min\{j\geq0:L^j x\in S\}
			\]
			satisfies $h(Lx)=\max(h(x)-1,0)$. Also $|L^tQ|$ is nonincreasing.
			For arbitrary finite $Q$ these are exact finite-time statements, not
			metric or stochastic stability assertions.
		\end{enumerate}
	\end{theorem}

	\begin{proof}
		\leavevmode
		\begin{enumerate}[label={(iii)},leftmargin=*,nosep]
			\item[(i)] $N$ is normal and $F(N)\subseteq N$. The history homomorphism
			$r^\infty$ has kernel $N$ and intertwines $F$ with the shift, so the
			induced isomorphism $Q\cong r^\infty(G)$ preserves the marked readout
			at time zero. Any equivariant surjection preserving $r$ has kernel
			contained in $N$ by iteration.
			The quotient homomorphism theorem gives the factor and its uniqueness,
			and the prefix kernel assertion is Theorem~\ref{dinv:thm:core}.
			\item[(ii)] A homomorphism kills all $x^{-1}L(x)$ exactly when $aL=a$, hence kills
			their normal closure. This proves the quotient assertion, including
			well-defined trivial induced dynamics since the quotient map equals
			its composite with $L$.
			\item[(iii)] The descending image sequence stabilizes by
			counting subgroup indices; on $S$ the map is surjective and finite,
			hence bijective. Any subset on which $L$ is a bijection is contained
			in every image. The hitting-time identity follows by forward invariance
			of $S$ and the minimality of the first entry time; if $x\notin S$,
			removing the first step reduces the first entry time by exactly one. \qedhere
		\end{enumerate}
	\end{proof}

	\subsection{Linear synchronization and exact emergence of conservation laws}
	\label{dinv:sec:conservation}

	The next results address a stronger question than existence of first
	integrals: which conservation laws of a synchronized source are not
	restrictions of any conservation law of the unsynchronized components?
	They hold over arbitrary fields and also for real continuous-time linear
	flows. For the discrete case let $A$ be a finite-dimensional vector
	space over $k$ with linear evolution $F$ (the product of component
	laws), and set $D=F-I$. For the real continuous-time case let $F$ instead
	denote the infinitesimal generator of $e^{tF}$ and set $D=F$.
	For any invariant subspace $W\leq A$ write
	\[
	\mathcal C(W)=\{\ell\in W^*: \ell D|_W=0\},
	\qquad \mathcal E(W)=\mathcal C(W)/\operatorname{res}_W\mathcal C(A).
	\]
	These are linear first integrals and their classes modulo inherited
	integrals. The discrete identity is $\ell F=\ell$; the continuous
	identity follows by differentiating at zero and conversely by the
	exponential series. In both cases $D$-invariance and $F$-invariance
	are equivalent. For product laws $\mathcal C(A)$ is the direct sum of
	local linear conservation spaces, so this definition of emergence
	compares against all linear conservation laws of the independent product.
	This classifies the zero-defect restriction problem in the sense of Definition~\ref{dyn45:eq-classification-tasks}, with an exact defect space and dimension. A conjugacy used in this comparison must also
	carry the marked subspace and boundary observations to their counterparts.

	\begin{theorem}[Linear conservation defect: dimension, exact sequence, and zero-defect criterion]
		\label{dinv:thm:emergence}
		\leavevmode
		\begin{enumerate}[label=(\alph*),leftmargin=*,nosep]
			\item \emph{Defect.} For every $D$-invariant $W$ there is a canonical isomorphism
			\begin{equation}\label{dinv:eq:emergence}
				\boxed{\mathcal E(W)\cong ((W\cap DA)/DW)^*.}
			\end{equation}
			In particular the number of independent emergent conservation laws is
			$\dim(W\cap DA)-\dim DW$.
			\item \emph{Observability.} For a linear boundary constraint $p(x)\in C$
			with a linear map $p:A\to P$ of vector spaces and a subspace $C\leq P$,
			write $b:A\to P/C$ for the composite with the quotient map.
			If $n=\dim A$, the maximal synchronized state space is
			$W=\ker O$ for $O=(b,bF,\ldots,bF^{n-1})$ (the empty matrix for $n=0$).
			In the continuous case $O$ is the generator observability matrix, not
			a list of sampled flow readouts; its kernel imposes the constraint at
			every $t\geq0$. Its emergent dimension is
			\begin{equation}\label{dinv:eq:rank}
				\dim\mathcal E(W)=\operatorname{rank}O-\operatorname{rank}(OD)
					-\operatorname{rank}(O|_{\ker D}).
			\end{equation}
			Equivalently it is $\dim\ker(D|_W)-\dim\ker D+
			\operatorname{rank}O-\operatorname{rank}(OD)$.
			\item \emph{Exact sequence.} Fix a $D$-invariant subspace $N$ and write $\widehat D$ for the induced
			endomorphism of $A/N$. In the observed application $N$ consists of
			states invisible to the declared complete boundary histories. For any invariant $W\supseteq N$ put $C_D(W)=(W\cap DA)/DW$, and write $\widehat W=W/N$. There is a canonical short exact sequence
			\begin{equation}\label{dinv:eq:defect-sequence}
				0\longrightarrow (N\cap DA)/(N\cap DW)
				\longrightarrow C_D(W)\longrightarrow C_{\widehat D}(\widehat W)
				\longrightarrow 0,
			\end{equation}
			where $C_{\widehat D}(\widehat W)=(\widehat W\cap\operatorname{im}\widehat D)/\widehat D\widehat W$. Consequently
			\begin{equation}\label{dinv:eq:defect-dimension}
				\begin{split}
					\dim\mathcal E(W)={}&\dim\mathcal E(N)+\dim\mathcal E_{\widehat D}(\widehat W)\\
					&-\dim((N\cap DW)/DN).
				\end{split}
			\end{equation}
			Here $\mathcal E_{\widehat D}$ is defined on $A/N$ with the same conservation convention. The subtracted term measures the invisible defect already absorbed into the image of $D|_W$; there is no assertion of a canonical direct-sum splitting.
			\item \emph{Zero defect.} Then the following are equivalent:
			\begin{enumerate}[label=(\roman*)]
				\item Every invariant $W$ containing $N$ has no emergent linear
				conservation law ($\mathcal E(W)=0$).
				\item Both $N\cap DA=DN$ and
				$\ker\widehat D\cap\operatorname{im}\widehat D=0$.
			\end{enumerate}
			The second condition says that the zero primary component of $\widehat D$
			is semisimple (no nontrivial Jordan chain at zero). In discrete time it
			concerns eigenvalue one of the observed update. When $N$ is the full
			boundary-history kernel, every such $W$ is realizable by a linear
			constraint on a complete finite observation image (a discrete prefix,
			or in continuous time a finite set of flow observations or the generator
			observation matrix). Thus the criterion is exact for this enriched
			attachment language, not for an arbitrary fixed menu of instantaneous
			constraints.
		\end{enumerate}
	\end{theorem}

	\begin{proof}
		\leavevmode
		\begin{enumerate}[label={(iii)},leftmargin=*,nosep]
			\item[(a)] Conservation on $W$ is precisely annihilation of $DW$. Restrictions
			of integrals on $A$ are precisely those functionals on $W$ annihilating
			$W\cap DA$: if a functional vanishes there, extend it to $(W+DA)/DA$
			and then to $A/DA$. Restriction of $\mathcal C(W)$ to $W\cap DA$ is
			surjective onto the dual of $(W\cap DA)/DW$, by extending a functional
			on a subspace of $W/DW$. Its kernel consists exactly of inherited laws.
			This proves \eqref{dinv:eq:emergence} and canonicity.
			\item[(b)] Cayley--Hamilton gives the finite prefix formula: annihilating the
			first $n$ powers implies annihilating all powers. In continuous time
			annihilation of all powers is equivalent to annihilation of the flow
			by differentiation at zero and the convergent exponential series. Next,
			$\dim(W\cap DA)=\operatorname{rank}D-\operatorname{rank}(OD)$ and
			$\dim DW=\dim A-\operatorname{rank}O-\dim(W\cap\ker D)$.
			Subtract and use rank-nullity for $D$ to obtain the simplified formula.
			\item[(c)] Let $\pi:A\to A/N$. Inclusion induces the injection on the left in \eqref{dinv:eq:defect-sequence}. Projection induces the right-hand map because $\pi(DW)=\widehat D\widehat W$. Any $\widehat z\in\widehat W\cap\operatorname{im}\widehat D$ has the form $\pi(Da)$, and $Da\in W$ because $W$ contains $N$; hence that map is onto. If $z\in W\cap DA$ projects to $\widehat D\widehat W$, choose $w\in W$ with $\pi z=\pi Dw$. Then $z-Dw\in N\cap DA$, so the kernel is the indicated image of the injection. The dimension identity follows by subtracting the dimensions of $DN\subseteq N\cap DW\subseteq N\cap DA$, using \eqref{dinv:eq:emergence}. This is ordinary exact-sequence dimension accounting; no choice of complement is needed.
			\item[(d)] Assume (i). Taking $W=N$ yields the first equality in (ii). If a nonzero $v$ belongs to
			$\ker\widehat D\cap\operatorname{im}\widehat D$, choose $a\in A$
			with $\pi(Da)=v$. For $W=\pi^{-1}(kv)$ we have $Da\in W\cap DA$,
			while $\pi(DW)=0$, so $Da\notin DW$. This contradicts (i). Conversely
			assume (ii). The zero intersection and rank-nullity give the direct
			sum $A/N=\ker\widehat D\oplus\operatorname{im}\widehat D$.
			For any invariant $W\supseteq N$, the inverse of $\widehat D$ on its
			image preserves $\pi(W)\cap\operatorname{im}\widehat D$: this is an
			invariant subspace of a finite-dimensional invertible restriction and
			hence its restriction is invertible. For $z\in W\cap DA$, the class
			$\pi z$ lies there, so choose $w\in W$ with $\pi(Dw)=\pi z$.
			Then $z-Dw\in N\cap DA=DN\subseteq DW$. Thus $z\in DW$.
			This proves (i). Finally, when the boundary-history kernel is $N$, the
			finite observation map identifies $A/N$ with its image. In discrete
			time use the prefix map by Cayley--Hamilton; in continuous time a finite
			subfamily of the kernels of the flow readouts has the same intersection
			as all of them, by finite dimension (select a new readout whenever the
			current intersection is too large). Alternatively use the generator
			observation matrix. Each $W\supseteq N$ is the inverse image of the
			subspace that is its image under this map. Imposing this constraint
			at every time does not further shrink $W$ because it is invariant.
			This realizes all such constraints in the stated enriched language. \qedhere
		\end{enumerate}
	\end{proof}

	\begin{corollary}[Sharp unrestricted criterion and arbitrarily large defects]
		\label{dinv:cor:sharp}
		For a fully observed $A$ ($N=0$), absence of emergent linear laws on
		all invariant synchronized subspaces is equivalent to
		$\ker D\cap DA=0$. For $k^r\oplus k^r$ take the discrete law
		$F(x,y)=(x+y,y)$ or the continuous generator $F(x,y)=(y,0)$ and
		impose the linear interface constraint $y=0$. The synchronized source
		has $r$ independent emergent linear conservation laws, whereas every
		inherited integral is zero there. Thus the defect is unbounded even
		for invertible discrete updates, and even in the absence of cyclic
		spatial matching constraints.
	\end{corollary}

	\begin{proof}
		Take $N=0$ in the theorem. In the examples $D(x,y)=(y,0)$, so
		$DA=W=k^r\oplus0$ and $DW=0$. Formula~\eqref{dinv:eq:emergence}
		gives dimension $r$. The interface constraint may be implemented by
		matching the readout $y$ to a zero-state fragment. In discrete time
		the inverse law is $(x,y)\mapsto(x-y,y)$.
	\end{proof}

	\paragraph{Observable subspaces and conserved mass.}
	The invariant-core iteration is the usual unobservable-subspace
	construction in linear systems theory \cite{dinv:Sontag1998}; the
	finite-group bound uses subgroup indices. Group-valued behaviors and
	graphical group-code realizations have substantial established theories
	\cite{ForneyTrott1993,dinv:ForneyTrott2004,dinv:ForneyGL2012}. The temporal record theorem
	extends the static joint-boundary continuation record by initialized
	all-time matching and finite trajectory probes. Its proof uses elementary
	homomorphism and coset arguments. The linear defect formula is an
	annihilator/quotient calculation, closely related to the standard exact
	sequence for kernels and cokernels. The observable classification gives
	a precise condition for the emergence question under this declared
	constraint language. 
	In a nonlinear two-population exchange system, a conserved total mass restricts trajectories to mass leaves. Stability within a leaf and stability transverse to it are separate questions. Nonlinear conservation, controlled evolution and collective phase behavior provide further directions.

	\subsection{Monotone cascades, protective shields, and intervention costs}
	\label{subsec:monotone-cascades}

	The following finite exposure model applies the state-set fixed-point and viability principles of Section~\ref{dyn45:fp-safety}. An exposure package contributes a loss to a member once every other member has failed. This specifies the interactions and timing of the cascade. Theorem~\ref{thm:shield} characterizes the greatest surviving shield, and Theorem~\ref{thm:rescue} computes the minimum cost of preserving a specified function.

	Let $V$ be finite and $\mathcal H$ a finite indexed family of exposure packages (subsets of $V$). Repeated occurrences contribute separately. Fix finite nonnegative real buffers $b_i$, shocks $s_i$, and loss coefficients $L_{iI}$. Let $D_0=D_{\rm seed}\subseteq V$ be any prescribed set of initially failed units. With no initial failures take $D_{\rm seed}=\varnothing$.
	Define	
	\[
	\ell_i(D)=\sum_{I\in\HH:\,i\in I}L_{iI}\,\mathbf 1_{I\setminus\{i\}\subseteq D},\qquad
	D_{t+1}=D_t\cup\{i:s_i+\ell_i(D_t)>b_i\}.
	\]
	
	Write $\Phi(D_t)$ for the right-hand side of the update and $D_\infty$ for the stabilized failed set. The map $\Phi$ on $\mathcal P(V)$ is monotone and inflationary; its fixed-point structure is the finite powerset case of Tarski's theorem~\cite{Tarski1955}. From $D_{\rm seed}$ it reaches the least fixed point containing $D_{\rm seed}$ after at most $|V\setminus D_{\rm seed}|$ strict changes. Indeed every fixed point containing $D_{\rm seed}$ contains each iterate, and the stabilized iterate is fixed. A \emph{protective shield} is a subset $K\subseteq V\setminus D_{\rm seed}$ such that
	
	\[
	s_i+\ell_i(V\setminus K)\le b_i\quad\text{for every }i\in K.
	\]
	\begin{theorem}[Greatest protective shield and functional survival criterion]\label{thm:shield}
		The surviving set $K^*=V\setminus D_\infty$ is the largest protective shield. Let $g_P:\mathcal P(V)\to\{0,1\}$ be a specified nondecreasing functional-collapse predicate for $P$, with $1$ denoting failure. Then $P$ survives, meaning $g_P(D_\infty)=0$, precisely when there is some protective shield $K$ such that $g_P(V\setminus K)=0$.
	\end{theorem}

	\begin{proof}
		For a shield $K$, $D_0\subseteq V\setminus K$. If $D_t\subseteq V\setminus K$, monotonicity gives $s_i+\ell_i(D_t)\le s_i+\ell_i(V\setminus K)\le b_i$ on $K$. Thus $D_{t+1}$ also misses $K$ and $K\subseteq K^*$. At stabilization each survivor satisfies its buffer inequality against $D_\infty$ and lies outside the seed. Hence $K^*$ is a shield. Write $g_P(D)=1$ for functional failure, with $g_P$ nondecreasing under inclusion. If any shield satisfies $g_P(V\setminus K)=0$, then $D_\infty\subseteq V\setminus K$ implies $g_P(D_\infty)=0$. The converse uses $K=K^*$.
		
	\end{proof}

	\paragraph{Minimum-cost intervention before a cascade.}

	Suppose a pre-cascade intervention adds an unrestricted nonnegative real increment $u_i$ to buffer $b_i$ at finite unit cost $c_i\ge0$, with total cost $\sum_i c_i u_i$. The initial failure set, exposure law and collapse predicate are unchanged. For any proposed shield, the least sufficient increments are obtained coordinate by coordinate.

	\begin{theorem}[Exact minimum-cost pre-cascade shield rescue]\label{thm:rescue}
		The minimum total buffer-injection cost that prevents collapse of $P$ is
		\[
		\mathcal C_P^*=\min_{\substack{K\subseteq V\setminus D_{\rm seed}\\g_P(V\setminus K)=0}}
		\sum_{i\in K}c_i\,[s_i+\ell_i(V\setminus K)-b_i]_+,
		\]
		with value $+\infty$ when no functional shield is possible.
	\end{theorem}

	\begin{proof}
		Any successful intervention has a final survivor shield $K$ preserving $P$. For each $i\in K$, intervention must cover at least $[s_i+\ell_i(V\setminus K)-b_i]_+$, giving the lower bound because all costs are nonnegative. Conversely, injecting precisely those amounts on any candidate $K$ and zero elsewhere makes its shield inequalities hold. The greatest-shield theorem (Theorem~\ref{thm:shield}) then keeps all of $K$ surviving and preserves $P$. There are finitely many candidate shields, so when the family is nonempty its minimum is attained. If the family is empty, the necessity argument rules out every intervention.
	\end{proof}

	\section{Discussion and further directions}
	\label{sec:discussion}

	The framework developed here connects native formation,
	structural admission, temporal evolution and observed behavior.
	Its recovery criteria, realization theorems and conservation
	formulas establish a common mathematical foundation for
	the study of dynamical packaging across different
	temporal-law classes.

	\paragraph{Intrinsic dynamical recovery.}
	The sharp classwise recovery thresholds share the spatial
	obstruction rank $\nu(\mathcal H)$ but arise from distinct
	semantic mechanisms. Time expansion preserves the spatial
	obstruction rank, while temporal witness defects depend
	on the local laws. A natural direction is the development
	of a unified intrinsic dynamical recovery rank extending
	the static original-chart rank $\operatorname{ipr}$,
	together with its behavior under composition and
	structural reconfiguration.

	\paragraph{Finite autonomous realization and admission.}
	The geometric classifications establish sharp conditions
	for universal finite autonomous realization and
	block-preserving globalization. Further questions concern
	the characterization and decidability of finite realizability
	for individual probability-extendible cyclic atlases,
	the effect of structural admission, and the minimum
	realization resources under specified markings.

	\paragraph{Invariants and structural evolution.}
	The temporal records and conservation formulas establish
	compositional invariants for synchronized systems.
	Further directions concern conservation under nonlinear
	and nonhomomorphic observations, and coherence of
	temporal descent along structural reconfiguration paths.

	\paragraph{Collective extension.}
	The finite-system interface of
	Section~\ref{dyn46:subsec-collective-interface}
	provides a foundation for the study of interacting,
	size-indexed packaged systems, macroscopic evolution,
	and collective invariants.

	\section{AI use declaration}

	The author used OpenAI's ChatGPT and Anthropic's Claude in preparing this work.
	ChatGPT assisted in developing the following material:
	the observation-depth and word-distinguishing bounds (Sections~\ref{sec:fixed}, \ref{sec:local-dynamical-atlases}, and~\ref{sec:width});
	the capacity obstructions and temporal recovery results	(Section~\ref{sec:global-theorems});
	the finite-driver comparison and six-state coupling	(Section~\ref{sec:Finite-utonomous-globalization});
	block-preserving globalization (Theorem~\ref{dyn53:thm:block-preserving});
	binary stationary coding (Corollary~\ref{dyn54:cor:binary-stationary});
	and the cascade formulas (Section~\ref{subsec:monotone-cascades}).
	Claude assisted with assessing dependencies and originality, checking references, and editorial revision.
	The author takes full responsibility for the mathematical results, proofs, citations, and submitted text.

	\phantomsection

\end{document}